\documentclass[10pt,a4paper]{article}
\usepackage[utf8]{inputenc}
\usepackage{amsmath,amssymb,amsthm,mathtools}
\usepackage{esint}
\usepackage{enumitem}
\usepackage{tikz}
\usetikzlibrary{decorations.pathmorphing,arrows.meta,calc}
\setlist[itemize]{label=$\bullet$}
\usepackage[left=2.5cm,right=2.5cm,top=2.5cm,bottom=2.5cm]{geometry}
\usepackage{hyperref}

\numberwithin{equation}{section}
\newtheorem{theo}{Theorem}[section]
\newtheorem{lem}[theo]{Lemma}
\newtheorem{prop}[theo]{Proposition}
\newtheorem{cor}[theo]{Corollary}
\theoremstyle{definition}
\newtheorem{defin}[theo]{Definition}
\theoremstyle{remark}
\newtheorem{remark}[theo]{Remark}

\newcommand{\TT}{\mathbb{T}}
\newcommand{\RR}{\mathbb{R}}
\newcommand{\ZZ}{\mathbb{Z}}
\newcommand{\NN}{\mathbb{N}}
\newcommand{\Rs}{\mathcal R}
\newcommand{\PP}{\mathbb P}
\newcommand{\dd}{\,\mathrm{d}}
\newcommand{\Vc}{\mathcal V}
\DeclareMathOperator{\D}{div}
\DeclareMathOperator{\Supp}{supp}
\DeclareMathOperator{\dist}{dist}
\newcommand{\Sz}{\mathcal S_0^{d\times d}}
\newcommand{\Bs}{\mathcal B}
\newcommand{\QQ}{\mathbb Q}
\newcommand{\ot}{\mathbin{\mathring\otimes}}
\newcommand{\Yz}{\mathcal Y_0}
\newcommand{\Yo}{\mathcal Y_1}
\newcommand{\Rv}{{R^{\mathrm v}}}
\newcommand{\proofstep}[1]{\par\smallskip\noindent\textit{#1.}\ }

\title{Sharp nonuniqueness for the Vlasov--Navier--Stokes system}
\date{}
\author{Xingyu Li}

\begin{document}
\maketitle

\begin{abstract}
We consider the Vlasov--Navier--Stokes system on $\TT^d\times\RR^d$, $d\ge2$. For smooth data and every $p<2$, $q<\infty$, we
construct infinitely many weak solutions with velocity in $L^p_tL^\infty\cap L^1_tW^{1,q}$ and bounded, renormalized distribution
function; infinitely many of them violate the energy inequality. For smooth velocity data, weak solutions with velocity in
$L^2_tL^\infty$ are unique, so the result is sharp at the endpoint of the Ladyzhenskaya--Prodi--Serrin scale. The fluid carries the
nonuniqueness, through the convex integration scheme of Cheskidov and Luo, and the kinetic equation is solved exactly. The main
new difficulty is the memory of the particles: the kinetic state at time $t$ depends on the whole past of the velocity, so a
modification on a short interval changes the system at all later times. We overcome it by a causal gluing, carried out from left
to right in time, which re-solves the system on every interval of regularity. This yields wild solutions that are strong outside a
closed set of times of arbitrarily small Hausdorff dimension. A causal mollification extends the construction to rough data,
including every finite-energy velocity datum in dimensions two and three; in two dimensions the energy class is thus exactly the
borderline of uniqueness. Finally, in three dimensions we construct Leray solutions satisfying the strong energy inequality and
show that they are strong outside a closed set of times of vanishing $\mathcal H^{1/2}$ measure.
\end{abstract}

\section{Introduction}\label{sec:intro}

\subsection{The Vlasov--Navier--Stokes system}

We consider the Vlasov--Navier--Stokes system (VNS) on $\TT^d\times\RR^d$, $d\ge2$,
\begin{equation}\label{eq:VNSintro}
\left\{
\begin{aligned}
&\partial_tf+v\cdot\nabla_xf+\D_v\big[(u-v)f\big]=0,\\
&\partial_tu-\Delta u+\D(u\otimes u)+\nabla p=\int_{\RR^d}(v-u)f\dd v,\qquad\D u=0 .
\end{aligned}\right.
\end{equation}
It describes a dispersed phase of small particles, with distribution function $f(t,x,v)\ge0$, transported by an
incompressible viscous fluid with velocity $u(t,x)$. The phases interact through the drag (Brinkman) force: the particles
are accelerated by $u-v$, and the fluid feels the opposite force $\int(v-u)f\dd v=j_f-\rho_fu$, where $\rho_f=\int f\dd v$ and
$j_f=\int vf\dd v$. All physical constants are set to one. Smooth solutions satisfy the energy identity
\[
\frac{\mathrm d}{\mathrm dt}\Big(\frac12\|u\|_{L^2}^2+\frac12\iint|v|^2f\dd x\dd v\Big)+\|\nabla u\|_{L^2}^2+\iint|v-u|^2f\dd x\dd v=0,
\]
and conserve the total momentum $\int u\dd x+\iint vf\dd x\dd v$; the spatial mean of $u$ alone is not conserved.

Global weak solutions satisfying the energy inequality exist in various settings \cite{ABMB97,BDGM09,Yu13,WY15,BGM17}. In two
dimensions, Han-Kwan, Miot, Moussa and Moyano proved that weak solutions in the energy class are unique, for kinetic data
satisfying a mild decay condition in velocity \cite{HMMM20}. The long-time behavior \cite{HMM20} and hydrodynamic limits
\cite{GJV04} have also been studied extensively. In three dimensions, uniqueness of energy solutions is open, as for the
Navier--Stokes equations; a weak--strong uniqueness principle beyond the Osgood class was obtained in \cite{HMMM25}. In view of
the nonuniqueness of Leray solutions of the forced Navier--Stokes equations \cite{ABC22}, uniqueness is not to be expected in
general. To our knowledge, nonuniqueness of weak solutions of \eqref{eq:VNSintro} has not been established. We prove it in a class that is
sharp with respect to the Ladyzhenskaya--Prodi--Serrin scaling, first for smooth data and then for rough data.

\subsection{Convex integration and the memory of the particles}\label{sec:memory}

Since the work of De Lellis and Sz\'ekelyhidi \cite{DLS09,DLS13}, convex integration has produced nonunique weak solutions
of many equations of fluid dynamics; see Isett's resolution of Onsager's conjecture \cite{Ise18}. For the Navier--Stokes
equations, Buckmaster and Vicol \cite{BV19} constructed nonunique weak solutions in $C_tL^2$, and Buckmaster, Colombo and Vicol
\cite{BCV22} obtained wild solutions whose singular set in time has Hausdorff dimension less than one. For the transport
equation, see \cite{MSz18,CL21}. Cheskidov and Luo \cite{CL22} proved that weak solutions of the Navier--Stokes equations are
nonunique in $L^p_tL^\infty\cap L^1_tW^{1,q}$ for every $p<2$, $q<\infty$, while solutions in $L^2_tL^\infty$ are unique, a classical fact
(see \cite[Theorem~1.3 and Appendix~A]{CL22}). Their wild solutions are smooth outside a closed set of times of arbitrarily small
Hausdorff dimension; this rests on a concentration step that glues local exact solutions onto the approximate one
\cite[Section~3]{CL22}. They also obtain arbitrary $L^2$ data for $2\le d\le4$, by modifying a Leray--Hopf solution after a time at
which it is regular \cite[Theorem~1.5 and Remark~1.7(3)]{CL22}.

The drag force makes \eqref{eq:VNSintro} a Navier--Stokes system with a self-consistent forcing, and in our construction the nonuniqueness
comes from the fluid. The new feature is that this forcing has \emph{memory}. The distribution function at time $t$ is the
push-forward of $f_0$ along the phase-space flow of $u$ on $[0,t]$, so it depends on the whole past of the velocity. A
modification of $u$ on a short time interval changes the particle state, hence the drag force, at all later times. This causes
no difficulty for nonuniqueness alone, where every perturbation is small globally in time, but it affects all arguments that are
local in time: in \cite[Section~3]{CL22} the approximate solution is replaced by local exact solutions on good intervals and left
untouched elsewhere, and for \eqref{eq:VNSintro} a local exact solution computed with the old particle state is no longer a solution
once anything earlier has changed. The main technical contribution of this paper is a \emph{causal} version of the gluing,
carried out from left to right in time, in which every good interval is re-solved with the particle state produced by
everything done before it. It gives intervals of regularity (Theorem~\ref{thm:C}). Causality enters in a different form in the
extension to rough data (Theorem~\ref{thm:D}): there the background is mollified causally, using only its past, so that the particle
state on an initial interval is left untouched.

\subsection{Main results}

Weak solutions are defined in Definition~\ref{def:weak}. In particular, $u\in L^2_{t,x}\cap L^1_tL^\infty$, and
$f\in L^\infty$ is nonnegative and compactly supported in $v$. Solutions of Leray--Hopf type are those satisfying the
energy inequality (Definition~\ref{def:LH}).

\begin{theo}[Nonuniqueness]\label{thm:A}
Let $d\ge2$, $1\le p<2$ and $1\le q<\infty$. Let $u_{\rm in}\in C^\infty(\TT^d)$ be divergence-free and let
$f_0\in C^\infty_c(\TT^d\times\RR^d)$ be nonnegative. Assume that \eqref{eq:VNSintro} has a smooth solution on some interval
$[0,t_1]$ with data $(u_{\rm in},f_0)$. Then there are infinitely many weak solutions $(u,f)$ of \eqref{eq:VNSintro} on $[0,1]$
with these data, all coinciding with the smooth solution on $[0,t_1/4]$, such that
\begin{enumerate}[label=(\roman*)]
\item $u\in L^2_{t,x}\cap L^p_tL^\infty\cap L^1_tW^{1,q}$;
\item $f$ is bounded, has uniformly bounded velocity support, belongs to $C_tL^s$ for every $s<\infty$, and is the unique
bounded (hence renormalized) solution of the kinetic equation with velocity field $u$;
\item the total momentum is conserved.
\end{enumerate}
Infinitely many of these solutions are not of Leray--Hopf type.
\end{theo}

Theorem~\ref{thm:A} follows from Theorem~\ref{thm:main} and Corollary~\ref{cor:nonunique}. Theorem~\ref{thm:main} is a
flexibility statement: modulo spatial means, the velocity can be prescribed up to an arbitrarily small error in
$L^p_tL^\infty\cap L^1_tW^{1,q}$. The endpoint $p=2$ cannot be reached.

\begin{theo}[Uniqueness in the critical class]\label{thm:B}
Let $d\ge2$ and $T\in(0,1]$. Let $u_{\rm in}\in C^\infty(\TT^d)$ be divergence-free and let $f_0\in L^\infty$ be nonnegative and
compactly supported in $v$. Then there is at most one weak solution of \eqref{eq:VNSintro} on $[0,T]$ with data
$(u_{\rm in},f_0)$ and $u\in L^2(0,T;L^\infty)$. Such a solution satisfies $\nabla u\in L^2(0,T;L^\infty)$, its distribution
function is transported by the classical phase-space flow, and it satisfies the energy equality. In particular, on any
interval on which a smooth solution exists, the smooth solution is the only weak solution with $u\in L^2_tL^\infty$, while for
every $p<2$ there are infinitely many weak solutions with $u\in L^p_tL^\infty$.
\end{theo}

Theorem~\ref{thm:B} is Theorem~\ref{thm:sharp} together with Corollary~\ref{cor:dichotomy}; Corollary~\ref{cor:roughtrace} extends
it to $u_{\rm in}\in B^1_{q,2}$, $q>d$. The causal gluing gives the analogue of the intervals of regularity of \cite{CL22}.

\begin{theo}[Intervals of regularity]\label{thm:C}
Under the assumptions of Theorem~\ref{thm:A}, and for every $\varepsilon>0$, the solutions can be chosen so that, in
addition, there is a closed set $\mathcal S\subset(0,1)$ with $\dim_{\mathcal H}\mathcal S\le\varepsilon$ such that, on every compact
interval $\mathcal K\subset(0,1]\setminus\mathcal S$, $u\in L^2(\mathcal K;W^{1,\infty})\cap C(\mathcal K;W^{1,s})\cap L^2(\mathcal K;W^{2,s})$ for every $s<\infty$. On such $\mathcal K$ the
solution is the only weak solution with velocity in $L^2(\mathcal K;L^\infty)$ starting from its own state at $\inf \mathcal K$.
\end{theo}

Theorem~\ref{thm:C} is Theorem~\ref{thm:regint}. Off the singular set the solutions are strong, but in general not $C^\infty$: the
distribution function inherits the limited regularity produced by earlier perturbations, and the drag force passes it on to the
velocity. This
is a direct effect of the memory, and it is why uniqueness on $\mathcal K$ is stated from the state at $\inf \mathcal K$ rather than from the
initial data.

The next result removes the smoothness assumptions on the data. The kinetic datum only needs to be bounded, and the
velocity datum only needs to generate a weak solution that becomes Lipschitz after some time.

\begin{theo}[Rough data]\label{thm:D}
Let $d\ge2$, $1\le p<2$, $1\le q<\infty$. Let $f_0\in L^\infty(\TT^d\times\RR^d)$ be nonnegative and vanish for large $|v|$, and let
$u_{\rm in}$ be divergence-free with
\begin{enumerate}[label=(\alph*)]
\item $u_{\rm in}\in B^1_{q_0,2}$ for some $q_0>d$ (for instance $u_{\rm in}\in H^s$, $s>\frac d2$), if $d\ge2$; or
\item $u_{\rm in}\in H^1$, if $d\in\{2,3\}$; or
\item $u_{\rm in}\in L^2$, if $d\in\{2,3\}$.
\end{enumerate}
Then there are infinitely many weak solutions of \eqref{eq:VNSintro} on $[0,1]$ with data $(u_{\rm in},f_0)$ such that
$u\in L^2_{t,x}\cap L^1_tL^\infty$; $f$ is bounded, has bounded velocity support, belongs to $C_tL^s$ for every $s<\infty$, and is
the unique bounded (hence renormalized) solution of the kinetic equation with velocity field $u$; and the total momentum is
conserved. Moreover $u\in L^p(t_1,1;L^\infty)\cap L^1(t_1,1;W^{1,q})$ for some $t_1\in[0,\frac12)$, with $t_1=0$ in cases (a), (b), and
$u\in L^p(0,1;L^\infty)$ unless $d=3$ in case (c). Infinitely many of them are not of Leray--Hopf type. For $d=2$ in case (c), at
most one of them belongs to $L^\infty(0,1;L^2)\cap L^2(0,1;H^1)$; in particular at most one is of Leray--Hopf type.
\end{theo}

Theorem~\ref{thm:D} follows from Theorems~\ref{thm:rough} and~\ref{thm:admissible} and Remark~\ref{rem:rough}. For $d=2$ in
case (c), together with the uniqueness theorem of \cite{HMMM20}, whose decay condition on $f_0$ holds here, it shows that the energy
class is exactly the borderline of uniqueness for \emph{every} finite-energy datum: uniqueness holds in $L^\infty_tL^2\cap L^2_tH^1$,
even without the energy inequality, and fails as soon as one leaves it, even for velocities in $L^p_tL^\infty$, $p<2$, and
bounded particle densities.

The tools developed for rough data also give a Leray-type structure theorem in three dimensions.

\begin{theo}[Structure of Leray solutions, $d=3$]\label{thm:E}
Let $u_{\rm in}\in L^2(\TT^3)$ be divergence-free, and let $f_0\in L^\infty$ be nonnegative and vanish for large $|v|$. Every Leray
solution of \eqref{eq:VNSintro} satisfying the strong energy inequality is strong outside a relatively closed set of times
$\mathcal S\subset(0,1)$ with $\mathcal H^{1/2}(\mathcal S)=0$. More precisely, the conclusions of Theorem~\ref{thm:structure} hold.
Such solutions exist for every such datum (Theorem~\ref{thm:SEI}).
\end{theo}

Together with Theorem~\ref{thm:D}(c), this gives the following picture for $d=3$ and every $u_{\rm in}\in L^2$: Leray solutions
with the strong energy inequality exist and are strong on an open set of full measure, while there are infinitely many weak solutions
with velocity in $L^1_tL^\infty$ that coincide with one of them on an initial interval and fail the energy inequality.

\subsection{Main ideas}\label{sec:ideas}

\paragraph{Perturbing the fluid only.} One could try to let the particles produce the nonuniqueness, as the density
does for the inhomogeneous Navier--Stokes equations \cite{Li26}. There, small density weights $a_k$ are moved from the fluid velocity
$u$ to concentrated profiles transported with large velocities $V_k$, the fluid being set equal to $V_k$ on their support; the
mismatch $a_k(V_k-u)$ and $a_k(V_k\otimes V_k-u\otimes u)$ of the mass and momentum fluxes cancels both defects, because the density enters
the single momentum equation through $\rho u\otimes u$. This fails for VNS for two structural reasons. Particles moving with the local
fluid velocity exert no drag, and the momentum flux of particles that solve the kinetic equation is balanced by their own
momentum, so they cannot cancel a Reynolds stress of the fluid. And in phase space the flux of the particles at a given velocity is
that velocity times their mass: a flux defect at bounded velocities can be changed only by changing the distribution function
there by a comparable amount, while the force, the only field that can be perturbed, produces a flux in $v$ whose oscillatory part
is not balanced by the transport of the profiles. We therefore perturb the fluid as in \cite{CL22} and solve the kinetic equation
exactly at every step: $f_n=f[u_n]$ is the push-forward of $f_0$ along the phase-space flow of $u_n$, and the approximate
solutions solve a VNS--Reynolds system (Definition~\ref{def:VNSR}) in which only the fluid equation carries a defect. Since
$u\in L^1_tW^{1,q}$, the DiPerna--Lions theory \cite{DPL89,Amb04} applies to $(v,u-v)$, and the distribution function is the
renormalized solution determined by $u$. The situation is thus opposite to that of \cite{Li26}, where the density is not
renormalized and the velocity is not bounded.

Adding a fluid perturbation $w$ changes $f[u]$ to $f[u+w]$. The main estimate is a dual stability bound
(Lemma~\ref{lem:dual}): for every test function $\psi$,
\[
\Big|\iint\big(f[u+w]-f[u]\big)(t)\,\psi\Big|\le e^{1+\|\nabla u\|_{L^1_tL^\infty}}[\psi]_{C^1}\,M_0\,\|w\|_{L^1_tL^\infty}.
\]
Only the Lipschitz norm of the old field enters; the gradient of $w$, which is very large, is not used. The constant
depends only on the previous step, so it can be compensated by choosing $\delta_n$ after $u_n$, while the temporal intermittency of
\cite{CL22} makes $\|w\|_{L^1_tL^\infty}$ arbitrarily small. The change in the drag force is converted into a small Reynolds stress by
an antidivergence estimate in a weak norm (Lemma~\ref{lem:weakR}), and a spatially constant correction, obtained by a Volterra fixed
point, restores the exact momentum balance (Lemma~\ref{lem:mean}). This reduces the iteration to a Navier--Stokes step
(Proposition~\ref{prop:NSstep}), which is the step of \cite{CL22} with minor adjustments (Section~\ref{sec:NSdiff}; full proof in
Appendices~\ref{sec:CI} and~\ref{sec:estimates}).

\paragraph{Memory and causal gluing.} For Theorem~\ref{thm:C}, the stress of an approximate solution is to be pushed
onto few short intervals, the \emph{bad} set, and the solution on the long \emph{good} intervals between them is to be replaced by an
exact solution. By the memory effect described in Section~\ref{sec:memory}, the replacement on a good interval $[a,b]$ depends on
every modification made before $a$, including the convex integration perturbations on earlier bad intervals, so the construction
must be causal (Section~\ref{sec:proofW} and Figure~\ref{fig:memory}). It maintains a current global solution and sweeps $[0,1]$ from
left to right in two passes. In the first pass, a \emph{segment corrector} (Proposition~\ref{prop:segcorr}) re-solves each good
interval with the current particle state and glues the result in thin layers, so that the stress is concentrated on the bad set. In
the second pass, a localized convex integration step is applied at each junction, and the next good interval is re-solved with the
particle state that this step has produced. The two passes cannot be merged (see Section~\ref{sec:proofW}). In each corrector, the
drag difference created by all earlier modifications enters as a forcing, not as an error to be absorbed later.

The scheme closes because of the following four points.
\begin{itemize}
\item \emph{A linear kinetic duality} in $(W^{1,1})^*$, which measures the perturbation in $L^1_tL^q_x$ rather than $L^1_tL^\infty_x$
(Lemma~\ref{lem:B2}). The instantaneous part of the drag difference, which is not controlled by the past of the perturbation, is kept
separate (Remark~\ref{rem:shift}).
\item \emph{A two-norm fixed point} for the corrector: bounded in maximal regularity, small in a space with one derivative less, and
hence, by interpolation, small in $L^2_tW^{1,\infty}$.
\item \emph{Constants fixed before the scales.} Kinetic comparisons are started at the beginning $a$ of each segment from the common
particle state, so the size of the base is measured only after $a$. Earlier convex integration perturbations, which are small in
$L^1_tW^{1,q}$ but large in $L^\infty_tW^{1,q}$, vanish there and do not enter. This allows us to fix all upper bounds in two levels,
one per pass, before the scales are chosen.
\item \emph{Mean-preserving gluing} (Proposition~\ref{prop:glue}). The spatial mean of the fluid is driven by the momentum exchange with
the particles, and so it too has memory. Cutting off the mean of a corrector would create a force with nonzero mean, which is not
the divergence of a stress. We cut off only the mean-zero part; the mean is carried through the gluing layers and follows the exact
momentum balance.
\end{itemize}

\paragraph{Rough data.} The memory also governs the extension to rough data. The kinetic estimates depend on the velocity only through
its Lipschitz norm integrated in time, and extend to fields in $L^1_tW^{1,\infty}$ and bounded kinetic data via the Carath\'eodory flow
of $(v,u-v)$ (Lemma~\ref{lem:kinstrong}). The convex integration step needs a smooth background; we run it on a \emph{causal} space--time
mollification of the background, which only looks at the past, and put the mollification error into the new stress
(Lemma~\ref{lem:NSrough}). Causality keeps the stress zero on an initial interval and leaves the particle state there untouched. The
iteration starts from a weak solution that is Lipschitz only after some time $t_1$ (Definition~\ref{def:strongcont}); all kinetic
comparisons start at $t_1$ from the common particle state (Lemma~\ref{lem:concat}), so before $t_1$ only $u\in L^1_tW^{1,1}\cap L^1_tL^\infty$
is used. Such strong continuations are obtained by maximal regularity (Proposition~\ref{prop:LWPX}), by an $H^1$ theory for $d\le3$
whose existence time depends only on the $H^1$ norm and on bounds for $f$ (Proposition~\ref{prop:LWPH1}), and, for $L^2$ data, from an
energy solution after a time at which it is in $H^1$. Before that time the velocity support of $f$ must stay bounded, which requires
$u\in L^1_tL^\infty$: for $d=2$ we prove a periodic version of \cite[Lemma~3]{HMMM20} for energy solutions (Lemma~\ref{lem:L2Linf2D}), and
for $d=3$ we propagate the third and sixth velocity moments using only the energy inequality (Lemma~\ref{lem:L1Linf}). In both
dimensions the energy solution is the one of Theorem~\ref{thm:SEI}, so the only result of \cite{HMMM20} used in this paper is their
uniqueness theorem \cite[Theorem~1]{HMMM20}, applied on $[0,1]$; its proof is local in time.

\paragraph{Sharpness.} For Theorem~\ref{thm:B}, a duality argument as in \cite[Appendix~A]{CL22}, after \cite{LM01}, identifies any
$L^2_{t,x}$ weak solution of the Oseen problem with drift in $L^2_tL^\infty$ with its strong solution. Since the drag force of a
critical solution is only bounded in $x$, maximal $L^2_tL^q$ regularity is used to obtain $\nabla u\in L^2_tL^\infty$, which makes the
phase-space flow Lipschitz. Uniqueness then follows from a Gronwall estimate for
\[
\|u_1-u_2\|_{L^2}^2+\int f_0(z)\,|Z_1(t,z)-Z_2(t,z)|^2\dd z,
\]
where $Z_i$ are the phase-space flows, following Loeper's functional \cite{Loe06} as used in \cite{HMMM20}; a maximal-function
inequality replaces the $C^1$ norm of the test function.

\paragraph{Structure theorem.} The same Lagrangian argument gives weak--strong uniqueness against Lipschitz strong solutions
(Lemma~\ref{lem:wsT3}). Leray solutions with the strong energy inequality are constructed by a Galerkin scheme for the fluid in which
the kinetic equation is again solved exactly (Theorem~\ref{thm:SEI}); the approximate solutions then satisfy the energy equality on
every interval. With the velocity bound and the $H^1$ theory, Leray's covering argument gives Theorem~\ref{thm:E}.

\paragraph{The inviscid case.} The convex integration step does not use the Laplacian (Remark~\ref{rem:final}), and the kinetic estimates
of Section~\ref{sec:kinetic} are insensitive to viscosity; we expect the nonuniqueness part of Theorem~\ref{thm:main} to extend to the
Vlasov--Euler system, but have not checked the details. The results of Sections~\ref{sec:sharp} and~\ref{sec:regint} use the viscosity in
an essential way.

\subsection{Further questions}

We list some problems suggested by this work.
\begin{enumerate}
\item \emph{Particle-driven nonuniqueness.} Can the particles, rather than the fluid, carry the nonuniqueness? This would
yield distribution functions that are not renormalized, and velocities with more Sobolev regularity in space and continuity in
time, as in \cite{Li26}, where the velocity is in $C_tH^s$, $s<\frac32$, but is not bounded.
The obstructions described above suggest that a new building block is needed.
\item \emph{Kinetic transport with rough forces.} For the linear equation $\partial_tf+v\cdot\nabla_xf+\D_v(Ef)=0$,
Champagnat and Jabin proved well-posedness of the Lagrangian flow for $E\in H^{3/4}$ \cite{CJ10}. Whether non-Lagrangian
solutions exist for Sobolev force fields below this threshold is, as far as we know, open. Convex integration in phase
space is constrained by the fact that the $x$-component of the field is $v$ and cannot be perturbed.
\item \emph{$C^\infty$ regularity} off the singular set in Theorem~\ref{thm:C}, and nonuniqueness for arbitrary $L^2$ data when
$d\ge4$ (the cases $d\le3$ are covered by Theorem~\ref{thm:D}). Already for the heat flow of a general $L^2$ datum,
$\|e^{t\Delta}u_{\rm in}\|_{L^\infty}\sim t^{-d/4}$ is not integrable near $t=0$ when $d\ge4$.
\item \emph{Critical data} such as $u_{\rm in}\in H^{1/2}(\TT^3)$; a Fujita--Kato type solution with bounded drag force should provide a
strong continuation.
\item \emph{The strong energy inequality.} Theorem~\ref{thm:SEI} provides Leray solutions satisfying it. Whether every energy
solution satisfies it, for instance those constructed in \cite{BDGM09}, is open, as for the Navier--Stokes equations.
\item \emph{Other critical spaces} $X^{p,q}$, $\frac2p+\frac dq=1$, $q<\infty$, in Theorem~\ref{thm:B}, and kinetic models with
diffusion in velocity (Vlasov--Fokker--Planck--Navier--Stokes).
\end{enumerate}

\subsection{Structure of the paper}

Section~\ref{sec:iteration} introduces weak solutions and the VNS--Reynolds system, states the main theorem in its precise
form, and reduces it to the Navier--Stokes step of Proposition~\ref{prop:NSstep}. Section~\ref{sec:kinetic} contains the
kinetic lemmas: characteristics, dual stability, the coupling error, the kinetic limit and the mean-momentum correction.
Section~\ref{sec:NSdiff} recalls the Navier--Stokes step of \cite{CL22} and lists the differences; its full proof is in
Appendices~\ref{sec:CI} and~\ref{sec:estimates}.
Section~\ref{sec:sharp} proves the uniqueness in the critical class. Section~\ref{sec:regint} proves the intervals of
regularity. Section~\ref{sec:statements} defines strong continuations and states the results for rough data;
Section~\ref{sec:rough} extends the construction to Lipschitz backgrounds and proves Theorem~\ref{thm:rough}, and
Section~\ref{sec:local} constructs strong continuations. Section~\ref{sec:structure} proves the structure theorem, and
Section~\ref{sec:HMMMtorus} contains the results on energy solutions used there: a velocity bound, weak--strong uniqueness, and the
existence of Leray solutions with the strong energy inequality.

\section{Weak solutions, main results and the iteration}\label{sec:iteration}

\subsection{Notation}\label{sec:notation}

Throughout, $d\ge2$, $\TT^d=\RR^d/\ZZ^d$ with $|\TT^d|=1$, and the time interval is $[0,1]$. We write
$z=(x,v)\in\TT^d\times\RR^d$ and $\fint g=\int_{\TT^d}g\dd x$. Norms such as $L^p_tL^q$ refer to
$L^p(0,1;L^q(\TT^d))$. For $a,b\in\RR^d$ let $a\ot b=a\otimes b-\frac1d(a\cdot b)\mathrm{Id}$, and let $\Sz$ be
the space of symmetric traceless $d\times d$ matrices.

For a smooth vector field $g$ on $\TT^d$ we use the projection onto mean-zero divergence-free fields,
$\PP g=g-\nabla\Delta^{-1}\D g-\fint g$, and the antidivergence $\Rs$ of \cite[Definition~B.2]{CL22}. We recall
that $\Rs g$ is $\Sz$-valued, that
\begin{equation}\label{eq:R-basic}
\D\Rs g=g-\fint g,\qquad \Rs(\text{const})=0,\qquad \|\Rs g\|_{L^q}\le C_{d,q}\|g\|_{L^q}\quad(1\le q\le\infty),
\end{equation}
see \cite[Theorem~B.3]{CL22}, and that $\Rs$ acts in space only, so it preserves smoothness in time and
vanishing on time intervals. Moreover $\Rs g=K*g$ for a matrix-valued kernel $K$ that is smooth on
$\TT^d\setminus\{0\}$ and satisfies
\begin{equation}\label{eq:kernel}
|K(z)|\le C_d\,\dist(z,\ZZ^d)^{1-d},\qquad |\nabla K(z)|\le C_d\,\dist(z,\ZZ^d)^{-d}.
\end{equation}
Indeed, the entries of $K$ are linear combinations of first derivatives of the periodic Green function $G$
of $\Delta$ and third derivatives of $G_2=G*G$; near the origin $G$ and $G_2$ differ from the Euclidean
fundamental solutions of $\Delta$ and $\Delta^2$ by smooth functions, and the logarithms occurring for
$d\in\{2,4\}$ disappear after the indicated number of derivatives.

\paragraph{Kinetic notation.}
The kinetic datum $f_0\in C^\infty_c(\TT^d\times\RR^d)$ is fixed and nonnegative. We set
\[
\Supp f_0\subset\TT^d\times B_{\Rv_0},\qquad M_0=\|f_0\|_{L^1},\qquad M_\infty=\|f_0\|_{L^\infty},\qquad
J_0=\iint vf_0\dd z .
\]
Radii in velocity are denoted by $\Rv$, with subscripts ($\Rv_0,\Rv_1,\Rv_\infty,\dots$), to distinguish them from Reynolds stresses,
which are always denoted by $R$.
Let $\Vc=C([0,1];C^1(\TT^d;\RR^d))$. For $U\in\Vc$ we denote by $Z_U(t;s,z)$ the flow of
$\dot X=V$, $\dot V=U(t,X)-V$, and by $f[U]$ the solution of
\begin{equation}\label{eq:kinetic}
\partial_tf+v\cdot\nabla_xf+\D_v\big[(U-v)f\big]=0,\qquad f(0)=f_0 .
\end{equation}
We write
\[
\rho[U]=\int f[U]\dd v,\qquad j[U]=\int vf[U]\dd v,\qquad J[U](t)=\iint vf[U](t)\dd z,
\]
and, for a distribution function $f$ and a velocity field $U$, the drag force
\[
F[f,U](t,x)=\int_{\RR^d}\big(v-U(t,x)\big)f(t,x,v)\dd v .
\]
Basic properties of $f[U]$ are collected in Lemma~\ref{lem:char}.

\subsection{Weak solutions}

The Vlasov--Navier--Stokes (VNS) system reads
\begin{equation}\label{eq:VNS}
\left\{
\begin{aligned}
&\partial_tf+v\cdot\nabla_xf+\D_v\big[(u-v)f\big]=0,\\
&\partial_tu-\Delta u+\D(u\otimes u)+\nabla p=F[f,u],\qquad \D u=0,
\end{aligned}\right.
\qquad (f,u)(0)=(f_0,u_{\rm in}).
\end{equation}

\begin{defin}[Weak solutions]\label{def:weak}
Let $u_{\rm in}\in L^2(\TT^d)$ be weakly divergence-free. A pair $(u,f)$ is a weak solution of \eqref{eq:VNS}
on $[0,1]$ with data $(u_{\rm in},f_0)$ if:
\begin{enumerate}[label=(\roman*)]
\item $u\in L^2((0,1)\times\TT^d)\cap L^1(0,1;L^\infty(\TT^d))$ and $u(t)$ is weakly divergence-free for a.e.\ $t$;
\item $f\in L^\infty((0,1)\times\TT^d\times\RR^d)$, $f\ge0$, and $f=0$ for $|v|\ge \Rv$, for some $\Rv<\infty$;
\item for every $\varphi\in C^\infty_c([0,1)\times\TT^d\times\RR^d)$,
\[
\iint f_0\,\varphi(0)\dd z+\int_0^1\!\!\iint f\big(\partial_t\varphi+v\cdot\nabla_x\varphi+(u-v)\cdot\nabla_v\varphi\big)\dd z\dd t=0;
\]
\item for every divergence-free $\varphi\in C^\infty([0,1]\times\TT^d;\RR^d)$ vanishing near $t=1$,
\[
\int u_{\rm in}\cdot\varphi(0)\dd x+\int_0^1\!\!\int\Big(u\cdot(\partial_t\varphi+\Delta\varphi)+u\otimes u:\nabla\varphi
+F[f,u]\cdot\varphi\Big)\dd x\dd t=0 .
\]
\end{enumerate}
\end{defin}

By (i)--(ii), $F[f,u]=j_f-\rho_fu\in L^1((0,1)\times\TT^d)$, so (iv) makes sense. Constant test fields are
admissible in (iv), so the definition contains the equation $\frac{\mathrm d}{\mathrm dt}\fint u=\fint F[f,u]$ in
$\mathcal D'(0,1)$. The spatial mean of $u$ is not conserved; the conserved quantity is the total momentum
$\fint u+\iint vf\dd z$.

\begin{defin}[Leray--Hopf solutions]\label{def:LH}
A weak solution is of Leray--Hopf type if $u\in L^\infty(0,1;L^2)\cap L^2(0,1;H^1)$ and, for a.e.\ $t\in[0,1]$,
\[
\mathcal E(t)+\int_0^t\Big(\|\nabla u\|_{L^2}^2+\iint|v-u|^2f\dd z\Big)\dd s\le\mathcal E(0),\qquad
\mathcal E(t)=\tfrac12\|u(t)\|_{L^2}^2+\tfrac12\iint|v|^2f(t)\dd z .
\]
\end{defin}

\subsection{Main results}

\begin{theo}[Main theorem]\label{thm:main}
Let $d\ge2$, $1\le p<2$, $1\le q<\infty$ and $\varepsilon\in(0,1)$. Let $t_0\in(0,1]$ and let $(\bar u,\bar f)$ be a
smooth solution of \eqref{eq:VNS} on $[0,t_0]$ with data $(u_{\rm in},f_0)$. Let
$\tilde u\in C^\infty([0,1]\times\TT^d;\RR^d)$ be divergence-free with $\tilde u=\bar u$ on $[0,t_0]$. Then there is a
weak solution $(u,f)$ of \eqref{eq:VNS} on $[0,1]$ with data $(u_{\rm in},f_0)$ such that:
\begin{enumerate}[label=(\arabic*)]
\item $u\in L^2_{t,x}\cap L^p_tL^\infty\cap L^1_tW^{1,q}$;
\item $(u,f)=(\bar u,\bar f)$ on $[0,t_0/2]$;
\item $\big\|(u-\fint u)-(\tilde u-\fint\tilde u)\big\|_{L^p_tL^\infty\cap L^1_tW^{1,q}}\le\varepsilon$;
\item $0\le f\le e^dM_\infty$, $f$ vanishes for $|v|\ge \Rv_\infty$ for some $\Rv_\infty<\infty$,
$f\in C([0,1];L^s(\TT^d\times\RR^d))$ for every $s<\infty$, and $f$ is the unique bounded weak solution of the
kinetic equation in Definition~\ref{def:weak}(iii) with velocity field $u$; in particular $f$ is a renormalized
solution;
\item the total momentum is conserved: $\fint u(t)+\iint vf(t)\dd z=\fint u_{\rm in}+J_0$ for a.e.\ $t$.
\end{enumerate}
\end{theo}

The norm in (3) is taken modulo spatial means because the fluid mean is slaved to the particle momentum by (5);
it is corrected at each step of the construction (Lemma~\ref{lem:mean}).

\begin{cor}[Nonuniqueness of the Cauchy problem]\label{cor:nonunique}
Let $d\ge2$, $1\le p<2$, $1\le q<\infty$, and let $(u_{\rm in},f_0)$ be smooth data for which \eqref{eq:VNS} has a
smooth solution on $[0,t_1]$ for some $t_1\in(0,1]$. Then there are infinitely many distinct weak solutions of
\eqref{eq:VNS} on $[0,1]$ with data $(u_{\rm in},f_0)$, all in $L^2_{t,x}\cap L^p_tL^\infty\cap L^1_tW^{1,q}$, all coinciding
with the smooth solution on $[0,t_1/4]$. Infinitely many of them are not of Leray--Hopf type.
\end{cor}

\begin{remark}\label{rem:local}
For smooth data with $f_0$ compactly supported, smooth solutions of \eqref{eq:VNS} exist on a short time
interval; this is classical and we take it as given. The hypothesis of Corollary~\ref{cor:nonunique} is thus
satisfied by all such data.
\end{remark}

\subsection{The VNS--Reynolds system}

\begin{defin}\label{def:VNSR}
A smooth solution of the VNS--Reynolds system is a triple $(u,p,R)$ of smooth functions on $[0,1]\times\TT^d$,
with $u$ vector-valued, $\D u=0$, $p$ scalar and $R$ taking values in $\Sz$, such that
\begin{equation}\label{eq:VNSR}
\partial_tu-\Delta u+\D(u\otimes u)+\nabla p=F\big[f[u],u\big]+\D R .
\end{equation}
\end{defin}

The kinetic equation is solved exactly, through $f[u]$, and the only defect is the stress $R$ in the fluid
equation. Since $\D R$, $\nabla p$ and $\D(u\otimes u)$ have zero mean, \eqref{eq:VNSR} contains the exact mean
equation $\frac{\mathrm d}{\mathrm dt}\fint u=\fint F[f[u],u]$. By Lemma~\ref{lem:char}(d) this is equivalent to
\begin{equation}\label{eq:momentum}
\fint u(t)+J[u](t)=\fint u(0)+J_0\qquad(0\le t\le1).
\end{equation}

\begin{lem}[Initial approximate solution]\label{lem:start}
Under the hypotheses of Theorem~\ref{thm:main}, there is a function $m\in C^\infty([0,1];\RR^d)$ with $m=0$ on
$[0,t_0]$ such that $u_0:=\tilde u+m$, together with suitable $p_0,R_0$, is a smooth solution of the VNS--Reynolds system,
and $R_0=0$ on $[0,t_0]$.
\end{lem}

\begin{proof}
Apply Lemma~\ref{lem:mean}(a) to $U=\tilde u$, and let $m$ be the resulting correction. On $[0,t_0]$ the pair
$(\tilde u,f[\tilde u])=(\bar u,\bar f)$ solves \eqref{eq:VNS}, so \eqref{eq:momentum} holds there and
Lemma~\ref{lem:mean}(b) gives $m=0$ on $[0,t_0]$. Set $u_0=\tilde u+m$, which is smooth and divergence-free, and
\[
E:=\partial_tu_0-\Delta u_0+\D(u_0\otimes u_0)-F[f[u_0],u_0].
\]
By Lemma~\ref{lem:mean}(a), $\fint E=\frac{\mathrm d}{\mathrm dt}\fint u_0-\fint F[f[u_0],u_0]=0$. Put
$p_0:=-\Delta^{-1}\D E$ and $R_0:=\Rs\PP E$. Since $\fint E=0$ and $\fint\PP E=0$,
\[
E+\nabla p_0=E-\nabla\Delta^{-1}\D E=\PP E=\D R_0 ,
\]
which is \eqref{eq:VNSR}. On $[0,t_0]$ we have $u_0=\bar u$ and $f[u_0]=\bar f$ by causality
(Lemma~\ref{lem:char}(b)), hence $E=-\nabla\bar p$ is a gradient, $\PP E=0$ and $R_0=0$.
\end{proof}

\subsection{The iteration}

The convex integration step is a Navier--Stokes statement:

\begin{prop}[Navier--Stokes step, after \cite{CL22}]\label{prop:NSstep}
Let $d\ge2$, $1\le p<2$, $1\le q<\infty$. There are $M\ge1$ depending only on $d$, and $r\in(1,\frac d{d-1})$
depending only on $d,p,q$, with the following property. Let $u\in C^\infty([0,1]\times\TT^d;\RR^d)$ be divergence-free,
$R\in C^\infty([0,1]\times\TT^d;\Sz)$, and $0<s<t_*\le1$ with $R=0$ on $[0,t_*]$. For every $\delta>0$ there are smooth
$w$, $P$, $R^{\rm CL}$, with $\D w=0$ and $R^{\rm CL}$ taking values in $\Sz$, such that
\begin{enumerate}[label=(\roman*)]
\item $\partial_tw-\Delta w+\D(u\otimes w+w\otimes u+w\otimes w)+\D R+\nabla P=\D R^{\rm CL}$;
\item $\fint w(t)=0$ for all $t$, and $w=0$, $R^{\rm CL}=0$ on $[0,s]$;
\item $\|w\|_{L^2_{t,x}}\le M\|R\|_{L^1_{t,x}}^{1/2}$;
\item $\|w\|_{L^p_tL^\infty}+\|w\|_{L^1_tW^{1,q}}\le\delta$;
\item $\|R^{\rm CL}\|_{L^1(0,1;L^r)}\le\delta$.
\end{enumerate}
\end{prop}

\begin{remark}[On the proof of Proposition~\ref{prop:NSstep}]\label{rem:NSstep}
Proposition~\ref{prop:NSstep} is the convex integration step of \cite[Sections~4--5]{CL22}, without the concentration step and with
minor adjustments. Section~\ref{sec:NSdiff} recalls the construction and lists the differences; the full proof is given in
Appendices~\ref{sec:CI} and~\ref{sec:estimates}.
\end{remark}

\begin{prop}[VNS iteration step]\label{prop:VNSstep}
Let $d,p,q$ and $M,r$ be as in Proposition~\ref{prop:NSstep}. Let $(u,p_u,R)$ be a smooth solution of the
VNS--Reynolds system and $0<s<t_*\le1$ with $R=0$ on $[0,t_*]$. For every $\delta\in(0,1]$ there is a smooth solution
$(u_1,p_1,R_1)$ of the VNS--Reynolds system such that
\begin{enumerate}[label=(\roman*)]
\item $u_1=u$ and $R_1=0$ on $[0,s]$;
\item $\|u_1-u\|_{L^2_{t,x}}\le M\|R\|_{L^1_{t,x}}^{1/2}+\delta$;
\item $\|u_1-u\|_{L^p_tL^\infty}+\|u_1-u\|_{L^1_tW^{1,q}}\le\delta$;
\item $\|R_1\|_{L^1(0,1;L^r)}\le\delta$.
\end{enumerate}
\end{prop}

Proposition~\ref{prop:VNSstep} is proved in Section~\ref{sec:proofstep}, assuming Proposition~\ref{prop:NSstep}.

\subsection{Proof of Theorem~\ref{thm:main}}

\proofstep{Step 1: construction}
Let $(u_0,p_0,R_0)$ be given by Lemma~\ref{lem:start}, and set $s_n=t_0\big(\frac12+2^{-n-1}\big)$ for $n\ge0$, so that
$s_0=t_0$ and $s_n\downarrow t_0/2$. Suppose $(u_n,p_n,R_n)$ is a smooth solution of the VNS--Reynolds system with
$R_n=0$ on $[0,s_n]$. Put $\Lambda_n=\|\nabla u_n\|_{L^1_tL^\infty}<\infty$ and
\[
\delta_n=\min\Big\{\varepsilon4^{-n-1},\ 4^{-n},\ \frac{2^{-n}e^{-1-\Lambda_n}}{1+M_0}\Big\}.
\]
Proposition~\ref{prop:VNSstep}, applied with $(t_*,s)=(s_n,s_{n+1})$ and $\delta=\delta_n$, gives
$(u_{n+1},p_{n+1},R_{n+1})$ with $R_{n+1}=0$ on $[0,s_{n+1}]$. By induction, for all $n\ge0$:
\begin{enumerate}[label=(P\arabic*)]
\item $u_{n+1}=u_n$ on $[0,s_{n+1}]$, hence $u_n=u_0=\bar u$ on $[0,t_0/2]$;
\item $\|R_{n+1}\|_{L^1_{t,x}}\le\|R_{n+1}\|_{L^1_tL^r}\le\delta_n$;
\item $\|u_{n+1}-u_n\|_{L^p_tL^\infty}+\|u_{n+1}-u_n\|_{L^1_tW^{1,q}}\le\delta_n$;
\item $\|u_{n+1}-u_n\|_{L^2_{t,x}}\le M\|R_n\|_{L^1_{t,x}}^{1/2}+\delta_n\le M\delta_{n-1}^{1/2}+\delta_n\le M2^{1-n}+4^{-n}$ for $n\ge1$.
\end{enumerate}

\proofstep{Step 2: the velocity}
By (P3)--(P4), $(u_n)$ is Cauchy in $X:=L^2_{t,x}\cap L^p_tL^\infty\cap L^1_tW^{1,q}$; let $u$ be its limit. Then
$u(t)$ is weakly divergence-free for a.e.\ $t$, and (1) holds. By (P3),
$\|u-u_0\|_{L^p_tL^\infty\cap L^1_tW^{1,q}}\le\sum_n\delta_n\le\varepsilon/3$. Since $u_0=\tilde u+m$ with $m$ constant in
space, $(u-\fint u)-(\tilde u-\fint\tilde u)=(u-u_0)-\fint(u-u_0)$. Subtracting the mean at most doubles the norms of
$L^\infty$ and $W^{1,q}$, which gives (3). Since $p\ge1$,
\begin{equation}\label{eq:Kbound}
\sup_n\|u_n\|_{L^1_tL^\infty}\le\|u_0\|_{L^1_tL^\infty}+1=:K,\qquad
\sum_ne^{1+\Lambda_n}\|u_{n+1}-u_n\|_{L^1_tL^\infty}\le\sum_n2^{-n}<\infty .
\end{equation}

\proofstep{Step 3: the distribution function}
By \eqref{eq:Kbound} and $u_n\to u$ in $L^1_tL^\infty$, Lemma~\ref{lem:convergence} applies to $f_n:=f[u_n]$. It yields a
limit $f$ with $0\le f\le e^dM_\infty$, $f=0$ for $|v|\ge \Rv_\infty:=\Rv_0+K$, which solves the weak kinetic equation of
Definition~\ref{def:weak}(iii) with field $u$; moreover $F[f_n,u_n]\to F[f,u]$ in the sense of
Lemma~\ref{lem:convergence}(d), and $J[u_n](t)\to\iint vf(t)\dd z$ for every $t$.

\proofstep{Step 4: the fluid equation}
Let $\varphi$ be an admissible test field in Definition~\ref{def:weak}(iv). Testing \eqref{eq:VNSR} for $u_n$ with
$\varphi$ and integrating by parts, using $u_n(0)=u_{\rm in}$ by (P1),
\[
\int u_{\rm in}\cdot\varphi(0)+\int_0^1\!\!\int\Big(u_n\cdot(\partial_t\varphi+\Delta\varphi)+u_n\otimes u_n:\nabla\varphi
+F[f_n,u_n]\cdot\varphi-R_n:\nabla\varphi\Big)=0 .
\]
As $n\to\infty$: $u_n\to u$ in $L^2_{t,x}$ gives convergence of the first two terms; $R_n\to0$ in $L^1_{t,x}$ by (P2);
and $\iint F[f_n,u_n]\cdot\varphi\to\iint F[f,u]\cdot\varphi$ by Lemma~\ref{lem:convergence}(d). Hence (iv) holds, and
$(u,f)$ is a weak solution.

\proofstep{Step 5: properties (2), (4), (5)}
By (P1), $u=\bar u$ on $[0,t_0/2]$, and $f_n=f[\bar u]=\bar f$ on $[0,t_0/2]$ by causality, hence $f=\bar f$ there;
this is (2). Since $u\in L^1_tW^{1,q}\cap L^1_tL^\infty$, Lemma~\ref{lem:DPL} gives the uniqueness, renormalization and
continuity statements in (4). Finally, each $u_n$ satisfies \eqref{eq:momentum}, which reads
$\fint u_n(t)+J[u_n](t)=\fint u_{\rm in}+J_0$. Along a subsequence $\fint u_n(t)\to\fint u(t)$ for a.e.\ $t$, while
$J[u_n](t)\to\iint vf(t)\dd z$ for every $t$; this gives (5).
\qed

\begin{remark}[Exponents uniform in $p$ and $q$]\label{rem:diag}
Theorem~\ref{thm:main} is stated for fixed $p<2$ and $q<\infty$. A single solution lying in $L^p_tL^\infty\cap L^1_tW^{1,q}$ for all $p<2$ and
$q<\infty$ simultaneously can be obtained by a diagonal argument: choose $p_n\uparrow2$ and $q_n\uparrow\infty$, and apply
Proposition~\ref{prop:VNSstep} at step $n$ with the exponents $(p_n,q_n)$, the corresponding $r_n$, and $\delta_n$ as above. Since every
constant at step $n$ is chosen after $u_n$ is fixed, the construction is unchanged. For fixed $m$ and $n\ge m$, the bound
$\|u_{n+1}-u_n\|_{L^{p_n}L^\infty\cap L^1W^{1,q_n}}\le\delta_n$ implies the same bound in $L^{p_m}L^\infty\cap L^1W^{1,q_m}$, since $|\TT^d|=1$ and the time
interval has length one. Hence $(u_n)$ is Cauchy in every $L^{p_m}L^\infty\cap L^1W^{1,q_m}$, and the limit lies in all of them. The
convergence of the stresses, $\|R_n\|_{L^1_{t,x}}\le\|R_n\|_{L^1L^{r_n}}\le\delta_{n-1}$, is unaffected. We do not use this refinement
elsewhere, and all statements below are to be read with fixed exponents.
\end{remark}

\subsection{Proof of Corollary~\ref{cor:nonunique}}

Let $(\bar u,\bar f)$ be the smooth solution on $[0,t_1]$. Choose $\chi\in C^\infty(\RR;[0,1])$ with $\chi=1$ on
$(-\infty,t_1/2]$ and $\chi=0$ on $[3t_1/4,\infty)$, and choose a smooth divergence-free $\Psi$ on $[0,1]\times\TT^d$
with $\fint\Psi(t)=0$, $\Psi=0$ for $t\le t_1$, and $\|\Psi\|_{L^p(0,1;L^2)}\ge1$. For instance,
$\Psi=c\,\eta(t)(\sin2\pi x_2,0,\dots,0)$ with $\eta$ supported in $(t_1,1]$, if $t_1<1$; if $t_1=1$, replace $t_1$ by
$1/2$ throughout. For $k\in\NN$ put
\[
\tilde u_k=\chi\bar u+k\Psi ,
\]
which is smooth, divergence-free and equal to $\bar u$ on $[0,t_1/2]$. Theorem~\ref{thm:main} with $t_0=t_1/2$ and
$\varepsilon=\frac14$ gives weak solutions $(u_k,f_k)$ with data $(u_{\rm in},f_0)$, coinciding with $(\bar u,\bar f)$ on
$[0,t_1/4]$.

\proofstep{Distinctness}
For $k\ne l$, $(\tilde u_k-\fint\tilde u_k)-(\tilde u_l-\fint\tilde u_l)=(k-l)\Psi$, so by Theorem~\ref{thm:main}(3) and
$\|\Psi\|_{L^p_tL^\infty}\ge\|\Psi\|_{L^p_tL^2}\ge1$,
\[
\big\|(u_k-\fint u_k)-(u_l-\fint u_l)\big\|_{L^p_tL^\infty}\ge|k-l|-2\varepsilon\ge\tfrac12 .
\]

\proofstep{Failure of the energy inequality}
If $(u_k,f_k)$ were of Leray--Hopf type, then $\|u_k(t)\|_{L^2}^2\le2\mathcal E(0)$ for a.e.\ $t$, since $f_k\ge0$, and hence
$\|u_k\|_{L^p_tL^2}\le(2\mathcal E(0))^{1/2}$. On the other hand, by orthogonality of the mean and by
Theorem~\ref{thm:main}(3),
\[
\|u_k\|_{L^p_tL^2}\ge\big\|u_k-\fint u_k\big\|_{L^p_tL^2}
\ge k\|\Psi\|_{L^p_tL^2}-\big\|\chi\bar u-\fint\chi\bar u\big\|_{L^p_tL^2}-2\varepsilon\xrightarrow[k\to\infty]{}\infty .
\]
Hence $(u_k,f_k)$ is not of Leray--Hopf type for all large $k$.
\qed

\section{Kinetic lemmas}\label{sec:kinetic}

This section collects the kinetic ingredients and proves Proposition~\ref{prop:VNSstep}. The central estimate is
the dual stability Lemma~\ref{lem:dual}: the distance between $f[u]$ and $f[u+\tilde w]$ is controlled by
$\|\tilde w\|_{L^1_tL^\infty}$ times a constant that depends only on the Lipschitz norm of the old field $u$. In the
convex integration step $\nabla(u+\tilde w)$ is of size $\sigma\mu$ and cannot be allowed to enter. The constants
depending on $u$ are compensated by choosing $\delta$ after $u$.

For a function $\psi\in C^1(\TT^d\times\RR^d)$ and $R>0$ we write
\[
[\psi]_R=\sup_{\TT^d\times B_R}\max\{|\nabla_x\psi|,|\nabla_v\psi|\},\qquad
\|\psi\|_{C^1_b}=\|\psi\|_{L^\infty}+\|\nabla_x\psi\|_{L^\infty}+\|\nabla_v\psi\|_{L^\infty},
\]
and for a vector field $h$ on $\TT^d$,
\[
\|h\|_{(C^1)^*}=\sup\Big\{\Big|\int h\cdot\phi\dd x\Big|:\ \phi\in C^1(\TT^d;\RR^d),\ \|\phi\|_{L^\infty}+\|\nabla\phi\|_{L^\infty}\le1\Big\}.
\]

\subsection{Characteristics}

\begin{lem}\label{lem:char}
Let $U\in\Vc$ and $f=f[U]$.
\begin{enumerate}[label=(\alph*)]
\item The flow $Z_U$ is $C^1$ in all variables, $\det D_zZ_U(t;s,z)=e^{-d(t-s)}$, and
\begin{equation}\label{eq:pushforward}
f(t,Z_U(t;0,z))=e^{dt}f_0(z),\qquad \iint f(t)\psi\dd z=\iint f_0(z)\psi\big(Z_U(t;0,z)\big)\dd z .
\end{equation}
In particular $f\in C^1$ solves \eqref{eq:kinetic} classically, $0\le f\le e^dM_\infty$, and $\|f(t)\|_{L^1}=M_0$. If
$U$ is smooth on $[0,1]\times\TT^d$, so is $f$.
\item (Causality.) If $U_1=U_2$ on $[0,t_0]$, then $f[U_1]=f[U_2]$ on $[0,t_0]$.
\item (Velocity support.) For $t\ge s$, $|V(t;s,z)|\le|v|+\int_s^t\|U(\tau)\|_{L^\infty}\dd\tau$. Consequently
\[
\Supp f(t)\subset\TT^d\times B_{\Rv_U(t)},\quad \Rv_U(t):=\Rv_0+\int_0^t\|U\|_{L^\infty}\dd\tau,\qquad
\|\rho[U](t)\|_{L^\infty}\le e^dM_\infty|B_1|\Rv_U(t)^d .
\]
\item (Momentum.) $J[U]\in C^1([0,1];\RR^d)$, $J[U](0)=J_0$, and
\begin{equation}\label{eq:Jprime}
\frac{\mathrm d}{\mathrm dt}J[U]=\iint(U-v)f\dd z=-\fint F[f,U].
\end{equation}
\end{enumerate}
\end{lem}

\begin{proof}
(a) The field $b_U(t,z)=(v,U(t,x)-v)$ is continuous in $t$ and $C^1$ in $z$, with globally bounded $z$-derivative, so
the flow is global and $C^1$ in $(t,s,z)$. Liouville's formula and $\D_zb_U=-d$ give the Jacobian, and
\eqref{eq:pushforward} is the push-forward of $f_0$ written in Eulerian and Lagrangian form. The remaining statements
follow from \eqref{eq:pushforward} and the smooth dependence of the flow on parameters.

(b) The flows coincide on $[0,t_0]$.

(c) Where $V\ne0$, $\frac{\mathrm d}{\mathrm d\tau}|V|=\frac V{|V|}\cdot(U(\tau,X)-V)\le\|U(\tau)\|_{L^\infty}$; at zeros of $V$
use $(|V|^2+\eta^2)^{1/2}$ and let $\eta\to0$. The support bound follows from \eqref{eq:pushforward}, and the density
bound from $f\le e^dM_\infty$.

(d) By (c), $f(t)$ has compact support in $v$, uniformly in $t$. Multiplying \eqref{eq:kinetic} by $v$ and integrating by
parts, the term $v\cdot\nabla_xf$ integrates to zero on $\TT^d$, and
$-\iint v\,\D_v[(U-v)f]=\iint(U-v)f$. The last identity in \eqref{eq:Jprime} is the definition of $F$.
\end{proof}

\subsection{Dual stability}

\begin{lem}[Dual stability]\label{lem:dual}
Let $u,\tilde w\in\Vc$, $f=f[u]$, $f^+=f[u+\tilde w]$, $\rho^+=\int f^+\dd v$, and
\[
\Lambda:=\int_0^1\|\nabla u(\tau)\|_{L^\infty}\dd\tau,\qquad
\Rv_1:=\Rv_0+\int_0^1\big(2\|u\|_{L^\infty}+\|\tilde w\|_{L^\infty}\big)\dd\tau .
\]
Then for every $t\in[0,1]$ and $\psi\in C^1(\TT^d\times\RR^d)$,
\begin{equation}\label{eq:dual}
\Big|\iint\big(f^+-f\big)(t)\,\psi\dd z\Big|\le e^{1+\Lambda}[\psi]_{\Rv_1}\int_0^t\!\!\int_{\TT^d}|\tilde w|\,\rho^+\dd x\dd s
\le e^{1+\Lambda}[\psi]_{\Rv_1}M_0\|\tilde w\|_{L^1(0,t;L^\infty)} .
\end{equation}
\end{lem}

\begin{proof}
\proofstep{Step 1: equation for the difference}
Write $b=b_u$. Since $b_{u+\tilde w}=b+(0,\tilde w)$, the function $g=f^+-f\in C^1$ satisfies
\begin{equation}\label{eq:gdiff}
\partial_tg+\D_z(bg)=-\D_v(\tilde wf^+),\qquad g(0)=0,
\end{equation}
and, by Lemma~\ref{lem:char}(c), $f(s)$ and $f^+(s)$ are supported in $\TT^d\times B_{\Rv_1}$ for all $s$.

\proofstep{Step 2: duality}
Fix $t$, $\psi$ and set $\Phi(s,z)=\psi(Z_u(t;s,z))$ for $s\in[0,t]$. Then $\Phi\in C^1$ and
$\partial_s\Phi+b\cdot\nabla_z\Phi=0$, $\Phi(t)=\psi$. Since $g$ and $f^+$ have compact support in $v$, integration by
parts in $z$ gives
\[
\frac{\mathrm d}{\mathrm ds}\iint g\Phi\dd z=\iint\big(bg\cdot\nabla_z\Phi+f^+\tilde w\cdot\nabla_v\Phi\big)\dd z-\iint g\,b\cdot\nabla_z\Phi\dd z
=\iint f^+\tilde w\cdot\nabla_v\Phi\dd z,
\]
and integrating over $[0,t]$,
\begin{equation}\label{eq:duality}
\iint g(t)\psi\dd z=\int_0^t\!\!\iint f^+(s,z)\,\tilde w(s,x)\cdot\nabla_v\Phi(s,z)\dd z\dd s .
\end{equation}

\proofstep{Step 3: the linearized flow}
For $s\le\tau\le t$ and $\xi=(\xi_x,\xi_v)$, let $(\delta X,\delta V)(\tau)=D_zZ_u(\tau;s,z)\xi$. Then
$\delta\dot X=\delta V$ and $\delta\dot V=\nabla u(\tau,X)\delta X-\delta V$. With $a=|\delta X|$, $b=|\delta V|$ we have, for a.e.\ $\tau$,
$a'\le b$ and $b'\le\|\nabla u(\tau)\|_{L^\infty}a$, hence $(a+b)'\le(1+\|\nabla u(\tau)\|_{L^\infty})(a+b)$ and
\begin{equation}\label{eq:linflow}
|\delta X(t)|+|\delta V(t)|\le e^{1+\Lambda}\big(|\xi_x|+|\xi_v|\big).
\end{equation}

\proofstep{Step 4: conclusion}
For $\eta\in\RR^d$, $\nabla_v\Phi(s,z)\cdot\eta=\nabla_z\psi(Z_u(t;s,z))\cdot D_zZ_u(t;s,z)(0,\eta)$. If $z\in\Supp f^+(s)$, then
$|v|\le \Rv_0+\int_0^s\|u+\tilde w\|_{L^\infty}$ by Lemma~\ref{lem:char}(c) for the field $u+\tilde w$, and the same lemma for
the field $u$ on $[s,t]$ shows $Z_u(t;s,z)\in\TT^d\times B_{\Rv_1}$. Therefore \eqref{eq:linflow} gives
$|\nabla_v\Phi(s,z)|\le e^{1+\Lambda}[\psi]_{\Rv_1}$ on $\Supp f^+(s)$, and \eqref{eq:duality} yields the first inequality in
\eqref{eq:dual}. The second follows from $\int\rho^+(s)\dd x=M_0$.
\end{proof}

\subsection{The coupling error}

\begin{lem}[Antidivergence in a weak norm]\label{lem:weakR}
Let $1<r<\frac d{d-1}$, $a:=\frac dr-(d-1)\in(0,1)$ and $\theta:=\frac a{a+d+1}$. There is $C=C(d,r)$ such that for every
smooth vector field $h$ on $\TT^d$,
\begin{equation}\label{eq:weakR}
\|\Rs h\|_{L^r}\le C\Big(\|h\|_{L^1}^{1-\theta}\|h\|_{(C^1)^*}^\theta+\|h\|_{(C^1)^*}\Big).
\end{equation}
\end{lem}

\begin{proof}
Let $\eta$ be a standard mollifier, $\eta_\varepsilon=\varepsilon^{-d}\eta(\cdot/\varepsilon)$, $\varepsilon\in(0,1]$, and write
$\Rs h=(K-K*\eta_\varepsilon)*h+(K*\eta_\varepsilon)*h$. Let $\delta(z)=\dist(z,\ZZ^d)$.

For the first term, by Young's inequality it suffices to bound $\|K-K*\eta_\varepsilon\|_{L^r}$. On $\{\delta\le2\varepsilon\}$ we have
$|K*\eta_\varepsilon|\le(|K|\mathbf 1_{\{\delta\le3\varepsilon\}})*\eta_\varepsilon$, so by \eqref{eq:kernel}
\[
\|K-K*\eta_\varepsilon\|_{L^r(\{\delta\le2\varepsilon\})}\le2\|K\|_{L^r(\{\delta\le3\varepsilon\})}
\le C\Big(\int_0^{3\varepsilon}\varrho^{(1-d)r+d-1}\dd\varrho\Big)^{1/r}=C\varepsilon^a ,
\]
using $r<\frac d{d-1}$. On $\{\delta\ge2\varepsilon\}$, all points of $B(z,\varepsilon)$ are at distance $\ge\delta(z)/2$ from $\ZZ^d$, so
$|K-K*\eta_\varepsilon|(z)\le\varepsilon\sup_{B(z,\varepsilon)}|\nabla K|\le C\varepsilon\delta(z)^{-d}$ and, using $r>1$,
\[
\|K-K*\eta_\varepsilon\|_{L^r(\{\delta\ge2\varepsilon\})}\le C\varepsilon\Big(\int_{2\varepsilon}^\infty\varrho^{-dr+d-1}\dd\varrho\Big)^{1/r}=C\varepsilon^a .
\]

For the second term, each entry of $(K*\eta_\varepsilon)*h$ at $x$ equals $\int h\cdot\phi^x\dd y$ with $\phi^x_k(y)=(K_{ijk}*\eta_\varepsilon)(x-y)$, and
$\|\phi^x\|_{L^\infty}+\|\nabla\phi^x\|_{L^\infty}\le\|K\|_{L^1}(\|\eta_\varepsilon\|_{L^\infty}+\|\nabla\eta_\varepsilon\|_{L^\infty})\le C\varepsilon^{-d-1}$.
Hence $\|(K*\eta_\varepsilon)*h\|_{L^r}\le\|(K*\eta_\varepsilon)*h\|_{L^\infty}\le C\varepsilon^{-d-1}\|h\|_{(C^1)^*}$.

Finally, with $A=\|h\|_{L^1}$ and $B=\|h\|_{(C^1)^*}$ we have $\|\Rs h\|_{L^r}\le C(\varepsilon^aA+\varepsilon^{-d-1}B)$ for all
$\varepsilon\in(0,1]$. If $B<A$, take $\varepsilon=(B/A)^{1/(a+d+1)}$; both terms equal $A^{1-\theta}B^\theta$. If $B\ge A$, take
$\varepsilon=1$.
\end{proof}

\begin{prop}[Coupling error]\label{prop:coupling}
Let $u,\tilde w\in\Vc$ with $W:=\|\tilde w\|_{L^1_tL^\infty}\le1$, and let $1<r<\frac d{d-1}$, $\theta$ as in
Lemma~\ref{lem:weakR}. Put $u^+=u+\tilde w$ and
\[
h:=F\big[f[u^+],u^+\big]-F\big[f[u],u\big].
\]
If $h$ is smooth (e.g.\ if $u,\tilde w$ are smooth), then
\begin{equation}\label{eq:coupling}
\|\Rs h\|_{L^1(0,1;L^r)}\le C_\star\big(W^\theta+W\big),
\end{equation}
where $C_\star$ depends only on $d,r,M_0,M_\infty,\Rv_0$ and $L:=\|u\|_{L^\infty(0,1;W^{1,\infty})}$.
\end{prop}

\begin{proof}
Write $f=f[u]$, $f^+=f[u^+]$. Since $F[f^+,u^+]=\int(v-u)f^+\dd v-\rho^+\tilde w$, we have $h=h_1-\rho^+\tilde w$ with
$h_1=\int(v-u)(f^+-f)\dd v$. As $W\le1$, the radius in Lemma~\ref{lem:dual} satisfies $\Rv_1\le \Rv_0+2L+1$, and $\Lambda\le L$.

To bound $h_1$, fix $t$ and $\phi\in C^1(\TT^d;\RR^d)$ with $\|\phi\|_{L^\infty}+\|\nabla\phi\|_{L^\infty}\le1$, and take
$\psi(x,v)=(v-u(t,x))\cdot\phi(x)$. Then $\nabla_v\psi=\phi$ and $\nabla_x\psi=(\nabla\phi)^T(v-u)-(\nabla u)^T\phi$, so
$[\psi]_{\Rv_1}\le1+\Rv_1+2L$, and Lemma~\ref{lem:dual} gives $\|h_1(t)\|_{(C^1)^*}\le B_0W$ with $B_0=e^{1+L}(1+\Rv_1+2L)M_0$.
By Lemma~\ref{lem:char}(c), $\|h_1(t)\|_{L^1}\le(\Rv_1+L)\iint(f^++f)\dd z=2(\Rv_1+L)M_0=:A_0$. Lemma~\ref{lem:weakR} gives
$\|\Rs h_1(t)\|_{L^r}\le C(A_0^{1-\theta}B_0^\theta W^\theta+B_0W)$ for every $t$.

For the second term, by \eqref{eq:R-basic} and Lemma~\ref{lem:char}(c),
\[
\|\Rs(\rho^+\tilde w)(t)\|_{L^r}\le C\|\rho^+(t)\|_{L^\infty}\|\tilde w(t)\|_{L^\infty}\le Ce^dM_\infty|B_1|\Rv_1^d\,\|\tilde w(t)\|_{L^\infty}.
\]

Integrating both bounds in $t$ proves \eqref{eq:coupling}.
\end{proof}

\subsection{The kinetic limit}

\begin{lem}[Convergence of distribution functions]\label{lem:convergence}
Let $u_n\in C^\infty([0,1]\times\TT^d;\RR^d)$, $n\ge0$, and $u\in L^1(0,1;L^\infty(\TT^d))$ satisfy
\[
K:=\sup_n\|u_n\|_{L^1_tL^\infty}<\infty,\qquad
\sum_ne^{1+\Lambda_n}\|u_{n+1}-u_n\|_{L^1_tL^\infty}<\infty,\qquad
u_n\to u\ \text{in }L^1_tL^\infty,
\]
where $\Lambda_n=\|\nabla u_n\|_{L^1_tL^\infty}$.
Let $f_n=f[u_n]$ and $\Rv_\infty=\Rv_0+K$. Then:
\begin{enumerate}[label=(\alph*)]
\item $0\le f_n\le e^dM_\infty$ and $f_n=0$ for $|v|\ge \Rv_\infty$;
\item there is $f\in L^\infty((0,1)\times\TT^d\times\RR^d)$, with a representative $t\mapsto f(t)$ defined for every $t$, such
that $\sup_t\sup_{\|\psi\|_{C^1_b}\le1}\big|\iint(f_n-f)(t)\psi\big|\to0$; $f_n(t)\rightharpoonup^*f(t)$ in $L^\infty$ for every $t$;
$f_n\rightharpoonup^*f$ in $L^\infty((0,1)\times\TT^d\times\RR^d)$; and $f$ satisfies the bounds of (a), with
$\iint f(t)\dd z=M_0$ and $f(0)=f_0$;
\item $f$ satisfies Definition~\ref{def:weak}(iii) with velocity field $u$;
\item for every $\Phi\in C([0,1]\times\TT^d;\RR^d)$, $\iint F[f_n,u_n]\cdot\Phi\dd x\dd t\to\iint F[f,u]\cdot\Phi\dd x\dd t$, and
$J[u_n](t)\to\iint vf(t)\dd z$ for every $t$.
\end{enumerate}
\end{lem}

\begin{proof}
(a) is Lemma~\ref{lem:char}(a),(c).

(b) Since $[\psi]_R\le\|\psi\|_{C^1_b}$, Lemma~\ref{lem:dual} with base field $u_k$ and perturbation $u_{k+1}-u_k$ gives, for
$m>n$ and all $t$,
\[
\sup_{\|\psi\|_{C^1_b}\le1}\Big|\iint(f_m-f_n)(t)\psi\Big|\le M_0\sum_{k\ge n}e^{1+\Lambda_k}\|u_{k+1}-u_k\|_{L^1_tL^\infty}\xrightarrow[n\to\infty]{}0 .
\]
Fix $t$. By (a), $(f_n(t))$ is bounded in $L^\infty$ and supported in a fixed compact set. Any two weak-$*$ cluster points
coincide on the dense subspace $C^1_c$ of $L^1$, by the Cauchy property, so $f_n(t)\rightharpoonup^*f(t)$ for a unique $f(t)$,
which inherits the bounds of (a). The convergence is uniform in $t$ in the dual $C^1_b$ norm. Testing with
$\psi=\chi(v)$, $\chi\in C^1_c$ equal to $1$ on $B_{\Rv_\infty}$, gives the mass identity; $f(0)=f_0$ because $f_n(0)=f_0$.
For $G\in L^1((0,1)\times\TT^d\times\RR^d)$, the functions $t\mapsto\iint f_n(t)G(t)$ are measurable, bounded by
$e^dM_\infty\|G(t)\|_{L^1}$ and converge pointwise to $\iint f(t)G(t)$; by dominated convergence,
$\iiint f_nG\to\int_0^1\iint f(t)G(t)\dd t$. Hence every weak-$*$ cluster point of $(f_n)$ in $L^\infty((0,1)\times\TT^d\times\RR^d)$
is represented by $t\mapsto f(t)$, and the whole sequence converges.

(c) Each $f_n$ is a classical, hence weak, solution with field $u_n$ and datum $f_0$. For
$\varphi\in C^\infty_c([0,1)\times\TT^d\times\RR^d)$, the terms $\iiint f_n(\partial_t\varphi+v\cdot\nabla_x\varphi-v\cdot\nabla_v\varphi)$
and $\iiint f_n\,u\cdot\nabla_v\varphi$ converge by (b), since the multipliers belong to $L^1$. The remaining term is bounded by
\[
\Big|\iiint f_n(u_n-u)\cdot\nabla_v\varphi\Big|\le e^dM_\infty\int_0^1\|u_n-u\|_{L^\infty}\|\nabla_v\varphi(t)\|_{L^1_z}\dd t\le C_\varphi\|u_n-u\|_{L^1_tL^\infty}\to0 .
\]

(d) Let $\chi$ be as in (b). Then $\iint F[f_n,u_n]\cdot\Phi=\iiint f_n\chi(v)v\cdot\Phi-\iiint f_n\chi(v)u_n\cdot\Phi$. The first term
and $\iiint f_n\chi\,u\cdot\Phi$ converge by (b), since $\chi v\cdot\Phi$ and $\chi u\cdot\Phi$ are integrable; the difference
$\iiint f_n\chi(u_n-u)\cdot\Phi$ is bounded by $e^dM_\infty|B_{\Rv_\infty}|\|\Phi\|_{L^\infty}\|u_n-u\|_{L^1_tL^\infty}$. The convergence of
$J[u_n](t)$ follows from (b) with $\psi=\chi(v)v_i$.
\end{proof}

\begin{lem}[Identification of the limit]\label{lem:DPL}
Let $u\in L^1(0,1;W^{1,1}(\TT^d))\cap L^1(0,1;L^\infty(\TT^d))$. Then the kinetic equation of
Definition~\ref{def:weak}(iii) with velocity field $u$ and datum $f_0$ has at most one weak solution in
$L^\infty((0,1)\times\TT^d\times\RR^d)$. The same holds for any bounded datum vanishing for large $|v|$, in place of $f_0$; the smoothness of $f_0$ is
not used, and this general form is the one applied in Sections~\ref{sec:sharp} and~\ref{sec:rough}--\ref{sec:HMMMtorus}. This solution, when it exists, is renormalized, i.e.
\[
\partial_t\beta(f)+\D_z\big(b\,\beta(f)\big)-d\big(f\beta'(f)-\beta(f)\big)=0\quad\text{in }\mathcal D'
\]
for every $\beta\in C^1(\RR)$, where $b=(v,u-v)$. If it vanishes for $|v|\ge \Rv$, it belongs to
$C([0,1];L^s(\TT^d\times\RR^d))$ for every $s<\infty$.
\end{lem}

\begin{proof}
The equation is the continuity equation $\partial_tf+\D_z(bf)=0$ on $\TT^d\times\RR^d$ with $b(t,x,v)=(v,u(t,x)-v)$. We verify
the hypotheses of the DiPerna--Lions theory \cite[Section~II]{DPL89} (see also \cite{Amb04}) for bounded solutions, in the form used
in \cite[Section~1.3]{HMMM25} for the same phase-space field:
$\nabla_xb=(0,\nabla u)$ and $\nabla_vb=(\mathrm{Id},-\mathrm{Id})$, so $b\in L^1(0,1;W^{1,1}_{\rm loc})$; $\D_zb=-d$ is bounded;
and $|b(t,z)|/(1+|z|)\le2+\|u(t)\|_{L^\infty}\in L^1(0,1)$, since $|b|\le2|v|+|u|$. The commutator lemma for $b\in L^1_tW^{1,1}_{\rm loc}$ and
$f\in L^\infty$ yields renormalization; uniqueness follows from \cite[Theorem~II.2]{DPL89} and strong continuity in time as in
\cite[Section~II]{DPL89}. The Lagrangian flow is given by \cite[Theorem~III.2]{DPL89}. The
arguments are local in $z$ and apply to bounded solutions, so periodicity in $x$ plays no role.
\end{proof}

\begin{remark}
Lemma~\ref{lem:DPL} shows that, in the solutions of Theorem~\ref{thm:main}, the distribution function is determined by the
velocity, and the nonuniqueness comes from the fluid. This is in contrast with nonuniqueness mechanisms driven by the
transported quantity, where the transported density fails to be renormalized.
\end{remark}

\subsection{Mean momentum correction}

The antidivergence $\Rs$ inverts the divergence only on mean-zero fields, while the drag force has a time-dependent mean.
The following lemma restores the exact mean equation by a spatially constant correction $c(t)$. Since $c$ changes the
kinetic solution, it is found by a fixed point.

\begin{lem}[Mean momentum correction]\label{lem:mean}
Let $U\in C^\infty([0,1]\times\TT^d;\RR^d)$ be divergence-free.
\begin{enumerate}[label=(\alph*)]
\item There is a unique $c\in C([0,1];\RR^d)$ such that
\begin{equation}\label{eq:fixed}
c(t)=\fint U(0)-\fint U(t)+J_0-J[U+c](t)\qquad(0\le t\le1).
\end{equation}
Moreover $c\in C^\infty([0,1];\RR^d)$, $c(0)=0$, and $U^\sharp:=U+c$ satisfies \eqref{eq:momentum} and
$\frac{\mathrm d}{\mathrm dt}\fint U^\sharp=\fint F[f[U^\sharp],U^\sharp]$.
\item If $U$ satisfies \eqref{eq:momentum} on $[0,t_0]$, then $c=0$ on $[0,t_0]$.
\item Suppose $U=u+w$, where $u$ is smooth, divergence-free and satisfies \eqref{eq:momentum} on $[0,1]$, and $w$ is smooth,
divergence-free with $\fint w(t)=0$ for all $t$. Then
\begin{equation}\label{eq:csize}
|c(t)|\le Ae^{At}\|w\|_{L^1(0,t;L^\infty)},\qquad A:=\sqrt d\,e^{1+\Lambda}M_0,\quad\Lambda=\|\nabla u\|_{L^1_tL^\infty}.
\end{equation}
\end{enumerate}
\end{lem}

\begin{proof}
\proofstep{Step 1: the map}
For $c\in C([0,1];\RR^d)$ we have $U+c\in\Vc$, so $J[U+c]$ is defined and continuous by Lemma~\ref{lem:char}. Let
$\mathcal T(c)$ denote the right-hand side of \eqref{eq:fixed}; then $\mathcal T(c)\in C([0,1];\RR^d)$, $\mathcal T(c)(0)=0$, and $\mathcal T$
is causal: $\mathcal T(c)(t)$ depends only on $c|_{[0,t]}$.

\proofstep{Step 2: existence and uniqueness}
Let $c_1,c_2\in C([0,1];\RR^d)$. Lemma~\ref{lem:dual} with base field $U+c_2\in\Vc$ (whose spatial gradient is $\nabla U$),
perturbation $c_1-c_2\in\Vc$, and test functions $\psi=v_i$ (for which $[\psi]_R=1$) gives
\[
|\mathcal T(c_1)(t)-\mathcal T(c_2)(t)|\le\beta_U\int_0^t|c_1-c_2|\dd s,\qquad \beta_U:=\sqrt d\,e^{1+\|\nabla U\|_{L^1_tL^\infty}}M_0 .
\]
In the norm $\|c\|_*=\sup_te^{-2\beta_Ut}|c(t)|$, $\mathcal T$ is a contraction with constant $\frac12$, and has a unique fixed
point. Applying the same argument on $[0,t_0]$ shows that the fixed point of the restricted problem is the restriction of $c$.

\proofstep{Step 3: regularity and the mean equation}
If $c\in C^k([0,1])$, then $U+c$ is $C^k$ in time with values in $C^\infty$, so the flow, $f[U+c]$ and
$F[f[U+c],U+c]$ are $C^k$ in time. By \eqref{eq:Jprime}, $J[U+c]\in C^{k+1}$, hence $c=\mathcal T(c)\in C^{k+1}$. By
induction $c\in C^\infty$, and $c(0)=0$. Since $\fint U^\sharp=\fint U+c$, \eqref{eq:fixed} is \eqref{eq:momentum} for
$U^\sharp$; differentiating and using \eqref{eq:Jprime} gives the mean equation.

\proofstep{Step 4: proof of (b)}
If $U$ satisfies \eqref{eq:momentum} on $[0,t_0]$, then $\mathcal T(0)=0$ on $[0,t_0]$, so $0$ is the fixed point of the problem
restricted to $[0,t_0]$, and $c=0$ there by Step~2.

\proofstep{Step 5: proof of (c)}
Since $\fint w=0$, we have $\fint U=\fint u$, and \eqref{eq:momentum} for $u$ turns \eqref{eq:fixed} into
\[
c(t)=J[u](t)-J[u+w+c](t).
\]
Lemma~\ref{lem:dual} with base field $u$ (not $U+c$), perturbation $w+c$, and $\psi=v_i$ gives
\[
|c(t)|\le A\int_0^t\big(\|w(s)\|_{L^\infty}+|c(s)|\big)\dd s ,
\]
and Gronwall's inequality yields \eqref{eq:csize}.
\end{proof}

\begin{remark}
The correction $c$ restores the conservation of total momentum. It is not compactly supported in time: after the support
of $w$, the particles carry a slightly different momentum, which the fluid mean has to compensate. This does not affect the
construction, since only $c=0$ before the support of $w$ is needed. For the intervals of regularity, where the concentration step of
\cite[Section~3]{CL22} has to be carried out causally in time, the mean is handled by the gluing of Proposition~\ref{prop:glue}; see
Section~\ref{sec:regint}.
\end{remark}

\subsection{Proof of Proposition~\ref{prop:VNSstep}}\label{sec:proofstep}

Let $(u,p_u,R)$, $0<s<t_*\le1$ and $\delta\in(0,1]$ be as in Proposition~\ref{prop:VNSstep}. Set
\[
L=\|u\|_{L^\infty(0,1;W^{1,\infty})},\qquad A=\sqrt d\,e^{1+L}M_0,\qquad \bar A=1+Ae^A,
\]
let $C_\star$ be the constant of Proposition~\ref{prop:coupling} for this $L$ and the exponent $r$ of
Proposition~\ref{prop:NSstep}, and choose $\delta'\in(0,\delta]$ so small that
\begin{equation}\label{eq:deltaprime}
\bar A\delta'\le1,\qquad 2\bar A\delta'\le\delta,\qquad
\delta'+C_\star\big((\bar A\delta')^\theta+\bar A\delta'\big)+4Ae^A(L+1)\delta'\le\delta .
\end{equation}
The choice depends only on $u$, $f_0$, $d$, $p$, $q$ and $\delta$.

\proofstep{Step 1: the Navier--Stokes step}
Apply Proposition~\ref{prop:NSstep} to $(u,R)$ with $s,t_*$ and $\delta'$, obtaining $w,P,R^{\rm CL}$. Since $p\ge1$,
\[
W:=\|w\|_{L^1_tL^\infty}\le\|w\|_{L^p_tL^\infty}\le\delta' .
\]

\proofstep{Step 2: the mean correction}
Since $(u,p_u,R)$ solves the VNS--Reynolds system, $u$ satisfies \eqref{eq:momentum}. Lemma~\ref{lem:mean} with $U=u+w$
gives $c\in C^\infty([0,1];\RR^d)$ with $|c(t)|\le Ae^AW\le Ae^A\delta'$, and $c=0$ on $[0,s]$ by \eqref{eq:csize}, since
$w=0$ there. Define
\[
u_1=u+w+c,\qquad h=F\big[f[u_1],u_1\big]-F\big[f[u],u\big].
\]
By Lemma~\ref{lem:mean}(a), $\frac{\mathrm d}{\mathrm dt}\fint u_1=\fint F[f[u_1],u_1]$, while
$\frac{\mathrm d}{\mathrm dt}\fint u=\fint F[f[u],u]$. Since $\fint u_1=\fint u+c$, this gives
\begin{equation}\label{eq:cprime}
c'=\fint h .
\end{equation}

\proofstep{Step 3: the new stress}
Define
\[
R_1=R^{\rm CL}+c\ot(u+w)+(u+w)\ot c-\Rs h,\qquad p_1=p_u+P-\tfrac2d\,c\cdot(u+w).
\]
These are smooth, $R_1$ is $\Sz$-valued, and $\D u_1=0$. Since $c$ depends only on $t$ and $\D(u+w)=0$, we have
$\D(c\otimes c)=0$ and
\[
\D\big((u+w)\otimes c+c\otimes(u+w)\big)=\D\big(c\ot(u+w)+(u+w)\ot c\big)+\tfrac2d\nabla\big(c\cdot(u+w)\big).
\]
Adding \eqref{eq:VNSR} for $u$ and Proposition~\ref{prop:NSstep}(i), and using $\partial_tu_1=\partial_tu+\partial_tw+c'$,
\[
\begin{aligned}
\partial_tu_1-\Delta u_1+\D(u_1\otimes u_1)+\nabla(p_u+P)&=F\big[f[u],u\big]+\D R^{\rm CL}+\D\big(c\ot(u+w)+(u+w)\ot c\big)\\
&\quad+\tfrac2d\nabla\big(c\cdot(u+w)\big)+c' .
\end{aligned}
\]
Moving the gradient term to the left-hand side gives $\nabla p_1$ with the minus sign in the definition of $p_1$.
By \eqref{eq:cprime} and \eqref{eq:R-basic}, $c'=\fint h=h-\D\Rs h$. The right-hand side therefore equals
$F[f[u_1],u_1]+\D R_1$, and $(u_1,p_1,R_1)$ is a smooth solution of the VNS--Reynolds system.

\proofstep{Step 4: the estimates}
(i) On $[0,s]$ we have $w=0$ and $c=0$, so $u_1=u$, $f[u_1]=f[u]$ by Lemma~\ref{lem:char}(b), and $h=0$. Together with
$R^{\rm CL}=0$ on $[0,s]$, this gives $R_1=0$ there.

(ii)--(iii) Since $c$ is constant in space and $|\TT^d|=1$, each of the norms $\|c\|_{L^2_{t,x}}$, $\|c\|_{L^p_tL^\infty}$,
$\|c\|_{L^1_tW^{1,q}}$ is at most $\sup_t|c(t)|\le Ae^A\delta'$. With Proposition~\ref{prop:NSstep}(iii),(iv) and
\eqref{eq:deltaprime} this gives (ii) and (iii).

(iv) The field $\tilde w=w+c$ is smooth with $\|\tilde w\|_{L^1_tL^\infty}\le\bar A\delta'\le1$, so Proposition~\ref{prop:coupling}
gives $\|\Rs h\|_{L^1_tL^r}\le C_\star((\bar A\delta')^\theta+\bar A\delta')$. Moreover
\[
\big\|c\ot(u+w)+(u+w)\ot c\big\|_{L^1_tL^r}\le4\sup_t|c|\,\big(\|u\|_{L^1_tL^\infty}+\|w\|_{L^1_tL^\infty}\big)\le4Ae^A(L+1)\delta',
\]
and $\|R^{\rm CL}\|_{L^1_tL^r}\le\delta'$. By \eqref{eq:deltaprime}, $\|R_1\|_{L^1_tL^r}\le\delta$.
\qed

\begin{remark}[Dependence of the constants]
The constants $A$ and $C_\star$ depend on $e^{L}$, i.e.\ on the Lipschitz norm of the current velocity, which grows along
the iteration. They do not multiply a quantity that depends on the new perturbation through anything other than
$\|w\|_{L^1_tL^\infty}$, and \eqref{eq:deltaprime} fixes $\delta'$ after $u$ is known. In the proof of Theorem~\ref{thm:main}
the same mechanism appears in the choice of $\delta_n$, which contains $e^{-\Lambda_n}$ so that the distribution functions
converge (Lemma~\ref{lem:convergence}).
\end{remark}

\section{The Navier--Stokes step: differences from \cite{CL22}}\label{sec:NSdiff}

Proposition~\ref{prop:NSstep} is the convex integration step of Cheskidov and Luo \cite[Sections~4--5]{CL22}, without their
concentration step. Neither the drag force nor the distribution function enters it. In this section we recall the
construction and list the changes. For completeness, a full proof in the notation used here is given in
Appendices~\ref{sec:CI} and~\ref{sec:estimates}.

\paragraph{The construction.} If $R\equiv0$ we take $w=0$. Otherwise the perturbation depends on four parameters: $\nu$ (temporal
oscillation), $\kappa$ (temporal concentration), $\sigma$ (spatial oscillation) and $\mu$ (spatial concentration). All four are powers of a
single frequency $\lambda$, see \eqref{eq:params}, with $\nu\ll\mu\ll\sigma\ll\kappa^{1/2}$.
\begin{itemize}
\item \emph{Spatial building blocks.} We use the stationary Mikado flows $W_k(\sigma x)$, $k\in\Xi$, of \cite{DS17,MSz18}, where $\Xi$ is given by the
geometric lemma. Each $W_k$ is concentrated in a tube of radius $\mu^{-1}$ around a periodic line in the direction $k$, is normalized in
$L^2$, and satisfies $\D(W_k\otimes W_k)=0$ and $W_k=\D\Omega_k$ with $\Omega_k$ skew-symmetric (Lemma~\ref{lem:mikado}).
\item \emph{Temporal intermittency.} Let $g_\kappa$ be the $1$-periodic extension of $\kappa^{1/2}g(\kappa t)$. Then $\int_0^1g_\kappa^2=1$, while $g_\kappa$ is
supported on a set of measure $\kappa^{-1}$ per period. This makes $\|w\|_{L^p_tL^\infty}$ small for $p<2$ while keeping
$\|w\|_{L^2_{t,x}}\lesssim\|R\|_{L^1}^{1/2}$.
\item \emph{Amplitudes and perturbation.} We set
$a_k=\theta(t)\,g_\kappa(\nu t)\,\rho^{1/2}\,\Gamma_k(\mathrm{Id}-R/\rho)$, where $\rho$ is a regularization of $4|R|$ \eqref{eq:rho} and $\theta$ is a cutoff in
time \eqref{eq:theta}. The perturbation is $w=w^{(p)}+w^{(c)}+w^{(t)}$:
  \begin{itemize}
  \item the principal part $w^{(p)}=\sum_ka_kW_k(\sigma x)$;
  \item a divergence corrector $w^{(c)}$, chosen so that $w^{(p)}+w^{(c)}$ is the divergence of a skew-symmetric field \eqref{eq:divform};
  \item the temporal corrector $w^{(t)}=\nu^{-1}h_\kappa(\nu t)\,(\D R-\nabla\Delta^{-1}\D\D R)$ \eqref{eq:wt}, which absorbs the part of $g_\kappa^2(\nu t)R$ that
  oscillates in time.
  \end{itemize}
\item \emph{New stress.} We have $R^{\rm CL}=R_{\rm lin}+R_{\rm cor}+R_{\rm osc}$ \eqref{eq:RCL}. The parameters \eqref{eq:gamma}--\eqref{eq:params} give the
scaling inequalities \eqref{eq:S1}--\eqref{eq:S5}, from which (iii)--(v) follow with $r=d/(d-\gamma)$ close to $1$.
\end{itemize}

\paragraph{Differences from \cite{CL22}.}
\begin{enumerate}
\item The concentration step of \cite[Section~3]{CL22} is omitted. The two-sided cutoff \cite[(4.12)]{CL22} is replaced by a
one-sided cutoff $\theta$ with $\theta=0$ on $(-\infty,s]$ and $\theta=1$ on $[t_*,\infty)$. Since $R=0$ on $[0,t_*]$, we still have
$\theta^2R=R$ \eqref{eq:thetaR}, which is all that the identity of \cite[Lemma~4.4]{CL22} uses. At $t=1$ no cutoff is needed. In
\cite{CL22} the stress vanishes near $t=0$ and $t=1$, so the functions of time to which \cite[Lemma~B.1]{CL22} is applied are periodic.
Here $R$ need not vanish near $t=1$, and the use of \cite[Lemma~B.1]{CL22} in time is replaced by an elementary decorrelation
estimate (Lemma~\ref{lem:temporal}(c)).
\item As in \cite{CL22}, the perturbation has zero mean: $w^{(p)}+w^{(c)}$ is the divergence of a skew-symmetric field, and $w^{(t)}$ is
built from $\D R$ and a gradient. We record it in (ii) because here the background $u$ need not have zero mean (item 4 below), and the
spatial mean of the velocity is changed only by the constant correction of Lemma~\ref{lem:mean}.
\item Every term of $R^{\rm CL}$ contains $a_k$, $w$ or $R$, or one of their time derivatives, and all of these vanish on $[0,s]$; this gives
(ii).
\item Nothing in the construction uses $\fint u=0$. The identity (i) is algebraic and does not require $(u,R)$ to solve any equation.
\item The $L^2$ bound (iii) follows from \cite[Propositions~5.3--5.5]{CL22} and the proof of \cite[Lemma~5.7]{CL22}, for large
$\lambda$. The statements of \cite[Propositions~2.2(2) and~4.1(2)]{CL22} are written with $\|R\|_{L^1}$; the proofs there give
the bound with $\|R\|_{L^1}^{1/2}$, and only the summability of the $L^2$ norms of the perturbations is used in \cite{CL22}.
\item The exponent $r$ of \cite[(5.2)]{CL22} satisfies $d-\frac dr\le\gamma<1$, hence $r<\frac d{d-1}$.
\item Three minor points of \cite[Section~4]{CL22} are adjusted:
  \begin{itemize}
  \item the directions are normalized to primitive vectors, $\Xi$ is treated as an indexed family in which two indices may carry the
  same vector, and the base points of the Mikado lines are chosen so that parallel lines are disjoint \eqref{eq:lines};
  \item the far-field interaction is replaced by its traceless part, the trace going into the pressure;
  \item the last term of the pressure appears with the opposite sign in \cite[(4.25)]{CL22}; this does not affect the estimates
  (Lemma~\ref{lem:osc}).
  \end{itemize}
\end{enumerate}

\paragraph{Use in Section~\ref{sec:regint}.} Outside Section~\ref{sec:proofstep}, Proposition~\ref{prop:NSstep} is used as a black box, except in the
localized step of Section~\ref{sec:regint} (Lemma~\ref{lem:CIloc}). That step uses one feature of the proof: the algebra of
Appendix~\ref{sec:CI} only uses $\theta^2R=R$, and the estimates of Appendix~\ref{sec:estimates} only use the smoothness of $\theta$ and
$0\le\theta\le1$. Moreover, the time supports of $a_k$, $w^{(p)}$ and $w^{(c)}$ lie in $\Supp\theta$, and those of $w^{(t)}$ and of the terms of
$R^{\rm CL}$ not containing $\theta$ lie in $\Supp_tR\subset\{\theta=1\}$, by $\theta^2R=R$; hence $\Supp_tw,\Supp_tR^{\rm CL}\subset\Supp\theta$. So $\theta$ may be
replaced by a two-sided cutoff. The part of the stress that is not cancelled passes unchanged into the new stress by the linearity of
the identity (i) in the added term $\D S_s$. Elsewhere in Section~\ref{sec:regint}, the exponent $r$ of Appendix~\ref{sec:parameters} and the
relations \eqref{eq:gamma}--\eqref{eq:r} are used (see the remark after Lemma~\ref{lem:B2}).
The Laplacian plays no role in the construction (Remark~\ref{rem:final}).

\section{Sharpness: uniqueness in the critical class}\label{sec:sharp}

The solutions of Theorem~\ref{thm:main} have $u\in L^p_tL^\infty$ for every $p<2$. In this section we show that the
endpoint $p=2$ cannot be reached: weak solutions with $u\in L^2_tL^\infty$ are unique, regular and of Leray--Hopf
type. This is the analogue for VNS of the Ladyzhenskaya--Prodi--Serrin endpoint in \cite[Theorem~1.3 and
Appendix~A]{CL22}. The proof has three steps.
\begin{enumerate}[label=(\arabic*)]
\item \emph{Duality.} Any $L^2_{t,x}$ weak solution of the Oseen problem with an $L^2_tL^\infty$ drift coincides with
its strong solution (Section~\ref{sec:oseen}). This is the argument of \cite[Appendix~A]{CL22}, after
\cite{LM01}.
\item \emph{Regularity.} A critical VNS solution solves an Oseen problem with drift $u$ and force
$F[f,u]\in L^2_tL^\infty$. Maximal $L^2_tL^q$ regularity then gives $\nabla u\in L^2_tL^\infty$ and the energy
equality (Section~\ref{sec:critreg}).
\item \emph{Lagrangian stability.} The difference of two solutions is controlled by the fluid energy of
$u_1-u_2$ together with an $L^2(f_0)$ distance between the two phase-space flows (Section~\ref{sec:lagstab}). This is
a Lagrangian version of the dual stability Lemma~\ref{lem:dual}. A maximal-function estimate replaces the
$C^1$ norm of the test function there, which is not available here. This step follows the method of Han-Kwan, Miot,
Moussa and Moyano \cite{HMMM20}, who combined Loeper's optimal-transport functional \cite{Loe06} with the maximal
function to prove uniqueness of energy solutions for $d=2$; see also \cite[Section~2]{HMMM25} for the three-dimensional
weak--strong setting.
\end{enumerate}

\paragraph{Conventions.} Throughout this section $T\in(0,1]$. Definitions~\ref{def:weak} and~\ref{def:LH} are used on
$[0,T]$ in place of $[0,1]$, with test functions vanishing near $t=T$. Let
$\QQ g:=g-\nabla\Delta^{-1}\D g$ be the Leray projector. Unlike $\PP$, it keeps the spatial mean. It is self-adjoint on
$L^2$, bounded on $L^q$ and $W^{k,q}$ for $1<q<\infty$, commutes with $\Delta$, satisfies $\QQ g=g$ if $\D g=0$, and
$g-\QQ g$ is a gradient. We write $\mathcal D_T$ for the set of divergence-free
$\varphi\in C^\infty([0,T]\times\TT^d;\RR^d)$ vanishing near $t=T$.

\subsection{Main result}

\begin{theo}[Uniqueness and regularity in the critical class]\label{thm:sharp}
Let $d\ge2$ and $T\in(0,1]$. Let $u_{\rm in}\in C^\infty(\TT^d;\RR^d)$ be divergence-free, and let
$f_0\in L^\infty(\TT^d\times\RR^d)$ be nonnegative and vanish for $|v|\ge \Rv_0$. Then:
\begin{enumerate}[label=(\alph*)]
\item every weak solution $(u,f)$ of \eqref{eq:VNS} on $[0,T]$ with data $(u_{\rm in},f_0)$ and
\[
u\in L^2(0,T;L^\infty(\TT^d))
\]
satisfies $\partial_tu\in L^2(0,T;L^q)$ and $u\in L^2(0,T;W^{2,q})\cap C([0,T];W^{1,q})$ for every $q<\infty$; in particular
$\nabla u\in L^2(0,T;L^\infty)$. Moreover $f$ is the push-forward of $f_0$ under the classical phase-space flow of $u$,
and $(u,f)$ satisfies the energy \emph{equality}
\begin{equation}\label{eq:energy}
\mathcal E(t)+\int_0^t\Big(\|\nabla u\|_{L^2}^2+\iint|v-u|^2f\dd z\Big)\dd s=\mathcal E(0)\qquad(0\le t\le T);
\end{equation}
in particular $(u,f)$ is of Leray--Hopf type;
\item there is at most one such weak solution.
\end{enumerate}
\end{theo}

\begin{cor}[Sharp dichotomy]\label{cor:dichotomy}
Let $d\ge2$, $1\le p<2$, $1\le q<\infty$. Let $(u_{\rm in},f_0)$ be smooth data, with $f_0\ge0$ compactly supported, and
let $(\bar u,\bar f)$ be a smooth solution of \eqref{eq:VNS} on $[0,t_1]$ with these data. Let $T\in(0,t_1]$. Then:
\begin{enumerate}[label=(\roman*)]
\item $(\bar u,\bar f)$ is the only weak solution on $[0,T]$ with $u\in L^2(0,T;L^\infty)$;
\item there are infinitely many weak solutions on $[0,T]$ with $u\in L^p(0,T;L^\infty)\cap L^1(0,T;W^{1,q})$. All of
them except possibly one satisfy $\|u\|_{L^2(0,T;L^\infty)}=\infty$, and infinitely many of them are not of Leray--Hopf type.
\end{enumerate}
\end{cor}

\begin{proof}
(i) follows from Theorem~\ref{thm:sharp}(b), since $\bar u\in L^2(0,T;L^\infty)$.

(ii) Argue as in the proof of Corollary~\ref{cor:nonunique}, with $t_1$ replaced by $T$. Choose $\chi\in C^\infty(\RR;[0,1])$
with $\chi=1$ on $[0,T/2]$ and $\chi=0$ on $[3T/4,\infty)$. Choose a smooth divergence-free $\Psi$ with $\fint\Psi=0$,
supported in $(T/2,T)\times\TT^d$, and with $\|\Psi\|_{L^p(0,T;L^2)}\ge1$. Apply Theorem~\ref{thm:main} with $t_0=T/2$,
$\varepsilon=\frac14$ and $\tilde u_k=\chi\bar u+k\Psi$, and restrict the resulting solutions to $[0,T]$. A weak solution
on $[0,1]$ restricts to a weak solution on $[0,T]$, because a test function in $\mathcal D_T$ extends by zero to an
admissible test function on $[0,1]$. Since the error terms in Theorem~\ref{thm:main}(3) are measured on $[0,1]\supset[0,T]$,
\[
\big\|(u_k-\textstyle\fint u_k)-(u_l-\fint u_l)\big\|_{L^p(0,T;L^\infty)}\ge|k-l|\,\|\Psi\|_{L^p(0,T;L^2)}-2\varepsilon\ge\tfrac12
\qquad(k\ne l).
\]
So the restrictions are pairwise distinct, and by (i) at most one of them lies in $L^2(0,T;L^\infty)$. The failure of
the Leray--Hopf property for all large $k$ follows as in Corollary~\ref{cor:nonunique}. Indeed, a Leray--Hopf solution
satisfies $\|u_k\|_{L^p(0,T;L^2)}\le(2\mathcal E(0))^{1/2}$. On the other hand,
$\|u_k\|_{L^p(0,T;L^2)}\ge k\|\Psi\|_{L^p(0,T;L^2)}-\|\chi\bar u-\fint\chi\bar u\|_{L^p(0,T;L^2)}-2\varepsilon\to\infty$.
\end{proof}

\begin{remark}\label{rem:sharp}
\begin{enumerate}
\item In Theorem~\ref{thm:sharp} the smoothness of $u_{\rm in}$ is used only through $u_{\rm in}\in B^1_{q,2}(\TT^d)$ for
one $q>d$; see Corollary~\ref{cor:roughtrace} for the precise statement. No regularity of $f_0$ is needed.
\item The class $L^2_tL^\infty$ is the endpoint $q=\infty$ of the scaling-critical family $X^{p,q}$, $\frac2p+\frac dq=1$.
The other critical spaces could be treated with the interpolation inequality \cite[Proposition~A.4]{CL22} in
Lemma~\ref{lem:adjoint}. We do not pursue this, since $L^2_tL^\infty$ is the endpoint matching Theorem~\ref{thm:main}.
\item Since the force $F[f,u]=j_f-\rho_fu$ is only bounded in $x$, one cannot expect more than $\nabla^2u\in L^2_tL^q$,
$q<\infty$, from the fluid equation. Energy estimates alone would give $\nabla^2u\in L^2_tL^2$, which does not make the
phase-space flow Lipschitz in $d\ge2$. For this reason we use $L^q$-maximal regularity, which yields
$\nabla u\in L^2_tL^\infty$.
\end{enumerate}
\end{remark}

\subsection{Linear tools}

For $2\le q<\infty$ let $X_q:=B^1_{q,2}(\TT^d)$. For an interval $I=[a,b]\subset[0,1]$ let
\[
\mathrm{MR}_q(I):=\big\{v\in L^2(I;W^{2,q}(\TT^d)):\ \partial_tv\in L^2(I;L^q(\TT^d))\big\},\qquad
\|v\|_{\mathrm{MR}_q(I)}:=\|v\|_{L^2(I;W^{2,q})}+\|\partial_tv\|_{L^2(I;L^q)} .
\]

\begin{lem}[Maximal regularity]\label{lem:MR}
Let $2\le q<\infty$. There is $C_q\ge1$, independent of the interval $I=[a,b]\subset[0,1]$, such that:
\begin{enumerate}[label=(\alph*)]
\item for $v_0\in X_q$ and $g\in L^2(I;L^q)$ there is a unique $v\in\mathrm{MR}_q(I)$ with $\partial_tv-\Delta v=g$ a.e.\ and $v(a)=v_0$, and
$\|v\|_{\mathrm{MR}_q(I)}\le C_q\big(\|v_0\|_{X_q}+\|g\|_{L^2(I;L^q)}\big)$;
\item every $v\in\mathrm{MR}_q(I)$ belongs to $C(I;X_q)$, with
\[
\sup_I\|v\|_{X_q}\le C_q\big(\|v(a)\|_{X_q}+\|v\|_{\mathrm{MR}_q(I)}\big),
\]
and $\sup_I\|v\|_{X_q}\le C_q\|v\|_{\mathrm{MR}_q(I)}$ if $v(a)=0$;
\item $C^\infty(\TT^d)\subset X_q\hookrightarrow W^{1,q}$, $\mathrm{MR}_q(I)\hookrightarrow\mathrm{MR}_2(I)$, and
$\mathrm{MR}_2(I)\hookrightarrow C(I;H^1)$;
\item if $v_0$ and $g$ are divergence-free, so is $v$.
\end{enumerate}
\end{lem}

\begin{proof}[References]
On $L^q(\TT^d)$, $1<q<\infty$, the Laplacian has maximal $L^2$-regularity on bounded intervals. The trace space of
$\mathrm{MR}_q$ is the real interpolation space $(L^q,W^{2,q})_{1/2,2}=B^1_{q,2}$; see \cite[Chapter~III]{Ama95} and
\cite[Chapter~3]{PS16}. In particular this yields (a) and (b) on $[0,1]$. The constants are uniform in $I$ by the usual
extensions: in (a), translate to $a=0$ and extend $g$ by zero to $[0,1]$; in (b) with $v(a)=0$, extend $v$ by zero to
$[a-1,a]$ and use the trace embedding on intervals of length in $[1,2]$. The general case of (b) follows by subtracting
$e^{(t-a)\Delta}v(a)$. For (c): $B^1_{q,2}\hookrightarrow F^1_{q,2}=W^{1,q}$ for $q\ge2$ \cite{Tri83}; the
other embeddings hold since $|\TT^d|=1$, and the last one is the Lions--Magenes lemma. For (d): $\QQ$ commutes with the heat
flow, so $\QQ v$ solves the same problem as $v$.
\end{proof}

\begin{lem}[Maximal function]\label{lem:maximal}
Let $g\in H^1(\TT^d)$, extended periodically to $\RR^d$, and let $\mathcal M$ be the Hardy--Littlewood maximal operator on $\RR^d$.
Then $\mathcal M|\nabla g|$ is periodic,
\begin{equation}\label{eq:maxL2}
\|\mathcal M|\nabla g|\|_{L^2(\TT^d)}\le C_M\|\nabla g\|_{L^2(\TT^d)},
\end{equation}
and there is a null set $N\subset\RR^d$ such that
\begin{equation}\label{eq:maxpw}
|g(x)-g(y)|\le C_d|x-y|\big(\mathcal M|\nabla g|(x)+\mathcal M|\nabla g|(y)\big)\qquad\text{for all }x,y\in\RR^d\setminus N .
\end{equation}
The constants $C_M$ and $C_d$ depend only on $d$.
\end{lem}

\begin{proof}
\eqref{eq:maxpw} is the classical estimate for $W^{1,1}_{\rm loc}$ functions, with $N$ the complement of the Lebesgue points of
$g$; see \cite[Lemma~A.3]{CDL08}. It is stated there with the maximal function truncated at a radius $R\ge|x-y|$, which
is bounded by $\mathcal M$. For \eqref{eq:maxL2}, let $h=|\nabla g|$ and $x\in[0,1]^d$. Balls of radius $r<1$ around $x$ lie in
$B:=[-1,2]^d$, so the supremum over $r<1$ is at most $\mathcal M(h\mathbf 1_B)(x)$. For $r\ge1$, the ball $B(x,r)$ meets at most
$C(1+r)^d$ unit cells, so $\fint_{B(x,r)}h\le C\|h\|_{L^1(\TT^d)}$. Hence
$\|\mathcal Mh\|_{L^2([0,1]^d)}\le\|\mathcal M(h\mathbf 1_B)\|_{L^2(\RR^d)}+C\|h\|_{L^1(\TT^d)}\le C\|h\|_{L^2(\TT^d)}$.
\end{proof}

\subsection{The Oseen problem: strong solutions and duality}\label{sec:oseen}

Let $b\in L^2(0,T;L^\infty(\TT^d;\RR^d))$ be weakly divergence-free for a.e.\ $t$. For $G\in L^1((0,T)\times\TT^d)$ and
$v_0\in L^2$, we call $v\in L^2((0,T)\times\TT^d)$, weakly divergence-free for a.e.\ $t$, a \emph{weak solution of the Oseen
problem}
\begin{equation}\label{eq:oseen}
\partial_tv-\Delta v+\QQ(b\cdot\nabla v)=\QQ G,\qquad \D v=0,\qquad v(0)=v_0,
\end{equation}
if for every $\varphi\in\mathcal D_T$
\begin{equation}\label{eq:oseenweak}
\int v_0\cdot\varphi(0)\dd x+\int_0^T\!\!\int\Big(v\cdot\big(\partial_t\varphi+\Delta\varphi+b\cdot\nabla\varphi\big)+G\cdot\varphi\Big)\dd x\dd t=0 .
\end{equation}
Here $(b\cdot\nabla\varphi)_i=b_j\partial_j\varphi_i$. This formulation is meaningful because $b\cdot\nabla\varphi\in L^2_{t,x}$.

\begin{prop}[Strong solutions]\label{prop:oseenstrong}
Let $2\le q<\infty$, $G\in L^2(0,T;L^q)$, and let $v_0\in X_q$ be divergence-free. There is a unique divergence-free
$v\in\mathrm{MR}_q([0,T])$ solving \eqref{eq:oseen} a.e.\ in $(0,T)\times\TT^d$. It is a weak solution in the sense of
\eqref{eq:oseenweak}. If the hypotheses hold for several $q$, the solutions coincide.
\end{prop}

\begin{proof}
For $v\in\mathrm{MR}_q(I)$, Lemma~\ref{lem:MR}(b),(c) gives $\nabla v\in L^\infty(I;L^q)$, hence
\begin{equation}\label{eq:drift}
\|\QQ(b\cdot\nabla v)\|_{L^2(I;L^q)}\le C\|b\|_{L^2(I;L^\infty)}\sup_{t\in I}\|\nabla v(t)\|_{L^q}<\infty .
\end{equation}
Since $\|b(\cdot)\|_{L^\infty}^2\in L^1(0,T)$, there is a partition $0=\tau_0<\tau_1<\dots<\tau_n=T$ with
$C\,C_q^2\|b\|_{L^2(I_j;L^\infty)}\le\frac12$ on $I_j=[\tau_{j-1},\tau_j]$, where $C$ is the constant in \eqref{eq:drift}.

We construct $v$ interval by interval. Suppose $v(\tau_{j-1})=:v_{j-1}\in X_q$ is divergence-free. On
$\mathcal X_j:=\{v\in\mathrm{MR}_q(I_j):v(\tau_{j-1})=v_{j-1}\}$, let $\Theta(v)$ be the solution of Lemma~\ref{lem:MR}(a) with
$g=\QQ G-\QQ(b\cdot\nabla v)$. For $v,v'\in\mathcal X_j$ the difference $\Theta(v)-\Theta(v')$ has zero initial value, so by
Lemma~\ref{lem:MR}(a),(b) and \eqref{eq:drift},
\[
\|\Theta(v)-\Theta(v')\|_{\mathrm{MR}_q(I_j)}\le C_q\,C\|b\|_{L^2(I_j;L^\infty)}\,C_q\|v-v'\|_{\mathrm{MR}_q(I_j)}\le\tfrac12\|v-v'\|_{\mathrm{MR}_q(I_j)} .
\]
The unique fixed point solves \eqref{eq:oseen} on $I_j$. It is divergence-free by Lemma~\ref{lem:MR}(d), and
$v(\tau_j)\in X_q$ by Lemma~\ref{lem:MR}(b). The concatenation is continuous at the $\tau_j$, hence lies in
$\mathrm{MR}_q([0,T])$. Uniqueness follows from the contraction on each $I_j$.

The solution is weak: test \eqref{eq:oseen} with $\varphi\in\mathcal D_T$ and integrate by parts in $t$ and $x$. Since
$\D\varphi=0$, $\int\QQ g\cdot\varphi=\int g\cdot\varphi$. Moreover, $\int(b\cdot\nabla v)\cdot\varphi=-\int v\cdot(b\cdot\nabla\varphi)$
because $\int b\cdot\nabla(v\cdot\varphi)=0$ for the weakly divergence-free $b\in L^\infty$ and $v\cdot\varphi\in W^{1,1}$.

For $q\ge q'\ge2$ we have $\mathrm{MR}_q\subset\mathrm{MR}_{q'}$ and $X_q\subset X_{q'}$. Both solutions solve the problem in
$\mathrm{MR}_{q'}$, where the solution is unique.
\end{proof}

\begin{lem}[Adjoint problem]\label{lem:adjoint}
Let $H\in C^\infty([0,T]\times\TT^d;\RR^d)$. There is a divergence-free $\Phi$ with
\[
\Phi\in L^2(0,T;H^2)\cap C([0,T];H^1),\qquad \partial_t\Phi\in L^2_{t,x},\qquad b\cdot\nabla\Phi\in L^2_{t,x},\qquad \Phi(T)=0,
\]
solving $-\partial_t\Phi-\Delta\Phi-\QQ(b\cdot\nabla\Phi)=\QQ H$ a.e.\ in $(0,T)\times\TT^d$.
\end{lem}

\begin{proof}
\proofstep{Galerkin scheme}
Let $V_n$ be the space of real divergence-free trigonometric polynomials of degree at most $n$, and $\Pi_n$ the
$L^2$-orthogonal projection onto $V_n$. Note that $\Pi_n=\QQ$ followed by truncation of Fourier modes, so it commutes with
$\Delta$. Let $\Phi_n\in W^{1,1}(0,T;V_n)$ solve the linear backward ODE
\[
-\partial_t\Phi_n-\Delta\Phi_n-\Pi_n(b\cdot\nabla\Phi_n)=\Pi_nH,\qquad\Phi_n(T)=0 .
\]
Its coefficients are integrable in time, since $b\in L^1(0,T;L^\infty)$, so the solution exists and is unique.

\proofstep{Energy estimates}
Since $\D b=0$, $\int(b\cdot\nabla\Phi_n)\cdot\Phi_n=0$, so $-\frac12\frac{\mathrm d}{\mathrm dt}\|\Phi_n\|_2^2+\|\nabla\Phi_n\|_2^2=\int H\cdot\Phi_n$. Testing
with $-\Delta\Phi_n\in V_n$,
\[
-\tfrac12\tfrac{\mathrm d}{\mathrm dt}\|\nabla\Phi_n\|_2^2+\|\Delta\Phi_n\|_2^2=-\int H\cdot\Delta\Phi_n-\int(b\cdot\nabla\Phi_n)\cdot\Delta\Phi_n
\le\tfrac12\|\Delta\Phi_n\|_2^2+\|H\|_2^2+\|b\|_{L^\infty}^2\|\nabla\Phi_n\|_2^2 .
\]
Gronwall's inequality backward from $t=T$ gives, uniformly in $n$,
\[
\sup_t\big(\|\Phi_n(t)\|_2^2+\|\nabla\Phi_n(t)\|_2^2\big)+\int_0^T\|\Delta\Phi_n\|_2^2\dd t\le C\|H\|_{L^2_{t,x}}^2\exp\big(2\|b\|_{L^2(0,T;L^\infty)}^2\big).
\]
Hence $\|b\cdot\nabla\Phi_n\|_{L^2_{t,x}}\le\|b\|_{L^2_tL^\infty}\|\nabla\Phi_n\|_{L^\infty_tL^2}$ and, from the equation,
$\|\partial_t\Phi_n\|_{L^2_{t,x}}$ are bounded uniformly in $n$.

\proofstep{Limit}
Along a subsequence, $\Phi_n\rightharpoonup^*\Phi$ in $L^\infty(0,T;H^1)$ and weakly in $L^2(0,T;H^2)$, and $\partial_t\Phi_n\rightharpoonup\partial_t\Phi$
in $L^2_{t,x}$. For $g\in L^2_{t,x}$ we have $\iint(b\cdot\nabla\Phi_n)\cdot g=\iint\nabla\Phi_n:(g\otimes b)$ with $g\otimes b\in L^1(0,T;L^2)$, so
$b\cdot\nabla\Phi_n\rightharpoonup b\cdot\nabla\Phi$ in $L^2_{t,x}$. Moreover $\iint\Pi_n(b\cdot\nabla\Phi_n)\cdot g=\iint(b\cdot\nabla\Phi_n)\cdot\Pi_ng$,
and $\Pi_ng\to\QQ g$ strongly in $L^2_{t,x}$. Hence the equation passes to the limit. $\Phi(T)=0$ because
$\Phi_n(T)\rightharpoonup\Phi(T)$, $\Phi$ is divergence-free, and $\Phi\in C([0,T];H^1)$ by Lemma~\ref{lem:MR}(c).
\end{proof}

\begin{lem}[Admissibility]\label{lem:admissible}
Let $\Phi$ be as in Lemma~\ref{lem:adjoint}. There are $\varphi_k\in\mathcal D_T$ with
$\partial_t\varphi_k+\Delta\varphi_k+b\cdot\nabla\varphi_k\to\partial_t\Phi+\Delta\Phi+b\cdot\nabla\Phi$ in $L^2((0,T)\times\TT^d)$.
\end{lem}

\begin{proof}
Extend $\Phi$ by $0$ to $(T,2T]$. Since $\Phi(T)=0$ and $\Phi\in C([0,T];H^1)$, the extension $\tilde\Phi$ lies in
$H^1(0,2T;L^2)\cap L^2(0,2T;H^2)$. For small $h>0$ let
$\varphi_h:=\eta_h*_{t,x}\tilde\Phi(\cdot+2h)$, restricted to $[0,T]$, where $\eta_h$ is a space--time mollifier supported in
$\{|t|,|x|\le h\}$. Then $\varphi_h(t)$ only involves $\tilde\Phi$ on $[t+h,t+3h]\subset(0,2T]$; it is smooth and divergence-free, and it vanishes for
$t\ge T-h$. By continuity of
translations and mollification,
\[
\varphi_h\to\Phi\quad\text{in }H^1(0,T;L^2)\cap L^2(0,T;H^2),\qquad\text{hence in }C([0,T];H^1).
\]
Therefore $\partial_t\varphi_h+\Delta\varphi_h\to\partial_t\Phi+\Delta\Phi$ in $L^2_{t,x}$, and
$\|b\cdot\nabla(\varphi_h-\Phi)\|_{L^2_{t,x}}\le\|b\|_{L^2_tL^\infty}\|\nabla(\varphi_h-\Phi)\|_{L^\infty_tL^2}\to0$.
\end{proof}

\begin{prop}[Uniqueness by duality]\label{prop:duality}
Let $w\in L^2((0,T)\times\TT^d)$ be weakly divergence-free for a.e.\ $t$ and satisfy
\[
\int_0^T\!\!\int w\cdot\big(\partial_t\varphi+\Delta\varphi+b\cdot\nabla\varphi\big)\dd x\dd t=0\qquad\text{for all }\varphi\in\mathcal D_T .
\]
Then $w=0$.
\end{prop}

\begin{proof}
Let $H$ be smooth and let $\Phi$ be given by Lemma~\ref{lem:adjoint}. By Lemma~\ref{lem:admissible}, the identity extends to
$\varphi=\Phi$. By the equation of $\Phi$,
\[
\partial_t\Phi+\Delta\Phi+b\cdot\nabla\Phi=-H+(\mathrm{Id}-\QQ)H+(\mathrm{Id}-\QQ)(b\cdot\nabla\Phi),
\]
and the last two terms are gradients of $H^1$ functions for a.e.\ $t$, which are orthogonal to $w(t)$. Hence $\iint w\cdot H=0$
for all smooth $H$, so $w=0$.
\end{proof}

\subsection{Regularity of critical solutions}\label{sec:critreg}

\begin{prop}\label{prop:critreg}
Under the hypotheses of Theorem~\ref{thm:sharp}, let $(u,f)$ be a weak solution on $[0,T]$ with $u\in L^2(0,T;L^\infty)$.
Let $\Rv_*:=\Rv_0+\|u\|_{L^1(0,T;L^\infty)}$. Then:
\begin{enumerate}[label=(\alph*)]
\item $F:=F[f,u]\in L^2(0,T;L^\infty)$;
\item $u\in\mathrm{MR}_q([0,T])$ for every $q\in[2,\infty)$, and $\partial_tu-\Delta u+\QQ(u\cdot\nabla u)=\QQ F$ a.e.; in particular
$\nabla u\in L^2(0,T;L^\infty)$;
\item the field $b_u(t,x,v)=(v,u(t,x)-v)$ has a unique flow $Z_u(t;s,z)$ in the Carath\'eodory sense. Each $Z_u(t;s,\cdot)$
is bi-Lipschitz with $\det D_zZ_u(t;s,z)=e^{-d(t-s)}$. Moreover $f(t)=Z_u(t;0,\cdot)_\#f_0$ for every $t$ (after
modification on a null set of times), $0\le f\le e^d\|f_0\|_{L^\infty}$, and $f(t)=0$ for $|v|\ge \Rv_*$;
\item the energy equality \eqref{eq:energy} holds.
\end{enumerate}
\end{prop}

\begin{proof}
(a) Let $\Rv$ be such that $f=0$ for $|v|\ge \Rv$ (Definition~\ref{def:weak}(ii)). Then
$\|\rho_f\|_{L^\infty}+\|j_f\|_{L^\infty}\le(1+\Rv)|B_{\Rv}|\|f\|_{L^\infty}$, and $F=j_f-\rho_fu$.

(b) Fix $q\in(d,\infty)$. Let $v\in\mathrm{MR}_q([0,T])$ be the strong solution of Proposition~\ref{prop:oseenstrong} with $b:=u$,
$G:=F$ and $v_0:=u_{\rm in}$; it is a weak solution in the sense of \eqref{eq:oseenweak}. Since
$u\otimes u:\nabla\varphi=u\cdot(u\cdot\nabla\varphi)$, Definition~\ref{def:weak}(iv) says that $u$ is also a weak solution of the same
Oseen problem. The difference $w=u-v\in L^2_{t,x}$ is weakly divergence-free and satisfies the homogeneous identity of
Proposition~\ref{prop:duality}, so $u=v$. Hence $u\in\mathrm{MR}_q$ for all $q\in[2,\infty)$ by Proposition~\ref{prop:oseenstrong},
and $\nabla u\in L^2(0,T;W^{1,q})\subset L^2(0,T;L^\infty)$ for $q>d$.

(c) By (b), $b_u$ is measurable in $t$, Lipschitz in $z$ with constant $1+\|\nabla u(t)\|_{L^\infty}\in L^1(0,T)$, and grows at most
linearly. Carath\'eodory theory gives a unique flow, bi-Lipschitz in $z$. The Jacobian identity follows from Liouville's
formula and $\D_zb_u=-d$. The proof of Lemma~\ref{lem:char}(a),(c) applies verbatim, so the push-forward
$\tilde f(t)=Z_u(t;0,\cdot)_\#f_0$ is a bounded weak solution with the stated bounds. Since
$u\in L^1(0,T;W^{1,1})\cap L^1(0,T;L^\infty)$, Lemma~\ref{lem:DPL} gives $f=\tilde f$.

(d) By (b) and Lemma~\ref{lem:MR}(c), $u\in C([0,T];H^1)$ and $\partial_tu\in L^2_{t,x}$, so $t\mapsto\frac12\|u(t)\|_2^2$ is absolutely
continuous with derivative $\int\partial_tu\cdot u$. Since $\QQ$ is self-adjoint and $\QQ u=u$, and since
$\int(u\cdot\nabla u)\cdot u=\int u\cdot\nabla\frac{|u|^2}2=0$ (as $|u|^2\in W^{1,1}$ for a.e.\ $t$),
\[
\frac{\mathrm d}{\mathrm dt}\frac12\|u\|_2^2+\|\nabla u\|_2^2=\int F\cdot u=\int j_f\cdot u-\int\rho_f|u|^2 .
\]
For the particles, by (c), $\frac12\iint|v|^2f(t)=\frac12\int f_0(z)|V(t;0,z)|^2\dd z$. The integrand is absolutely continuous
in $t$ with derivative $V\cdot(u(t,X)-V)$, bounded by an integrable function on $\Supp f_0$. Hence
\[
\frac{\mathrm d}{\mathrm dt}\frac12\iint|v|^2f\dd z=\int f_0\,V\cdot(u(t,X)-V)\dd z=\int j_f\cdot u\dd x-\iint|v|^2f\dd z .
\]
Adding the two identities and using $|v-u|^2=|v|^2-2v\cdot u+|u|^2$ gives
$\frac{\mathrm d}{\mathrm dt}\mathcal E+\|\nabla u\|_2^2+\iint|v-u|^2f=0$, which integrates to \eqref{eq:energy}.
\end{proof}

\begin{cor}[Rough initial velocity]\label{cor:roughtrace}
Let $d\ge2$, $T\in(0,1]$, $q_0\in(d,\infty)$, let $u_{\rm in}\in B^1_{q_0,2}(\TT^d)$ be divergence-free, and let $f_0\in L^\infty$ be nonnegative and
vanish for $|v|\ge \Rv_0$. Then every weak solution $(u,f)$ on $[0,T]$ with these data and $u\in L^2(0,T;L^\infty)$ satisfies:
\begin{enumerate}[label=(\alph*)]
\item $u\in\mathrm{MR}_{q_0}([0,T])$, hence $\nabla u\in L^2(0,T;L^\infty)$, and $u\in\mathrm{MR}_s([t_*,T])$ for every $s<\infty$ and $t_*\in(0,T)$;
\item $f$ is the push-forward of $f_0$ under the Lipschitz phase-space flow of $u$, and the energy equality \eqref{eq:energy} holds;
\item there is at most one such weak solution.
\end{enumerate}
\end{cor}

\begin{proof}
The proof of Proposition~\ref{prop:critreg}(b) applies with $q=q_0$, since Proposition~\ref{prop:oseenstrong} only needs $v_0\in X_{q_0}$.
This gives $u\in\mathrm{MR}_{q_0}$ and, as $q_0>d$, $\nabla u\in L^2(0,T;L^\infty)$. Parts (c) and (d) of that proposition then follow
verbatim, which gives (b). For the second half of (a), let $\eta\in C^\infty(\RR;[0,1])$ vanish near $0$ and equal $1$ on $[t_*,\infty)$. Then
$v:=\eta u$ is a weak solution of the Oseen problem \eqref{eq:oseen} with drift $u$, zero initial value and forcing
$\eta F+\eta'u\in L^2(0,T;L^s)$ for every $s<\infty$, since $F$ is bounded and $u\in L^2(0,T;L^\infty)$. By Propositions~\ref{prop:oseenstrong}
and~\ref{prop:duality}, $v$ is the strong solution, so $v\in\mathrm{MR}_s([0,T])$ and $u\in\mathrm{MR}_s([t_*,T])$. Finally, the proof of
Theorem~\ref{thm:sharp}(b) in Section~\ref{sec:lagstab} uses only (a) with $q=q_0$ and (b), which gives (c).
\end{proof}

\subsection{Lagrangian stability and proof of Theorem~\ref{thm:sharp}}\label{sec:lagstab}

Part (a) of Theorem~\ref{thm:sharp} is Proposition~\ref{prop:critreg}. Part (b) adapts \cite[Section~4.3]{HMMM20}
to dimension $d$ and to the class $L^2_tL^\infty$. In two dimensions, the bound $u\in L^2_tL^\infty$ holds for every weak solution in the energy class
\cite[Lemma~3]{HMMM20}, as a consequence of the energy identity \cite[Lemma~2]{HMMM20}, and the velocity is log-Lipschitz, so an Osgood argument is used there. In our
setting $\nabla u_i\in L^2_tL^\infty$ by Proposition~\ref{prop:critreg}, and Gronwall's lemma suffices. For (b), let $(u_1,f_1)$ and $(u_2,f_2)$ be two weak
solutions as in the theorem, and let $Z_i=(X_i,V_i)$ be the flows of Proposition~\ref{prop:critreg}(c), lifted to
$\RR^d\times\RR^d$ (with $u_i$ extended periodically) and started at time $0$. Set
\[
w:=u_1-u_2,\qquad \delta Z:=(\delta X,\delta V):=Z_1-Z_2,\qquad
D(t)^2:=\int f_0(z)\,|\delta Z(t,z)|^2\dd z,\qquad Y:=\|w\|_{L^2}^2+D^2 .
\]
Let $\Rv_*$ bound the velocity supports of $f_1,f_2$, and let $\varrho_*:=\max_i\|\rho_{f_i}\|_{L^\infty}\le e^d\|f_0\|_{L^\infty}|B_{\Rv_*}|$. By
Proposition~\ref{prop:critreg}(c), for every periodic measurable $h\ge0$,
\begin{equation}\label{eq:push}
\int f_0(z)\,h\big(X_i(t,z)\big)\dd z=\int_{\TT^d}\rho_{f_i}(t,x)h(x)\dd x\le\varrho_*\|h\|_{L^1(\TT^d)} .
\end{equation}

\proofstep{Step 1: the flows}
For $f_0$-a.e.\ $z$, the map $t\mapsto\delta Z(t,z)$ is absolutely continuous with $\delta Z(0,z)=0$, $\frac{\mathrm d}{\mathrm dt}\delta X=\delta V$ and
\[
\frac{\mathrm d}{\mathrm dt}\delta V=\big(u_1(t,X_1)-u_1(t,X_2)\big)+w(t,X_2)-\delta V .
\]
Since $|u_1(t,X_1)-u_1(t,X_2)|\le\|\nabla u_1(t)\|_{L^\infty}|\delta X|$, Young's inequality gives
\[
\frac{\mathrm d}{\mathrm dt}|\delta Z|^2\le2\big(1+\|\nabla u_1(t)\|_{L^\infty}\big)|\delta Z|^2+|w(t,X_2)|^2 .
\]
Integrating against $f_0$ and using
\eqref{eq:push} with $h=|w(t)|^2$,
\begin{equation}\label{eq:D}
\frac{\mathrm d}{\mathrm dt}D^2\le2\big(1+\|\nabla u_1(t)\|_{L^\infty}\big)D^2+\varrho_*\|w(t)\|_{L^2}^2 .
\end{equation}
This is justified by Fubini, since $|\delta Z|$ is bounded on $[0,T]\times\Supp f_0$.

\proofstep{Step 2: the fluid}
By Proposition~\ref{prop:critreg}(b), $w\in\mathrm{MR}_2$ solves
$\partial_tw-\Delta w+\QQ(u_1\cdot\nabla w+w\cdot\nabla u_2)=\QQ(F_1-F_2)$ with $F_i=F[f_i,u_i]$. Testing with $w$, we use
$\int(u_1\cdot\nabla w)\cdot w=0$ and $\int(w\cdot\nabla u_2)\cdot w=-\int(w\cdot\nabla w)\cdot u_2$. Since
$F_1-F_2=\int(v-u_1)(f_1-f_2)\dd v-\rho_{f_2}w$,
\begin{equation}\label{eq:fluid}
\frac12\frac{\mathrm d}{\mathrm dt}\|w\|_2^2+\|\nabla w\|_2^2\le\|u_2\|_{L^\infty}\|w\|_2\|\nabla w\|_2+I(t),\qquad
I(t):=\iint(f_1-f_2)(t)\,\psi_t\dd z,
\end{equation}
where $\psi_t(x,v)=(v-u_1(t,x))\cdot w(t,x)$ and we dropped $-\int\rho_{f_2}|w|^2\le0$.

\proofstep{Step 3: the coupling term}
By Proposition~\ref{prop:critreg}(c), $I(t)=\int f_0(z)\big(\psi_t(Z_1)-\psi_t(Z_2)\big)\dd z$, and
\[
\psi_t(Z_1)-\psi_t(Z_2)=\delta V\cdot w(X_1)-\big(u_1(X_1)-u_1(X_2)\big)\cdot w(X_1)+\big(V_2-u_1(X_2)\big)\cdot\big(w(X_1)-w(X_2)\big),
\]
with all functions evaluated at time $t$. We bound the three terms after integration against $f_0$.
\begin{itemize}
\item By Cauchy--Schwarz and \eqref{eq:push}, the first is at most $D\,\varrho_*^{1/2}\|w\|_2$.
\item Likewise, the second is at most $\|\nabla u_1\|_{L^\infty}D\,\varrho_*^{1/2}\|w\|_2$.
\item For the third, note that $|V_2-u_1(X_2)|\le \Rv_*+\|u_1(t)\|_{L^\infty}$ on $\Supp f_0$. For a.e.\ $t$ we have $w(t)\in H^1$. Let
$N_t$ be the null set of Lemma~\ref{lem:maximal} for $g=w(t)$. The sets $\{z:X_i(t,z)\in N_t\}$ are null, because $Z_i(t;0,\cdot)$ is
bi-Lipschitz. By \eqref{eq:maxpw}, Cauchy--Schwarz, \eqref{eq:push} and \eqref{eq:maxL2}, the third term is at most
\[
\big(\Rv_*+\|u_1\|_{L^\infty}\big)C_d\sum_{i=1,2}\int f_0|\delta X|\,\mathcal M|\nabla w|(X_i)\dd z
\le2C_dC_M\varrho_*^{1/2}\big(\Rv_*+\|u_1\|_{L^\infty}\big)D\,\|\nabla w\|_2 .
\]
\end{itemize}

\proofstep{Step 4: Gronwall}
Insert Step~3 into \eqref{eq:fluid} and absorb $\|\nabla w\|_2$ by Young's inequality. Combined with \eqref{eq:D}, this gives,
for a.e.\ $t$,
\[
\frac{\mathrm d}{\mathrm dt}Y\le C\,a(t)\,Y,\qquad a(t):=1+\|u_1(t)\|_{L^\infty}^2+\|u_2(t)\|_{L^\infty}^2+\|\nabla u_1(t)\|_{L^\infty},
\]
with $C$ depending only on $d$, $\Rv_*$ and $\varrho_*$. Here $a\in L^1(0,T)$ by the hypothesis $u_i\in L^2_tL^\infty$ and by
Proposition~\ref{prop:critreg}(b). Since $Y$ is absolutely continuous and $Y(0)=0$, we get $Y\equiv0$. Thus $u_1=u_2$, and then
$f_1=f_2$ by Lemma~\ref{lem:DPL}, or directly since $\delta Z=0$ on $\Supp f_0$. This proves Theorem~\ref{thm:sharp}(b).
\qed

\begin{remark}\label{rem:weakstrong}
Step~4 uses the Lipschitz bound only for $u_1$. The argument is therefore a weak--strong uniqueness statement for the
flows, combined with the fact that $L^2_tL^\infty$ solutions are automatically strong (Proposition~\ref{prop:critreg}). The
endpoint enters in two places: in the integrability of $\|u_i\|_{L^\infty}^2$ in $a(t)$, and in the absorption of the drift
term in the maximal regularity argument (Proposition~\ref{prop:oseenstrong}). Both fail for $u\in L^p_tL^\infty$ with $p<2$, in
agreement with Theorem~\ref{thm:main}.
\end{remark}

\section{Intervals of regularity}\label{sec:regint}

Cheskidov and Luo show that their wild solutions are smooth outside a closed set of times of Hausdorff dimension at
most $\varepsilon$ \cite[Theorem~1.8]{CL22}. The tool is a concentration step that glues local exact solutions, so that the stress
of each approximate solution lives on few short intervals \cite[Section~3]{CL22}. In this section we prove the analogue
for VNS.

\begin{figure}[t]
\centering
\resizebox{\textwidth}{!}{%
\begin{tikzpicture}[x=1cm,y=1cm,>=Stealth,font=\small]
\definecolor{cipert}{RGB}{200,60,40}
\definecolor{kin}{RGB}{40,100,170}
\definecolor{bad}{RGB}{225,225,225}
\definecolor{corr}{RGB}{60,140,80}
\node[anchor=west,font=\small\bfseries] at (-0.3,5.15) {(a) Memory of the particles};
\draw[->] (0,4) -- (14.6,4) node[right] {$t$};
\node[anchor=east,align=right,font=\footnotesize] at (-0.1,4.35) {perturbation\\of $u$};
\fill[cipert!70] (5,4) .. controls (5.3,4.85) and (5.9,4.85) .. (6.2,4) -- cycle;
\node[cipert,font=\footnotesize,anchor=west] at (6.5,4.45) {$w$, supported in $[a,b]$};
\foreach \x/\l in {0/0,5/a,6.2/b,14/1}{\draw (\x,3.93)--(\x,4.07); \node[below,font=\footnotesize] at (\x,3.93) {$\l$};}
\draw[->] (0,2) -- (14.6,2) node[right] {$t$};
\node[anchor=east,align=right,font=\footnotesize] at (-0.1,2.35) {change of\\$f[u]$};
\fill[kin!25] (5,2) .. controls (5.6,2.1) and (5.8,2.7) .. (6.2,2.75)
      .. controls (8,2.85) and (10,2.6) .. (12,2.7) .. controls (13,2.75) and (13.6,2.65) .. (14,2.68) -- (14,2) -- cycle;
\draw[kin,thick] (0,2.02) -- (5,2.02) .. controls (5.6,2.1) and (5.8,2.7) .. (6.2,2.75)
      .. controls (8,2.85) and (10,2.6) .. (12,2.7) .. controls (13,2.75) and (13.6,2.65) .. (14,2.68);
\node[kin,font=\footnotesize,align=center] at (10.4,3.25) {persists for $t>b$};
\node[font=\footnotesize,align=center] at (2.5,2.35) {unchanged for $t<a$};
\draw[dashed,gray] (5,4) -- (5,2); \draw[dashed,gray] (6.2,4) -- (6.2,2);
\node[anchor=west,font=\small\bfseries] at (-0.3,0.65) {(b) Causal gluing in two passes};
\foreach \a in {4.2,9.2}{\fill[bad] (\a,-4.75) rectangle (\a+0.6,0.2);}
\node[font=\footnotesize] at (4.5,-5.0) {$J_1$};
\node[font=\footnotesize] at (9.5,-5.0) {$J_2$};
\node[font=\footnotesize] at (2.1,-5.0) {$G_0$};
\node[font=\footnotesize] at (7.0,-5.0) {$G_1$};
\node[font=\footnotesize] at (11.9,-5.0) {$G_2$};
\draw[->] (0,-4.75) -- (14.6,-4.75) node[right] {$t$};
\foreach \x/\l in {0/0,14/1}{\draw (\x,-4.82)--(\x,-4.68);}
\node[anchor=east,font=\footnotesize,align=right] at (-0.1,-0.4) {pass 1\\(concentrate)};
\foreach \a/\b/\n in {0.15/4.05/0,4.95/9.05/1,9.95/13.85/2}{
  \draw[corr,thick,rounded corners=2pt,fill=corr!12] (\a,-0.75) rectangle (\b,-0.05);
  \node[font=\footnotesize] at ({(\a+\b)/2},-0.4) {re-solve $G_{\n}$ with current $f$};}
\draw[->,kin,thick] (4.05,-0.4) -- (4.95,-0.4);
\draw[->,kin,thick] (9.05,-0.4) -- (9.95,-0.4);
\node[anchor=east,font=\footnotesize,align=right] at (-0.1,-1.45) {stress after\\pass 1};
\draw[gray] (0,-1.7) -- (14,-1.7);
\foreach \a in {4.2,9.2}{\fill[cipert!60] (\a,-1.7) rectangle (\a+0.6,-1.2);}
\node[font=\footnotesize,gray!60!black] at (2.1,-1.45) {$0$};
\node[font=\footnotesize,gray!60!black] at (7.0,-1.45) {$0$};
\node[font=\footnotesize,gray!60!black] at (11.9,-1.45) {$0$};
\node[anchor=east,font=\footnotesize,align=right] at (-0.1,-2.65) {pass 2\\(convex\\integration)};
\foreach \a in {4.2,9.2}{
  \draw[cipert,thick,decorate,decoration={zigzag,segment length=1.6mm,amplitude=0.9mm}] (\a+0.05,-2.65) -- (\a+0.55,-2.65);
  \node[cipert,font=\footnotesize] at (\a+0.3,-2.2) {CI};}
\foreach \a/\b/\n in {4.95/9.05/1,9.95/13.85/2}{
  \draw[corr,thick,rounded corners=2pt,fill=corr!12] (\a,-3.0) rectangle (\b,-2.3);
  \node[font=\footnotesize] at ({(\a+\b)/2},-2.65) {re-solve $G_{\n}$ again};}
\draw[->,kin,thick] (4.8,-2.65) -- (4.95,-2.65);
\draw[->,kin,thick] (9.8,-2.65) -- (9.95,-2.65);
\draw[->,kin,thick] (4.5,-3.0) .. controls (4.6,-3.75) and (5.2,-3.75) .. (5.35,-3.05);
\node[kin,font=\footnotesize,anchor=west] at (5.55,-3.4)
  {entry layer on $G_1$ depends on the CI step at $J_1$};
\node[font=\footnotesize,anchor=west] at (0.15,-4.3) {order of operations};
\draw[->,thick] (3.0,-4.3) -- (14,-4.3) node[midway,above,font=\footnotesize] {left to right in time};
\end{tikzpicture}}
\caption{(a) A perturbation $w$ of the velocity on a short interval $[a,b]$ changes the particle state, and hence the drag force,
at all later times. (b) The causal construction of Section~\ref{sec:proofW}. In pass~1, each good interval $G_i$ is re-solved with the
particle state produced by everything before it, and the stress is concentrated on the junctions $J_i$. In pass~2, a localized convex
integration step (CI) is applied at each junction, and the next good interval is re-solved again. The corrector on $G_1$ depends,
through the particles, on the CI perturbation at $J_1$, which is built to cancel the stress left by pass~1, so the two passes
cannot be merged.}
\label{fig:memory}
\end{figure}
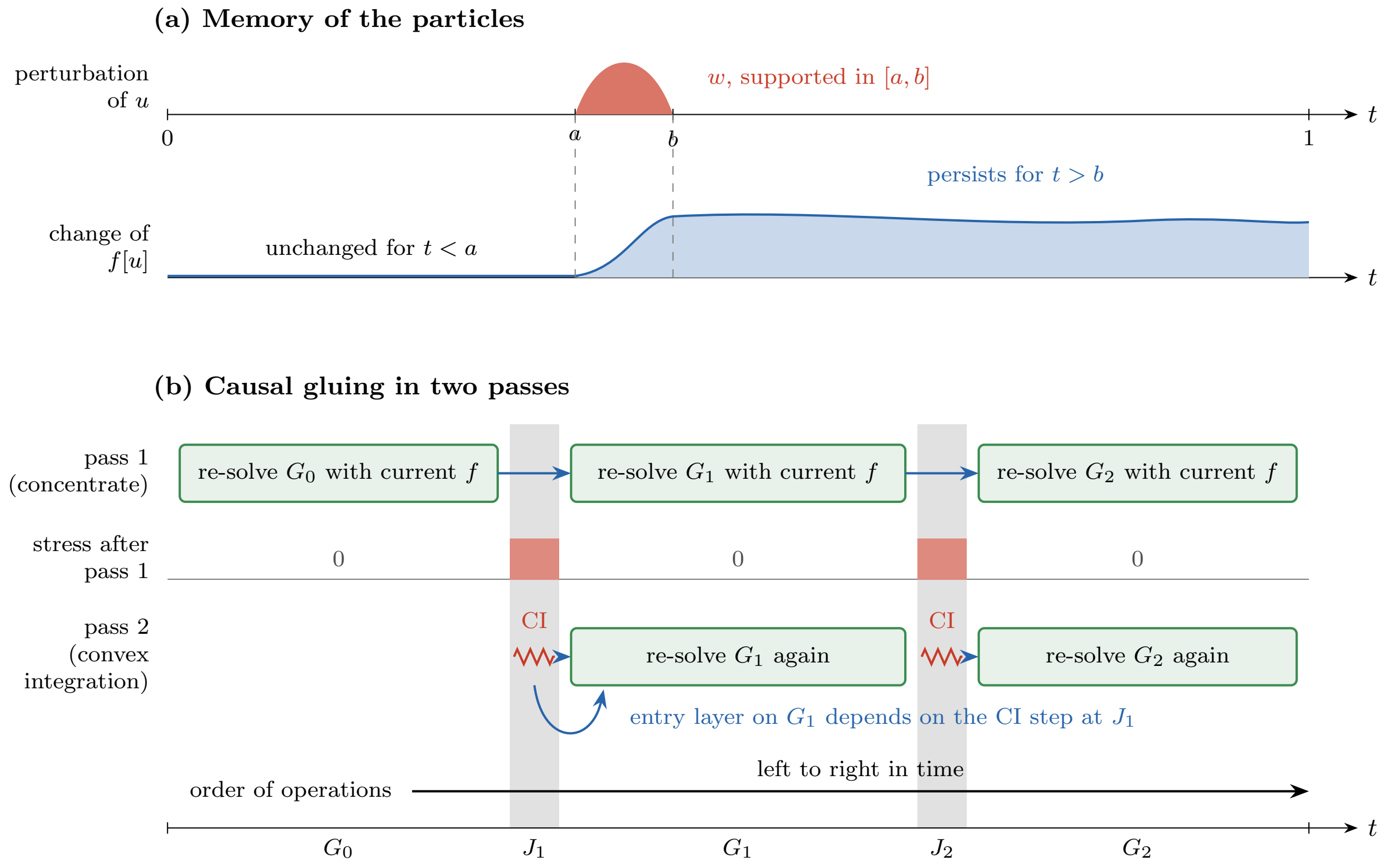

Because of the memory of the particles (Section~\ref{sec:memory}), the solution can no longer be left untouched outside the
bad set: on good intervals it must be \emph{re-solved} with the new kinetic state, and every gluing must be \emph{causal}, that is,
carried out from left to right in time. The regularity obtained on good intervals is that of a strong solution, in the sense of
Section~\ref{sec:sharp}, and not $C^\infty$ in general (see the comment after Theorem~\ref{thm:C}).
Figure~\ref{fig:memory} illustrates the memory effect and the resulting two-pass construction.

\subsection{Statement}

Since $\|\cdot\|_{W^{1,q}}\le\|\cdot\|_{W^{1,q'}}$ on $\TT^d$ for $q\le q'$, it suffices to prove Theorem~\ref{thm:regint} below for
large $q$. Throughout this section, $1\le p<2$, $2\le q<\infty$ with $q>d$, $\beta\in(0,1-\frac dq)$, $\varepsilon\in(0,1)$, and $r=r(d,p,q)$ is the
exponent of Appendix~\ref{sec:parameters}. All VNS--Reynolds solutions (Definition~\ref{def:VNSR}) carry the kinetic state $f[U]$,
the push-forward of $f_0$ along the whole history of $U$.

\paragraph{Notation.} In this section $\varsigma$ denotes the size of a base stress, $\beta\in(0,1-\frac dq)$ an interpolation
exponent, and $N_{\rm c}$ the number of components of $\tilde I$. As in Section~\ref{sec:kinetic}, radii in velocity are denoted by
$\Rv_*,\Rv_1$, and $R$ denotes a Reynolds stress.

\begin{defin}[Well-prepared]\label{def:wp}
A smooth VNS--Reynolds solution $(U,P,S)$ on $[0,1]$ is \emph{well-prepared} for a set $I$ and a scale $\tau>0$ if:
\begin{itemize}
\item $I$ is a union of at most $\tau^{-\varepsilon}$ pairwise disjoint closed intervals of length $5\tau$, with $0,1\notin I$;
\item $S(t)=0$ whenever $\dist(t,I^c)\le\tau$.
\end{itemize}
\end{defin}

\begin{theo}[Intervals of regularity]\label{thm:regint}
Let $d\ge2$, $1\le p<2$, $q<\infty$ and $\varepsilon>0$. Let $(u_{\rm in},f_0)$ be smooth data, with $f_0\ge0$ compactly supported, admitting
a smooth VNS solution on some $[0,t_1]$. Then there are infinitely many weak solutions with these data, all coinciding
with the smooth solution near $t=0$, such that:
\begin{itemize}
\item $u\in L^2_{t,x}\cap L^p_tL^\infty\cap L^1_tW^{1,q}$;
\item there is a closed set $\mathcal S\subset(0,1)$ with $\dim_{\mathcal H}\mathcal S\le\varepsilon$ such that, on every compact interval
$\mathcal K\subset(0,1]\setminus\mathcal S$ and for every $s<\infty$,
\[
u\in L^2(\mathcal K;W^{1,\infty})\cap L^2(\mathcal K;W^{2,s})\cap C(\mathcal K;W^{1,s}),
\]
and $f$ is transported on $\mathcal K$ by a Lipschitz phase-space flow. In particular, by Corollary~\ref{cor:roughtrace}, the
solution is the unique weak solution on $\mathcal K$ with velocity in $L^2(\mathcal K;L^\infty)$ and initial state $(u,f)(\inf \mathcal K)$.
\end{itemize}
\end{theo}

The theorem follows from iterating the following proposition.

\begin{prop}[Well-prepared iteration]\label{prop:W}
There are $M_W$ and $C_r$, depending only on $d$ and $r$, with the following property. Let $(u,p_u,R)$ be well-prepared for
$(\tilde I,\tilde\tau)$, and let $\delta>0$. Then there is $(u_1,p_1,R_1)$, well-prepared for some $(I,\tau)$ with $I\subset\tilde I$ and
$\tau\le\tilde\tau/2$, such that:
\begin{enumerate}[label=(W\arabic*)]
\item $\|R_1\|_{L^1(0,1;L^r)}\le\delta$;
\item $\|u_1-u\|_{L^2_{t,x}}\le M_W\|R\|_{L^1(0,1;L^r)}^{1/2}+\delta$;
\item $\|u_1-u\|_{L^pL^\infty}+\|u_1-u\|_{L^1W^{1,q}}\le\delta$;
\item $u_1=u$ on $[0,\inf\tilde I]$;
\item $\|u_1-u\|_{L^2(\tilde I^c;W^{1,\infty})}\le\delta$.
\end{enumerate}
\end{prop}

(W1)--(W4) are, up to the square root in (W2) (see Section~\ref{sec:NSdiff}) and the weaker support condition in (W4), contained in
\cite[Proposition~2.2]{CL22}. (W5) is new. Since the kinetic state on $\tilde I^c$ depends on the
modifications made in $\tilde I$, the velocity on $\tilde I^c$ must be corrected. (W5) states that this correction is small in the
Lipschitz norm, which is what makes the limit strong off $\mathcal S$.

\subsection{Two analytic tools}

For $s\in\RR$ let $J^s=(1-\Delta)^{s/2}$ and $H^{s,q}=J^{-s}L^q(\TT^d)$; recall that $H^{1,q}=W^{1,q}$ for $1<q<\infty$. Let $X_q=B^1_{q,2}$ and let
$\mathrm{MR}_q(I)$ be as in Lemma~\ref{lem:MR}. Set
\[
\Yz(I):=\mathrm{MR}_q(I),\qquad \Yo(I):=J\,\mathrm{MR}_q(I)=L^2(I;H^{1,q})\cap H^1(I;H^{-1,q}),\qquad \|v\|_{\Yo(I)}:=\|J^{-1}v\|_{\mathrm{MR}_q(I)} .
\]

\begin{lem}[Shifted maximal regularity]\label{lem:MRshift}
There is $C_q\ge1$, independent of $I\subset[0,1]$, such that:
\begin{enumerate}[label=(\alph*)]
\item Lemma~\ref{lem:MR}(a), (b) hold in $\Yo$, with data in $JX_q=B^0_{q,2}$ and forcing in $L^2(I;H^{-1,q})$;
\item $B^0_{q,2}\hookrightarrow L^q$, $X_q\hookrightarrow W^{1,q}\hookrightarrow L^\infty$, $\Yo\hookrightarrow L^2L^\infty$ and $\Yz\hookrightarrow L^2W^{1,\infty}$;
\item $\|v\|_{L^2H^{2-\beta,q}}\le C\|v\|_{L^2H^{2,q}}^{1-\beta}\|v\|_{L^2H^{1,q}}^\beta$ and $H^{2-\beta,q}\hookrightarrow W^{1,\infty}$;
\item $\|v(t)\|_{B^{1-\beta}_{q,2}}\le C\|v(t)\|_{X_q}^{1-\beta}\|v(t)\|_{B^0_{q,2}}^\beta$ and $B^{1-\beta}_{q,2}\hookrightarrow L^\infty$;
\item for $1<s<\infty$, the antidivergence satisfies $\|\Rs g\|_{L^s}\le C_s\|g\|_{H^{-1,s}}$ and $\|\Rs g\|_{W^{1,s}}\le C_s\|g\|_{L^s}$.
\end{enumerate}
\end{lem}

\begin{proof}
(a): conjugate Lemma~\ref{lem:MR} by $J^{-1}$, which commutes with $\partial_t$, $\Delta$ and the heat semigroup.
(b), (c), (d): Sobolev and Besov embeddings and real or complex interpolation \cite{Tri83}, using $q>d$ and $1-\beta>\frac dq$.
(e): by \cite[Definition~B.2]{CL22}, $\Rs$ is a linear combination of $\nabla\Delta^{-1}$ and $\nabla^3\Delta^{-2}$, applied to
mean-free functions. These are Fourier multipliers of order $-1$, so $\nabla\Rs$ is a Calder\'on--Zygmund operator and
$\Rs J$ is bounded on $L^s$.
\end{proof}

The next lemma sharpens Lemma~\ref{lem:dual}: the perturbation is measured in $L^1_tL^s_x$ for any $s\in[1,\infty]$ instead of
$L^1_tL^\infty_x$, and the bound is linear in the perturbation.

\begin{lem}[Kinetic duality in $(W^{1,1})^*$]\label{lem:B2}
Let $0\le t_0<t\le1$. Let $b,U$ be fields on $[t_0,1]\times\TT^d$, Lipschitz in $x$ for a.e.\ time, with
$\nabla b,\nabla U\in L^1(t_0,1;L^\infty)$ and $b,U\in L^1(t_0,1;L^\infty)$. Set $\Lambda_b:=\int_{t_0}^t\|\nabla b\|_\infty$. Let $f_*\ge0$ be bounded
and vanish for $|v|>\Rv_*$, and let $f^U$, $f^b$ be its push-forwards from $t_0$ along $b_U$ and $b_b$. Put
\[
\Rv_1:=\Rv_*+\|U\|_{L^1L^\infty}+\|b\|_{L^1L^\infty},\qquad K:=e^{1+\Lambda_b+d}\|f_*\|_\infty|B_{\Rv_1}| .
\]
\begin{enumerate}[label=(\alph*)]
\item Let $1\le s\le\infty$, $\psi\in C^1(\TT^d\times\RR^d)$ and $0\le k\in L^{s'}(\TT^d)$ with
$\max\{|\nabla_x\psi(x,v)|,|\nabla_v\psi(x,v)|\}\le k(x)$ for a.e.\ $x$ and all $|v|\le \Rv_1$. Then
\[
\Big|\iint(f^U-f^b)(t)\psi\Big|\le K\|k\|_{L^{s'}(\TT^d)}\int_{t_0}^t\|U-b\|_{L^s}\dd\tau .
\]
\item Let $d<q<\infty$, $q'<s\le q$, $W\in W^{1,q}(\TT^d)$, and $h:=\int(v-W)(f^U-f^b)(t)\dd v$. Then
\[
\|h\|_{H^{-1,s}}\le C_{s,q}K\big(1+\Rv_1+\|W\|_{W^{1,q}}\big)\int_{t_0}^t\|U-b\|_{L^s}\dd\tau .
\]
\end{enumerate}
\end{lem}

\begin{remark}
The restriction $s>q'$ in (b) is needed: the term $|\nabla W||\varphi|$ below is in $L^{s'}$ only if $s'<q$. In particular (b) does
not apply with $s=r$, since $r'=d/\gamma>q$ by \eqref{eq:gamma}--\eqref{eq:r}. Below, (b) is only used with $s=q$.
\end{remark}

\begin{proof}
\proofstep{Step 1: smooth fields}
Assume first that $b$ and $U$ are smooth on $[t_0,1]\times\TT^d$; $f_*$ is only bounded, so we argue in Lagrangian form. Let
$\Phi(\tau,y):=\psi(Z_b(t;\tau,y))$, which is $C^1$ and satisfies $\partial_\tau\Phi+b_b\cdot\nabla_y\Phi=0$. For fixed $\zeta$ the function
$\gamma(\tau):=\Phi\big(\tau,Z_U(\tau;t_0,\zeta)\big)$ is $C^1$ on $[t_0,t]$, with $\gamma(t_0)=\psi(Z_b(t;t_0,\zeta))$, $\gamma(t)=\psi(Z_U(t;t_0,\zeta))$ and
\[
\gamma'(\tau)=\big(\partial_\tau\Phi+b_U\cdot\nabla_y\Phi\big)(\tau,Z_U)=\big((b_U-b_b)\cdot\nabla_y\Phi\big)(\tau,Z_U)=(U-b)(\tau,X_U)\cdot\nabla_v\Phi(\tau,Z_U),
\]
since $b_U-b_b=(0,U-b)$. Integrating in $\tau$, multiplying by $f_*(\zeta)$, integrating in $\zeta$, and using the push-forward identities
$\iint f^U(t)\psi=\int f_*\psi(Z_U(t;t_0,\cdot))$, $\iint f^b(t)\psi=\int f_*\psi(Z_b(t;t_0,\cdot))$ and $f^U(\tau)=Z_U(\tau;t_0,\cdot)_\#f_*$, we get
\begin{equation}\label{eq:B2id}
\iint(f^U-f^b)(t)\psi=\int_{t_0}^t\iint f^U(\tau,y)\,(U-b)(\tau,x)\cdot\nabla_v\Phi(\tau,y)\dd y\dd\tau,\qquad y=(x,v).
\end{equation}
No regularity of $f_*$ is used. If $y\in\Supp f^U(\tau)$, its velocity is at most $\Rv_*+\|U\|_{L^1L^\infty}$, and it grows by at most
$\|b\|_{L^1L^\infty}$ along $Z_b$, so $Z_b(t;\tau,y)\in\{|v|\le \Rv_1\}$. Step~3 of the proof of Lemma~\ref{lem:dual} then gives
$|\nabla_v\Phi(\tau,y)|\le e^{1+\Lambda_b}k(X_b(t;\tau,y))$ on $\Supp f^U(\tau)$. Let $\varrho:=e^d\|f_*\|_\infty|B_{\Rv_1}|$. Then
$\rho^U:=\int f^U\dd v\le\varrho$. The push-forward of the measure $f^U(\tau)\dd y$ by $Z_b(t;\tau,\cdot)$ has density
$e^{d(t-\tau)}f^U(\tau,Z_b(\tau;t,\cdot))\le e^d\|f_*\|_\infty$ and is supported in $\{|v|\le \Rv_1\}$. By H\"older's inequality in $y$,
\[
\begin{aligned}
\iint f^U|U-b|(x)\,k(X_b(t;\tau,y))\dd y&\le\Big(\iint f^U|U-b|^s\Big)^{1/s}\Big(\iint f^Uk(X_b)^{s'}\Big)^{1/s'}\\
&\le\varrho^{1/s}\|(U-b)(\tau)\|_{L^s}\cdot\varrho^{1/s'}\|k\|_{L^{s'}} ,
\end{aligned}
\]
with the obvious changes for $s\in\{1,\infty\}$. This proves (a) for smooth fields.

\proofstep{Step 2: general fields}
The chain rule along $\tau\mapsto Z_b(t;\tau,Z_U(\tau;t_0,\zeta))$ is not available directly for Carath\'eodory flows, so we pass to the
limit in the inequality rather than in \eqref{eq:B2id}. Extend $b,U$ by $0$ to $\tau<t_0$, and let $b_\varepsilon,U_\varepsilon$ be their
mollifications in $x$ and, causally, in time (averages over $(\tau-\varepsilon,\tau)$). Then $\int_{t_0}^t\|\nabla b_\varepsilon\|_\infty\le\Lambda_b$,
$\|b_\varepsilon\|_{L^1L^\infty}\le\|b\|_{L^1L^\infty}$, the same holds for $U_\varepsilon$, and $\int_{t_0}^t\|U_\varepsilon-b_\varepsilon\|_{L^s}\le\int_{t_0}^t\|U-b\|_{L^s}$ by Minkowski's
inequality. Hence Step~1 applies to $(b_\varepsilon,U_\varepsilon)$ with constants $\Rv_{1,\varepsilon}\le \Rv_1$ and $K_\varepsilon\le K$. Since $b(\tau)$ is Lipschitz,
$b\in L^1(t_0,1;C(\TT^d))$ as a Bochner function, so $b_\varepsilon\to b$ in $L^1(t_0,1;L^\infty)$, and likewise for $U$. Since the phase-space fields are Lipschitz in $z$ with constant
$1+\|\nabla b(\tau)\|_{L^\infty}$, Gronwall's inequality gives $|Z_{b_\varepsilon}(t;t_0,\zeta)-Z_b(t;t_0,\zeta)|\le e^{1+\Lambda_b}\int_{t_0}^t\|b_\varepsilon-b\|_{L^\infty}$, and
likewise for $U$; hence the flows from $t_0$ converge uniformly on $\Supp f_*$. Therefore
$\iint(f^{U_\varepsilon}-f^{b_\varepsilon})(t)\psi=\iint f_*(\zeta)[\psi(Z_{U_\varepsilon}(t;t_0,\zeta))-\psi(Z_{b_\varepsilon}(t;t_0,\zeta))]\dd\zeta$ converges to the corresponding
expression for $(b,U)$ by dominated convergence. This proves (a).

\proofstep{Step 3: part (b)}
First let $W$ be smooth. Take $\psi=(v-W)\cdot\varphi$ with $\varphi\in C^1$. Then $\nabla_v\psi=\varphi$ and $\nabla_x\psi=(\nabla\varphi)^T(v-W)-(\nabla W)^T\varphi$, so on
$\{|v|\le \Rv_1\}$ we may take $k=(1+\Rv_1+\|W\|_\infty)(|\varphi|+|\nabla\varphi|)+|\nabla W||\varphi|$. Since $q'<s\le q$, we have $s'<q$; define $p_*$ by
$\frac1{p_*}=\frac1{s'}-\frac1q\in(0,1)$. As $q>d$, $\frac1{p_*}>\frac1{s'}-\frac1d$, so $W^{1,s'}\hookrightarrow L^{p_*}$. Hence, using $W^{1,q}\hookrightarrow L^\infty$,
\[
\|k\|_{L^{s'}}\le(1+\Rv_1+\|W\|_\infty)\|\varphi\|_{W^{1,s'}}+\|\nabla W\|_{L^q}\|\varphi\|_{L^{p_*}}\le C_{s,q}\big(1+\Rv_1+\|W\|_{W^{1,q}}\big)\|\varphi\|_{W^{1,s'}} .
\]
Apply (a). For general $W\in W^{1,q}$, take smooth $W_n\to W$ in $W^{1,q}$, hence uniformly; since $f^U-f^b$ is bounded with compact
velocity support, $\int h_n\cdot\varphi\to\int h\cdot\varphi$. Conclude by density, since $H^{-1,s}=(W^{1,s'})^*$.
\end{proof}

\begin{remark}[Starting at $t_0>0$; the full drag difference]\label{rem:shift}
Lemmas~\ref{lem:dual} and~\ref{lem:mean} and Proposition~\ref{prop:coupling} hold when started at a time $t_0>0$ from a common kinetic
state $f_*$. More precisely, the statements hold on $[t_0,1]$, which has length at most one, with the following replacements: $f_0$ by $f_*$;
$M_0$ by $\|f_*\|_{L^1}=M_0$; $M_\infty$ by $\|f_*\|_\infty\le e^dM_\infty$; $\Rv_0$ by the radius $\Rv_0+\int_0^{t_0}\|U\|_{L^\infty}$ of the velocity support of
$f_*$; the integrals $\int_0^t$ and the norms $L^1_tL^\infty$, $L^1_tW^{1,\infty}$ by the same quantities on $[t_0,t]$ and $[t_0,1]$; and, in
Lemma~\ref{lem:mean}, \eqref{eq:fixed} by $c(t)=\fint U(t_0)-\fint U(t)+J(t_0)-J[U+c](t)$ for $t\ge t_0$, with $c(t_0)=0$, where $J(t_0)$ is the
common momentum of the particles at $t_0$.
In this section $f_*=f[U](t_0)$ for a smooth field $U$ and smooth $f_0$, so $f_*$ is smooth and the proofs apply verbatim after a
translation in time.
For a merely bounded $f_*$ the Eulerian computation in Step~2 of the proof of Lemma~\ref{lem:dual} is not available as written; one uses
instead the Lagrangian identity \eqref{eq:B2id}, which only involves the flows, or an approximation of $f_*$ in $L^1$ by smooth functions.

For two fields $U$, $b$ as in Lemma~\ref{lem:B2}, write $a:=U-b$ and
\[
h:=F[f^U,U]-F[f^b,b]=\int(v-b)(f^U-f^b)\dd v-\rho^U\,a .
\]
The first term is a \emph{history} term, controlled by Lemma~\ref{lem:B2}(b). The second is \emph{instantaneous}: it is of the size
of $a(t)$ itself, and is not controlled by the past of $a$. By Lemma~\ref{lem:MRshift}(e), for $d<q<\infty$,
\begin{equation}\label{eq:fulldrag}
\|\Rs h(t)\|_{L^q}\le C_qK\big(1+\Rv_1+\|b(t)\|_{W^{1,q}}\big)\int_{t_0}^t\|a\|_{L^q}\dd\tau+C_q\varrho\|a(t)\|_{L^q},\qquad\varrho=e^d\|f_*\|_\infty|B_{\Rv_1}|.
\end{equation}
This replaces the H\"older bound of Proposition~\ref{prop:coupling}. In the applications below, the instantaneous term is
controlled for two reasons. On the intervals where $a$ is spatially constant, it is bounded by the momentum defect, which is a
history term. Where $a$ contains a convex integration perturbation $w$, the term $\Rs(\rho^Uw)$ is small in $L^1_tL^q$, since
$\|w\|_{L^1_tL^\infty}$ is small, and it only enters the new stress on the support of $w$.
\end{remark}

\subsection{Segment correctors and mean-preserving gluing}\label{sec:segcorrsub}

Let $(U,P,S)$ be a smooth VNS--Reynolds solution on $[0,1]$, called the \emph{base}, and let $a\in[0,1)$ be the time at which a
segment (or, in Lemma~\ref{lem:CIloc}, a convex integration zone) starts. The size of the base is measured \emph{after $a$}:
\begin{equation}\label{eq:basenorms}
B_U:=\sup_{t\in[a,1]}\|U(t)\|_{W^{1,q}},\qquad \Lambda_U:=\int_a^1\|\nabla U\|_{L^\infty}\dd t,\qquad L_U:=\int_0^a\|U\|_{L^\infty}\dd t .
\end{equation}
Note that $\|U\|_{L^\infty((a,1)\times\TT^d)}\le C_EB_U$ by Lemma~\ref{lem:MRshift}(b). The Lipschitz norm of $U$ enters only through $\Lambda_U$, an
$L^1$ norm in time. Only the base on $[a,1]$ enters the kinetic comparisons, since they are started at $a$ from the common
kinetic state $f[U](a)$, whose velocity support has radius at most $\Rv_0+L_U$ by Lemma~\ref{lem:char}(c). This is used
in the second pass of Section~\ref{sec:proofW}, where the base contains earlier convex integration perturbations $w_i$, which are small in
$L^1_tW^{1,q}$ but whose $\sup_t\|w_i(t)\|_{W^{1,q}}$ depends on their frequency; they vanish on $[a,1]$.

A \emph{segment} is an interval $\Sigma=[a,b]\subset[0,1]$ with a cutoff $\chi\in C^\infty(\RR;[0,1])$ that vanishes near $a$ and either
vanishes near $b$ (\emph{two-sided}) or equals $1$ near $b=1$ (\emph{one-sided}). For a field $z$ on $\Sigma$, write $c_z:=\fint z$,
$\mathring z:=z-c_z$, and define the \emph{glued field} and its kinetic state
\[
\hat U_z:=U+\chi\mathring z+c_z\quad\text{on }\Sigma,\qquad \hat f_z:=f\big[U|_{[0,a]}\oplus\hat U_z\big].
\]
The cutoff acts only on the mean-zero part. The \emph{segment corrector} is the solution of
\begin{equation}\label{eq:segcorr}
\partial_tz-\Delta z+\QQ\D(U\otimes z+z\otimes U+z\otimes z)=\QQ\big[-\D S+g_z\big],\quad z(a)=0,\quad g_z:=F[\hat f_z,\hat U_z]-F[f[U],U].
\end{equation}
Where $\chi=1$ we have $\hat U_z=U+z$, and \eqref{eq:segcorr} says that $(U+z,\hat f_z)$ solves VNS exactly.

Put $\varsigma:=\|S\|_{L^\infty(\Sigma;L^q)}$ and $N_S:=\|S\|_{L^\infty(\Sigma;W^{1,q})}$. Let $M_\infty=\|f_0\|_\infty$, let $\Rv_0$ bound the velocity support
of $f_0$, and set
\[
\Rv_1:=\Rv_0+L_U+2C_EB_U+8,\qquad \varrho:=e^dM_\infty|B_{\Rv_1}|,\qquad C_h:=2(\Rv_1+C_EB_U)\varrho .
\]
The term $L_U$ accounts for the velocity that the particles gain before the segment starts. Every glued field $\hat U$
considered below satisfies $\|\hat U-U\|_{L^1(a,1;L^\infty)}\le4$: the contribution of $\Sigma$ is at most $3$, by \eqref{eq:smallLip} in (a)
and by \eqref{eq:adm} in (b), and after an exit by the mean the contribution of $[b,1]$ is at most $1$, by the continuity argument in
the proof of (c). Hence $\|\hat U\|_{L^1(a,1;L^\infty)}\le C_EB_U+4$. By Lemma~\ref{lem:char}(c), all kinetic states are supported in
$\{|v|\le \Rv_0+L_U+C_EB_U+4\}\subset\{|v|\le \Rv_1\}$, and the radius of Lemma~\ref{lem:B2}, applied with two such fields (or with
one such field and $U$) in the roles of $U$ and $b$, is at most $\Rv_0+L_U+2(C_EB_U+4)=\Rv_1$.

\begin{prop}[Segment correctors]\label{prop:segcorr}
In this proposition, the constants $C$ and $K$ depend only on $d,q,r,\beta,\Rv_0,M_\infty$ and on upper bounds for $B_U$,
$\Lambda_U$ and $L_U$ (see \eqref{eq:basenorms}). Primed constants, $C'$, $\ell_0'$ and $c_*'$, depend in addition on an upper bound for $N_S$.
\begin{enumerate}[label=(\alph*)]
\item (Short segments.) If $|\Sigma|\le\ell_0'$, then \eqref{eq:segcorr} has a solution $z\in\Yz(\Sigma)$, unique among solutions with
$\|z\|_{\Yz(\Sigma)}\le1$, and for $t\in\Sigma$:
\begin{enumerate}[label=(\roman*)]
\item $\|z\|_{\Yz(\Sigma)}+\sup_\Sigma\|z\|_{X_q}\le C'|\Sigma|^{1/2}$;
\item $\|z\|_{\Yo(\Sigma)}+\sup_\Sigma\|z\|_{L^q}\le C|\Sigma|^{1/2}\varsigma$;
\item $\|z\|_{L^2(\Sigma;W^{1,\infty})}+\sup_\Sigma\|z\|_{L^\infty}\le C'|\Sigma|^{1/2}\varsigma^\beta$;
\item $\|\Rs z(t)\|_{L^r}\le C_r\int_a^t\|S\|_{L^r}+C'|\Sigma|^{3/2}\varsigma$.
\end{enumerate}
\item (Long segments.) There are $\varepsilon_*'>0$ and $c_*'\in(0,\varepsilon_*']$ such that, if $\varsigma\le c_*'$, then \eqref{eq:segcorr} has a solution
$z\in\Yz(\Sigma)$, unique among solutions with $\|z\|_{\Yo(\Sigma)}\le\varepsilon_*'$ and $\|z\|_{L^1(\Sigma;W^{1,\infty})}\le1$, and (i)--(iv) hold with
$|\Sigma|^{1/2}$ and $|\Sigma|^{3/2}$ replaced by $1$; the constants, $\varepsilon_*'$ and $c_*'$ do not depend on $|\Sigma|\le1$.
\item (Momentum and mean.) In both cases, $c_z'=\fint g_z$ and $c_z=J[f[U]]-J[\hat f_z]$ on $\Sigma$. Write $a_z:=\hat U_z-U=\chi\mathring z+c_z$.
For $t\in\Sigma$,
\begin{equation}\label{eq:cg}
|c_z(t)|\le K\int_a^t\|a_z\|_{L^q},\qquad \|\Rs g_z(t)\|_{L^q}\le K\Big(\int_a^t\|a_z\|_{L^q}+\|a_z(t)\|_{L^q}\Big),
\end{equation}
so that $\sup_\Sigma\big(|c_z|+\|\Rs g_z\|_{L^q}\big)\le9K\sup_\Sigma\|z\|_{L^q}$. If $\chi$ is two-sided and the glued field exits by the mean
(Proposition~\ref{prop:glue}), then $a_z=c_z$ on $[b,1]$, \eqref{eq:cg} holds for $t\in[b,1]$, and
\begin{equation}\label{eq:persist}
\sup_{[b,1]}\big(|c_z|+\|\Rs g_z\|_{L^q}\big)\le C|\Sigma|\sup_\Sigma\|z\|_{L^q} .
\end{equation}
\item (Persistent stress.) In the situation of \eqref{eq:persist}, the stress $\Pi_z:=U\ot c_z+c_z\ot U-\Rs g_z$ satisfies
\[
\sup_{[b,1]}\|\Pi_z\|_{L^q}\le C|\Sigma|\sup_\Sigma\|z\|_{L^q},\qquad \sup_{[b,1]}\|\Pi_z\|_{W^{1,q}}\le C .
\]
\end{enumerate}
In both cases $z$ is smooth on $[a,b]\times\TT^d$, with one-sided derivatives at $t=a$ and $t=b$.
\end{prop}

Estimate (ii) gains a factor $\varsigma$ in the weaker norm $\Yo$; the factor $|\Sigma|^{1/2}$ in (i)--(iv) makes the
bounds summable over many short segments. The persistent stress in (d) is bounded in $W^{1,q}$ by a constant because
$\Rs$ gains one derivative; no derivative of the kinetic state is needed.

\begin{proof}
\proofstep{Kinetic bounds}
For $\zeta\in\Yz(\Sigma)$ with $\zeta(a)=0$ and $\|\zeta\|_{\Yz}\le1$, the glued field $\hat U_\zeta$ is Lipschitz in $x$ with an integrable
constant, since $\Yz\hookrightarrow L^2W^{1,\infty}$. Its kinetic state is defined by the Carath\'eodory flow, as in
Proposition~\ref{prop:critreg}(c). We always take $\ell_0'\le C_E^{-2}$, so that for $|\Sigma|\le\ell_0'$
\begin{equation}\label{eq:smallLip}
\|\zeta\|_{L^1(\Sigma;W^{1,\infty})}\le|\Sigma|^{1/2}\|\zeta\|_{L^2W^{1,\infty}}\le C_E|\Sigma|^{1/2}\|\zeta\|_{\Yz}\le1 .
\end{equation}
In (b) the bound $\|\zeta\|_{L^1W^{1,\infty}}\le1$ is part of the admissibility condition \eqref{eq:adm} below and is
verified there. Hence $\|\chi\mathring\zeta+c_\zeta\|_{L^1L^\infty}\le3\|\zeta\|_{L^1L^\infty}\le3$, so by Lemma~\ref{lem:char}(c),
$\hat f_\zeta\le e^dM_\infty$ is supported in $\{|v|\le \Rv_1\}$ and $\rho_{\hat f_\zeta}\le\varrho$. Write
\[
g_\zeta=h_\zeta-\rho_{\hat f_\zeta}a_\zeta,\qquad h_\zeta:=\int(v-U)(\hat f_\zeta-f[U])\dd v,\qquad a_\zeta:=\chi\mathring\zeta+c_\zeta .
\]
Then $|h_\zeta|\le C_h$ and $\|a_\zeta(t)\|_{L^s}\le3\|\zeta(t)\|_{L^s}$ for every $s$. The kinetic states agree at $a$. Lemma~\ref{lem:B2}(b), with
$t_0=a$, base $U$, $W=U(t)$ and $s=q$, gives
\begin{equation}\label{eq:wkB}
\|h_\zeta(t)\|_{H^{-1,q}}\le K_1\int_a^t\|a_\zeta\|_{L^q}\le3K_1\int_a^t\|\zeta\|_{L^q}.
\end{equation}
Here $K_1$ depends on $B_U$ through $\Rv_1$ and $\|U(t)\|_{W^{1,q}}$, and on $\Lambda_U$. For two candidates $\zeta,\zeta'$ that coincide up to
a time $s_0$,
\[
g_\zeta-g_{\zeta'}=\int(v-\hat U_{\zeta'})(\hat f_\zeta-\hat f_{\zeta'})\dd v-\rho_{\hat f_\zeta}(a_\zeta-a_{\zeta'}).
\]
Lemma~\ref{lem:B2}(b), with $t_0=s_0$, base $\hat U_{\zeta'}$, whose Lipschitz norm is at most $\Lambda_U+3$ in $L^1_t$ by \eqref{eq:smallLip}, gives
\begin{equation}\label{eq:difB}
\|(g_\zeta-g_{\zeta'})(t)\|_{H^{-1,q}}\le K_2\big(1+\|\zeta'(t)\|_{W^{1,q}}\big)\int_{s_0}^t\|\zeta-\zeta'\|_{L^q}+3\varrho\|(\zeta-\zeta')(t)\|_{L^q}.
\end{equation}

\proofstep{The map}
For an interval $I=[s_0,s_1]\subset\Sigma$ and a given solution on $[a,s_0]$, let $\Theta\zeta$ solve, on $I$,
$\partial_t\vartheta-\Delta\vartheta=\QQ[-\D S+g_\zeta-\D(U\otimes\zeta+\zeta\otimes U+\zeta\otimes\zeta)]$, with the given value at $s_0$. The forcing lies in
$L^2(I;L^q)$, by the bound $|h|\le C_h$. Let $C_T$ and $C_E$ denote the trace and embedding constants of Lemma~\ref{lem:MRshift}. We
use repeatedly the trace estimates of Lemmas~\ref{lem:MR}(b) and~\ref{lem:MRshift}(a) on $I=[s_0,s_1]$,
\[
\sup_I\|\zeta\|_{X_q}\le C_T\big(\|\zeta(s_0)\|_{X_q}+\|\zeta\|_{\Yz(I)}\big),\qquad \sup_I\|\zeta\|_{B^0_{q,2}}\le C_T\big(\|\zeta(s_0)\|_{B^0_{q,2}}+\|\zeta\|_{\Yo(I)}\big),
\]
and $\|\cdot\|_{L^q}\le C\|\cdot\|_{B^0_{q,2}}$, and
\[
\|\D(U\otimes\zeta+\zeta\otimes U)\|_{L^q}\le C_EB_U\|\zeta\|_{W^{1,q}}+\|\nabla U\|_{L^q}\|\zeta\|_{L^\infty}\le2C_EB_U\|\zeta\|_{X_q}.
\]
The initial traces cannot be dropped when $s_0>a$, even though $\zeta(a)=0$. They vanish in case~(a), where $I=\Sigma$ and $s_0=a$, and
for the difference of two candidates with the same value at $s_0$; in case~(b) they are bounded by the a priori bounds
\eqref{eq:aprioriY1} and \eqref{eq:aprioriY0} at $s_0$.

\proofstep{(a) Short segments}
Take $I=\Sigma$ and $E:=\{\zeta\in\Yz(\Sigma):\zeta(a)=0,\ \|\zeta\|_{\Yz}\le A'|\Sigma|^{1/2}\}$, with $A':=2C_q(CN_S+C_h)$ and the metric of $\Yo$. By
Lemma~\ref{lem:MR},
\[
\|\Theta\zeta\|_{\Yz}\le C_q\Big(|\Sigma|^{1/2}(CN_S+C_h)+|\Sigma|^{1/2}(3\varrho+2C_EB_U)C_TA'|\Sigma|^{1/2}+|\Sigma|^{1/2}C_EC_T^2(A'|\Sigma|^{1/2})^2\Big)\le A'|\Sigma|^{1/2}
\]
if $|\Sigma|\le\ell_0'$ is small. By Lemma~\ref{lem:MRshift}(a), \eqref{eq:difB} and $\|\D(\zeta\otimes\zeta-\zeta'\otimes\zeta')\|_{H^{-1,q}}\le C(\|\zeta\|_\infty+\|\zeta'\|_\infty)\|\zeta-\zeta'\|_{L^q}$,
\[
\begin{aligned}
\|\Theta\zeta-\Theta\zeta'\|_{\Yo}&\le C_q\Big(K_2\big(1+C_TA'\big)|\Sigma|^{3/2}C_T+(3\varrho+2C_EB_U)|\Sigma|^{1/2}C_T+2CC_EC_T^2A'|\Sigma|\Big)\|\zeta-\zeta'\|_{\Yo}\\
&\le\tfrac12\|\zeta-\zeta'\|_{\Yo}
\end{aligned}
\]
for $|\Sigma|\le\ell_0'$. The set $E$ is complete for the $\Yo$ metric: a $\Yo$-limit of elements of $E$ is a weak limit in the reflexive space
$\Yz$, so the bound and the initial value pass to the limit. The fixed point $z$ exists and is unique in $E$. For uniqueness on the larger set $\{\|\zeta\|_{\Yz(\Sigma)}\le1\}$, let $\tilde z$ be a second solution
with $\|\tilde z\|_{\Yz(\Sigma)}\le1$. The kinetic bounds above and \eqref{eq:difB} hold for $\tilde z$ by \eqref{eq:smallLip}. In the contraction
estimate the bound $\|\zeta\|_{\Yz}\le A'|\Sigma|^{1/2}$ entered in two places. In the kinetic term it entered only through
$1+\|\zeta'(t)\|_{W^{1,q}}\le1+C_TA'$, and $1+C_T$ suffices there, since the factor $|\Sigma|^{3/2}$ is explicit. In the quadratic term it
entered through $\|\zeta\|_{L^2(\Sigma;L^\infty)}$, and merely replacing $A'|\Sigma|^{1/2}$ by $1$ would lose the small factor. Instead, since
$\zeta(a)=0$, Lemma~\ref{lem:MR}(b) and $X_q\hookrightarrow L^\infty$ give
\[
\|\zeta\|_{L^2(\Sigma;L^\infty)}\le|\Sigma|^{1/2}\sup_\Sigma\|\zeta\|_{L^\infty}\le C_EC_T|\Sigma|^{1/2}\|\zeta\|_{\Yz(\Sigma)}\le C_EC_T|\Sigma|^{1/2},
\]
and the same for $\zeta'$. Hence, for $z=\Theta z$ and $\tilde z=\Theta\tilde z$,
\[
\|z-\tilde z\|_{\Yo}\le C_q\Big(K_2(1+C_T)|\Sigma|^{3/2}C_T+(3\varrho+2C_EB_U)|\Sigma|^{1/2}C_T+2CC_E C_T^2|\Sigma|^{1/2}\Big)\|z-\tilde z\|_{\Yo}
\le\tfrac12\|z-\tilde z\|_{\Yo}
\]
after decreasing $\ell_0'$, so $\tilde z=z$. This proves (i).
The gain $|\Sigma|^{1/2}$ does not require $\varsigma$ or $N_S$ to be small: every forcing term is bounded in $L^\infty(\Sigma;L^q)$, and the
factor comes from the $L^2(\Sigma)$ norm in time.

For (ii), estimate $z=\Theta z$ in $\Yo$ with zero data. By Lemma~\ref{lem:MRshift}(a), \eqref{eq:wkB} with $s=q$, and
$\|\D S\|_{H^{-1,q}}\le C\|S\|_{L^q}$,
\[
\|z\|_{\Yo}\le C_q\Big(C|\Sigma|^{1/2}\varsigma+3K_1|\Sigma|^{3/2}C_T\|z\|_{\Yo}+(3\varrho+2C_EB_U)|\Sigma|^{1/2}C_T\|z\|_{\Yo}+C_EC_TA'|\Sigma|^{1/2}C_T\|z\|_{\Yo}\Big),
\]
where the last term comes from $\|\D(z\otimes z)\|_{H^{-1,q}}\le\|z\|_\infty\|z\|_{L^q}$. For $|\Sigma|\le\ell_0'$ the terms containing $\|z\|_{\Yo}$ are
absorbed, and $\|z\|_{\Yo}\le2C_qC|\Sigma|^{1/2}\varsigma$.

(iii) follows from (i), (ii) and Lemma~\ref{lem:MRshift}(c), (d): $\|z\|_{L^2W^{1,\infty}}\le C\|z\|_{\Yz}^{1-\beta}\|z\|_{\Yo}^\beta$ and
$\|z(t)\|_{L^\infty}\le C\|z(t)\|_{X_q}^{1-\beta}\|z(t)\|_{B^0_{q,2}}^\beta$.

For (iv): since $\Rs$ and $\QQ$ are Fourier multipliers and the heat semigroup is an $L^r$ contraction,
$\|\Rs z(t)\|_{L^r}\le\int_a^t\|\Rs\QQ(\cdots)\|_{L^r}$. The integrand is bounded as follows, using Lemma~\ref{lem:MRshift}(e) and $r\le2\le q$:
\begin{itemize}
\item $\|\Rs\QQ\D S\|_{L^r}\le C_r\|S\|_{L^r}$;
\item $\|\Rs\QQ h_z\|_{L^r}\le\|\Rs\QQ h_z\|_{L^q}\le C\|h_z\|_{H^{-1,q}}\le3CK_1|\Sigma|\sup_\Sigma\|z\|_{L^q}$, by $|\TT^d|=1$, $r\le q$,
Lemma~\ref{lem:MRshift}(e) and \eqref{eq:wkB} (Lemma~\ref{lem:B2}(b) is not available with $s=r$);
\item $\|\Rs\QQ(\rho a_z)\|_{L^r}\le3C\varrho\sup_\Sigma\|z\|_{L^q}$;
\item the terms with $U\otimes z$ and $z\otimes z$ are at most $C(C_EB_U+\sup\|z\|_\infty)\sup_\Sigma\|z\|_{L^q}$.
\end{itemize}
Integrating over $[a,t]$ and using (ii) gives (iv). The constant is primed because of the factor $\sup\|z\|_\infty\le C'|\Sigma|^{1/2}\varsigma^\beta$
from (iii).

\proofstep{(b) Long segments}
Here $|\Sigma|$ is not small, so the forcing $h$ is not small in $L^2L^q$ and the nonlinear term is controlled by the
smallness of $\varsigma$. The argument is the one of (a) on pieces, with a priori bounds between the pieces. Call a solution $z$ of
\eqref{eq:segcorr} on $[a,t']\subset\Sigma$ \emph{admissible} if
\begin{equation}\label{eq:adm}
\|z\|_{\Yo(a,t')}\le\varepsilon_0,\qquad \|z\|_{L^1(a,t';W^{1,\infty})}\le1,
\end{equation}
where $\varepsilon_0\in(0,1]$ is fixed below. For admissible $z$ the kinetic bounds of the first paragraph hold on $[a,t']$: the glued field has
$L^1L^\infty$ norm at most $C_EB_U+3$ and Lipschitz norm at most $\Lambda_U+3$ in $L^1_t$. So \eqref{eq:wkB} and \eqref{eq:difB} hold with the
constants $K_1,K_2$ of (a), and $|h_z|\le C_h$. The constants are chosen in the following order; each depends only on the
quantities listed in the statement and on those chosen before it:
\[
\varepsilon_0,\ \ell_1\ \longrightarrow\ C_B\ \longrightarrow\ A,\ A_0\ \longrightarrow\ A_2,\ \ell_2\ \longrightarrow\ c_*'\le\varepsilon_*':=\varepsilon_0 .
\]
The parameters $\varepsilon_0$ and $\ell_1$ are fixed once, at the beginning, so that both absorption conditions below (in the
a priori bounds in $\Yo$ and in $\Yz$) hold; both involve only unprimed constants. In particular $\ell_1$, $N_1$ and $C_B$ are
unprimed and are not changed afterwards; $N_S$ enters from $A$ on.

\emph{A priori bound in $\Yo$.} Let $z$ be admissible on $[a,t']$, and cover $[a,t']$ by consecutive intervals $I_j=[s_j,s_{j+1}]$, $0\le j<N_1$, of
length at most $\ell_1$, with $N_1\le1/\ell_1+1$. On $I_j$ we use Lemma~\ref{lem:MRshift}(a) with data $z(s_j)$. Since $z(a)=0$,
$\sup_{[a,s]}\|z\|_{L^q}\le C_T\|z\|_{\Yo(a,s)}$; on $I_j$, $\sup_{I_j}\|z\|_{L^q}\le C_T(\|z(s_j)\|_{B^0_{q,2}}+\|z\|_{\Yo(I_j)})$; and
$\|\D(z\otimes z)\|_{L^2(I_j;H^{-1,q})}\le C\|z\|_{L^2(I_j;L^\infty)}\sup_{I_j}\|z\|_{L^q}\le CC_E\varepsilon_0\sup_{I_j}\|z\|_{L^q}$. Hence
\[
\begin{aligned}
\|z\|_{\Yo(I_j)}\le C_q\Big(&\|z(s_j)\|_{B^0_{q,2}}+\ell_1^{1/2}\big(C\varsigma+3K_1C_T\|z\|_{\Yo(a,s_{j+1})}\big)\\
&+\big((3\varrho+2C_EB_U)\ell_1^{1/2}+CC_E\varepsilon_0\big)C_T\big(\|z(s_j)\|_{B^0_{q,2}}+\|z\|_{\Yo(I_j)}\big)\Big).
\end{aligned}
\]
We choose $\varepsilon_0$ and then $\ell_1$ so small that the coefficients of $\|z\|_{\Yo(I_j)}$ on the right are at most $\frac1{4}$ in total,
and that $C_q\ell_1^{1/2}(3\varrho+2C_EB_U)C_T\le\frac14$, which is the absorption condition used in the a priori bound in $\Yz$ below. With
$S_j:=\|z(s_j)\|_{B^0_{q,2}}+\|z\|_{\Yo(a,s_j)}$, $S_0=0$, and $\|z(s_{j+1})\|_{B^0_{q,2}}\le C_T(\|z(s_j)\|_{B^0_{q,2}}+\|z\|_{\Yo(I_j)})$, this gives
$S_{j+1}\le C''(S_j+\varsigma)$, with $C''$ unprimed. Hence
\begin{equation}\label{eq:aprioriY1}
\|z\|_{\Yo(a,t')}+\sup_{[a,t']}\|z\|_{B^0_{q,2}}\le C_B\varsigma,\qquad C_B:=(2C'')^{N_1}.
\end{equation}

\emph{A priori bound in $\Yz$.} Let $z$ be admissible on $[a,t']$, and assume in addition
\begin{equation}\label{eq:boot}
\sup_{[a,t']}\|z\|_{L^\infty}\le\tfrac1{8C_q}.
\end{equation}
On $I_j$, Lemma~\ref{lem:MR}(a), (b), the bounds $\|\D S\|_{L^q}\le CN_S$, $\|g_z\|_{L^q}\le C_h+3\varrho\|z\|_{L^q}$,
$\|\D(U\otimes z+z\otimes U)\|_{L^q}\le2C_EB_U\|z\|_{X_q}$ and $\|\D(z\otimes z)\|_{L^q}\le2\|z\|_\infty\|\nabla z\|_{L^q}$ give
\[
\|z\|_{\Yz(I_j)}\le C_q\Big(\|z(s_j)\|_{X_q}+\ell_1^{1/2}\big(CN_S+C_h+(3\varrho+2C_EB_U)C_T(\|z(s_j)\|_{X_q}+\|z\|_{\Yz(I_j)})\big)\Big)+\tfrac14\|z\|_{\Yz(I_j)} .
\]
By the choice of $\ell_1$, the coefficient of $\|z\|_{\Yz(I_j)}$ coming from the linear terms is at most $\frac14$, so it is absorbed
together with the last term, and $X_{j+1}:=\sup_{[a,s_{j+1}]}\|z\|_{X_q}$ satisfies $X_{j+1}\le C'''(X_j+1)$, $X_0=0$, where $C'''$
depends on $N_S$. Hence, with $A:=(2C''')^{N_1}$ and $A_0:=2C_q N_1(C_TA+CN_S+C_h+1)$,
\begin{equation}\label{eq:aprioriY0}
\sup_{[a,t']}\|z\|_{X_q}\le A,\qquad \|z\|_{\Yz(a,t')}\le A_0 .
\end{equation}
By Lemma~\ref{lem:MRshift}(c), (d), applied on each $I_j$, and by \eqref{eq:aprioriY1}--\eqref{eq:aprioriY0},
\begin{equation}\label{eq:aprioriLip}
\sup_{[a,t']}\|z\|_{L^\infty}\le CA^{1-\beta}(C_B\varsigma)^\beta,\qquad
\|z\|_{L^1(a,t';W^{1,\infty})}\le\|z\|_{L^2(a,t';W^{1,\infty})}\le CN_1^{1/2}A_0^{1-\beta}(C_B\varsigma)^\beta .
\end{equation}
We require of $c_*'$ that, for $\varsigma\le c_*'$, the right-hand sides of \eqref{eq:aprioriLip} are at most $\frac1{16C_q}$ and $\frac12$, and that
$C_qC_B\varsigma\le\varepsilon_0/4$; in particular $C_B\varsigma\le\varepsilon_0/2$, since $C_q\ge1$.
Then \eqref{eq:boot} is not needed as a hypothesis: $t\mapsto\sup_{[a,t]}\|z\|_\infty$ is continuous, since
$z\in C([a,t'];X_q)$ and $X_q\hookrightarrow L^\infty$, and it vanishes at $t=a$; by \eqref{eq:aprioriLip} applied on $[a,t]$ it can never reach
$\frac1{8C_q}$. Consequently, for $\varsigma\le c_*'$, every admissible solution on $[a,t']$ satisfies the \emph{strict} bounds
\begin{equation}\label{eq:strict}
\|z\|_{\Yo(a,t')}\le\tfrac{\varepsilon_0}2,\qquad \|z\|_{L^1(a,t';W^{1,\infty})}\le\tfrac12,\qquad \sup_{[a,t']}\|z\|_{X_q}\le A,\qquad \|z\|_{\Yz(a,t')}\le A_0 .
\end{equation}

\emph{Existence.} Let $A_2:=2C_q(A+1)+A_0$, and fix a partition $a=s_0<s_1<\dots<s_{N_2}=b$ of mesh $\ell_2\le\ell_1$. Suppose an admissible solution $z$
is constructed on $[a,s_i]$ (for $i=0$, $z(a)=0$). On $I=[s_i,s_{i+1}]$ consider
\[
E_I:=\big\{\zeta\in\Yz(I):\ \zeta(s_i)=z(s_i),\ \|\zeta\|_{\Yz(I)}\le A_2,\ \|\zeta\|_{\Yo(I)}\le\tfrac{\varepsilon_0}2\big\},
\]
with the metric of $\Yo(I)$, each $\zeta$ being extended by $z$ on $[a,s_i]$. For $\zeta\in E_I$ the extension is admissible on $[a,s_{i+1}]$, because
\eqref{eq:strict} on $[a,s_i]$ and $\|\zeta\|_{L^1(I;W^{1,\infty})}\le\ell_2^{1/2}C_EA_2\le\frac12$ for $\ell_2$ small; so the kinetic bounds are available.
Let $\Theta$ be the map of (a) on $I$ with data $z(s_i)$. For $\zeta\in E_I$ the trace estimates include the initial value: by
\eqref{eq:strict} at $s_i$, $\sup_I\|\zeta\|_{X_q}\le C_T(A+A_2)$ and $\sup_I\|\zeta\|_{L^q}\le C(C_B\varsigma+\varepsilon_0)$, and these bounds enter the
constant $C(A_2)$ below, which therefore also depends on $A$; this is admissible, since $A$ is fixed before $A_2$. In the contraction
estimate the difference of two candidates has zero value at $s_i$, so the zero-trace form applies to it, while the factor
$1+\|\zeta'(t)\|_{W^{1,q}}$ in \eqref{eq:difB} is bounded by $1+C_T(A+A_2)$. As in (a), $\|\Theta\zeta\|_{\Yz(I)}\le C_q(\|z(s_i)\|_{X_q}+\ell_2^{1/2}C(A_2))\le A_2$ and
$\|\Theta\zeta\|_{\Yo(I)}\le C_q(\|z(s_i)\|_{B^0_{q,2}}+\ell_2^{1/2}C(A_2))\le C_qC_B\varsigma+C_q\ell_2^{1/2}C(A_2)\le\frac{\varepsilon_0}2$, since $C_qC_B\varsigma\le\frac{\varepsilon_0}4$ and after decreasing $\ell_2$ so that
$C_q\ell_2^{1/2}C(A_2)\le\frac{\varepsilon_0}4$. The
contraction estimate of (a), in which the history term \eqref{eq:difB} starts at $s_i$ since two candidates coincide on $[a,s_i]$, has factor
at most $\frac12$ for $\ell_2$ small. $E_I$ is complete as in (a). The fixed point extends $z$ to an admissible solution on $[a,s_{i+1}]$,
which satisfies \eqref{eq:strict} again. After $N_2$ steps $z$ is defined on $\Sigma$.

\emph{Uniqueness.} Let $z,\tilde z$ be solutions on $\Sigma$ with $\|\cdot\|_{\Yo(\Sigma)}\le\varepsilon_*'=\varepsilon_0$ and
$\|\cdot\|_{L^1(\Sigma;W^{1,\infty})}\le1$. Both are admissible, hence satisfy \eqref{eq:strict}. Let $s_*:=\sup\{s:\ z=\tilde z\text{ on }[a,s]\}$ and
suppose $s_*<b$. On $I=[s_*,\min\{s_*+\ell_2,b\}]$ both restrictions have $\Yz(I)$ norm at most $A_0\le A_2$ and $\Yo(I)$ norm at most $\frac{\varepsilon_0}2$,
and the contraction estimate gives $\|z-\tilde z\|_{\Yo(I)}\le\frac12\|z-\tilde z\|_{\Yo(I)}$, so $z=\tilde z$ on $I$, a contradiction.

\emph{Bounds (i)--(iv).} (i) and (ii) are \eqref{eq:aprioriY0} and \eqref{eq:aprioriY1}, (iii) is \eqref{eq:aprioriLip}, and (iv) follows as in (a),
since every integrand there is bounded by a constant times $\sup_\Sigma\|z\|_{L^q}\le C_B\varsigma$, or by $C_r\|S\|_{L^r}$, and $|\Sigma|\le1$.

\proofstep{(c) Momentum and mean}
Take the mean of \eqref{eq:segcorr}: $c_z'=\fint g_z$. By \eqref{eq:Jprime}, which holds a.e.\ for fields in $L^1W^{1,\infty}$ by
Lemma~\ref{lem:kinstrong}(a), $\frac{\mathrm d}{\mathrm dt}(J[\hat f_z]-J[f[U]])=-\fint g_z$, since the kinetic fields are
$\hat U_z$ and $U$. Both quantities vanish at $a$, which gives the identity. The bound on $c_z$ is Lemma~\ref{lem:B2}(a) with
$\psi=v_i$, $k\equiv1$ and $s=q$. The bound on $\Rs g_z$ is \eqref{eq:fulldrag}, i.e.\ \eqref{eq:wkB} for the history term and
Lemma~\ref{lem:MRshift}(e) for the instantaneous term $\rho a_z$. Since $\|a_z\|_{L^q}\le3\|z\|_{L^q}$ and $|\Sigma|\le1$, the supremum bound
follows.

After an exit by the mean, the glued field on $[b,1]$ is $U+c_z$, where $c_z$ is continued by \eqref{eq:voltW}; this is the equation
of Lemma~\ref{lem:mean} started at $b$, so $c_z$ is defined and continuous on all of $[b,1]$. Thus $a_z=c_z$ there, and the kinetic
states still agree at $a$. The constants of Lemma~\ref{lem:B2} with radius $\Rv_1$ are available on $[a,t]$ as long as $\int_b^t|c_z|\le1$; we
verify this by a continuity argument. Let $t^*:=\sup\{t\in[b,1]:\int_b^t|c_z|\le1\}$. For $t\le t^*$, Lemma~\ref{lem:B2} applies on $[a,t]$
with the base $U$, and the first bound in \eqref{eq:cg} gives
\[
|c_z(t)|\le K\Big(3|\Sigma|\sup_\Sigma\|z\|_{L^q}+\int_b^t|c_z|\Big),
\]
so Gronwall's inequality yields $\sup_{[b,t^*]}|c_z|\le3Ke^K|\Sigma|\sup_\Sigma\|z\|_{L^q}$. We require
\begin{equation}\label{eq:exitsmall}
3Ke^K|\Sigma|\sup_\Sigma\|z\|_{L^q}\le\tfrac12 .
\end{equation}
In case (a), $\sup_\Sigma\|z\|_{L^q}\le C|\Sigma|^{1/2}\varsigma$ by (ii) and $\varsigma\le CN_S$, so \eqref{eq:exitsmall} holds after decreasing $\ell_0'$. In case
(b), $\sup_\Sigma\|z\|_{L^q}\le C_B\varsigma$ by \eqref{eq:aprioriY1} and $|\Sigma|\le1$, so \eqref{eq:exitsmall} holds after decreasing $c_*'$, which is
chosen last. Then $\int_b^{t^*}|c_z|\le\frac12$, so $t^*=1$ by continuity, and
$\sup_{[b,1]}|c_z|\le3Ke^K|\Sigma|\sup_\Sigma\|z\|_{L^q}$. In particular the contribution of $[b,1]$ to $\|\hat U-U\|_{L^1(a,1;L^\infty)}$ is at
most $\frac12\le1$, as used in the definition of $\Rv_1$. Inserting this bound in the second bound of \eqref{eq:cg} gives \eqref{eq:persist}.

\proofstep{(d) Persistent stress}
By \eqref{eq:persist}, $\|U\ot c_z+c_z\ot U\|_{W^{1,q}}\le2\|U\|_{W^{1,q}}|c_z|\le2B_U|c_z|$. By Lemma~\ref{lem:MRshift}(e),
$\|\Rs g_z\|_{W^{1,q}}\le C\|g_z\|_{L^q}\le C(C_h+\varrho|c_z|)$, since $|h_z|\le C_h$ and $|\rho a_z|\le\varrho|c_z|$ on $[b,1]$. The $L^q$ bound
follows from \eqref{eq:persist} and $\|U\|_\infty\le C_EB_U$.

\proofstep{Smoothness}
Recall that in this section $f_0$ and the base are smooth, so $f_\Sigma:=f[U](a)$ is smooth with compact velocity support, and
$\hat f_z(t)=e^{d(t-a)}f_\Sigma\circ Z_{\hat U_z}(a;t,\cdot)$ for $t\in\Sigma$. We show that $z\in C^\infty(\Sigma\times\TT^d)$, with one-sided derivatives at
$t=a$ and $t=b$. Then $c_z$, $\hat U_z$ and $\hat f_z$ are smooth on $\Sigma$ as well. (The extension of $z$ by zero to $t<a$ is in general not smooth,
since $\partial_tz(a)=-\QQ\D S(a)$; this is why the cutoff $\chi$ vanishes near $a$.)

\emph{Space.} We prove by induction on $k\ge0$ that $J^kz\in\mathrm{MR}_q(\Sigma)$, that is, $z\in L^2(\Sigma;H^{k+2,q})$ and
$\partial_tz\in L^2(\Sigma;H^{k,q})$. By Lemma~\ref{lem:MR}(b),(c), applied to $J^kz$, this implies $z\in C(\Sigma;W^{k+1,q})$. The case $k=0$ is
$z\in\Yz(\Sigma)$. Assume the claim for some $k$, and let $G$ be the right-hand side of the heat equation solved by $z$,
\[
\partial_tz-\Delta z=G:=\QQ\big[-\D S+g_z-\D(U\otimes z+z\otimes U+z\otimes z)\big].
\]
Since $H^{k+2,q}\hookrightarrow C^{k+1}$ for $q>d$, the derivatives in $y=(x,v)$ of order $1,\dots,k+1$ of the phase-space field of
$\hat U_z=U+\chi\mathring z+c_z$ are bounded by a function in $L^1(\Sigma)$. By differentiating the equation of the flow, the maps
$Z_{\hat U_z}(a;t,\cdot)$ are of class $C^{k+1}$, with derivatives bounded on compact sets uniformly in $t\in\Sigma$. Hence $\hat f_z(t)$ is
bounded in $C^{k+1}$, uniformly in $t$, with velocity support in a fixed ball, and so are its density and momentum in $C^{k+1}(\TT^d)$.
Since $\hat U_z\in L^\infty(\Sigma;W^{k+1,q})$, we get $g_z\in L^\infty(\Sigma;W^{k+1,q})$. Moreover $H^{k+1,q}$ is an algebra, because $q>d$, so
$\D(z\otimes z)=z\cdot\nabla z\in L^2(\Sigma;H^{k+1,q})$; the terms containing $U$ and $S$ belong to the same space. Thus $G\in L^2(\Sigma;H^{k+1,q})$.
By Lemma~\ref{lem:MR}(a), the heat equation with forcing $J^{k+1}G$ and zero datum at $a$ has a solution $\zeta\in\mathrm{MR}_q(\Sigma)$, and
$J^{-(k+1)}\zeta=z$ by uniqueness. This is the claim for $k+1$.

\emph{Time.} By the first part, $z$, and hence $\hat U_z$, belong to $C(\Sigma;W^{m,q})$ for every $m$. Therefore the flow of $\hat U_z$ and
all its derivatives in $y$ are continuous in $(t,y)$, and $\hat f_z\in C(\Sigma;C^m)$ for every $m$. We prove by induction on $j\ge0$ that
$z\in C^j(\Sigma;W^{m,q})$ for every $m$. If this holds for some $j$, then $\hat U_z\in C^j(\Sigma;W^{m,q})$ for every $m$. Since
$\partial_t\hat f_z=-\D_y\big(b_{\hat U_z}\hat f_z\big)$, a second induction shows that $\hat f_z\in C^{j+1}(\Sigma;C^m)$ for every $m$. Hence
$g_z\in C^j(\Sigma;W^{m,q})$, so $G\in C^j(\Sigma;W^{m,q})$, and $\partial_tz=\Delta z+G\in C^j(\Sigma;W^{m,q})$ for every $m$. This is the claim for
$j+1$.
\end{proof}

\begin{prop}[Mean-preserving gluing]\label{prop:glue}
Let $z$ be a segment corrector, $c=c_z$, $\hat U=\hat U_z$, $\hat f=\hat f_z$, and let
$q_z:=\Delta^{-1}\D\big(-\D S+g_z-\D(U\otimes z+z\otimes U+z\otimes z)\big)$, so that $\partial_tz-\Delta z+\D(U\otimes z+z\otimes U+z\otimes z)+\nabla q_z=-\D S+g_z$
exactly; $q_z$ is smooth on $\Sigma$. On $\Sigma$,
\begin{equation}\label{eq:glueid}
\partial_t\hat U-\Delta\hat U+\D(\hat U\otimes\hat U)+\nabla(P+\chi q_z)-F[\hat f,\hat U]=\chi'\mathring z+(\chi^2-\chi)\D(\mathring z\otimes\mathring z)+(1-\chi)\big[\D S+c\cdot\nabla U+c'-g_z\big].
\end{equation}
Since $c'=\fint g_z$, the right-hand side has zero mean. Hence $(\hat U,\hat f)$ is a smooth VNS--Reynolds solution on $\Sigma$, with stress
\[
\hat S=\chi'\Rs z+(\chi^2-\chi)\,\mathring z\ot\mathring z+(1-\chi)\big(S+U\ot c+c\ot U-\Rs g_z\big),
\]
and pressure $P+\chi q_z-\frac1d(\chi^2-\chi)|\mathring z|^2-\frac2d(1-\chi)U\cdot c$. In particular $\hat S=0$ where $\chi=1$.

If $\chi$ is two-sided, define the \emph{exit by the mean}: after $b$, set $\hat U:=U+c(t)$, where $c$ solves the Volterra equation
\begin{equation}\label{eq:voltW}
c'=\fint\Big(F\big[f[U|_{[0,a]}\oplus\hat U],\hat U\big]-F[f[U],U]\Big),
\end{equation}
continuing $c_z$; here $\hat U=\hat U_z$ on $\Sigma$ and $\hat U=U+c$ on $[b,1]$, so the kinetic state is $\hat f_z$ continued. Then \eqref{eq:glueid} holds with $\chi=0$ on $[b,1]$,
where $g_z$ is understood as the continued drag difference $F[\hat f,\hat U]-F[f[U],U]$, and the glued field is a smooth VNS--Reynolds solution on
$[0,1]$.
\end{prop}

\begin{proof}
Write $\hat U=U+a$ with $a=\chi\mathring z+c$, and $(b\cdot\nabla g)_i=b_j\partial_jg_i$. Since $\D U=\D\mathring z=0$ and $c$ depends only on $t$,
\[
\D(U\otimes a+a\otimes U)=\chi\big(\mathring z\cdot\nabla U+U\cdot\nabla\mathring z\big)+c\cdot\nabla U,\qquad
\D(a\otimes a)=\chi^2\,\mathring z\cdot\nabla\mathring z+\chi\,c\cdot\nabla\mathring z .
\]
Hence, by the equation of the base,
\[
\begin{aligned}
&\partial_t\hat U-\Delta\hat U+\D(\hat U\otimes\hat U)+\nabla P-F[f[U],U]-\D S\\
&\qquad=\chi'\mathring z+c'+c\cdot\nabla U+\chi\Big[\partial_t\mathring z-\Delta\mathring z+\mathring z\cdot\nabla U+U\cdot\nabla\mathring z+c\cdot\nabla\mathring z\Big]+\chi^2\,\mathring z\cdot\nabla\mathring z .
\end{aligned}
\]
On the other hand, \eqref{eq:segcorr}, written with its pressure $q_z$ and with $z=\mathring z+c$, reads
\[
\partial_t\mathring z-\Delta\mathring z+\mathring z\cdot\nabla U+U\cdot\nabla\mathring z+c\cdot\nabla\mathring z+\nabla q_z=-\D S+g_z-c'-c\cdot\nabla U-\mathring z\cdot\nabla\mathring z .
\]
Inserting this into the bracket, and using $F[\hat f,\hat U]=F[f[U],U]+g_z$ and $\mathring z\cdot\nabla\mathring z=\D(\mathring z\otimes\mathring z)$, we obtain
\eqref{eq:glueid}: the terms $c'+c\cdot\nabla U$, $\D S$ and $g_z$ appear with the factor $1-\chi$, and the quadratic term with the factor
$\chi^2-\chi$.

The stress form follows from $\chi'\mathring z=\D(\chi'\Rs z)$, $\D(\mathring z\otimes\mathring z)=\D(\mathring z\ot\mathring z)+\frac1d\nabla|\mathring z|^2$,
$c\cdot\nabla U=\D(U\ot c+c\ot U)+\frac2d\nabla(U\cdot c)$, and $c'-g_z=-\D\Rs g_z$. The gradient terms are moved to the left-hand side,
which produces the minus signs in the pressure.

Near $a$, where $\chi=0$, we have $\hat U_z=U+c_z$, so $c_z$ solves the Volterra equation $c'=\fint\big(F[f[U|_{[0,a]}\oplus(U+c)],U+c]-F[f[U],U]\big)$
with $c(a)=0$. The function $c\equiv0$ solves it, so $c_z=0$ there by uniqueness (Lemma~\ref{lem:mean}, Step~2). Hence $\hat U=U$, $\hat f=f[U]$ and
$\hat S=S$ near $a$, and the glued solution is smooth across $a$.

After $b$, \eqref{eq:voltW} is the equation of Lemma~\ref{lem:mean} started at $b$, written for $c$ alone, and it has a unique smooth
solution. On the zone where $\chi=0$ near $b$, $\hat U$ depends on $z$ only through $c_z$, which satisfies the same equation. So the
continuation is smooth.
\end{proof}

\begin{lem}[Junctions]\label{lem:junction}
Suppose a segment exits by the mean. Suppose a later segment or convex integration step (Lemma~\ref{lem:CIloc}), whose
modifications vanish near its start, begins in the zone where the first cutoff vanishes. Then, on the overlap, the glued velocity is $B+c_{\rm tot}(t)$. Here
$B$ is the base before the first segment, and $c_{\rm tot}$ is the unique solution of \eqref{eq:voltW} relative to $B$. In particular,
the glued solution is smooth across the junction.
\end{lem}

\begin{proof}
After the second object starts, its mean satisfies $c_{\rm new}'=\fint(F[\hat f,U+c_{\rm new}]-F[f[U],U])$, where $U=B+c_{\rm prev}$. Adding the
equation of $c_{\rm prev}$, the middle terms cancel, so $c_{\rm prev}+c_{\rm new}$ solves \eqref{eq:voltW} relative to $B$ with the true kinetic
state. Before the start, $c_{\rm tot}=c_{\rm prev}$ solves the same equation. Conclude by uniqueness (Lemma~\ref{lem:mean}, Step~2).
\end{proof}

\begin{cor}[Layer stress]\label{cor:layer}
Let $\chi$ be monotone on each transition layer, with total transition length $\lambda_\chi$, so that $\|\chi'\|_{L^1}\le2$ and
$\|\chi-\chi^2\|_{L^1}+\|1-\chi\|_{L^1(\Sigma)}\le\lambda_\chi$ on the transition and zero zones inside $\Sigma$. Then
\[
\|\hat S\|_{L^1(\Sigma;L^r)}\le2\sup_\Sigma\|\Rs z\|_{L^r}+\lambda_\chi\Big(\sup_\Sigma\|S\|_{L^r}+\sup_\Sigma\|\mathring z\|_{L^\infty}^2+2C_EB_U\sup_\Sigma|c|+\sup_\Sigma\|\Rs g_z\|_{L^q}\Big).
\]
\end{cor}

\begin{proof}
Bound each term of $\hat S$ by the $L^1$ norm in time of its cutoff factor times the supremum of the remaining factor.
\end{proof}

The persistent stress of a single segment is bounded in $W^{1,q}$ by a constant, but not by a small one: the drag difference
$g_z$ is bounded pointwise, while the kinetic states differ by an amount that is small only in a weak norm. These bounds are not
summable over many segments. Instead we use that the drag differences telescope:

\begin{lem}[Accumulated persistent stress]\label{lem:telescope}
Let $(U_0,P_0,S_0)$ be a smooth VNS--Reynolds solution, and let $(U_l,P_l,S_l)$, $l=1,\dots,n$, be obtained by successive
applications of Propositions~\ref{prop:segcorr} and~\ref{prop:glue} with exit by the mean, or of Lemma~\ref{lem:CIloc}, each starting
after the previous one has exited (after $\sup Z_1$ for Lemma~\ref{lem:CIloc}). Let $t$ be a time after the last exit, and let
$\Pi_l$ be the persistent stress of the $l$-th operation (Proposition~\ref{prop:segcorr}(d), respectively $c\ot U+U\ot c-\Rs h$ in
Lemma~\ref{lem:CIloc}). Then $U_n(t)=U_0(t)+C_n(t)$ with $C_n=\sum_{l\le n}c_l$ spatially constant, and
\[
\sum_{l\le n}\Pi_l(t)=U_0\ot C_n+C_n\ot U_0-\Rs\big(F[f[U_n],U_n]-F[f[U_0],U_0]\big)(t)+M_n(t),
\]
where $M_n(t)$ is a spatially constant traceless matrix with $|M_n(t)|\le2\max_{l\le n}|C_l(t)|\sum_{l\le n}|c_l(t)|$. In particular, if $\Rv_1$
bounds the velocity supports of all kinetic states involved and $\varrho,C_h$ are as in Section~\ref{sec:segcorrsub}, then
\[
\Big\|\sum_{l\le n}\Pi_l(t)\Big\|_{W^{1,q}}\le2\|U_0(t)\|_{W^{1,q}}|C_n(t)|+C\big(C_h+\varrho|C_n(t)|\big)+|M_n(t)|,
\]
with $C$ depending only on $d,q$, and independent of $n$. The stress is not modified: $M_n$ is kept, and since $\D M_n=0$ it plays no
role in the Reynolds equation; it only enters norms through the last term.
\end{lem}

\begin{proof}
After the $l$-th exit the base of the $(l+1)$-st operation is the glued field $U_l=U_{l-1}+c_l$ with its kinetic state
$f[U_l]$. Hence the drag difference of the $l$-th operation is $g_l=F[f[U_l],U_l]-F[f[U_{l-1}],U_{l-1}]$ at time $t$, and
$\sum_lg_l$ telescopes. For the cross terms, $U_{l-1}=U_0+C_{l-1}$, so
$\sum_l(U_{l-1}\ot c_l+c_l\ot U_{l-1})=U_0\ot C_n+C_n\ot U_0+M_n$ with $M_n=\sum_l(C_{l-1}\ot c_l+c_l\ot C_{l-1})$, and $|a\ot b|\le|a||b|$
gives the bound on $|M_n|$. The bound
follows from Lemma~\ref{lem:MRshift}(e), $\|\Rs g\|_{W^{1,q}}\le C\|g\|_{L^q}$, and the pointwise bound
$|F[f[U_n],U_n]-F[f[U_0],U_0]|\le C_h+\varrho|C_n|$, as in Proposition~\ref{prop:segcorr}(d).
\end{proof}

\subsection{The localized convex integration step}

\begin{lem}[Localized CI step]\label{lem:CIloc}
Let $(U,P,S)$ be a smooth VNS--Reynolds solution with $B_U,\Lambda_U<\infty$. Let $Z_0\subset Z_1'\subset Z_1\subset[0,1]$ be intervals, let
$\theta_0\in C^\infty_c(Z_1';[0,1])$ equal $1$ on $Z_0$, and let $\theta_1\in C^\infty_c(Z_1;[0,1])$ equal $1$ on $Z_1'$. Set $S_c:=\theta_0S$ and
$S_s:=(1-\theta_0)S$. For every $\delta'>0$ there are smooth $w$ and $c$ with the following properties:
\begin{itemize}
\item $\Supp_tw\subset Z_1$, $\D w=0$, $\fint w=0$;
\item $c$ is spatially constant, with $c=0$ before $\inf Z_1$;
\item $U^\sharp=U+w+c$ is a smooth VNS--Reynolds solution with stress
\[
S^\sharp=S^{\rm CL}+S_s+c\ot(U+w)+(U+w)\ot c-\Rs h,\qquad h=F[f[U^\sharp],U^\sharp]-F[f[U],U];
\]
\item $\Supp_tS^{\rm CL}\subset Z_1$, $\|S^{\rm CL}\|_{L^1L^r}\le\delta'$, $\|w\|_{L^pL^\infty\cap L^1W^{1,q}}\le\delta'$, and $\|w\|_{L^2}\le M\|S_c\|_{L^1_{t,x}}^{1/2}$;
\item $\sup_t|c(t)|\le K^\sharp\delta'$; $\|\Rs h\|_{L^1(0,1;L^q)}\le K^\sharp\delta'$; and for $t\notin Z_1$, $\|\Rs h(t)\|_{L^q}\le K^\sharp\delta'$ and
$\|c\ot U+U\ot c-\Rs h\|_{W^{1,q}}\le C^\sharp$. The constants $K^\sharp,C^\sharp$ depend only on $d,q,\Rv_0,M_\infty$ and upper bounds for $B_U,\Lambda_U,L_U$, defined by
\eqref{eq:basenorms} with $a=\inf Z_1$.
\end{itemize}
\end{lem}

\begin{proof}
Apply Proposition~\ref{prop:NSstep} to $(U,S_c)$, with the two-sided cutoff $\theta_1$ in place of \eqref{eq:theta}. The algebra of
Appendix~\ref{sec:CI} uses only $\theta_1^2S_c=S_c$, see \eqref{eq:thetaR}, which holds since $\theta_1=1$ on $\Supp S_c$. The estimates of
Appendix~\ref{sec:estimates} use only the smoothness of $\theta_1$ and $0\le\theta_1\le1$; the constants $C_u$ now depend on $\theta_1$, while $M$
does not.
Since $\|\partial_t^k\theta_1\|_\infty\lesssim_k\tau^{-k}$, these constants depend on $\tau$; this is admissible, since the frequency $\lambda$ of the
step is chosen after $\tau$ (see the order of choices at the end of the proof of Proposition~\ref{prop:W}). The stress $S_s$ is carried
along unchanged: Proposition~\ref{prop:NSstep}(i) holds with $S_c$ in place of $R$, and adding $\D S_s$ to both sides shows that the new
stress is $S^{\rm CL}+S_s$, up to the kinetic terms treated below. Then proceed
as in Section~\ref{sec:proofstep}, starting Lemma~\ref{lem:mean} at $t_0:=\inf Z_1$ from the kinetic state $f[U](t_0)$.

For the kinetic bounds, apply Remark~\ref{rem:shift} with $b=U$, $a=w+c$. The momentum identity and Lemma~\ref{lem:B2}(a) with
$\psi=v_i$ give $|c(t)|\le K(\|w\|_{L^1(t_0,t;L^q)}+\int_{t_0}^t|c|)$, hence $\sup|c|\le Ke^K\delta'$ by Gronwall, using
$\|w\|_{L^1L^q}\le\|w\|_{L^pL^\infty}\le\delta'$. Here $K$ is the constant of Lemma~\ref{lem:B2} for the radius
$\Rv_1=\Rv_0+L_U+2(C_EB_U+2)$, which is admissible as long as $\|w\|_{L^1L^\infty}+\int_{t_0}^t|c|\le2$. As in the proof of
Proposition~\ref{prop:segcorr}(c), a continuity argument in $t$ shows that this holds on $[t_0,1]$ if $\delta'\le1$ and $Ke^K\delta'\le1$, which we
assume. By \eqref{eq:fulldrag},
\[
\|\Rs h(t)\|_{L^q}\le C_qK\big(1+\Rv_1+\|U(t)\|_{W^{1,q}}\big)\Big(\int_{t_0}^t\|w\|_{L^q}+\int_{t_0}^t|c|\Big)+C_q\varrho\big(\|w(t)\|_{L^q}+|c(t)|\big).
\]
The instantaneous term $\varrho\|w(t)\|_{L^q}$ is small only after integration in time, and it vanishes for $t\notin Z_1$. This gives the
stated bounds. The $W^{1,q}$ bound off $Z_1$ follows as in Proposition~\ref{prop:segcorr}(d), since $|h|\le C_h+\varrho|c|$ there.
\end{proof}

\subsection{Time geometry}\label{sec:geom}

Let $\tilde I=\bigcup_{k=1}^{N_{\rm c}}[\alpha_k,\alpha_k+5\tilde\tau]$, with $\alpha_1<\dots<\alpha_{N_{\rm c}}$. Since these intervals are pairwise
disjoint, $\alpha_{k+1}>\alpha_k+5\tilde\tau$. By well-preparedness, $R=0$ on $[\alpha_k,\alpha_k+\tilde\tau]\cup[\alpha_k+4\tilde\tau,\alpha_k+5\tilde\tau]$. Choose
$\tau\le\tilde\tau/2$ and $T=\tau^\varepsilon$ with
\[
T<\tfrac{\tilde\tau}2,\qquad5\tau<T,\qquad\tfrac52\tau<\tfrac{\tilde\tau}2 .
\]
As in \cite[Section~3.1]{CL22}, junctions lie on the global grid $T\ZZ$. The chain of component $k$ consists of the grid points
$t_{k,0}<\dots<t_{k,N_k}$ in $[\alpha_k+\frac{\tilde\tau}2,\alpha_k+\frac{9\tilde\tau}2]$. Then $t_{k,0}\le\alpha_k+\tilde\tau$ and $t_{k,N_k}\ge\alpha_k+4\tilde\tau$, so the
windows $[t_{k,i-1},t_{k,i}]$ cover the old stress support. Around a junction $t$ define
\[
Z_0=[t\pm\tfrac\tau2],\quad Z_1'=[t\pm\tfrac{3\tau}4],\quad Z_1=[t\pm\tau],\quad Z_2=[t\pm\tfrac{5\tau}4],\quad W=[t\pm\tfrac{5\tau}2],
\]
where $[t\pm\eta]:=[t-\eta,t+\eta]$, and set $I:=\bigcup W(t_{k,i})$.

\begin{lem}[Geometry]\label{lem:geom}
\begin{enumerate}[label=(\alph*)]
\item $I\subset\tilde I$.
\item The windows $W$ are pairwise disjoint, and $0,1\notin I$.
\item $Z_2(t)\subset\{s:\dist(s,I^c)\ge\tau\}$.
\item There are fewer than $\tau^{-\varepsilon}$ windows.
\item Consecutive junctions are more than $5\tau$ apart.
\end{enumerate}
\end{lem}

\begin{proof}
Interval arithmetic. For (b) and (e), junctions of the same chain are distinct points of $T\ZZ$, hence at distance at least
$T>5\tau$. For junctions of different chains $k<k'$, the disjointness of the components gives $\alpha_{k'}>\alpha_k+5\tilde\tau$, hence
\[
t_{k',0}-t_{k,N_k}\ge\Big(\alpha_{k'}+\frac{\tilde\tau}2\Big)-\Big(\alpha_k+\frac{9\tilde\tau}2\Big)>\tilde\tau>5\tau,
\]
using $\frac52\tau<\frac{\tilde\tau}2$. Since the windows have length $5\tau$, they are pairwise disjoint; in particular the new set $I$ is again a union
of pairwise disjoint intervals, as required in Definition~\ref{def:wp}. For (d), distinct junctions are distinct grid points in $(0,1)$,
and there are fewer than $T^{-1}$ of these.
\end{proof}

\subsection{Proof of Proposition~\ref{prop:W}}\label{sec:proofW}

The construction runs from left to right in time and maintains a \emph{current global solution} $U$: a smooth VNS--Reynolds
solution on $[0,1]$ with stress $S_U$. Initially $(U,S_U)=(u,R)$. It has two passes. The second pass cannot be merged with
the first. The entry layer of a corrector to the right of a junction depends on the kinetic state, and therefore on the
convex integration perturbation at that junction, which is built to cancel the stress there.

\proofstep{Uniformity of constants}
The constants of Proposition~\ref{prop:segcorr} and Lemma~\ref{lem:CIloc} depend on the base only through upper bounds for
$B_U,\Lambda_U,L_U$ of \eqref{eq:basenorms} and for $N_S$. We fix these upper bounds \emph{before} $T$, in two levels, and then check, for
$T$, $\tau$ and $\delta_\star$ small, that every base used in the construction satisfies them. The bounds of
level~2 depend on the constants of level~1, but not conversely, because a base of pass~1 never contains a corrector of pass~1 on
$[a,1]$, and a base of pass~2 never contains a corrector of pass~2 on $[a,1]$.
\begin{itemize}
\item \emph{Level 1 (pass 1).} Put $\bar B_1:=B_u+1$, $\bar\Lambda_1:=\Lambda_u+1$, $\bar L_1:=\|u\|_{L^1L^\infty}+1$, where $B_u,\Lambda_u$ are taken with $a=0$, and
\[
\bar N_1:=\|R\|_{L^\infty W^{1,q}}+2\bar B_1+C_{d,q}\big(C_h(\bar B_1,\bar L_1)+\varrho(\bar B_1,\bar L_1)\big)+2 .
\]
Let $C_1',\ell_{0,1}',c_{*,1}'$ and $K_1,C_1$ be the constants of Proposition~\ref{prop:segcorr} for these bounds. Every operation of pass~1
starts after all earlier operations of pass~1 have exited, and its segment and cutoff vanish near its start. Hence on $[a,1]$ its base is
$u+C_{\rm tot}$ with $C_{\rm tot}$ spatially constant, so $B_U\le B_u+|C_{\rm tot}|$ and $\Lambda_U\le\Lambda_u$. Moreover $L_U\le\|u\|_{L^1L^\infty}+\sum\|z_i\|_{L^1L^\infty}+|C_{\rm tot}|$,
where the window correctors contribute at most $T^{-1}\cdot T\cdot C_1'T^{1/2}$ and the long correctors at most $N_{\rm c}C_1'(CT^{1/2})^\beta$, by
Proposition~\ref{prop:segcorr}(a)(i), (b)(iii) and the bound $\varsigma\le CT^{1/2}$ on good segments shown below. The base stress is $R$ plus
the accumulated persistent stress, whose $W^{1,q}$ norm is bounded by Lemma~\ref{lem:telescope} with $U_0=u$. So all level-1 bounds hold
as soon as $|C_{\rm tot}|\le\frac12$, $\sum_l\sup|c_l|\le\frac12$ and $T$ is small.
\item \emph{Level 2 (pass 2).} Let $\bar u$ be the result of pass~1. The pass-1 correctors have disjoint supports in time, so
$\sup_t\|\bar u-u-C_{\rm tot}\|_{W^{1,q}}\le C\max(\sup\|z_i\|_{X_q})\le CC_1'$. Also $\|\nabla(\bar u-u)\|_{L^1L^\infty}\le T^{-1}\cdot T^{1/2}\cdot C_1'T^{1/2}\bar\varsigma_1^\beta+N_{\rm c}$,
with $\bar\varsigma_1:=\|R\|_{L^\infty L^q}+C$ the size of the base stresses of the windows. Put
\[
\bar B_2:=B_u+CC_1'+1,\qquad \bar\Lambda_2:=\Lambda_u+C_1'\bar\varsigma_1^\beta+N_{\rm c}+1,\qquad \bar L_2:=\bar L_1+1,
\]
\[
\bar N_2:=\bar N_1+2\bar B_2+C_{d,q}\big(C_h(\bar B_2,\bar L_2)+\varrho(\bar B_2,\bar L_2)\big)+2 ,
\]
and let $C_2',\ell_{0,2}',c_{*,2}',\varepsilon_{*,2}'$, $K_2^\sharp,C_2^\sharp$ be the constants of Proposition~\ref{prop:segcorr} and Lemma~\ref{lem:CIloc} for these
bounds. Every operation of pass~2 starts after all earlier operations of pass~2 have exited, and every $w_j$ vanishes after $\sup Z_1(t_j)$,
which precedes the start of every later operation (Lemma~\ref{lem:geom}(e)). Hence on $[a,1]$ the base of a pass-2 operation is
$\bar u+C_{\rm tot}$, and the $w_j$ enter only through $L_U$, with $\sum_j\|w_j\|_{L^1L^\infty}\le N_J\delta_\star\le\frac12$. The base stress on a
re-solved segment is $\bar R$ plus the accumulated persistent stresses of pass~2; the former vanishes there, since it is supported in
$\bigcup Z_0(t_j)$ (see the paragraph following the description of pass~1), and the latter is bounded by
Lemma~\ref{lem:telescope} with $U_0=\bar u$. So all level-2 bounds hold as soon as $|C_{\rm tot}|\le\frac12$, $\sum_l\sup|c_l|\le\frac12$ and
$N_J\delta_\star\le\frac12$.
\end{itemize}
Summing the bounds of Proposition~\ref{prop:segcorr}(d) termwise, instead of using Lemma~\ref{lem:telescope}, would only give $N_S\le CT^{-1}$,
which would make $\ell_0'$ depend on $T$. From now on $C'$, $\ell_0'$, $c_*'$, $\varepsilon_*'$ denote the larger, respectively smaller, of the two levels,
and $C$ denotes a constant depending only on $d,p,q,r,\beta,\varepsilon,f_0$, on $(u,R,\tilde I,\tilde\tau)$ and on $C_1',C_2'$; it does not depend on
$T,\tau$ or the parameters $\delta_j'$. We take $T\le\ell_0'$. The conditions $|C_{\rm tot}|\le\frac12$ and $\sum_l\sup|c_l|\le\frac12$ are part of
\eqref{eq:tol}: by Proposition~\ref{prop:segcorr}(c) and Lemma~\ref{lem:CIloc}, $\sum_l\sup|c_l|\le CT^{1/2}+CN_J\delta_\star+C(N_{\rm c}+1)\bar\varsigma$.

\proofstep{Pass 1: causal concentration}
On $[0,t_{1,0}]$ nothing is done. For each component $k$, in order:
\begin{itemize}
\item for $i=1,\dots,N_k$, apply Proposition~\ref{prop:segcorr}(a) on the window $\Sigma_{k,i}=[t_{k,i-1},t_{k,i}]$ with base the current $U$. Use a
two-sided cutoff vanishing on $[t_{k,i-1},t_{k,i-1}+\frac\tau4]\cup[t_{k,i}-\frac\tau4,t_{k,i}]$ and equal to $1$ on $[t_{k,i-1}+\frac\tau2,t_{k,i}-\frac\tau2]$.
Glue by Proposition~\ref{prop:glue}, exiting by the mean;
\item then apply Proposition~\ref{prop:segcorr}(b) on the good segment $[t_{k,N_k},t_{k+1,0}]$, or $[t_{N_{\rm c},N_{N_{\rm c}}},1]$ with a one-sided cutoff, with
layers of length $\frac\tau4$ inside the $Z_0$ of the endpoints.
\end{itemize}
The base stress on a window is $R$ plus the accumulated persistent stress of earlier windows, so $\varsigma\le\|R\|_{L^\infty L^q}+C$ and,
by Lemma~\ref{lem:telescope}, $N_S\le C$; it need not be small, and (a) does not require it. After a window with corrector $z_i$ exits, the base stress for all
later times acquires the \emph{persistent} term $\Pi_i$. By Proposition~\ref{prop:segcorr}(a)(ii) and (d), with $|\Sigma_{k,i}|=T$,
\[
\sup_{t\ge t_{k,i}}\|\Pi_i(t)\|_{L^q}\le CT\sup_{\Sigma_{k,i}}\|z_i\|_{L^q}\le CT^{3/2}.
\]
Since there are fewer than $T^{-1}$ windows, $\sum_i\sup\|\Pi_i\|_{L^q}\le CT^{1/2}$. On the good segments $R=0$, so the base stress
there consists of the persistent terms of all earlier windows and of all earlier good segments. Let $\varsigma_k$ be its size on the
$k$-th good segment. By \eqref{eq:persist} and Proposition~\ref{prop:segcorr}(b)(ii), the $k'$-th good segment leaves a persistent term
of size at most $CC_B\varsigma_{k'}$. Hence $\varsigma_k\le CT^{1/2}+C\sum_{k'<k}\varsigma_{k'}$, and by induction
$\varsigma_k\le(1+C)^{N_{\rm c}}CT^{1/2}\le c_*'$ for $T$ small, as required by (b). The factor $(1+C)^{N_{\rm c}}$ depends only on $\tilde I$,
which is fixed before $T$. The same persistent terms enter the base stress of the later windows; they are bounded by $CT^{1/2}$ in
$L^q$ and, by Lemma~\ref{lem:telescope}, by a constant in $W^{1,q}$, so the bounds $\varsigma\le\|R\|_{L^\infty L^q}+C$ and $N_S\le C$ for the
windows are unaffected. In what follows, constants $C$ in pass~1 may depend on $N_{\rm c}$.

Let $(\bar u,\bar p,\bar R)$ be the result. The persistent terms do not survive outside the layers. Indeed, the corrector of each
segment solves \eqref{eq:segcorr}, whose forcing contains $-\D S$ for the full base stress $S$, which includes all persistent terms
created earlier; hence, by Proposition~\ref{prop:glue}, the glued stress vanishes wherever the cutoff of that segment equals $1$.
After the exit of a segment the base stress, together with the new persistent term, reappears, but only until the cutoff of the
next segment reaches $1$. The windows and good segments cover $[t_{1,0},1]$, their cutoffs equal $1$ outside $\bigcup Z_0(t_{k,i})$,
the last cutoff equals $1$ near $t=1$, and $R=0$ on $[0,t_{1,0}]$. Therefore $\bar R$ is supported in $\bigcup Z_0(t_{k,i})$, where it
consists of $R$, the layer terms and the persistent terms. By Corollary~\ref{cor:layer} and Proposition~\ref{prop:segcorr}(a), the layer stress of the window $\Sigma_{k,i}$ is at
most
\[
2C_r\int_{\Sigma_{k,i}}\|R\|_{L^r}+CT^{3/2}+C\tau\big(\|R\|_{L^\infty L^r}+C+T\big),
\]
and summing over windows, as in \cite[Proposition~3.3]{CL22},
\begin{equation}\label{eq:pass1}
\|\bar R\|_{L^1L^r}\le2C_r\|R\|_{L^1L^r}+C\big(T^{1/2}+\tau^{1-\varepsilon}\big).
\end{equation}
Moreover, by Proposition~\ref{prop:segcorr}(a)(i) on windows (as in (W3) below), and (b)(i), (iii) on good segments:
\begin{itemize}
\item $\|\bar u-u\|_{L^\infty L^\infty\cap L^1W^{1,q}}\le CT^{1/2}+C(T^{1/2})^\beta$;
\item $\|\bar u-u\|_{L^2(\tilde I^c;W^{1,\infty})}\le C(T^{1/2})^\beta$, since $\tilde I^c$ meets only the good segments;
\item $\bar u=u$ on $[0,t_{1,0}]$.
\end{itemize}

\proofstep{Pass 2: convex integration and re-solving}
Enumerate the junctions $t_1<\dots<t_{N_J}$, with $N_J<\tau^{-\varepsilon}$. Start from $U=\bar u$. For $j=1,\dots,N_J$:
\begin{itemize}
\item \emph{(CI)} apply Lemma~\ref{lem:CIloc} at $t_j$ with the zones $Z_0,Z_1',Z_1$ of $t_j$, the current stress $S$, and a parameter $\delta_j'$.
The cancelled part $S_c=\theta_0S$ contains $\bar R$ on $Z_0(t_j)$;
\item \emph{(re-solve)} on $\Sigma_j:=[t_j+\tau,t_{j+1}-\frac{9\tau}8]$, or on $\Sigma_{N_J}:=[t_{N_J}+\tau,1]$ with a one-sided cutoff, with base the current
$U$: apply Proposition~\ref{prop:segcorr}(a) if $t_j$ and $t_{j+1}$ belong to the same chain, and Proposition~\ref{prop:segcorr}(b)
otherwise, i.e.\ if $t_j=t_{k,N_k}$ is the last junction of a chain or $j=N_J$. In the latter case $|\Sigma_j|$ is comparable to the
gap between two components of $\tilde I$ and need not be smaller than $\ell_0'$. Let $\mathcal L$ be the set of these indices;
$\#\mathcal L\le N_{\rm c}$. The layers are $[t_j+\tau,t_j+\frac{9\tau}8]$ and $[t_{j+1}-\frac{5\tau}4,t_{j+1}-\frac{9\tau}8]$, and the segment
starts after $\Supp w_j$.
\end{itemize}
On the re-solved segment $\bar R$ vanishes, so the base stress consists of the persistent terms of pass~2: those of the convex
integration steps $i\le j$ (the terms $c_i\ot U+U\ot c_i-\Rs h_i$ off $Z_1(t_i)$), and those of the exits of earlier re-solved segments.
They are the mismatch between the pass-1 solution $\bar u$, which is exact there, and the kinetic state modified by pass~2.

\proofstep{Supports}
After pass~2, the stress vanishes on the bulks $[t_j+\frac{5\tau}4,t_{j+1}-\frac{5\tau}4]$. In $Z_2(t_j)$ it consists of:
\begin{itemize}
\item $S^{\rm CL}_j$, the cross terms $c_j\ot(U+w_j)+(U+w_j)\ot c_j$ and $-\Rs h_j$, on $Z_1(t_j)$;
\item the layers of the two adjacent segments;
\item the persistent terms, alive between the exit of one segment and the entry of the next.
\end{itemize}
Before $t_1-\frac{5\tau}4$ nothing has changed. By Lemma~\ref{lem:geom}(b)--(d), $u_1$ is well-prepared for $(I,\tau)$, and $I\subset\tilde I$ by Lemma~\ref{lem:geom}(a).

\proofstep{Smallness}
Let $\varsigma_j$ be the size, in $L^\infty L^q$, of the base stress on the $j$-th re-solved segment $\Sigma_j$. By the bounds above,
Lemma~\ref{lem:CIloc} and \eqref{eq:persist},
\[
\varsigma_j\le CT^{1/2}+C\sum_{i\le j}\delta_i'+C\sum_{i<j}|\Sigma_i|\sup_{\Sigma_i}\|z_i\|_{L^q}
\le CT^{1/2}+C\sum_{i\le j}\delta_i'+C\sum_{i<j,\,i\notin\mathcal L}|\Sigma_i|\varsigma_i+C\sum_{i<j,\,i\in\mathcal L}\varsigma_i ,
\]
using Proposition~\ref{prop:segcorr}(a)(ii) for $i\notin\mathcal L$ and (b)(ii) for $i\in\mathcal L$. By induction on $j$,
$\varsigma_j\le A\prod_{i<j,\,i\notin\mathcal L}(1+C|\Sigma_i|)\prod_{i<j,\,i\in\mathcal L}(1+C)$ with $A=CT^{1/2}+CN_J\max_j\delta_j'$, and since
$\sum_i|\Sigma_i|\le1$ and $\#\mathcal L\le N_{\rm c}$,
\begin{equation}\label{eq:sigmabar}
\max_j\varsigma_j\le\bar\varsigma:=e^C(1+C)^{N_{\rm c}}\Big(CT^{1/2}+CN_J\max_j\delta_j'\Big).
\end{equation}
The application of (b) on $\Sigma_i$, $i\in\mathcal L$, requires $\varsigma_i\le c_*'$, which is part of the induction and follows from
$\bar\varsigma\le c_*'$. All bounds in this recursion are linear in the stress, because the persistent terms are controlled by the time
integral of the perturbation, through \eqref{eq:persist}; no factor exponential in $N_J$ appears. The factor $(1+C)^{N_{\rm c}}$ depends only
on $\tilde I$, which is fixed before $T$. The H\"older power $\varsigma^\beta$ appears only in
the outputs below. We choose $\delta_j'\le\delta_\star$, with $\delta_\star$ and $T$ so small that
\begin{equation}\label{eq:tol}
\begin{aligned}
&\bar\varsigma\le c_*',\qquad C(N_{\rm c}+1)\bar\varsigma+CN_J\delta_\star\le\tfrac\delta{10},\qquad (N_{\rm c}+1)^{1/2}C'\bar\varsigma^\beta\le\tfrac\delta{10},\\
&CT^{1/2}+CN_J\delta_\star+C(N_{\rm c}+1)\bar\varsigma\le\tfrac12,\qquad N_J\delta_\star\le\tfrac12,
\end{aligned}
\end{equation}
where $C'$ is the larger of the constants of Proposition~\ref{prop:segcorr}(a)(iii) and (b)(iii). These conditions are met by
choosing first $T$ (hence $\tau$) small and then $\delta_\star$ small in terms of $N_J$. Note that $N_J\le\tau^{-\varepsilon}$ is fixed before
$\delta_\star$.

\proofstep{Conclusion}
\begin{itemize}
\item (W1): on $Z_1(t_j)$, $\|S^{\rm CL}_j\|_{L^1L^r}+\|c_j\ot(U+w_j)+(U+w_j)\ot c_j\|_{L^1L^r}+\|\Rs h_j\|_{L^1L^q}\le C\delta_j'$ by Lemma~\ref{lem:CIloc}.
The layer stress of $\Sigma_j$ is bounded by Corollary~\ref{cor:layer} and Proposition~\ref{prop:segcorr}(a)(ii)--(iv):
\[
2\big(C_r|\Sigma_j|\varsigma_j+C|\Sigma_j|^{3/2}\varsigma_j\big)+C\tau\big(\varsigma_j+C'^2|\Sigma_j|\varsigma_j^{2\beta}+|\Sigma_j|^{3/2}\varsigma_j+|\Sigma_j|^{1/2}\varsigma_j\big)\le C|\Sigma_j|\bar\varsigma+C\tau\bar\varsigma^\beta
\]
for $j\notin\mathcal L$, and by (b)(ii)--(iv) at most $C\bar\varsigma$ for $j\in\mathcal L$. In the bracket, $|\Sigma_j|^{3/2}\varsigma_j$ bounds
$2C_EB_U\sup|c|$, and $|\Sigma_j|^{1/2}\varsigma_j$ bounds $\sup_{\Sigma_j}\|\Rs g_z\|_{L^q}$, whose instantaneous part is of the size of
$\sup_{\Sigma_j}\|z\|_{L^q}\le C|\Sigma_j|^{1/2}\varsigma_j$; we used $\varsigma_j\le\bar\varsigma\le1$. The persistent terms on $Z_2(t_j)\setminus Z_0(t_j)$,
which include the part $S_s=(1-\theta_0)S$ carried along in Lemma~\ref{lem:CIloc}, contribute at most $C\tau\bar\varsigma$. Summing over
$j\le N_J<\tau^{-\varepsilon}$ and using $\sum|\Sigma_j|\le1$ and $\#\mathcal L\le N_{\rm c}$,
\[
\|R_1\|_{L^1L^r}\le CN_J\delta_\star+C(N_{\rm c}+1)\bar\varsigma+C\tau^{1-\varepsilon}\bar\varsigma^\beta\le\delta
\]
by \eqref{eq:tol}, after decreasing $\tau$.
\item (W2): the $w_j$ have disjoint supports, and the $S_{c,j}=\theta_0S_j$ are supported in the disjoint intervals $Z_1'(t_j)$. On $Z_1'(t_j)$
the current stress $S_j$ is $\bar R$ plus the accumulated persistent stresses of pass~2, and the latter are bounded in $L^\infty L^q$ by
$\bar\varsigma$, as in the recursion for $\varsigma_j$. Since $|Z_1'(t_j)|=\frac{3\tau}2$ and $N_J<\tau^{-\varepsilon}$,
\[
\Big\|\sum_jw_j\Big\|_{L^2}^2\le M^2\sum_j\|S_{c,j}\|_{L^1_{t,x}}\le M^2\Big(\|\bar R\|_{L^1L^r}+\tfrac32\tau^{1-\varepsilon}\bar\varsigma\Big)\le M^2\big(2C_r\|R\|_{L^1L^r}+o(1)\big),
\]
by \eqref{eq:pass1}, where $o(1)\to0$ as $T,\tau,\delta_\star\to0$. All other contributions are small. The cut-off window correctors $\chi_i\mathring z_i$
of both passes have time supports in the windows, and at any time at most two of them (one from each pass) are nonzero. By
Proposition~\ref{prop:segcorr}(a)(i), $\sup\|z_i\|_{L^2}\le\sup\|z_i\|_{X_q}\le C'T^{1/2}$, so, adding squares over fewer than $2T^{-1}$ windows of length at
most $T$,
\[
\Big\|\sum_i\chi_i\mathring z_i\Big\|_{L^2_{t,x}}^2\le2\sum_i|\Sigma_i|\sup_{\Sigma_i}\|z_i\|_{L^2}^2\le2\cdot2T^{-1}\cdot T\cdot C'^2T=4C'^2T .
\]
Here the disjointness of the supports is needed: the sum of the norms, $\sum_i|\Sigma_i|^{1/2}\sup\|z_i\|_{L^2}$, is only bounded
by $C'$. The long segment correctors contribute $C(N_{\rm c}+1)\bar\varsigma^\beta$, and the means $CN_J\delta_\star+C\bar\varsigma$. Hence (W2) holds with $M_W:=M\sqrt{2C_r}$.
\item (W3): the window correctors contribute at most $CT^{1/2}$. Here (a)(iii) is not enough: it gives $\|z_i\|_{L^1(\Sigma_i;W^{1,\infty})}\le T^{1/2}\|z_i\|_{L^2(\Sigma_i;W^{1,\infty})}\le C'T$,
and the sum over fewer than $2T^{-1}$ windows is only bounded. Instead we use (a)(i): $\sup_{\Sigma_i}\|z_i\|_{W^{1,q}}+\sup_{\Sigma_i}\|z_i\|_{L^\infty}\le C\sup_{\Sigma_i}\|z_i\|_{X_q}\le C'T^{1/2}$.
At any time at most two cut-off window correctors are nonzero, and $\|\chi_i\mathring z_i\|\le2\|z_i\|$ in these norms, so
$\big\|\sum_i\chi_i\mathring z_i\big\|_{L^\infty(0,1;L^\infty\cap W^{1,q})}\le4C'T^{1/2}$, which bounds the $L^pL^\infty$ and $L^1W^{1,q}$ norms; the re-solved segments contribute at most
$C'\bar\varsigma^\beta$ for the short ones, for the same reason, and $N_{\rm c}C'\bar\varsigma^\beta$ for the long ones; the $w_j$ contribute $N_J\delta_\star$; the means contribute $CN_J\delta_\star+C\bar\varsigma$.
\item (W5): $\tilde I^c$ meets only the good segments of pass~1 and the re-solved segments. The former contribute $C(T^{1/2})^\beta$. For
the latter, by (a)(iii), (b)(iii) and since $L^2$ norms over disjoint intervals add in squares,
\[
\Big\|\sum_jz_j\Big\|_{L^2(\tilde I^c;W^{1,\infty})}^2\le\sum_{j\notin\mathcal L}\big(C'|\Sigma_j|^{1/2}\varsigma_j^\beta\big)^2+\sum_{j\in\mathcal L}\big(C'\varsigma_j^\beta\big)^2\le(N_{\rm c}+1)(C'\bar\varsigma^\beta)^2.
\]
Together with the means, this is at most $\delta$ by \eqref{eq:tol}.
\item (W4) holds since nothing changes before $t_{1,0}-\frac{5\tau}4>\inf\tilde I$.
\end{itemize}
The order of choices is: $u,R,\tilde I,\tilde\tau,\delta$ (hence $N_{\rm c}$); then the level-1 bounds and constants $C_1',\ell_{0,1}',c_{*,1}'$;
then the level-2 bounds and constants $C_2',\ell_{0,2}',c_{*,2}',K_2^\sharp,C_2^\sharp$; then $T$ and $\tau$; then pass~1, which fixes
$\bar u$ and $\bar R$; then $N_J$ and $\delta_\star$; then the frequencies of the convex integration steps, in time order. Every constant
depends only on objects fixed earlier. The large gradients of the convex integration perturbations do not enter, since all
kinetic estimates use only the flow of the base and the $L^1_tL^q$ norm of the perturbation.
\qed

\subsection{The first step}

The initial approximate solution of Lemma~\ref{lem:start} is exact on $[0,t_0]$ but its stress may reach $t=1$, so it is not
well-prepared. The following variant of Proposition~\ref{prop:W} produces a well-prepared solution.

\begin{prop}[Preparation]\label{prop:prep}
Let $(u_0,p_0,R_0)$ be a smooth VNS--Reynolds solution on $[0,1]$ with $R_0=0$ on $[0,t_0]$, $t_0\in(0,1)$, and let $\delta>0$. Then there is a
smooth VNS--Reynolds solution $(u_1,p_1,R_1)$, well-prepared for some $(I_1,\tau_1)$ with $I_1\subset(\frac{t_0}2,1)$, such that (W1)--(W3) of
Proposition~\ref{prop:W} hold with $(u,R)=(u_0,R_0)$, and $u_1=u_0$ on $[0,\frac{t_0}2]$.
\end{prop}

\begin{proof}
We run the proof of Proposition~\ref{prop:W} with a single chain and no good segments. Only the changes are listed.

\emph{Geometry.} Choose $\tau$ and $T=\tau^\varepsilon$ with $5\tau<T<\frac{t_0}4$ and $2T\le\ell_0'$, where $\ell_0'$ is fixed as in the paragraph
``Uniformity of constants'' with $(u,R)=(u_0,R_0)$ and $N_{\rm c}:=0$. Let $t_1<\dots<t_N$ be the points of $T\ZZ$ in $[\frac{t_0}2+T,1-T]$.
Then $t_1\le\frac{t_0}2+2T<t_0$ and $t_N>1-2T$. The windows are $\Sigma_i=[t_{i-1},t_i]$, $2\le i\le N$, and the \emph{last window}
$\Sigma_{N+1}=[t_N,1]$, of length less than $2T\le\ell_0'$; together they cover $[t_1,1]\supset\Supp_tR_0$. With the zones of
Section~\ref{sec:geom}, put $I_1:=\bigcup_iW(t_i)$ and $\tau_1:=\tau$. Since $t_1-\frac{5\tau}2>\frac{t_0}2$ and $t_N+\frac{5\tau}2<1$, we have
$I_1\subset(\frac{t_0}2,1)$, and Lemma~\ref{lem:geom}(b)--(e) hold with the same proofs.

\emph{Pass 1.} For $i=2,\dots,N$ apply Proposition~\ref{prop:segcorr}(a) on $\Sigma_i$ with the two-sided cutoff of Proposition~\ref{prop:W}, and
exit by the mean. On $\Sigma_{N+1}$ apply Proposition~\ref{prop:segcorr}(a) with a one-sided cutoff vanishing on $[t_N,t_N+\frac\tau4]$ and equal to
$1$ on $[t_N+\frac\tau2,1]$; there is no exit. The base stress on each window is $R_0$ plus the accumulated persistent stress, so
$\varsigma\le\|R_0\|_{L^\infty L^q}+C$ and $N_S\le C$ by Lemma~\ref{lem:telescope}, as in pass~1 of Proposition~\ref{prop:W}. The result $(\bar u,\bar R)$
has $\bar R$ supported in $\bigcup_iZ_0(t_i)$, since the stress vanishes where the cutoffs equal $1$, in particular on
$[t_N+\frac\tau2,1]$, and \eqref{eq:pass1} holds with $R=R_0$, by the same summation over fewer than $T^{-1}+1$ windows. Nothing is changed
before $t_1$.

\emph{Pass 2.} For $j=1,\dots,N$: apply Lemma~\ref{lem:CIloc} at $t_j$, and then Proposition~\ref{prop:segcorr}(a) on
$\Sigma_j':=[t_j+\tau,t_{j+1}-\frac{9\tau}8]$ for $j<N$, and on $\Sigma_N':=[t_N+\tau,1]$ with a one-sided cutoff. All these segments have length at
most $2T\le\ell_0'$, so (b) is never used and $\mathcal L=\emptyset$. The paragraphs ``Supports'', ``Smallness'' and ``Conclusion'' of the proof of
Proposition~\ref{prop:W} apply verbatim with $N_{\rm c}=0$ and give (W1)--(W3). (W5) is not claimed.

\emph{Well-preparedness and (W4).} After pass~2 the stress vanishes on the bulks of the re-solved segments, on $[t_N+\frac{9\tau}8,1]$ (where the
last cutoff equals $1$), and on $[0,t_1-\frac{5\tau}4]$, where nothing has changed and $R_0=0$ because $t_1<t_0$. Elsewhere it is supported in
$\bigcup_jZ_2(t_j)$, so $(u_1,p_1,R_1)$ is well-prepared for $(I_1,\tau_1)$ by Lemma~\ref{lem:geom}(b), (c). The number of intervals in $I_1$ is
$N<T^{-1}=\tau_1^{-\varepsilon}$, and $0,1\notin I_1$. Finally $u_1=u_0$ on $[0,t_1-\frac{5\tau}4]\supset[0,\frac{t_0}2]$.
\end{proof}

\subsection{Proof of Theorem~\ref{thm:regint}}

\proofstep{First step}
Let $(u_0,p_0,R_0)$ be the initial approximate solution of Lemma~\ref{lem:start}, exact on $[0,t_0]$. Proposition~\ref{prop:prep} with $\delta=\delta_0$
gives a solution well-prepared for some $(I_1,\tau_1)$ with $I_1\subset(\frac{t_0}2,1)$, equal to $u_0$ on $[0,\frac{t_0}2]$.

\proofstep{Iteration}
Apply Proposition~\ref{prop:W} with $\delta_n$ chosen as in the proof of Theorem~\ref{thm:main}, including the factor $e^{-\Lambda_n}$ needed
for Lemma~\ref{lem:convergence}. The limit is a weak solution, as in Section~\ref{sec:iteration}. Let $\mathcal S:=\bigcap_nI_n$. Each $I_n$ is
covered by $\tau_n^{-\varepsilon}$ intervals of length $5\tau_n$, with $\tau_n\to0$, so $\dim_{\mathcal H}\mathcal S\le\varepsilon$ \cite[Section~2.6]{CL22}.

\proofstep{Regularity off $\mathcal S$}
Let $\mathcal K\subset(0,1]\setminus\mathcal S$ be a compact interval. Since $\mathcal S$ is closed and $\inf \mathcal K>0$, there is $\eta>0$ such that
$\mathcal K':=[\inf \mathcal K-\eta,\sup \mathcal K]\subset(0,1]\setminus\mathcal S$. Then $\mathcal K'\subset I_m^c$ for all large $m$. By (W5),
$\sum_m\|u_{m+1}-u_m\|_{L^2(\mathcal K';W^{1,\infty})}<\infty$, so $u\in L^2(\mathcal K';W^{1,\infty})$ and the phase-space flow is Lipschitz on $\mathcal K'$; $f$ is its
push-forward by Lemma~\ref{lem:DPL}.

\emph{Trace.} Let $\eta_1\in C^\infty(\RR;[0,1])$ vanish on $(-\infty,\inf \mathcal K-\eta/2]$ and equal $1$ on $[\inf \mathcal K,\infty)$, and set $v:=\eta_1u$ on $\mathcal K'$.
Since $u$ is a weak solution, $v$ is a weak solution of the Oseen problem \eqref{eq:oseen} on $\mathcal K'$, with drift $b=u\in L^2(\mathcal K';L^\infty)$,
zero initial value, and forcing
\[
G=\eta_1F[f,u]+\eta_1'u\in L^2(\mathcal K';L^s)\qquad\text{for every }s<\infty,
\]
because $F[f,u]$ is bounded and $u\in L^2(\mathcal K';L^\infty)$. By Propositions~\ref{prop:oseenstrong} and~\ref{prop:duality}, $v$ is the strong
solution, so $v\in\mathrm{MR}_s(\mathcal K')$ and, by Lemma~\ref{lem:MR}(b), $u=v\in C(\mathcal K;X_s)\cap L^2(\mathcal K;W^{2,s})$ for every $s<\infty$. In particular
$u(\inf \mathcal K)\in X_s$ for every $s$.

\emph{Uniqueness.} Apply Corollary~\ref{cor:roughtrace} on $\mathcal K$, with the initial state $(u(\inf \mathcal K),f(\inf \mathcal K))$; the velocity datum lies
in $X_s$ for some $s>d$, and the kinetic datum is bounded and vanishes for $|v|$ large. This gives the uniqueness statement.
Infinitely many distinct solutions are obtained as in Corollary~\ref{cor:nonunique}.
\qed

\section{Strong continuations and precise statements}\label{sec:statements}

Theorem~\ref{thm:main} assumes smooth data and a smooth local solution. Here we remove both assumptions. The kinetic
datum is only bounded, and the velocity datum only needs to generate a solution that becomes strong after some time.
Only Lipschitz bounds on the velocity are ever used:
\begin{itemize}
\item the kinetic estimates of Section~\ref{sec:kinetic} depend on the base field through its Lipschitz norm;
\item the convex integration step can be run on a mollified background, at the price of a small error in the stress.
\end{itemize}
From now on, the kinetic datum $f_0\in L^\infty(\TT^d\times\RR^d)$ is only assumed nonnegative and to vanish for $|v|\ge \Rv_0$; it need
not be smooth. The quantities $M_0,M_\infty,J_0$ are defined as in Section~\ref{sec:notation}.

\begin{defin}[Strong continuation]\label{def:strongcont}
Let $u_{\rm in}\in L^2(\TT^d)$ be weakly divergence-free. A \emph{strong continuation} of the data $(u_{\rm in},f_0)$ is a weak
solution $(\bar u,\bar f)$ of \eqref{eq:VNS} on $[0,T_0]$, $0<T_0\le1$, together with times $0\le t_1\le a<T_0$ and an exponent
$q_0\in(d,\infty)\cap[2,\infty)$, such that:
\begin{enumerate}[label=(S\arabic*)]
\item $\bar u\in C_w([0,T_0];L^2)\cap L^1(0,T_0;W^{1,1})\cap L^1(0,T_0;L^\infty)$, and $\bar u\in C([t_1,T_0];L^2)\cap L^2(t_1,T_0;L^\infty)$;
\item $\nabla\bar u\in L^1(t_1,T_0;L^\infty)$;
\item $\bar u\in\mathrm{MR}_{q_0}([a,T_0])$.
\end{enumerate}
Data admitting a strong continuation are called \emph{admissible}.
\end{defin}

Condition (S2) asks for a Lipschitz velocity only after time $t_1$. Before $t_1$ the solution may be rough. By (S1) and
Lemma~\ref{lem:DPL}, $\bar f$ is the unique bounded solution of the kinetic equation with field $\bar u$, and it is continuous
with values in $L^1$.

\begin{theo}[Nonuniqueness for admissible data]\label{thm:rough}
Let $d\ge2$, $1\le p<2$, $1\le q<\infty$. Let $(u_{\rm in},f_0)$ be admissible, with strong continuation
$(\bar u,\bar f,t_1,a,T_0,q_0)$. Then there are infinitely many distinct weak solutions $(u,f)$ of \eqref{eq:VNS} on
$[0,1]$ with data $(u_{\rm in},f_0)$ such that:
\begin{enumerate}[label=(\arabic*)]
\item $u\in L^2_{t,x}\cap L^1(0,1;L^\infty)$ and $u\in L^p(t_1,1;L^\infty)\cap L^1(t_1,1;W^{1,q})$; if moreover
$\bar u\in L^p(0,t_1;L^\infty)$, then $u\in L^p(0,1;L^\infty)$;
\item $(u,f)=(\bar u,\bar f)$ on $[0,t_1+\tau]$, for some $\tau>0$ that does not depend on the solution;
\item $f$ is bounded, vanishes for $|v|$ large, belongs to $C([0,1];L^s)$ for every $s<\infty$, and is the unique bounded
(hence renormalized) solution of the kinetic equation with field $u$;
\item the total momentum is conserved.
\end{enumerate}
Infinitely many of these solutions are not of Leray--Hopf type.
\end{theo}

\begin{theo}[Admissible data]\label{thm:admissible}
Let $f_0\in L^\infty(\TT^d\times\RR^d)$ be nonnegative and vanish for $|v|\ge \Rv_0$, and let $u_{\rm in}$ be divergence-free.
The data $(u_{\rm in},f_0)$ are admissible in each of the following cases:
\begin{enumerate}[label=(\alph*)]
\item $d\ge2$ and $u_{\rm in}\in B^1_{q_0,2}(\TT^d)$ for some $q_0\in(d,\infty)\cap[2,\infty)$. This includes $H^s$ for $s>\frac d2$ and $C^{1,\alpha}$
for $\alpha>0$.
\item $d\in\{2,3\}$ and $u_{\rm in}\in H^1(\TT^d)$.
\item $d=2$ and $u_{\rm in}\in L^2(\TT^2)$;
\item $d=3$ and $u_{\rm in}\in L^2(\TT^3)$.
\end{enumerate}
In cases (a) and (b) one can take $t_1=0$. In cases (a)--(c), $\bar u\in L^2(0,T_0;L^\infty)$.
\end{theo}

\begin{remark}\label{rem:rough}
\begin{enumerate}[label=(\roman*)]
\item In case (c), Theorem~\ref{thm:rough} gives nonuniqueness for \emph{every} finite-energy velocity datum in two
dimensions. By \cite[Theorem~1]{HMMM20}, stated there for $\RR^2$ and $\TT^2$, there is at most one weak solution in the sense of
\cite[Definition~1]{HMMM20}; the proof gives uniqueness on every interval $[0,T]$. The admissibility condition (4) of
\cite[Definition~2]{HMMM20} holds: if $f_0$ is bounded and vanishes for $|v|\ge \Rv_0$, then the supremum of $f_0$ over $\TT^2\times B(e^tv,R)$ is at most
$\|f_0\|_{L^\infty}\mathbf 1_{\{e^t|v|\le R+\Rv_0\}}$, so $(1+|v|^2)\sup_{\TT^2\times B(e^tv,R)}f_0$ is bounded in $L^1(\RR^2_v)$ uniformly in $t\ge0$. Now let $(u,f)$ be
a weak solution in the sense of Definition~\ref{def:weak} with $u\in L^\infty(0,1;L^2)\cap L^2(0,1;H^1)$. Then $\rho_f$ and $j_f$ are bounded,
so $F[f,u]\in L^2(0,1;L^2)$, and $u\otimes u\in L^2(0,1;L^2)$ by Ladyzhenskaya's inequality. From Definition~\ref{def:weak}(iv),
$\partial_tu\in L^2(0,1;H^{-1}_{\rm div})$, so $u\in C([0,1];H^{-1}_{\rm div})$ and the weak formulation holds for every $t$ in the form required
by \cite[Definition~1]{HMMM20}; moreover $f\in L^\infty(0,1;L^1\cap L^\infty)$. Hence $(u,f)$ is a solution on $[0,1]$ in the sense of
\cite[Definition~1]{HMMM20}, that is, it satisfies the requirements of that definition, which is stated for global solutions, on
$[0,1]$. The proof of \cite[Theorem~1]{HMMM20} is local in time, as noted above, and applies to such solutions. No energy inequality
is needed. Consequently at most one of the solutions of Theorem~\ref{thm:rough} lies in
$L^\infty(0,1;L^2)\cap L^2(0,1;H^1)$, in particular at most one is of Leray--Hopf type.
\item Condition (S2) only concerns times after $t_1$. For $d=3$ we expect critical data, such as $u_{\rm in}\in H^{1/2}$, to be
admissible: a Fujita--Kato type solution with bounded drag force is smooth for $t>0$. We do not pursue this.
\item Case (d) rests on the bound $u\in L^1_{\rm loc}L^\infty$ for Leray solutions of the three-dimensional system, due to
Han-Kwan, Miot, Moussa and Moyano \cite[Lemma~B.1]{HMMM25}; Lemma~\ref{lem:L1Linf} gives a self-contained proof on $\TT^3$. It keeps the velocity support of $f$ bounded before $t_1$. For
$d\ge4$ this bound already fails for the heat flow of a general $L^2$ datum, since $\|e^{t\Delta}u_{\rm in}\|_{L^\infty}\sim t^{-d/4}$. The
$H^1$ theory of Proposition~\ref{prop:LWPH1} is also restricted to $d\le3$. Arbitrary $L^2$ data for $d\ge4$ remain open.
\item In case (d) the solutions of Theorem~\ref{thm:rough} only satisfy $u\in L^1_tL^\infty$ near $t=0$. This is the
integrability of the Leray solution they coincide with initially.
\end{enumerate}
\end{remark}

\section{The construction on Lipschitz backgrounds}\label{sec:rough}

\subsection{Kinetic lemmas for Lipschitz fields}

Let $\Vc_{\rm s}$ be the space of measurable $U:(0,1)\times\TT^d\to\RR^d$ such that $U(t)$ is Lipschitz for a.e.\ $t$ and
$\|U\|_{L^1W^{1,\infty}}:=\int_0^1\|U(t)\|_{W^{1,\infty}}\dd t<\infty$; we write $\Vc_{\rm s}=L^1(0,1;W^{1,\infty})$ with this meaning. Such a $U$ is a
Carath\'eodory function, so $t\mapsto U(t)\in C(\TT^d)$ is weakly measurable; since $C(\TT^d)$ is separable, Pettis' theorem shows that it is
strongly measurable, and $U\in L^1(0,1;C(\TT^d))$ as a Bochner function. For $U\in\Vc_{\rm s}$ the field $b_U=(v,U-v)$ is Lipschitz in $z$, with constant
$1+\|\nabla U(t)\|_{L^\infty}\in L^1(0,1)$, and grows linearly. So it has a unique Carath\'eodory flow $Z_U$, bi-Lipschitz in $z$,
with $\det D_zZ_U(t;s,\cdot)=e^{-d(t-s)}$; see the proof of Proposition~\ref{prop:critreg}(c). We set
$f[U](t):=Z_U(t;0,\cdot)_\#f_0$, and keep the notation $\rho[U],j[U],J[U]$ of Section~\ref{sec:notation}.

\begin{lem}\label{lem:kinstrong}
The following hold.
\begin{enumerate}[label=(\alph*)]
\item Lemma~\ref{lem:char} holds for $U\in\Vc_{\rm s}$, except for the statements about $C^1$ and smooth regularity of $f$.
In their place: $f[U]$ is the unique bounded weak solution of \eqref{eq:kinetic}, it belongs to $C([0,1];L^1)$, and
$J[U]\in W^{1,1}(0,1)$ satisfies \eqref{eq:Jprime} a.e.
\item The second inequality in \eqref{eq:dual} holds for $u,\tilde w\in\Vc_{\rm s}$, and for every $\psi$ that is continuous on $\TT^d\times\RR^d$ and
Lipschitz on $\TT^d\times B_{\Rv_1+1}$, with $[\psi]_{\Rv_1}$ replaced by $[\psi]_{\Rv_1+1}$, defined through the essential suprema of $|\nabla_x\psi|$
and $|\nabla_v\psi|$ on $\TT^d\times B_{\Rv_1+1}$.
\item Lemma~\ref{lem:weakR} holds for every $h\in L^1(\TT^d;\RR^d)$, with $\Rs h:=K*h\in L^r$.
\item Proposition~\ref{prop:coupling} holds for $u,\tilde w\in\Vc_{\rm s}$, without assuming that $h$ is smooth, with
$L:=\|u\|_{L^1(0,1;W^{1,\infty})}$.
\item Lemma~\ref{lem:convergence} holds for $u_n\in\Vc_{\rm s}$.
\item Let $U\in\Vc_{\rm s}\cap C([0,1];L^2)$ be divergence-free. Then Lemma~\ref{lem:mean} holds with the following change. In (a),
$c\in C([0,1];\RR^d)$, $\fint(U+c)\in W^{1,1}(0,1)$, and the mean equation holds a.e. If moreover $\fint U\in W^{1,1}(0,1)$, then
$c\in W^{1,1}(0,1)$. In (c), $u,w\in\Vc_{\rm s}\cap C([0,1];L^2)$.
\end{enumerate}
\end{lem}

\begin{remark}[Bounded data at a later time]\label{rem:kinbounded}
Lemmas~\ref{lem:DPL} and~\ref{lem:kinstrong}(a),(b) also hold, with the same proofs and with constants computed from the norms of
the datum, for every bounded nonnegative kinetic datum vanishing for $|v|\ge \Rv$, prescribed at any time $t_0\in[0,1)$, on any interval
$[t_0,T]\subset[0,1]$; only the $C^1$ statements are lost. Indeed, the push-forward argument does not use the smoothness of $f_0$, and for
a fixed flow both sides of the bound in (b) are linear in the datum and continuous with respect to its $L^1$ norm, so (b) follows by
approximating the datum by smooth ones. Translation in time and restriction to subintervals are used below without comment.
\end{remark}

\begin{proof}
\proofstep{Approximation}
Extend $U\in\Vc_{\rm s}$ by zero outside $[0,1]$, and let $U_\varepsilon$ be its mollification in $(t,x)$. Since $U\in L^1(0,1;C(\TT^d))$ as a
Bochner function, $U_\varepsilon\to U$ in $L^1(0,1;L^\infty)$. Moreover, since $\nabla U_\varepsilon=(\nabla U)_\varepsilon$ and by Minkowski's inequality,
\[
\|U_\varepsilon\|_{L^1L^\infty}\le\|U\|_{L^1L^\infty},\qquad \|\nabla U_\varepsilon\|_{L^1L^\infty}\le\|\nabla U\|_{L^1L^\infty}=:\Lambda .
\]
By Gronwall's inequality,
\[
|Z_{U_\varepsilon}(t;s,z)-Z_U(t;s,z)|\le e^{1+\Lambda}\int_s^t\|U_\varepsilon-U\|_{L^\infty}\to0,
\]
uniformly in $(t,s,z)$. Also let $f_0^\eta\in C^\infty_c$ satisfy $0\le f_0^\eta\le M_\infty$, $f_0^\eta\to f_0$ in $L^1$, and
$f_0^\eta=0$ for $|v|\ge \Rv_0+\eta$.

\proofstep{(a)}
The push-forward identities, the Jacobian and the support bound pass to the limit from $(U_\varepsilon,f_0^\eta)$; alternatively,
they follow directly from the Carath\'eodory flow, as in Proposition~\ref{prop:critreg}(c). The push-forward is a weak
solution: for $\varphi\in C^\infty_c([0,1)\times\TT^d\times\RR^d)$ and every $z$, the function $t\mapsto\varphi(t,Z_U(t;0,z))$ is absolutely continuous with
derivative $(\partial_t\varphi+v\cdot\nabla_x\varphi+(U-v)\cdot\nabla_v\varphi)(t,Z_U(t;0,z))$ for a.e.\ $t$, by the chain rule for the
Carath\'eodory solution $Z_U$; integrating in $t$ and against $f_0$, and using the push-forward identity, gives
Definition~\ref{def:weak}(iii). Uniqueness and $C([0,1];L^1)$ follow from Lemma~\ref{lem:DPL}, since
$U\in L^1W^{1,1}\cap L^1L^\infty$. For the momentum, $J[U](t)=\int f_0(z)V(t;0,z)\dd z$, where $V(\cdot;0,z)$ is absolutely continuous
with $\dot V=U(t,X)-V$, bounded by an integrable function on $\Supp f_0$. Hence $J[U]\in W^{1,1}$ and
\[
J[U]'=\int f_0\big(U(t,X)-V\big)=\iint(U-v)f[U] .
\]

\proofstep{(b)}
First let $\psi\in C^1$. Apply Lemma~\ref{lem:dual} to $(u_\varepsilon,\tilde w_\varepsilon)$ with datum $f_0^\eta$. The left side equals
$\int f_0^\eta\big(\psi(Z_{u_\varepsilon+\tilde w_\varepsilon})-\psi(Z_{u_\varepsilon})\big)$, which converges as $\varepsilon,\eta\to0$. On the right side,
$\Lambda_\varepsilon\le\Lambda$, $M_0^\eta\to M_0$, and $\Rv_{1,\varepsilon}\le \Rv_1+\eta$. Moreover $[\psi]_{\Rv_1+\eta}\to[\psi]_{\Rv_1}$, since $\nabla\psi$ is continuous.

Now let $\psi$ be continuous, and Lipschitz on $\TT^d\times B_{\Rv_1+1}$. Only the values of $\psi$ on $\TT^d\times B_{\Rv_1}$ enter the left side,
since all kinetic states involved are supported there. Let $\psi_\kappa:=\psi*\eta_\kappa$ be a mollification in $z$, $\kappa\in(0,1)$. Then $\psi_\kappa\in C^1$,
$\psi_\kappa\to\psi$ uniformly on $\TT^d\times B_{\Rv_1}$, and $[\psi_\kappa]_{\Rv_1}\le[\psi]_{\Rv_1+\kappa}\le[\psi]_{\Rv_1+1}$. Applying the $C^1$ case to $\psi_\kappa$ and letting
$\kappa\to0$ gives the claim.

\proofstep{(c)}
$K\in L^r$ by \eqref{eq:kernel} and $r<\frac d{d-1}$, so $\Rs h=K*h\in L^r$ by Young's inequality. Let $h_\varepsilon=h*\eta_\varepsilon$. Then
$\|h_\varepsilon\|_{L^1}\le\|h\|_{L^1}$ and $\|h_\varepsilon\|_{(C^1)^*}\le\|h\|_{(C^1)^*}$, because $\phi\mapsto\phi*\check\eta_\varepsilon$ preserves the
$C^1$ norm. Also $\Rs h_\varepsilon=(\Rs h)*\eta_\varepsilon\to\Rs h$ in $L^r$. Apply \eqref{eq:weakR} to $h_\varepsilon$ and let $\varepsilon\to0$.

\proofstep{(d)}
In the proof of Proposition~\ref{prop:coupling}, use (b) and (c) instead of Lemmas~\ref{lem:dual} and~\ref{lem:weakR}, and keep the
time dependence of the constants. Fix $t$ outside the null set where $u(t)$ fails to be Lipschitz; the bounds below then hold for
a.e.\ $t$, which suffices since they are integrated in time. The test function $\psi(x,v)=(v-u(t,x))\cdot\phi(x)$ is only Lipschitz, since $u(t)\in W^{1,\infty}$; it is
admissible in (b), and its gradient bounds below hold for a.e.\ $x$ and all $|v|\le \Rv_1+1$, so that $[\psi]_{\Rv_1+1}\le \Rv_1+1+2\|u(t)\|_{W^{1,\infty}}$; we absorb the $+1$
into $\Rv_1$. At time $t$:
\begin{itemize}
\item $[\psi]_{\Rv_1}\le \Rv_1+2\|u(t)\|_{W^{1,\infty}}$;
\item $\|h_1(t)\|_{(C^1)^*}\le e^{1+L}\big(\Rv_1+2\|u(t)\|_{W^{1,\infty}}\big)M_0W$;
\item $\|h_1(t)\|_{L^1}\le2\big(\Rv_1+\|u(t)\|_{L^\infty}\big)M_0$.
\end{itemize}
Integrate in $t$, using H\"older's inequality in $t$ for the product $A_0(t)^{1-\theta}B_0(t)^\theta$. This gives \eqref{eq:coupling} with a
constant depending on $L=\|u\|_{L^1W^{1,\infty}}$. Note that $h\in L^1(0,1;L^\infty)$.

\proofstep{(e)}
The proof of Lemma~\ref{lem:convergence} uses only (a), (b) and the weak formulation.

\proofstep{(f)}
Steps 1, 2, 4 and 5 of the proof of Lemma~\ref{lem:mean} use only (a) and (b). Continuity of $\mathcal T(c)$ uses that
$\fint U\in C([0,1])$, which holds since $U\in C([0,1];L^2)$. In Step~3, the fixed-point identity gives
$\fint(U+c)=\fint U(0)+J_0-J[U+c]\in W^{1,1}$. Differentiating with (a) gives the mean equation a.e.
\end{proof}

\begin{defin}[Strong VNS--Reynolds solutions]\label{def:sVNSR}
A strong solution of the VNS--Reynolds system on $[0,1]$ is a pair $(u,R)$ such that:
\begin{itemize}
\item $u\in C([0,1];L^2)\cap\Vc_{\rm s}$ is divergence-free;
\item $R\in L^1(0,1;L^1)$ takes values in $\Sz$;
\item for every divergence-free $\varphi\in C^\infty([0,1]\times\TT^d;\RR^d)$ vanishing near $t=1$,
\begin{equation}\label{eq:weak}
\int u(0)\cdot\varphi(0)\dd x+\int_0^1\!\!\int\Big(u\cdot(\partial_t\varphi+\Delta\varphi)+u\otimes u:\nabla\varphi+F\big[f[u],u\big]\cdot\varphi-R:\nabla\varphi\Big)\dd x\dd t=0 .
\end{equation}
\end{itemize}
\end{defin}

Every smooth solution of Definition~\ref{def:VNSR} is a strong solution. Constant test fields give the mean equation. By
Lemma~\ref{lem:kinstrong}(a), $\fint u\in W^{1,1}$ and \eqref{eq:momentum} holds.

\subsection{The convex integration step on a rough background}

\begin{lem}[Navier--Stokes step, rough background]\label{lem:NSrough}
Let $d,p,q,M,r$ be as in Proposition~\ref{prop:NSstep}. Let $u\in C([0,1];L^2)$ be divergence-free, let
$R\in L^1(0,1;L^r)$ take values in $\Sz$, and let $0<s<t_*\le1$ with $R=0$ on $[0,t_*]$. For every $\delta>0$ there are
smooth $w$, a distribution $P$ and $\tilde R\in L^1(0,1;L^r)$ with values in $\Sz$, such that:
\begin{enumerate}[label=(\roman*)]
\item $\partial_tw-\Delta w+\D(u\otimes w+w\otimes u+w\otimes w)+\D R+\nabla P=\D\tilde R$ in $\mathcal D'((0,1)\times\TT^d)$;
\item $\D w=0$, $\fint w(t)=0$ for all $t$, and $w=0$, $\tilde R=0$ on $[0,s]$;
\item $\|w\|_{L^2_{t,x}}\le M\|R\|_{L^1_{t,x}}^{1/2}$;
\item $\|w\|_{L^p_tL^\infty}+\|w\|_{L^1_tW^{1,q}}\le\delta$;
\item $\|\tilde R\|_{L^1(0,1;L^r)}\le\delta$.
\end{enumerate}
\end{lem}

\begin{proof}
Let $\eta\in C^\infty_c((0,1)\times B_1)$ be nonnegative with $\int\eta=1$, and $\eta_\ell(\tau,y)=\ell^{-1-d}\eta(\tau/\ell,y/\ell)$. Extend $u$ and $R$
by zero to $t<0$, and define the \emph{causal} mollifications
\[
u_\ell(t,x)=\iint u(t-\tau,x-y)\eta_\ell(\tau,y)\dd\tau\dd y,\qquad R_\ell(t,x)=\iint R(t-\tau,x-y)\eta_\ell(\tau,y)\dd\tau\dd y .
\]
Since $\eta_\ell$ is supported in $\tau\in(0,\ell)$, $R_\ell(t)$ depends only on $R$ on $(t-\ell,t)$. Hence $R_\ell=0$ on $[0,t_*]$
and $R_\ell$ is smooth with values in $\Sz$. Moreover $u_\ell$ is smooth and divergence-free, with
$\sup_t\|u_\ell(t)\|_{L^2}\le\sup_t\|u(t)\|_{L^2}=:N$, and $\|R_\ell\|_{L^1_{t,x}}\le\|R\|_{L^1_{t,x}}$, $R_\ell\to R$ in $L^1(0,1;L^r)$.
Fix $\ell$ with $\|R-R_\ell\|_{L^1L^r}\le\delta/3$. Apply Proposition~\ref{prop:NSstep} to $(u_\ell,R_\ell)$ with $s,t_*$ and
$\delta'':=\min\{\delta,\delta/(24(N+1))\}$, obtaining $w,P_\ell,R^{\rm CL}$. Set
\[
\tilde R:=R^{\rm CL}+(R-R_\ell)+(u-u_\ell)\ot w+w\ot(u-u_\ell),\qquad P:=P_\ell-\tfrac2d(u-u_\ell)\cdot w .
\]
Since $a\otimes b+b\otimes a=a\ot b+b\ot a+\frac2d(a\cdot b)\mathrm{Id}$, the trace part of $(u-u_\ell)\otimes w+w\otimes(u-u_\ell)$ is a gradient,
which is moved to the left-hand side; identity (i) then follows from Proposition~\ref{prop:NSstep}(i). The
matrix $\tilde R$ is symmetric and traceless.

(ii) and (iv) are inherited, since $R-R_\ell=0$ on $[0,t_*]$. (iii) follows from $\|R_\ell\|_{L^1_{t,x}}\le\|R\|_{L^1_{t,x}}$. For (v), since
$r\le2$ and $|\TT^d|=1$,
\[
\|(u-u_\ell)\ot w+w\ot(u-u_\ell)\|_{L^1L^r}\le4\sup_t\|u-u_\ell\|_{L^2}\|w\|_{L^1L^\infty}\le8N\delta''\le\delta/3 .
\]
Together with $\|R^{\rm CL}\|_{L^1L^r}\le\delta''\le\delta/3$, this gives (v).
\end{proof}

\begin{prop}[VNS step, strong solutions]\label{prop:VNSrough}
Let $d,p,q,M,r$ be as in Proposition~\ref{prop:NSstep}. Let $(u,R)$ be a strong solution of the VNS--Reynolds system with
$R\in L^1(0,1;L^r)$, and let $0<s<t_*\le1$ with $R=0$ on $[0,t_*]$. For every $\delta\in(0,1]$ there is a strong solution
$(u_1,R_1)$ with $R_1\in L^1(0,1;L^r)$ satisfying (i)--(iv) of Proposition~\ref{prop:VNSstep}.
\end{prop}

\begin{proof}
Follow Section~\ref{sec:proofstep} with these replacements:
\begin{itemize}
\item $L:=\|u\|_{L^1(0,1;W^{1,\infty})}$;
\item Lemma~\ref{lem:NSrough} instead of Proposition~\ref{prop:NSstep}, with $R^{\rm CL}$ replaced by $\tilde R$;
\item Lemma~\ref{lem:kinstrong}(a), (b), (d), (f) instead of Lemmas~\ref{lem:char}, \ref{lem:dual}, \ref{lem:mean} and Proposition~\ref{prop:coupling}.
\end{itemize}
The mean correction $c$ lies in $W^{1,1}(0,1)$, since $\fint(u+w)=\fint u\in W^{1,1}$. Hence $u_1=u+w+c\in C([0,1];L^2)\cap\Vc_{\rm s}$, and
\eqref{eq:cprime} holds a.e. Step~3 is an identity between distributions. Tested against divergence-free $\varphi$, it gives
\eqref{eq:weak} for $(u_1,R_1)$. Here we use that $\int c\otimes c:\nabla\varphi\dd x=0$, that the trace part of
$(u+w)\otimes c+c\otimes(u+w)$ is a gradient, and that $\int_0^1\!\!\int c\cdot\partial_t\varphi=-\int_0^1c'\cdot\fint\varphi$, since $c(0)=0$.
The initial values agree since $w=c=0$ near $t=0$. The estimates of Step~4 are unchanged, with
$\|c\ot(u+w)+(u+w)\ot c\|_{L^1L^r}\le4\sup|c|\,(L+1)$.
\end{proof}

\subsection{Concatenation}

\begin{lem}\label{lem:concat}
Let $0<t_1<T\le1$. Let $(u^1,f^1)$ be a weak solution on $[0,t_1]$ with data $(u_{\rm in},f_0)$, such that
$u^1\in C_w([0,t_1];L^2)\cap L^1(0,t_1;W^{1,1})$. Let $(u^2,f^2)$ be a weak solution on $[t_1,T]$ with data $(u^1(t_1),f^1(t_1))$,
in the sense of Definition~\ref{def:weak} translated in time. Then the concatenation is a weak solution on $[0,T]$ with data
$(u_{\rm in},f_0)$.
\end{lem}

\begin{proof}
The integrability requirements of Definition~\ref{def:weak} hold on each piece. By Lemma~\ref{lem:DPL}, the
function $f^1$ is continuous from $[0,t_1]$ to $L^1$.
Let $\varphi$ be a test function on $[0,T]$. Let $\zeta_n\in C^\infty(\RR;[0,1])$ be nonincreasing, equal to $1$ on
$(-\infty,t_1-\frac2n]$ and to $0$ on $[t_1-\frac1n,\infty)$. Test the equations of the first piece with $\zeta_n\varphi$, which is admissible
and divergence-free if $\varphi$ is. The terms with $\zeta_n'$ are
\[
\int\zeta_n'\iint f^1\varphi\dd z\dd t\to-\iint f^1(t_1)\varphi(t_1)\dd z,\qquad
\int\zeta_n'\int u^1\cdot\varphi\dd x\dd t\to-\int u^1(t_1)\cdot\varphi(t_1)\dd x,
\]
by continuity of $f^1$ and weak continuity of $u^1$ at $t_1$; for the latter, write
$u^1(t)\cdot\varphi(t)-u^1(t_1)\cdot\varphi(t_1)=u^1(t)\cdot(\varphi(t)-\varphi(t_1))+(u^1(t)-u^1(t_1))\cdot\varphi(t_1)$ and use
$\sup_t\|u^1(t)\|_{L^2}<\infty$. All other terms converge by dominated convergence. The resulting identities on
$[0,t_1]$ contain the boundary terms at $t_1$ with a minus sign. The identities on $[t_1,T]$ contain the same terms as
initial data. Adding the two gives Definition~\ref{def:weak}(iii) and (iv) on $[0,T]$.
\end{proof}

\subsection{Proof of Theorem~\ref{thm:rough}}

\proofstep{Step 0: reduction to $t_1=0$}
Let $(\bar u,\bar f,t_1,a,T_0,q_0)$ be the strong continuation. By (S2), the definition of weak solutions and
Lemma~\ref{lem:kinstrong}(a), on $[t_1,T_0]$ the function $\bar f$ is the push-forward of $\bar f(t_1)$ along the Carath\'eodory
flow of $\bar u$. Moreover $\bar f(t_1)\in L^\infty$ is nonnegative and vanishes for $|v|$ large. Translate time by $t_1$: the
data become $(\bar u(t_1),\bar f(t_1))$, and the translated solution lives on $[0,T_*]$ with $T_*:=T_0-t_1\le1-t_1$.

All results of Sections~\ref{sec:iteration}--\ref{sec:kinetic} and of Section~\ref{sec:rough} hold with the kinetic datum
$\bar f(t_1)$, by Lemma~\ref{lem:kinstrong}. We construct solutions on $[0,1]$, restrict them to $[0,1-t_1]$, translate back, and
concatenate them with $(\bar u,\bar f)|_{[0,t_1]}$ by Lemma~\ref{lem:concat}. When $t_1>0$, the hypotheses of Lemma~\ref{lem:concat} are (S1).

We may therefore assume $t_1=0$, and write $(\bar u,\bar f)$ on $[0,T_*]$ with
\[
\bar u\in C([0,T_*];L^2)\cap L^2(0,T_*;L^\infty)\cap L^1(0,T_*;W^{1,\infty})\cap\mathrm{MR}_{q_0}([a,T_*]),\qquad \bar f=f[\bar u].
\]

\proofstep{Step 1: initial approximate solutions}
Fix $a<a'<t_2'<t_2<t_3<T_*$. Choose $\chi\in C^\infty(\RR;[0,1])$ with $\chi=1$ on $(-\infty,t_2]$ and $\chi=0$ on $[t_3,\infty)$. Choose
$\Psi$ smooth, divergence-free, with $\fint\Psi=0$, supported in $t\in(t_3,1-t_1)$, and with
$\|\Psi\|_{L^p(0,1-t_1;L^2)}\ge1$. For $k\in\NN$ put $\tilde u_k:=\chi\bar u+k\Psi\in C([0,1];L^2)\cap\Vc_{\rm s}$.

Since $\bar u$ is a weak solution, it satisfies \eqref{eq:momentum} on $[0,T_*]$, and $\fint\bar u\in W^{1,1}$. Lemma~\ref{lem:kinstrong}(f)
gives a correction $m_k\in W^{1,1}$ with $m_k=0$ on $[0,t_2]$. Set $u_{0,k}:=\tilde u_k+m_k$. We claim that there is $R_{0,k}\in L^1(0,1;L^{q_0})$,
vanishing on $[0,t_2]$, such that $(u_{0,k},R_{0,k})$ is a strong solution of the VNS--Reynolds system.

On $(a,1)$ define the function
\[
G:=\partial_tu_{0,k}-\Delta u_{0,k}+u_{0,k}\cdot\nabla u_{0,k}-F\big[f[u_{0,k}],u_{0,k}\big]\in L^1(a,1;L^{q_0}) .
\]
It is integrable: on $(a,T_*)$ we have $\bar u\in\mathrm{MR}_{q_0}\subset C(X_{q_0})\subset L^\infty_{t,x}$, so
$\partial_t\bar u,\Delta\bar u,\bar u\cdot\nabla\bar u\in L^2L^{q_0}$; moreover $m_k'\in L^1$, $\Psi$ is smooth, and $F\in L^1L^\infty$ by
Lemma~\ref{lem:kinstrong}(a).
By causality, $u_{0,k}=\bar u$ and $f[u_{0,k}]=\bar f$ on $[0,t_2]$. Since $\bar u$ is a weak solution, $\int\!\!\int G\cdot\psi=0$ for every
divergence-free $\psi$ supported in $(a,t_2)$, so $\PP G=0$ a.e.\ on $(a,t_2)$. Moreover $\fint G=0$ a.e.\ on $(a,1)$ by
Lemma~\ref{lem:kinstrong}(f).

Let $\zeta\in C^\infty(\RR;[0,1])$ with $\zeta=1$ on $(-\infty,a']$ and $\zeta=0$ on $[t_2',\infty)$, and set $R_{0,k}:=(1-\zeta)\Rs\PP G$. Then $R_{0,k}=0$ on
$[0,t_2]$. Let $\varphi$ be a divergence-free test field vanishing near $t=1$, and split $\varphi=\zeta\varphi+(1-\zeta)\varphi$. Both parts are
divergence-free.
\begin{itemize}
\item For $\zeta\varphi$, which vanishes after $t_2'$, the terms of \eqref{eq:weak} without the stress vanish, since $u_{0,k}=\bar u$ and
$f[u_{0,k}]=\bar f$ on $[0,t_2]$, and $\bar u$ is a weak solution on $[0,T_*]$.
\item For $\psi:=(1-\zeta)\varphi$, which vanishes near $[0,a']$, integration by parts turns these terms into $-\iint G\cdot\psi$. Since $\psi$ is
divergence-free and $\fint G=0$, this equals $-\iint\PP G\cdot\psi=\iint\Rs\PP G:\nabla\psi=\iint R_{0,k}:\nabla\varphi$, because $\zeta$ depends
only on $t$.
\end{itemize}
These terms are linear in the test field, so for $\varphi$ they add up to $\iint R_{0,k}:\nabla\varphi$, and \eqref{eq:weak} holds with $R=R_{0,k}$.

\proofstep{Step 2: iteration}
Run the proof of Theorem~\ref{thm:main} with:
\begin{itemize}
\item $s_n:=t_2\big(\frac12+2^{-n-1}\big)$;
\item Proposition~\ref{prop:VNSrough} instead of Proposition~\ref{prop:VNSstep};
\item Lemma~\ref{lem:kinstrong}(e) instead of Lemma~\ref{lem:convergence}.
\end{itemize}
Every iterate is a strong solution, so $\Lambda_n<\infty$. Step~4 of that proof uses \eqref{eq:weak} instead of \eqref{eq:VNSR}. The
result is a weak solution $(u_k,f_k)$ on $[0,1]$ with data $(\bar u(0),\bar f(0))$. It satisfies properties (1), (3), (4), (5) of
Theorem~\ref{thm:main}, with $\tilde u=\tilde u_k$ and $\varepsilon=\frac14$, and coincides with $(\bar u,\bar f)$ on $[0,t_2/2]$. Property (1) uses
$u_{0,k}\in L^2(0,1;L^\infty)\cap L^1(0,1;W^{1,\infty})$.

\proofstep{Step 3: conclusion}
Restrict to $[0,1-t_1]$, translate back, and concatenate with $(\bar u,\bar f)|_{[0,t_1]}$ by Lemma~\ref{lem:concat}. This gives (2) with
$\tau=t_2/2$. Property (1) on $[0,t_1]$ follows from (S1). For (3), $u\in L^1(0,1;W^{1,1})\cap L^1(0,1;L^\infty)$, so Lemma~\ref{lem:DPL}
applies on $[0,1]$. Momentum is conserved on each piece.

Distinctness and the failure of the energy inequality follow as in the proof of Corollary~\ref{cor:nonunique}. The
fluctuations of $u_k-u_l$ are within $\frac12$ of $(k-l)\Psi$ in $L^p(0,1-t_1;L^\infty)$, after translation. A solution of Leray--Hopf type
satisfies $\|u_k\|_{L^p_tL^2}\le(2\mathcal E(0))^{1/2}$.
\qed

\section{Local strong solutions: proof of Theorem~\ref{thm:admissible}}\label{sec:local}

\begin{prop}[Data in $B^1_{q_0,2}$]\label{prop:LWPX}
Let $d\ge2$ and $q_0\in(d,\infty)\cap[2,\infty)$. Let $u_*\in X_{q_0}$ be divergence-free, and let $f_*\in L^\infty$ be nonnegative and vanish for
$|v|\ge \Rv_*$. There is $T>0$, depending only on $d$, $q_0$ and upper bounds for $\|u_*\|_{X_{q_0}}$, $\Rv_*$ and $\|f_*\|_{L^\infty}$, and there is
$u\in\mathrm{MR}_{q_0}([0,T])$ with $\D u=0$, such that
\begin{equation}\label{eq:mild}
\partial_tu-\Delta u+\QQ(u\cdot\nabla u)=\QQ F\big[f[u],u\big]\quad\text{a.e.},\qquad u(0)=u_* .
\end{equation}
Here $f[u]$ is the push-forward of $f_*$ along $b_u$. The pair $(u,f[u])$ is a weak solution on $[0,T]$ and a strong
continuation of $(u_*,f_*)$ with $t_1=a=0$.
\end{prop}

\begin{proof}
Write $q=q_0$, and let $C_q,C_T,C_E$ be the constants of Lemmas~\ref{lem:MR} and~\ref{lem:MRshift}. For $T\le1$ and $A\ge1$ let
\[
E_T:=\big\{u\in\mathrm{MR}_q([0,T]):\ \D u=0,\ u(0)=u_*,\ \|u\|_{\mathrm{MR}_q}\le A\big\},
\]
with the metric $\|u-u'\|_{\Yo(0,T)}$. The set $E_T$ is complete, as in the proof of Proposition~\ref{prop:segcorr}(a).

\proofstep{Bounds on $E_T$}
Let $u\in E_T$. By Lemma~\ref{lem:MR}(b) and Lemma~\ref{lem:MRshift}(b),
\[
\sup_t\|u(t)\|_{X_q}\le A_1:=C_q(\|u_*\|_{X_q}+A),\qquad \|u\|_{L^\infty_{t,x}}\le C_EA_1,\qquad \int_0^T\|\nabla u\|_{L^\infty}\le C_ET^{1/2}A .
\]
Take $T\le T_1$, with $T_1$ such that $T_1C_EA_1\le1$ and $C_ET_1^{1/2}A\le1$. Then, by Lemma~\ref{lem:kinstrong}(a), $f[u]$ vanishes for
$|v|\ge \Rv_*+1$, and $\rho[u]\le\varrho:=e^d\|f_*\|_\infty|B_{\Rv_*+1}|$. So $|F[f[u],u]|\le\varrho(\Rv_*+1+C_EA_1)$.

\proofstep{Invariance}
Let $\Theta u$ be the solution of $\partial_t\vartheta-\Delta\vartheta=\QQ\big(F[f[u],u]-u\cdot\nabla u\big)$ with $\vartheta(0)=u_*$. It is divergence-free by
Lemma~\ref{lem:MR}(d). Since $\|u\cdot\nabla u\|_{L^q}\le\|u\|_\infty\|\nabla u\|_{L^q}\le C_E^2A_1^2$,
\[
\|\Theta u\|_{\mathrm{MR}_q}\le C_q\|u_*\|_{X_q}+C_qT^{1/2}\big(C_E^2A_1^2+\varrho(\Rv_*+1+C_EA_1)\big).
\]
Choose $A:=2C_q\|u_*\|_{X_q}+1$, and then $T$ small; this gives $\Theta(E_T)\subset E_T$.

\proofstep{Contraction}
Let $u,u'\in E_T$ and $\delta u=u-u'$, so $\delta u(0)=0$. Write
\[
F[f[u],u]-F[f[u'],u']=\int(v-u')\big(f[u]-f[u']\big)\dd v-\rho[u]\,\delta u .
\]
Lemma~\ref{lem:B2}(b), with $t_0=0$, base $b=u'$, $U=u$ and $W=u'$, gives the following bound. Here $K'$ depends only on $\Rv_*$,
$\|f_*\|_\infty$ and $A_1$, since $\Lambda_{u'}\le1$, $\Rv_1\le\Rv_*+2$ and $\sup_t\|u'(t)\|_{W^{1,q}}\le\sup_t\|u'(t)\|_{X_q}\le A_1$:
\[
\Big\|\int(v-u')\big(f[u]-f[u']\big)(t)\dd v\Big\|_{H^{-1,q}}\le K'\|\delta u\|_{L^1(0,t;L^q)}\le K'\|\delta u\|_{L^1(0,t;L^\infty)} .
\]
By Lemma~\ref{lem:MRshift},
\[
\|\delta u\|_{L^1(0,T;L^\infty)}\le C_ET^{1/2}\|\delta u\|_{\Yo},\qquad \sup_t\|\delta u\|_{L^q}\le C_T\|\delta u\|_{\Yo}.
\]
Using $\|\D(u\otimes u-u'\otimes u')\|_{H^{-1,q}}\le C\|u\otimes u-u'\otimes u'\|_{L^q}$ and
Lemma~\ref{lem:MRshift}(a),
\[
\|\Theta u-\Theta u'\|_{\Yo}\le C_q\Big(T^{1/2}C\big(C_EA_1+\varrho\big)C_T+K'C_ET\Big)\|\delta u\|_{\Yo}\le\tfrac12\|\delta u\|_{\Yo}
\]
for $T$ small. The resulting $T$ depends only on the quantities listed.

\proofstep{Conclusion}
The fixed point $u$ solves \eqref{eq:mild}. Since $\QQ$ is self-adjoint and $\QQ\varphi=\varphi$ for divergence-free $\varphi$, testing
\eqref{eq:mild} gives Definition~\ref{def:weak}(iv). Lemma~\ref{lem:kinstrong}(a) gives (iii). Properties (S1)--(S3) hold:
\[
\mathrm{MR}_q\subset C(X_q)\cap L^2(0,T;W^{1,\infty}),\qquad X_q\subset W^{1,q}\subset L^\infty .
\]
\end{proof}

\begin{prop}[Data in $H^1$, $d\le3$]\label{prop:LWPH1}
Let $d\in\{2,3\}$. Let $u_*\in H^1$ be divergence-free, and let $f_*\in L^\infty$ be nonnegative and vanish for $|v|\ge \Rv_*$. There are
$T>0$ and a weak solution $(u,f[u])$ on $[0,T]$ with data $(u_*,f_*)$, such that
\[
u\in C([0,T];L^2)\cap L^\infty(0,T;H^1)\cap L^2(0,T;H^2),\qquad \partial_tu\in L^2(0,T;L^2),\qquad \nabla u\in L^1(0,T;L^\infty).
\]
Moreover $u(t)\in X_4$ for a.e.\ $t\in(0,T)$, and $(u,f[u])$ satisfies the energy inequality of Definition~\ref{def:LH} on
$[0,T]$. One can take
\begin{equation}\label{eq:Tlower}
T=c_*\big(1+\|u_*\|_{H^1}^2\big)^{-2},
\end{equation}
where $c_*>0$ depends only on $d$ and on upper bounds for $\|u_*\|_{L^2}$, $\Rv_*$, $\|f_*\|_{L^\infty}$ and $\|f_*\|_{L^1}$; in particular $c_*$
can be taken nonincreasing in each of these quantities.
\end{prop}

\begin{proof}
Constants $C$ depend only on $d$, $\|u_*\|_{H^1}$, $\Rv_*$, $\|f_*\|_{L^\infty}$ and $\|f_*\|_{L^1}$, unless indicated otherwise.

\proofstep{Step 1: approximation}
Let $u_*^\varepsilon$ be the Fourier truncations of $u_*$. They are smooth, divergence-free, satisfy $\|u_*^\varepsilon\|_{H^1}\le\|u_*\|_{H^1}$, and converge
to $u_*$ in $H^1$. Since $4>d$, Proposition~\ref{prop:LWPX} with $q_0=4$ gives solutions $u^\varepsilon\in\mathrm{MR}_4$ of \eqref{eq:mild} with data
$(u_*^\varepsilon,f_*)$ on some interval. We show that they exist on an interval $[0,T]$ independent of $\varepsilon$, with uniform bounds.

For a solution on $[0,t]$ write
\[
\Rv^\varepsilon(t):=\Rv_*+\int_0^t\|u^\varepsilon\|_{L^\infty},\qquad \varrho(t):=e^d\|f_*\|_\infty|B_{\Rv^\varepsilon(t)}| .
\]
Then $f[u^\varepsilon](t)$ vanishes for $|v|\ge \Rv^\varepsilon(t)$, $\rho[u^\varepsilon](t)\le\varrho(t)$, and $\|f[u^\varepsilon](t)\|_\infty\le e^d\|f_*\|_\infty$.

\proofstep{Step 2: uniform a priori bounds}
Let $u=u^\varepsilon$ be an $\mathrm{MR}_4$ solution on $[0,t]$ with $t\le T$ and $\Rv(t)\le \Rv_*+1$. Write $\varrho_1:=\varrho|_{\Rv=\Rv_*+1}$.
\begin{itemize}
\item \emph{Energy.} The computation in the proof of Proposition~\ref{prop:critreg}(d) gives the energy equality. Hence
$\|u(s)\|_{L^2}^2\le\|u_*\|_{L^2}^2+\Rv_*^2\|f_*\|_{L^1}\le C$.
\item \emph{$H^1$ bound.} Test \eqref{eq:mild} with $-\Delta u$, which is divergence-free. By the Gagliardo--Nirenberg inequalities in
$d\le3$, $|\int(u\cdot\nabla u)\cdot\Delta u|\le\frac14\|\Delta u\|_2^2+CY^3$, where $Y:=1+\|u\|_{H^1}^2$. Also $\|F\|_{L^2}\le\varrho_1(\Rv_*+1+\|u\|_{L^2})\le C$.
Hence
\[
Y'+\|\Delta u\|_{L^2}^2\le CY^3+C .
\]
The constant $C$ here depends only on $\|u_*\|_{L^2}$, $\Rv_*$, $\|f_*\|_{L^\infty}$ and $\|f_*\|_{L^1}$, not on $\|u_*\|_{H^1}$. Comparison
with the Riccati equation gives $Y\le2Y(0)$ for $s\le T_2:=c\,Y(0)^{-2}$. Integrating, $\int_0^s\|\Delta u\|_2^2\le Y(0)+CsY(0)^3\le2Y(0)$ for
$s\le T_2$, after decreasing $c$.
\item \emph{Velocity support.} $H^2\subset L^\infty$ for $d\le3$, so $\int_0^s\|u\|_{L^\infty}\le Cs^{1/2}\big(\int_0^s\|u\|_{H^2}^2\big)^{1/2}\le Cs^{1/2}Y(0)^{1/2}$.
Hence $\Rv(s)\le \Rv_*+\frac12$ for $s\le T_3:=c\,Y(0)^{-1}$.
\end{itemize}
Set $T:=\min\{T_2,T_3,1\}$. Thus every $\mathrm{MR}_4$ solution on $[0,t]$ with $t\le T$ and $\Rv(t)\le\Rv_*+1$ satisfies the bounds above, and in
particular $\Rv(t)\le\Rv_*+\frac12$. In addition:
\begin{itemize}
\item \emph{Lipschitz bound.} Write $G:=F-u\cdot\nabla u$. Then
$\|G\|_{L^6}\le\|u\|_\infty\|\nabla u\|_{L^6}+\|F\|_{L^6}\le C\|u\|_{H^2}^2+C(1+\|u\|_{H^2})$, so $\|G\|_{L^1(0,t;L^6)}\le C$. On $\TT^d$, for
$0<s\le1$,
\[
\|\nabla e^{s\Delta}g\|_\infty\le Cs^{-d/4}\|\nabla g\|_{L^2},\qquad \|\nabla e^{s\Delta}\QQ g\|_\infty\le Cs^{-\frac12-\frac d{12}}\|g\|_{L^6}.
\]
The exponents $\frac d4$ and $\frac12+\frac d{12}$ are less than $1$. From the Duhamel formula, $\int_0^t\|\nabla u\|_{L^\infty}\le C$.
\item \emph{$X_4$ bound (non-uniform).} Let $y(s):=\sup_{[0,s]}\|u\|_{X_4}^2$. By Lemma~\ref{lem:MR}(a), (b) and
$\|u\cdot\nabla u\|_{L^4}\le\|u\|_\infty\|\nabla u\|_{L^4}\le Cy^{1/2}\|u\|_{H^2}$,
\[
y(s)\le C\big(\|u_*^\varepsilon\|_{X_4}^2+1\big)+C\int_0^sy(\sigma)\|u(\sigma)\|_{H^2}^2\dd\sigma .
\]
Gronwall's inequality gives $y\le B_\varepsilon<\infty$ on $[0,t]$.
\end{itemize}

\proofstep{Step 3: existence on $[0,T]$}
Fix $\varepsilon$. Let $\mathcal T$ be the set of $t\in(0,T]$ such that \eqref{eq:mild} with data $(u_*^\varepsilon,f_*)$ has a solution
$u\in\mathrm{MR}_4([0,t])$ with $\Rv(t)\le\Rv_*+1$. Two such solutions coincide on their common interval: both are weak solutions with
velocity in $L^2L^\infty$ and datum $u_*^\varepsilon\in X_4$, $4>d$, so Corollary~\ref{cor:roughtrace}(c) applies. Hence the solutions for
$t\in\mathcal T$ are restrictions of a single function. By Proposition~\ref{prop:LWPX}, $\mathcal T\ne\emptyset$, since $\Rv$ is continuous with
$\Rv(0)=\Rv_*$. For $t_0\in\mathcal T$, Step~2 gives $\Rv(t_0)\le\Rv_*+\frac12$ and $\|u(t_0)\|_{X_4}^2\le B_\varepsilon$, where $B_\varepsilon$ does not depend
on $t_0$. Proposition~\ref{prop:LWPX}, applied from $t_0$ with data $(u(t_0),f[u](t_0))$, extends the solution to $[t_0,t_0+\tau_\varepsilon]$. Here
$\|f[u](t_0)\|_\infty\le e^d\|f_*\|_\infty$, the velocity support of $f[u](t_0)$ has radius at most $\Rv_*+1$, and $\tau_\varepsilon>0$ depends only on
$B_\varepsilon$, $\Rv_*$ and $\|f_*\|_\infty$. The kinetic state of the extension is the push-forward of $f_*$ from time $0$, by the flow
property, and the concatenation is in $\mathrm{MR}_4$, since the traces match (Lemma~\ref{lem:MR}(b)). On the extension
$\Rv\le\Rv_*+1$ persists up to $\min\{t_0+\tau_\varepsilon,T\}$. Otherwise there is a first time $t_\sharp$ with $\Rv(t_\sharp)=\Rv_*+1$; then $t_\sharp\in\mathcal T$,
and Step~2 gives $\Rv(t_\sharp)\le\Rv_*+\frac12$, a contradiction. So $\min\{t_0+\tau_\varepsilon,T\}\in\mathcal T$ for every $t_0\in\mathcal T$. Since
$\tau_\varepsilon$ does not depend on $t_0$, this gives $T\in\mathcal T$.

\proofstep{Step 4: limit}
By Step~2 and \eqref{eq:mild}, $(u^\varepsilon)$ is bounded in $L^\infty H^1\cap L^2H^2$ and $\partial_tu^\varepsilon$ is bounded in $L^2L^2$. By the
Aubin--Lions lemma, a subsequence converges in $L^2(0,T;H^{7/4})\cap C([0,T];L^2)$ to some $u$ with the stated regularity. Since
$H^{7/4}\subset L^\infty$, $u^\varepsilon\to u$ in $L^1L^\infty$. Along a further subsequence, $\nabla u^\varepsilon(t)\to\nabla u(t)$ in $L^2$ for a.e.\ $t$. By
lower semicontinuity of $\|\cdot\|_{L^\infty}$ under this convergence and Fatou's lemma, $\int_0^T\|\nabla u\|_\infty\le C$. Hence
$u\in\Vc_{\rm s}$ on $[0,T]$, and Lemma~\ref{lem:kinstrong}(b), with base $u$, gives
\[
\sup_t\sup_{\|\psi\|_{C^1_b}\le1}\Big|\iint\big(f[u^\varepsilon]-f[u]\big)(t)\psi\Big|\le e^{1+C}M_0\|u^\varepsilon-u\|_{L^1L^\infty}\to0 .
\]
The velocity supports are uniformly bounded. As in Lemma~\ref{lem:convergence}(d), the drag forces converge against continuous
test functions. Since $u^\varepsilon\to u$ in $L^2_{t,x}$, we may pass to the limit in the weak formulation of \eqref{eq:mild}, so $(u,f[u])$ is a
weak solution.

Finally, $(u,f[u])$ satisfies the energy equality, hence the energy inequality. The computation in the proof of
Proposition~\ref{prop:critreg}(d) applies: it only uses $\partial_tu\in L^2L^2$, $u\in L^2H^2$ and the Lipschitz phase-space flow.

\proofstep{Step 5: $X_4$ at a.e.\ time}
$u(t)\in H^2$ for a.e.\ $t$. Moreover $H^2=B^2_{2,2}\hookrightarrow B^{2-d/4}_{4,2}\hookrightarrow B^1_{4,2}=X_4$, since $d\le3$ \cite{Tri83}.
\end{proof}

\begin{proof}[Proof of Theorem~\ref{thm:admissible}]
(a) This is Proposition~\ref{prop:LWPX}, with $u_*=u_{\rm in}$ and $f_*=f_0$. The examples embed into $X_{q_0}$ for a suitable $q_0>d$:
$H^s=B^s_{2,2}\hookrightarrow B^{s-d/2+d/q_0}_{q_0,2}\hookrightarrow B^1_{q_0,2}$ if $s\ge1+\frac d2-\frac d{q_0}$, which can be arranged for
$s>\frac d2$ with $q_0$ close to $d$, and $C^{1,\alpha}\hookrightarrow B^{1+\alpha}_{\infty,\infty}\hookrightarrow B^{1}_{q_0,2}$.

(b) Let $(u,f)$ be given by Proposition~\ref{prop:LWPH1} on $[0,T]$; replacing $T$ by $\min\{T,\frac12\}$, we may assume $T\le\frac12$. Choose
$a\in(0,T)$ with $u(a)\in X_4$. Note that $f(a)$ is bounded
and vanishes for $|v|$ large. Proposition~\ref{prop:LWPX}, applied from time $a$ with data $(u(a),f(a))$, gives an $\mathrm{MR}_4$
solution on $[a,a+T']$. Concatenate by Lemma~\ref{lem:concat}, using $u\in C([0,a];L^2)\cap L^1W^{1,1}$. The concatenation is a
strong continuation with $t_1=0$, $q_0=4$ and $T_0=\min\{a+T',\frac12\}>a$. Restricting to a shorter final time preserves (S1)--(S3), and
$T_0\le\frac12$. Indeed, (S1) holds since $L^2H^2\subset L^2L^\infty$ for $d\le3$ and the pieces match continuously in $L^2$ at $a$;
(S2) holds on both pieces; and (S3) holds on $[a,T_0]$.

(c) and (d). Let $d\in\{2,3\}$. By Theorem~\ref{thm:SEI}, there is an energy solution $(u,f)$ on $[0,1]$ with data $(u_{\rm in},f_0)$ which
satisfies the strong energy inequality \eqref{eq:SEI}. By Lemma~\ref{lem:L1Linf} if $d=3$, and by Lemma~\ref{lem:L2Linf2D} if $d=2$,
$u\in L^1(0,1;L^\infty)$, $f$ is bounded and vanishes for $|v|\ge \Rv_0+\|u\|_{L^1L^\infty}$, and $(u,f)$ is a weak solution in the sense of
Definition~\ref{def:weak}; if $d=2$, moreover $u\in L^2(0,1;L^\infty)$. Since
$u\in C_w([0,1];L^2)\cap L^2(0,1;H^1)$, (S1) holds on $[0,t_1]$ for every $t_1$.

Choose $t_1\in(0,\frac12)$ such that $u(t_1)\in H^1$ and \eqref{eq:SEI} holds with $s=t_1$; almost every $t_1$ has these properties. Note
that $f(t_1)$ is bounded and vanishes for $|v|$ large. Apply (b) from time $t_1$ to the data $(u(t_1),f(t_1))$; it gives a solution on
$[t_1,t_1+T_0']$ with $T_0'\le\frac12$ (see the proof of (b)), so that $T_0:=t_1+T_0'<1$. Concatenate by Lemma~\ref{lem:concat}, which needs
$u\in C_w([0,t_1];L^2)\cap L^1(0,t_1;W^{1,1})$; this holds by (E1). The result is a strong continuation $(\bar u,\bar f)$ on $[0,T_0]$ with this
$t_1$. It
satisfies $\bar u\in C([t_1,T_0];L^2)\cap L^2(t_1,T_0;L^\infty)$ by the construction of (b): Proposition~\ref{prop:LWPH1} on the first piece
and Proposition~\ref{prop:LWPX} on the second. If $d=2$, then also
$\bar u=u\in L^2(0,t_1;L^\infty)$, so $\bar u\in L^2(0,T_0;L^\infty)$, as claimed for case (c).

In fact $(\bar u,\bar f)=(u,f)$ on $[0,T_0]$, so the strong continuation is a restriction of the Leray solution. Indeed, on
$[t_1,T_0]$ the velocity $\bar u$ satisfies hypothesis (a) of Lemma~\ref{lem:wsT3}: on the first piece of the construction of (b) by
Proposition~\ref{prop:LWPH1}, on the second piece because $\mathrm{MR}_4\subset\mathrm{MR}_2\cap L^2W^{1,\infty}$, and the two pieces match in
$H^1$. The restriction of $(u,f)$ to $[t_1,T_0]$ is a weak solution with data $(u(t_1),f(t_1))$, as in the proof of
Lemma~\ref{lem:concat}, and it satisfies hypothesis (b) of Lemma~\ref{lem:wsT3} by the choice of $t_1$.
\end{proof}

\begin{remark}
As explained in Section~\ref{sec:ideas}, the kinetic comparisons start at $t_1$ from the common kinetic state $\bar f(t_1)$, so that
no regularity of $f_0$ and of $u_{\rm in}$ beyond (S1) is needed.
\end{remark}

\section{A structure theorem in three dimensions}\label{sec:structure}

The same ingredients give a Leray-type structure theorem for energy solutions of the three-dimensional system. Let $d\in\{2,3\}$. We
say that $(u,f)$ is a \emph{Leray solution with the strong energy inequality} if it is an energy solution on $[0,1]$ in the sense of
Definition~\ref{def:energy} and, for $s=0$ and almost every $s\in(0,1)$, and for all $t\in[s,1]$,
\begin{equation}\label{eq:SEI}
\mathcal E(t)+\int_s^t\Big(\|\nabla u\|_{L^2}^2+\iint|v-u|^2f\dd z\Big)\dd\sigma\le\mathcal E(s).
\end{equation}
Since \eqref{eq:SEI} is required for every $t$, we fix representatives: $u$ is the representative in $C_w([0,1];L^2)$, given by (E1), and $f$
is the representative in $C([0,1];L^1)$ given by Lemma~\ref{lem:L1Linf}(c) or Lemma~\ref{lem:L2Linf2D}(c). By Lemma~\ref{lem:L1Linf}
if $d=3$, and by Lemma~\ref{lem:L2Linf2D} if $d=2$, such a solution is a weak solution in the sense of
Definition~\ref{def:weak}: $u\in L^1(0,1;L^\infty)$, $f$ has
bounded velocity support, and $f\in C([0,1];L^1)$ is the unique bounded, hence renormalized, solution of the kinetic equation with
field $u$. By Theorem~\ref{thm:SEI}, Leray solutions with the strong energy inequality exist for all data considered here.

\begin{theo}[Structure of Leray solutions, $d=3$]\label{thm:structure}
Let $d=3$, let $u_{\rm in}\in L^2(\TT^3)$ be divergence-free, and let $f_0\in L^\infty$ be nonnegative and vanish for $|v|\ge \Rv_0$. Let
$(u,f)$ be a Leray solution with the strong energy inequality and data $(u_{\rm in},f_0)$. Then there is an open set
$\mathcal O\subset(0,1)$ such that:
\begin{enumerate}[label=(\alph*)]
\item the singular set $\mathcal S:=(0,1)\setminus\mathcal O$ satisfies $\mathcal H^{1/2}(\mathcal S)=0$; in particular $|\mathcal S|=0$;
\item on every compact interval $\mathcal K\subset\mathcal O$, $u\in L^\infty(\mathcal K;H^1)\cap L^2(\mathcal K;H^2)$ and $\nabla u\in L^1(\mathcal K;L^\infty)$; $f$ is transported on $\mathcal K$ by a
bi-Lipschitz phase-space flow;
\item for almost every $t_0\in(0,1)$, the interval $(t_0,t_0+T(t_0))\cap(0,1)$ is contained in $\mathcal O$, where
$T(t_0)=c\,\big(1+\|u(t_0)\|_{H^1}^2\big)^{-2}$ and $c>0$ depends only on the data.
\end{enumerate}
\end{theo}

\begin{proof}
\proofstep{Step 1: uniform kinetic bounds}
By Lemma~\ref{lem:L1Linf}, $K_1:=\|u\|_{L^1(0,1;L^\infty)}<\infty$. For every $t$, $f(t)$ vanishes for
$|v|\ge \Rv_1:=\Rv_0+K_1$, $\|f(t)\|_{L^\infty}\le e^3\|f_0\|_{L^\infty}$, $\|f(t)\|_{L^1}=\|f_0\|_{L^1}$, and $\|u(t)\|_{L^2}^2\le2\mathcal E(0)$
by \eqref{eq:SEI} with $s=0$. Let $c>0$ be the constant $c_*$ of \eqref{eq:Tlower}, computed with the upper bounds $(2\mathcal E(0))^{1/2}$ for $\|u(t)\|_{L^2}$, $\Rv_1$,
$e^3\|f_0\|_\infty$ and $\|f_0\|_{L^1}$; since $c_*$ is nonincreasing in these bounds, it applies at every $t_0\in G$. It depends only on the data.

\proofstep{Step 2: local identification}
Let $G\subset(0,1)$ be the set of times $t_0$ such that $u(t_0)\in H^1$ and \eqref{eq:SEI} holds with $s=t_0$. It has full measure. Fix
$t_0\in G$ and put $T(t_0)=c(1+\|u(t_0)\|_{H^1}^2)^{-2}$. Proposition~\ref{prop:LWPH1}, applied from time $t_0$ with data $(u(t_0),f(t_0))$, gives
a solution $(\tilde u,\tilde f)$ on $J:=[t_0,t_0+T(t_0)]\cap[0,1]$. It satisfies $\tilde u\in L^\infty H^1\cap L^2H^2$,
$\partial_t\tilde u\in L^2L^2$ and $\nabla\tilde u\in L^1L^\infty$, hence $\tilde u\in L^2(J;L^\infty)\cap L^1(J;W^{1,\infty})$, and the energy
inequality from $t_0$.

In particular $(\tilde u,\tilde f)$ satisfies hypothesis (a) of Lemma~\ref{lem:wsT3} on $J$. The restriction $(u,f)|_J$ is a weak solution
on $J$ with data $(u(t_0),f(t_0))$, in the sense of Definition~\ref{def:weak} translated in time: this follows from the weak formulations
on $[0,1]$ by a cutoff in time, using the weak continuity of $u$ and the continuity of $f$ in $L^1$, as in the proof of
Lemma~\ref{lem:concat}. Since $t_0\in G$, it satisfies hypothesis (b) of Lemma~\ref{lem:wsT3}, with $u_*=u(t_0)$ and $f_*=f(t_0)$, which is bounded
and vanishes for $|v|\ge \Rv_1$. By Lemma~\ref{lem:wsT3}, $(u,f)=(\tilde u,\tilde f)$ on $J$.

\proofstep{Step 3: the open set}
Let $\mathcal O:=\bigcup_{t_0\in G}\big(t_0,t_0+T(t_0)\big)\cap(0,1)$. It is open, and (c) holds by construction. A compact interval $\mathcal K\subset\mathcal O$ is covered by
finitely many of these intervals, on each of which $u=\tilde u$ has the regularity of Step~2. Hence (b) holds, with the flow of
Lemma~\ref{lem:kinstrong}(a) started at $\inf \mathcal K$; this flow exists because $\nabla u\in L^1(\mathcal K;L^\infty)$ on the whole interval $\mathcal K$, and $f$ is
transported by it because, by the uniqueness part of Lemma~\ref{lem:kinstrong}(a) (see Remark~\ref{rem:kinbounded}), $f$ is the only bounded
weak solution on $\mathcal K$ with datum $f(\inf \mathcal K)$. For a
compact set $\mathcal K$ that is not an interval, the regularity of $u$ still holds on $\mathcal K$, but $[\inf \mathcal K,\sup \mathcal K]$ may meet $\mathcal S$.

\proofstep{Step 4: the singular set}
Let $t\in\mathcal S$. For every $s\in G\cap(0,t)$ we have $t\notin(s,s+T(s))$, so $T(s)\le t-s$, that is,
\[
\|\nabla u(s)\|_{L^2}^2\ge c^{1/2}(t-s)^{-1/2}-1-2\mathcal E(0)\qquad\text{for a.e. }s\in(0,t).
\]
Hence there is $\varepsilon_0>0$, depending only on $c$ and $\mathcal E(0)$, such that for $\varepsilon\le\min\{\varepsilon_0,t\}$,
\begin{equation}\label{eq:Leray}
\int_{t-\varepsilon}^t\|\nabla u\|_{L^2}^2\dd s\ge2c^{1/2}\varepsilon^{1/2}-(1+2\mathcal E(0))\varepsilon\ge c^{1/2}\varepsilon^{1/2}.
\end{equation}

Fix $\delta\in(0,\varepsilon_0)$, and let $\mathcal S_\delta$ be the $\delta$-neighborhood of $\mathcal S\cap[\delta,1)$. Consider the closed intervals $[t-\varepsilon,t]$ with
$t\in\mathcal S\cap[\delta,1)$ and $\varepsilon<\delta$; each point of $\mathcal S\cap[\delta,1)$ belongs to one of them (as its right endpoint), and their lengths are
bounded. By the Vitali covering lemma there are pairwise disjoint intervals $I_k=[t_k-\varepsilon_k,t_k]$ of this family, at most countably many,
such that the closed intervals $5I_k$, with the same centers and five times the lengths, cover $\mathcal S\cap[\delta,1)$. The bound
\eqref{eq:Leray} holds for the closed intervals $I_k$ as well, since endpoints have measure zero. Then, by \eqref{eq:Leray},
\[
\mathcal H^{1/2}_{10\delta}\big(\mathcal S\cap[\delta,1)\big)\le\sum_k(5\varepsilon_k)^{1/2}\le5^{1/2}c^{-1/2}\sum_k\int_{I_k}\|\nabla u\|_{L^2}^2\le C\int_{\mathcal S_\delta}\|\nabla u\|_{L^2}^2 .
\]
The right side is at most $C\int_0^1\|\nabla u\|_{L^2}^2<\infty$, so $\mathcal H^{1/2}(\mathcal S\cap[\delta,1))<\infty$ for every $\delta$, and hence $|\mathcal S|=0$.
Every point of $\mathcal S_\delta$ lies within $\delta$ of $\mathcal S$, so $\limsup_{\delta\to0}\mathbf 1_{\mathcal S_\delta}\le\mathbf 1_{\overline{\mathcal S}}$ pointwise, and
$|\overline{\mathcal S}|=|\mathcal S|=0$ since $\mathcal S$ is relatively closed in $(0,1)$. By dominated convergence, $\int_{\mathcal S_\delta}\|\nabla u\|_{L^2}^2\to0$. Now fix $\delta'>0$. For
$\delta\le\delta'$, $\mathcal H^{1/2}_{10\delta}(\mathcal S\cap[\delta',1))\le\mathcal H^{1/2}_{10\delta}(\mathcal S\cap[\delta,1))\le C\int_{\mathcal S_\delta}\|\nabla u\|_{L^2}^2$, and letting $\delta\to0$ gives
$\mathcal H^{1/2}(\mathcal S\cap[\delta',1))=0$. Taking the union over $\delta'=1/m$, $m\in\NN$, gives $\mathcal H^{1/2}(\mathcal S)=0$.
\end{proof}

\begin{remark}\label{rem:structure}
\begin{enumerate}[label=(\roman*)]
\item Theorem~\ref{thm:structure} is the analogue of Leray's structure theorem for the Navier--Stokes equations. The
coupling enters at two points. The local existence time of Proposition~\ref{prop:LWPH1} must be controlled uniformly in
time. This requires the bound $u\in L^1_tL^\infty$ of Lemma~\ref{lem:L1Linf}, after \cite[Lemma~B.1]{HMMM25}, which keeps the velocity
support, the density and the drag force bounded along the solution. And the identification needs the weak--strong
uniqueness of Lemma~\ref{lem:wsT3}; see \cite[Theorem~1.2]{HMMM25} for strong solutions beyond the Lipschitz class.
\item The strong energy inequality \eqref{eq:SEI} is needed to restart the comparison at a.e.\ $t_0$. The solutions of
Theorem~\ref{thm:SEI} satisfy it, because the approximate solutions of its Galerkin scheme satisfy the energy equality on every
interval $[s,t]$. We do not know whether every energy solution satisfies it.
\item For the consequences in $d=3$ in combination with Theorem~\ref{thm:admissible}(d), see the discussion after
Theorem~\ref{thm:E}.
\end{enumerate}
\end{remark}

\section{Energy solutions on the torus}\label{sec:HMMMtorus}

This section collects the facts on energy solutions used in Sections~\ref{sec:local} and~\ref{sec:structure}: a velocity
bound (Lemma~\ref{lem:L1Linf}), weak--strong uniqueness for Lipschitz strong solutions (Lemma~\ref{lem:wsT3}), the
existence of Leray solutions satisfying the strong energy inequality (Theorem~\ref{thm:SEI}), and a two-dimensional
velocity bound (Lemma~\ref{lem:L2Linf2D}). The first and the last are periodic versions, for energy solutions, of
\cite[Lemma~B.1]{HMMM25} and \cite[Lemma~3]{HMMM20}. The proofs are self-contained, apart from standard tools (maximal
regularity of the heat equation, the Aubin--Lions lemma, the maximal-function estimate of Lemma~\ref{lem:maximal}) and the
uniqueness of bounded solutions of the kinetic equation for velocities in $L^1_tW^{1,1}\cap L^1_tL^\infty$ (Lemma~\ref{lem:DPL}). The
kinetic lemmas of Sections~\ref{sec:kinetic} and~\ref{sec:rough} are used for bounded data and on subintervals $[t_0,T]$, as explained in
Remark~\ref{rem:kinbounded}. Throughout, $d\in\{2,3\}$ and $f_0$ is bounded,
nonnegative and vanishes for $|v|\ge \Rv_0$; Lemma~\ref{lem:L1Linf} is stated for $d=3$, and the last paragraph treats $d=2$.
Let $\Delta_j$, $j\ge0$, be the periodic Littlewood--Paley blocks,
localizing to $|n|\sim2^j$, and let $P_0$ be the projection on the zero mode.

\paragraph{Velocity moments.} For a nonnegative function $g$ on $\TT^d\times\RR^d$ we write
\[
m_kg(x):=\int_{\RR^d}g(x,v)\,|v|^k\dd v,\qquad M_k(g):=\iint g\,|v|^k\dd z .
\]
If $0\le g\le F$ on $\RR^d_v$ and $0\le k<m$, splitting the integral at $|v|=r$ and optimizing in $r$ gives
$\int g|v|^k\dd v\le C\,F^{\frac{m-k}{m+d}}\big(\int g|v|^m\dd v\big)^{\frac{k+d}{m+d}}$. Hence, if $0\le f\le F_\infty$ on $\TT^d\times\RR^d$,
\begin{equation}\label{eq:interp}
\|m_kf\|_{L^{\frac{m+d}{k+d}}(\TT^d)}\le C\,F_\infty^{\frac{m-k}{m+d}}M_m(f)^{\frac{k+d}{m+d}}\qquad(0\le k<m).
\end{equation}
Since $|\TT^d|=1$, we also have $\|g\|_{L^r}\le\|g\|_{L^s}$ for $r\le s$, and $H^1(\TT^d)\hookrightarrow L^6$ for $d\le3$.

\begin{defin}[Energy solutions]\label{def:energy}
Let $d\in\{2,3\}$, $0<T\le1$, and let $u_{\rm in}\in L^2(\TT^d)$ be weakly divergence-free. An \emph{energy solution} on $[0,T]$ with data
$(u_{\rm in},f_0)$ is a pair $(u,f)$ such that:
\begin{enumerate}[label=(E\arabic*)]
\item $u\in C_w([0,T];L^2)\cap L^2(0,T;H^1)$, and $u(t)$ is weakly divergence-free;
\item $f\in L^\infty((0,T)\times\TT^d\times\RR^d)\cap L^\infty(0,T;L^1(\TT^d\times\RR^d))$, $f\ge0$, and $t\mapsto M_2(f(t))$ is essentially
bounded;
\item the identities of Definition~\ref{def:weak}(iii) and (iv) hold on $[0,T]$, with test functions vanishing near $t=T$;
\item the energy inequality of Definition~\ref{def:LH} holds for a.e.\ $t\in[0,T]$, with $\mathcal E(0):=\frac12\|u_{\rm in}\|_{L^2}^2+\frac12M_2(f_0)$.
\end{enumerate}
\end{defin}

In contrast with Definition~\ref{def:weak}(i)--(ii), neither $u\in L^1_tL^\infty$ nor a bound on the velocity support of $f$ is
assumed; both are consequences (Lemmas~\ref{lem:L1Linf} and~\ref{lem:L2Linf2D}). No renormalization property is assumed either. The
identities in (E3) are meaningful: by (E2) and \eqref{eq:interp} with $m=2$, $\rho_f\in L^\infty(0,T;L^{(d+2)/d})$ and
$j_f\in L^\infty(0,T;L^{(d+2)/(d+1)})$, while $u\in L^2(0,T;L^6)$ by (E1); since $\frac{d+2}d\ge\frac65$, $F[f,u]=j_f-\rho_fu\in L^1((0,T)\times\TT^d)$,
and $f|u|$ is locally integrable.

\subsection{The velocity bound}

\begin{lem}[Velocity bound for energy solutions]\label{lem:L1Linf}
Let $d=3$, and let $(u,f)$ be an energy solution on $[0,T]$, $T\le1$, with data $(u_{\rm in},f_0)$. Then:
\begin{enumerate}[label=(\alph*)]
\item $u\in L^1(0,T;L^\infty)$, with a bound depending only on $T$, $\mathcal E(0)$, $\|f_0\|_{L^1}$, $\|f_0\|_{L^\infty}$ and $\Rv_0$;
\item $0\le f\le e^3\|f_0\|_{L^\infty}$, and $f(t)$ vanishes for $|v|\ge \Rv_0+\|u\|_{L^1(0,T;L^\infty)}$;
\item $(u,f)$ is a weak solution on $[0,T]$ in the sense of Definition~\ref{def:weak}; $f$ is the unique bounded (hence
renormalized) solution of the kinetic equation with field $u$, it belongs to $C([0,T];L^1)$, and $\iint f(t)\dd z=\|f_0\|_{L^1}$ for
every $t$.
\end{enumerate}
\end{lem}

\begin{proof}
Let $F_\infty:=\|f\|_{L^\infty}$. Until Step~8, constants $C$ depend only on $T$, $\mathcal E(0)$, $\|f_0\|_{L^1}$, $\Rv_0$ and $F_\infty$. In
Step~8 we show that $F_\infty\le e^3\|f_0\|_{L^\infty}$, so that in the end they depend only on the quantities listed in (a). By (E1) and
(E4), $\|u\|_{L^\infty L^2}+\|u\|_{L^2H^1}\le C$, hence $\|u\|_{L^2(0,T;L^6)}\le C$, and $M_2(f(t))\le2\mathcal E(0)$ for a.e.\ $t$. By
\eqref{eq:interp} with $m=2$ and $|\TT^3|=1$,
\begin{equation}\label{eq:rhoenergy}
\|\rho_f\|_{L^\infty(0,T;L^{6/5})}\le\|\rho_f\|_{L^\infty(0,T;L^{5/3})}\le C .
\end{equation}
The main
difficulty is the drag force: the pointwise bound $|F|^2\le\rho_f\int|v-u|^2f\dd v$ would require $\rho_f\in L^\infty$, that is, a bound on
the velocity support, which is what we want to prove. We therefore propagate velocity moments, first with the energy alone and
then with an improved integrability of $u$. This replaces \cite[Lemma~B.1]{HMMM25}, where the bound $F\in L^2_{\rm loc}L^2$ is obtained by
referring to \cite[Lemma~4.7]{EHKM21}. The moments are propagated directly from the weak formulation, with truncated weights, so
that no uniqueness or stability theorem for the kinetic equation with an unbounded velocity field is needed.

\proofstep{Step 1: test functions of $v$}
Let $\Phi\in C^1(\RR^3)$ be bounded, with $\nabla\Phi$ compactly supported. We claim that $N_\Phi(t):=\iint f(t)\,\Phi(v)\dd z$ has an
absolutely continuous representative on $[0,T]$ with
\begin{equation}\label{eq:NPhi}
N_\Phi(0)=\iint f_0\,\Phi\dd z,\qquad N_\Phi'(t)=\iint f(t)\,(u-v)\cdot\nabla_v\Phi\dd z\quad\text{for a.e.\ }t .
\end{equation}
By density, the identity of Definition~\ref{def:weak}(iii) holds for test functions in $C^1_c([0,T)\times\TT^3\times\RR^3)$. Let
$\chi\in C^\infty_c(B_2)$ with $\chi=1$ on $B_1$, let $\chi_L:=\chi(\cdot/L)$ for $L\ge1$, and let $\eta\in C^\infty_c([0,T))$. Testing with
$\varphi=\eta(t)\Phi(v)\chi_L(v)$, which does not depend on $x$, gives
\[
\eta(0)\iint f_0\,\Phi\chi_L\dd z+\int_0^T\!\!\iint f\,\Big(\eta'\,\Phi\chi_L+\eta\,\chi_L\,(u-v)\cdot\nabla_v\Phi+\eta\,\Phi\,(u-v)\cdot\nabla_v\chi_L\Big)\dd z\dd t=0 .
\]
As $L\to\infty$, the first two terms converge by dominated convergence, since $f\in L^\infty(0,T;L^1)$, and $\chi_L\nabla\Phi=\nabla\Phi$ for
large $L$. Since $(|u|+|v|)|\nabla\chi_L|\le C(|u|L^{-1}+1)\mathbf 1_{\{|v|\ge L\}}$, the last term is bounded by
\[
C\|\eta\|_{L^\infty}\|\Phi\|_{L^\infty}\Big(L^{-1}\int_0^T\|u\|_{L^6}\|\rho_f\|_{L^{6/5}}\dd t+L^{-2}\int_0^TM_2(f(t))\dd t\Big)\xrightarrow[L\to\infty]{}0,
\]
by \eqref{eq:rhoenergy}. Hence, for every $\eta\in C^\infty_c([0,T))$,
\[
-\int_0^T\eta'N_\Phi\dd t=\eta(0)\iint f_0\,\Phi\dd z+\int_0^T\eta\,G\dd t,\qquad G(t):=\iint f(t)\,(u-v)\cdot\nabla_v\Phi\dd z ,
\]
and $G\in L^1(0,T)$ because $\nabla\Phi$ is bounded with compact support. This is \eqref{eq:NPhi}. With $\Phi\equiv1$ we obtain
$\iint f(t)\dd z=\|f_0\|_{L^1}$ for a.e.\ $t$.

\proofstep{Step 2: truncated moments}
Let $k\in\{3,6\}$ and $p:=k+3$. Fix a nondecreasing $\theta\in C^\infty([0,\infty))$ with $\theta(s)=s^k$ for $s\le2$ and $\theta$ constant
on $[3,\infty)$, and let $\Phi_L(v):=L^k\theta(|v|/L)$ for $L\ge \Rv_0$. Then $\Phi_L\in C^1(\RR^3)$ is bounded and radially nondecreasing,
$\Phi_L(v)=|v|^k$ for $|v|\le2L$, and, with constants depending only on $k$ and $\theta$,
\[
c_k|v|^k\mathbf 1_{\{|v|\le3L\}}\le\Phi_L(v)\le C_k|v|^k,\qquad |\nabla\Phi_L(v)|\le C_k|v|^{k-1}\mathbf 1_{\{|v|\le3L\}} .
\]
Since $\nabla\Phi_L(v)$ is a nonnegative multiple of $v$, we have $(u-v)\cdot\nabla_v\Phi_L\le|u||\nabla\Phi_L|$. By Step~1,
$N_L(t):=\iint f(t)\Phi_L\dd z$ is absolutely continuous and
\[
N_L'\le C_k\int_{\TT^3}|u|\,m^L\dd x,\qquad m^L(t,x):=\int_{\{|v|\le3L\}}f\,|v|^{k-1}\dd v .
\]
Splitting the last integral at $|v|=r$ and using $|v|^k\mathbf 1_{\{|v|\le3L\}}\le c_k^{-1}\Phi_L$, we get
$m^L\le C\big(F_\infty r^{k+2}+r^{-1}\int f\,\Phi_L\dd v\big)$ for every $r>0$, and after optimizing in $r$,
\[
m^L\le C\,F_\infty^{1/p}\Big(\int f\,\Phi_L\dd v\Big)^{\frac{p-1}p},\qquad\text{hence}\qquad
\|m^L(t)\|_{L^{\frac p{p-1}}(\TT^3)}\le C\,F_\infty^{1/p}N_L(t)^{\frac{p-1}p}.
\]
By H\"older's inequality, $N_L'\le C\|u(t)\|_{L^p}N_L^{1-1/p}$ for a.e.\ $t$. Thus, for every $\epsilon>0$, the function
$(N_L+\epsilon)^{1/p}$ is absolutely continuous with derivative at most $C\|u(t)\|_{L^p}$. Since
$N_L(0)\le C_kM_k(f_0)\le C_k\Rv_0^k\|f_0\|_{L^1}$, we get $N_L(t)^{1/p}\le C\big(1+\|u\|_{L^1(0,T;L^p)}\big)$. Since $\Phi_L\ge0$ and $\Phi_L\to|v|^k$ pointwise as
$L\to\infty$, Fatou's lemma, applied along the integers outside a common null set of times, gives
\begin{equation}\label{eq:trunc}
M_k(f(t))\le C\big(1+\|u\|_{L^1(0,T;L^{k+3})}\big)^{k+3}\qquad\text{for a.e.\ }t,\quad k\in\{3,6\},
\end{equation}
whenever the right-hand side is finite.

\proofstep{Step 3: the third moment}
Since $\|u\|_{L^1(0,T;L^6)}\le C$, \eqref{eq:trunc} with $k=3$ gives $M_3(f(t))\le C$ for a.e.\ $t$. By \eqref{eq:interp} with $m=3$ and
$k=0,1$,
\[
\|\rho_f\|_{L^\infty(0,T;L^2)}+\|j_f\|_{L^\infty(0,T;L^{3/2})}\le C,\qquad\text{hence}\qquad
\|F\|_{L^2(0,T;L^{3/2})}\le C+\|\rho_f\|_{L^\infty L^2}\|u\|_{L^2L^6}\le C ,
\]
where $F:=F[f,u]$.

\proofstep{Step 4: an $L^9$ bound for $u$}
The mean $P_0u$ is bounded by $\|u\|_{L^\infty L^2}$. The mean-free part satisfies the heat equation with forcing $\QQ(F-\D(u\otimes u))$
in the weak sense of Definition~\ref{def:weak}(iv), with $u\in L^\infty L^2\cap L^2H^1$. Comparing Fourier coefficients as in the
Decomposition step of the proof of Lemma~\ref{lem:L2Linf2D}, with $e$ ranging over an orthonormal basis of $n^\perp$, we obtain the Duhamel
formula for a.e.\ $t\in[0,T]$ (the integral $B(t)$ below converges for a.e.\ $t$, which suffices for the $L^1_t$ bounds that follow):
\begin{equation}\label{eq:duhamel}
(1-P_0)u(t)=e^{t\Delta}(1-P_0)u_{\rm in}+\int_0^te^{(t-s)\Delta}(1-P_0)\QQ F\dd s-\int_0^te^{(t-s)\Delta}(1-P_0)\QQ\D(u\otimes u)\dd s=:A+D+B .
\end{equation}
On mean-free functions on $\TT^3$, for $0<\tau\le1$, $1<p<\infty$ and $p\le s\le\infty$, the heat semigroup satisfies
\[
\|e^{\tau\Delta}\QQ g\|_{L^s}\le C\tau^{-\frac32(\frac1p-\frac1s)}\|g\|_{L^p},\qquad \|e^{\tau\Delta}\QQ\D g\|_{L^s}\le C\tau^{-\frac12-\frac32(\frac1p-\frac1s)}\|g\|_{L^p},
\]
and the first bound also holds for $p=2$ without $\QQ$. The restriction $p<\infty$ is needed: $\QQ$ is not bounded on $L^\infty$, so the
first bound fails for $p=s=\infty$ as $\tau\to0$. All applications below use finite $p$. With $s=9$, and with $p=2$, $p=\frac32$ and $p=3$ respectively,
\begin{gather*}
\|A(t)\|_{L^9}\le Ct^{-7/12}\|u_{\rm in}\|_{L^2},\\
\|D(t)\|_{L^9}\le C\int_0^t(t-s)^{-5/6}\|F(s)\|_{L^{3/2}}\dd s,\qquad
\|B(t)\|_{L^9}\le C\int_0^t(t-s)^{-5/6}\|u(s)\|_{L^6}^2\dd s .
\end{gather*}
The kernels are integrable, so Young's inequality in time and Step~3 give $\|u\|_{L^1(0,T;L^9)}\le C$.

\proofstep{Step 5: the sixth moment and the drag force}
By Step~4 and \eqref{eq:trunc} with $k=6$,
\begin{equation}\label{eq:M6}
M_6(f(t))\le C\qquad\text{for a.e.\ }t\in[0,T] .
\end{equation}
By \eqref{eq:interp} with $m=6$ and $k=0,1$,
\[
\|\rho_f(t)\|_{L^3}\le C F_\infty^{2/3}M_6(f(t))^{1/3},\qquad \|j_f(t)\|_{L^{9/4}}\le C F_\infty^{5/9}M_6(f(t))^{4/9}.
\]
Since $|\TT^3|=1$ and $H^1(\TT^3)\hookrightarrow L^6$,
\[
\|F(t)\|_{L^2}\le\|j_f(t)\|_{L^{9/4}}+\|\rho_f(t)\|_{L^3}\|u(t)\|_{L^6}\le C\big(1+\|u(t)\|_{H^1}\big),
\]
and $F\in L^2(0,T;L^2)$ follows from $u\in L^2(0,T;H^1)$.

\proofstep{Step 6: the linear terms in $L^1L^\infty$}
We return to \eqref{eq:duhamel}. By the bounds of Step~4 with $s=\infty$ and $p=2$,
\[
\|A(t)\|_{L^\infty}\le Ct^{-3/4}\|u_{\rm in}\|_{L^2},\qquad \int_0^T\|D\|_{L^\infty}\le C\int_0^T\!\!\int_0^t(t-s)^{-3/4}\|F(s)\|_{L^2}\dd s\dd t\le CT^{1/4}\|F\|_{L^1L^2}.
\]

\proofstep{Step 7: the bilinear term in $L^1L^\infty$}
This replaces the Besov lemma of Chemin used in \cite{HMMM25}. We first claim that, for $a,b\in H^1(\TT^3)$,
\begin{equation}\label{eq:prod}
\sum_{j\ge0}2^{j/2}\|\Delta_j(ab)\|_{L^2}\le C\big(\|\nabla a\|_{L^2}+|P_0a|\big)\big(\|\nabla b\|_{L^2}+|P_0b|\big).
\end{equation}
Write $a=P_0a+\mathring a$. Products involving $P_0a$ or $P_0b$ are bounded by
$|P_0a|\sum_j2^{j/2}\|\Delta_j\mathring b\|_{L^2}\le C|P_0a|\|\nabla b\|_{L^2}$, using $2^{j/2}=2^{-j/2}2^j$ and Cauchy--Schwarz. For $\mathring a\mathring b$, use Bony's
decomposition and set $a_k=2^k\|\Delta_k\mathring a\|_{L^2}$, $b_k=2^k\|\Delta_k\mathring b\|_{L^2}$, so that $\|(a_k)\|_{\ell^2}\simeq\|\nabla a\|_{L^2}$.
\begin{itemize}
\item \emph{Paraproduct.} By Bernstein, $\|S_{k-1}\mathring a\|_{L^\infty}\le C\sum_{k'<k}2^{k'/2}a_{k'}$. Hence
\[
\sum_j2^{j/2}\|\Delta_j(T_{\mathring a}\mathring b)\|_{L^2}\le C\sum_kb_k\sum_{k'<k}2^{(k'-k)/2}a_{k'}\le C\|(a_k)\|_{\ell^2}\|(b_k)\|_{\ell^2},
\]
by Young's inequality for sequences; here we used that $\Delta_j(S_{k-1}\mathring a\,\Delta_k\mathring b)=0$ unless $|j-k|\le2$. The
paraproduct $T_{\mathring b}\mathring a$ is bounded in the same way.
\item \emph{Remainder.} By Bernstein from $L^1$ to $L^2$ on the block $j$,
\[
2^{j/2}\|\Delta_j(\Delta_k\mathring a\,\tilde\Delta_k\mathring b)\|_{L^2}\le C2^{2j}\|\Delta_k\mathring a\|_{L^2}\|\tilde\Delta_k\mathring b\|_{L^2}=C2^{2(j-k)}a_k\tilde b_k .
\]
Summing over $j\le k+3$, then over $k$, gives $C\sum_ka_k\tilde b_k$.
\end{itemize}
This proves \eqref{eq:prod}.

For $j\ge0$, $\|\Delta_je^{\tau\Delta}g\|_{L^2}\le e^{-c\tau4^j}\|\Delta_jg\|_{L^2}$ and $\|\Delta_j\QQ\D g\|_{L^2}\le C2^j\|\Delta_jg\|_{L^2}$. Hence
\[
\int_0^T\|\Delta_jB(t)\|_{L^2}\dd t\le C2^j\int_0^T\!\!\int_0^te^{-c(t-s)4^j}\|\Delta_j(u\otimes u)(s)\|_{L^2}\dd s\dd t\le C2^{-j}\int_0^T\|\Delta_j(u\otimes u)\|_{L^2}\dd s .
\]
By Bernstein, $\|\Delta_jB\|_{L^\infty}\le C2^{3j/2}\|\Delta_jB\|_{L^2}$. Summing over $j$ and using \eqref{eq:prod},
\[
\int_0^T\|B\|_{L^\infty}\le C\int_0^T\sum_j2^{j/2}\|\Delta_j(u\otimes u)\|_{L^2}\le C\int_0^T\big(\|\nabla u\|_{L^2}+\|u\|_{L^2}\big)^2<\infty .
\]
Together with $|P_0u|\le\|u\|_{L^\infty L^2}$ and Step~6, this gives $\|u\|_{L^1(0,T;L^\infty)}\le C$.

\proofstep{Step 8: velocity support, the bound on $f$, and uniqueness}
By Step~7 and (E1), $u\in L^1(0,T;L^\infty)\cap L^1(0,T;W^{1,1})$. Extend $u$ by zero outside $[0,T]$ and let $u_\varepsilon$ be its
mollification in $(t,x)$. Then $u_\varepsilon$ is smooth, $\|u_\varepsilon\|_{L^1(0,T;L^\infty)}\le\|u\|_{L^1(0,T;L^\infty)}$ by Minkowski's inequality, and
$u_\varepsilon\to u$ in $L^1((0,T)\times\TT^3)$. Let $f_\varepsilon:=f[u_\varepsilon]$ be the push-forward of $f_0$ along the flow of $(v,u_\varepsilon-v)$. By
Lemma~\ref{lem:kinstrong}(a), $f_\varepsilon$ is a weak solution of the kinetic equation with field $u_\varepsilon$ and datum $f_0$,
\[
0\le f_\varepsilon\le e^3\|f_0\|_{L^\infty},\qquad f_\varepsilon=0\ \text{ for }|v|\ge \Rv_0+\|u\|_{L^1(0,T;L^\infty)} .
\]
Along a subsequence, $f_\varepsilon\rightharpoonup^*g$ in $L^\infty((0,T)\times\TT^3\times\RR^3)$, and $g$ satisfies the same two bounds. For a test function
$\varphi$,
\[
\Big|\iiint f_\varepsilon\,(u_\varepsilon-u)\cdot\nabla_v\varphi\Big|\le e^3\|f_0\|_{L^\infty}\,C_\varphi\,\|u_\varepsilon-u\|_{L^1((0,T)\times\TT^3)}\to0,
\]
while the other terms of the weak formulation pass to the limit because their multipliers are integrable. Hence $g$ is a bounded
weak solution of the kinetic equation with field $u$ and datum $f_0$, and so is $f$, by (E2)--(E3). By Lemma~\ref{lem:DPL},
applied on $[0,T]$, $f=g$. This proves (b), in particular $F_\infty\le e^3\|f_0\|_{L^\infty}$, and it shows that $(u,f)$ satisfies
Definition~\ref{def:weak}(i)--(ii). The remaining statements of (c) follow from Lemma~\ref{lem:DPL} and from the conservation of mass
of Step~1, which extends to every $t$ by the continuity of $f$ in $L^1$.
\end{proof}

\subsection{Weak--strong uniqueness}

The next lemma compares an energy solution with a strong solution whose velocity is Lipschitz in space. It is a special case of the
situation treated in \cite[Theorem~1.2]{HMMM25}, where the strong solution may be much rougher. In the Lipschitz case the Lagrangian
argument of Section~\ref{sec:lagstab} applies directly, on the torus, and we give the proof.

\begin{lem}[Weak--strong uniqueness]\label{lem:wsT3}
Let $d\in\{2,3\}$ and $0\le t_0<T\le1$. Let $u_*\in L^2(\TT^d)$ be divergence-free, and let $f_*\in L^\infty(\TT^d\times\RR^d)$ be nonnegative and vanish for
$|v|\ge \Rv_*$. Let $(u_1,f_1)$ and $(u_2,f_2)$ be weak solutions of \eqref{eq:VNS} on $[t_0,T]$ with data $(u_*,f_*)$ at time $t_0$, in
the sense of Definition~\ref{def:weak} translated in time. Assume that:
\begin{enumerate}[label=(\alph*)]
\item $u_1\in L^\infty(t_0,T;H^1)\cap L^2(t_0,T;H^2)$, $\partial_tu_1\in L^2(t_0,T;L^2)$ and $\nabla u_1\in L^1(t_0,T;L^\infty)$;
\item $u_2\in C_w([t_0,T];L^2)\cap L^2(t_0,T;H^1)$ and, for every $t\in[t_0,T]$,
\[
\mathcal E_2(t)+\int_{t_0}^t\Big(\|\nabla u_2\|_{L^2}^2+\iint|v-u_2|^2f_2\dd z\Big)\dd\sigma\le\tfrac12\|u_*\|_{L^2}^2+\tfrac12M_2(f_*),
\]
where $\mathcal E_2(t)=\frac12\|u_2(t)\|_{L^2}^2+\frac12M_2(f_2(t))$.
\end{enumerate}
Then $(u_1,f_1)=(u_2,f_2)$ on $[t_0,T]$.
\end{lem}

\begin{proof}
We follow Section~\ref{sec:lagstab}. There both velocities are Lipschitz; here only $u_1$ is, and the flow of $u_2$ is replaced by the
flows of mollified fields. Let $w:=u_1-u_2$, $F_i:=F[f_i,u_i]$,
\begin{gather*}
\Rv_\sharp:=\Rv_*+\|u_1\|_{L^1(t_0,T;L^\infty)}+\|u_2\|_{L^1(t_0,T;L^\infty)},\qquad \varrho_\sharp:=e^d\|f_*\|_{L^\infty}|B_{\Rv_\sharp}|,\\
a(t):=1+\|u_1(t)\|_{L^\infty}^2+\|\nabla u_1(t)\|_{L^\infty}.
\end{gather*}
By (a) and $H^2(\TT^d)\hookrightarrow L^\infty$ ($d\le3$), $u_1\in L^2(t_0,T;L^\infty)$; hence $a\in L^1(t_0,T)$ and $u_1\in L^1(t_0,T;W^{1,\infty})$. Moreover
$u_1\in C([t_0,T];H^1)$ by Lemma~\ref{lem:MR}(c).

\proofstep{Step 1: the kinetic states}
Both velocities belong to $L^1(t_0,T;W^{1,1})\cap L^1(t_0,T;L^\infty)$, so by Lemma~\ref{lem:DPL} each $f_i$ is the unique bounded solution
of its kinetic equation and belongs to $C([t_0,T];L^1)$. By Lemma~\ref{lem:kinstrong}(a), started at $t_0$ from $f_*$, $f_1(t)$ is the
push-forward of $f_*$ under the Carath\'eodory flow $Z_1(t,\cdot):=Z_{u_1}(t;t_0,\cdot)$. For $f_2$ we argue as in Step~8 of the proof of
Lemma~\ref{lem:L1Linf}: extend $u_2$ by zero outside $[t_0,T]$, let $u_2^\varepsilon$ be its mollification in $(t,x)$, let
$Z_2^\varepsilon(t,\cdot)$ be the flow of $(v,u_2^\varepsilon-v)$ started at $t_0$, and let $f_2^\varepsilon(t)$ be the push-forward of $f_*$ under
$Z_2^\varepsilon(t,\cdot)$. Then $\|u_2^\varepsilon\|_{L^1(t_0,T;L^\infty)}\le\|u_2\|_{L^1(t_0,T;L^\infty)}$, $u_2^\varepsilon\to u_2$ in $L^2((t_0,T)\times\TT^d)$, and
\begin{equation}\label{eq:ws-approx}
f_2^\varepsilon\rightharpoonup^*f_2\quad\text{in }L^\infty((t_0,T)\times\TT^d\times\RR^d)\qquad(\varepsilon\to0),
\end{equation}
for the whole family, since every weak-$*$ cluster point is a bounded weak solution with field $u_2$ and datum $f_*$, and hence equals $f_2$ by Lemma~\ref{lem:DPL}. The functions
$f_1$, $f_2$, $f_2^\varepsilon$ are bounded by $e^d\|f_*\|_{L^\infty}$ and vanish for $|v|\ge \Rv_\sharp$, so their densities are bounded by
$\varrho_\sharp$. Consequently, for every $t$ and every periodic measurable $h\ge0$,
\begin{equation}\label{eq:ws-push}
\int f_*(z)\,h\big(X_1(t,z)\big)\dd z=\int_{\TT^d}\rho_{f_1}(t)\,h\dd x\le\varrho_\sharp\|h\|_{L^1(\TT^d)},\qquad
\int f_*(z)\,h\big(X_2^\varepsilon(t,z)\big)\dd z\le\varrho_\sharp\|h\|_{L^1(\TT^d)},
\end{equation}
where $Z_1=(X_1,V_1)$ and $Z_2^\varepsilon=(X_2^\varepsilon,V_2^\varepsilon)$. In particular, compositions of almost everywhere defined functions
of $x$ with $X_1(t,\cdot)$ or $X_2^\varepsilon(t,\cdot)$ are well defined $f_*$-almost everywhere.

\proofstep{Step 2: the energy inequality of the fluid}
Let $\chi\in C^\infty_c(\RR^d)$ be equal to $1$ on $B_{\Rv_\sharp}$. Testing Definition~\ref{def:weak}(iii) for $f_2$ with
$\varphi=\eta(t)\frac12|v|^2\chi(v)$ shows that $\frac{\mathrm d}{\mathrm dt}\frac12M_2(f_2)=\int j_{f_2}\cdot u_2\dd x-M_2(f_2)$ in the sense of
distributions, with initial value $\frac12M_2(f_*)$. Since $t\mapsto M_2(f_2(t))$ is continuous, by Step~1,
\[
\tfrac12M_2(f_2(t))+\int_{t_0}^tM_2(f_2)\dd\sigma=\tfrac12M_2(f_*)+\int_{t_0}^t\!\!\int j_{f_2}\cdot u_2\dd x\dd\sigma\qquad(t_0\le t\le T).
\]
Since $\iint|v-u_2|^2f_2\dd z=M_2(f_2)-2\int j_{f_2}\cdot u_2\dd x+\int\rho_{f_2}|u_2|^2\dd x$ and $F_2\cdot u_2=j_{f_2}\cdot u_2-\rho_{f_2}|u_2|^2$,
subtracting this identity from (b) gives
\begin{equation}\label{eq:ws-fluid}
\tfrac12\|u_2(t)\|_{L^2}^2+\int_{t_0}^t\|\nabla u_2\|_{L^2}^2\dd\sigma\le\tfrac12\|u_*\|_{L^2}^2+\int_{t_0}^t\!\!\int F_2\cdot u_2\dd x\dd\sigma\qquad(t_0\le t\le T).
\end{equation}

\proofstep{Step 3: the relative energy}
We claim that, for every $t\in[t_0,T]$,
\begin{equation}\label{eq:ws-rel}
\tfrac12\|w(t)\|_{L^2}^2+\int_{t_0}^t\|\nabla w\|_{L^2}^2\dd\sigma\le\int_{t_0}^t\Big(\|\nabla u_1\|_{L^\infty}\|w\|_{L^2}^2+I(\sigma)\Big)\dd\sigma,\qquad
I(\sigma):=\iint(f_1-f_2)(\sigma)\,\psi_\sigma\dd z,
\end{equation}
where $\psi_\sigma(x,v):=(v-u_1(\sigma,x))\cdot w(\sigma,x)$.

By (a) and Step~1, $u_1\cdot\nabla u_1$ and $F_1$ belong to $L^2(t_0,T;L^2)$. Hence Definition~\ref{def:weak}(iv) gives
$\partial_tu_1-\Delta u_1+\QQ(u_1\cdot\nabla u_1)=\QQ F_1$ a.e., $u_1(t_0)=u_*$, and, as in the proof of Proposition~\ref{prop:critreg}(d),
\[
\tfrac12\|u_1(t)\|_{L^2}^2+\int_{t_0}^t\|\nabla u_1\|_{L^2}^2\dd\sigma=\tfrac12\|u_*\|_{L^2}^2+\int_{t_0}^t\!\!\int F_1\cdot u_1\dd x\dd\sigma .
\]
Next, as in the proof of Lemma~\ref{lem:concat}, Definition~\ref{def:weak}(iv) for $u_2$ and the weak continuity of $u_2$ give, for every
$t\in[t_0,T]$ and every smooth divergence-free $\varphi$ on $[t_0,T]\times\TT^d$,
\[
\int u_2(t)\cdot\varphi(t)\dd x-\int u_*\cdot\varphi(t_0)\dd x=\int_{t_0}^t\!\!\int\Big(u_2\cdot(\partial_t\varphi+\Delta\varphi)+u_2\otimes u_2:\nabla\varphi+F_2\cdot\varphi\Big)\dd x\dd\sigma .
\]
Smooth divergence-free fields are dense in the space of divergence-free fields in $H^1(t_0,T;L^2)\cap L^2(t_0,T;H^2)$ (reflect in time and
mollify), which is embedded in $C([t_0,T];H^1)$. Every term of the identity is continuous in this topology, because $u_2\in L^2_{t,x}$,
$u_2\otimes u_2\in L^1(t_0,T;L^3)$ while $\nabla\varphi$ is controlled in $C([t_0,T];L^{3/2})$, and $F_2\in L^2_{t,x}$. Hence the identity holds
with $\varphi=u_1$. Inserting the equation of $u_1$, and using $\int u_2\cdot\QQ g=\int u_2\cdot g$ and $u_2\otimes u_2:\nabla u_1=(u_2\cdot\nabla u_1)\cdot u_2$,
\[
\int u_1(t)\cdot u_2(t)\dd x=\|u_*\|_{L^2}^2+\int_{t_0}^t\!\!\int\Big(-2\nabla u_1:\nabla u_2+(u_2\cdot\nabla u_1)\cdot u_2-(u_1\cdot\nabla u_1)\cdot u_2+F_1\cdot u_2+F_2\cdot u_1\Big)\dd x\dd\sigma .
\]
The cubic terms equal $-(w\cdot\nabla u_1)\cdot u_2=(w\cdot\nabla u_1)\cdot w-w\cdot\nabla\frac{|u_1|^2}2$, and the last term integrates to zero
because $w$ is divergence-free. Adding the energy equality of $u_1$ and \eqref{eq:ws-fluid}, and subtracting the last identity, we obtain
\[
\tfrac12\|w(t)\|_{L^2}^2+\int_{t_0}^t\|\nabla w\|_{L^2}^2\dd\sigma\le\int_{t_0}^t\!\!\int\Big((F_1-F_2)\cdot w-(w\cdot\nabla u_1)\cdot w\Big)\dd x\dd\sigma .
\]
Since $F_1-F_2=\int(v-u_1)(f_1-f_2)\dd v-\rho_{f_2}w$ and $\rho_{f_2}\ge0$, this implies \eqref{eq:ws-rel}.

\proofstep{Step 4: the coupling term}
Lift $Z_1$ and $Z_2^\varepsilon$ to $\RR^d\times\RR^d$, with the fields extended periodically, started at $t_0$ from the same point, and set
\[
\delta Z^\varepsilon:=(\delta X^\varepsilon,\delta V^\varepsilon):=Z_1-Z_2^\varepsilon,\qquad D_\varepsilon(t)^2:=\int f_*(z)\,|\delta Z^\varepsilon(t,z)|^2\dd z,\qquad w_\varepsilon:=u_1-u_2^\varepsilon .
\]
As in Step~1 of Section~\ref{sec:lagstab},
$\frac{\mathrm d}{\mathrm dt}|\delta Z^\varepsilon|^2\le2(1+\|\nabla u_1(t)\|_{L^\infty})|\delta Z^\varepsilon|^2+|w_\varepsilon(t,X_2^\varepsilon)|^2$. Integrating against
$f_*$, and using \eqref{eq:ws-push} and Gronwall's inequality,
\[
D_\varepsilon(t)^2\le C_1\int_{t_0}^t\|w_\varepsilon(\sigma)\|_{L^2}^2\dd\sigma,\qquad C_1:=\varrho_\sharp\exp\big(2+2\|\nabla u_1\|_{L^1(t_0,T;L^\infty)}\big).
\]
Since $w_\varepsilon\to w$ in $L^2((t_0,T)\times\TT^d)$, the functions $D_\varepsilon$ are bounded uniformly in $\varepsilon$ and $t$, and
\begin{equation}\label{eq:ws-D}
\limsup_{\varepsilon\to0}D_\varepsilon(t)^2\le C_1\int_{t_0}^t\|w(\sigma)\|_{L^2}^2\dd\sigma\qquad(t_0\le t\le T).
\end{equation}

Let $I_\varepsilon(\sigma):=\iint(f_1-f_2^\varepsilon)(\sigma)\,\psi_\sigma\dd z$. For a.e.\ $\sigma$ we have $u_1(\sigma)\in W^{1,\infty}$ and $w(\sigma)\in H^1$, and
\[
I_\varepsilon=\int f_*\,\Big(\delta V^\varepsilon\cdot w(X_1)-\big(u_1(X_1)-u_1(X_2^\varepsilon)\big)\cdot w(X_1)
+\big(V_2^\varepsilon-u_1(X_2^\varepsilon)\big)\cdot\big(w(X_1)-w(X_2^\varepsilon)\big)\Big)\dd z,
\]
with all functions evaluated at time $\sigma$. The three terms are bounded as in Step~3 of Section~\ref{sec:lagstab}. By the
Cauchy--Schwarz inequality and \eqref{eq:ws-push}, the first two are at most $(1+\|\nabla u_1\|_{L^\infty})\varrho_\sharp^{1/2}D_\varepsilon\|w\|_{L^2}$. For
the third, $|V_2^\varepsilon-u_1(X_2^\varepsilon)|\le \Rv_\sharp+\|u_1\|_{L^\infty}$ on the support of $f_*$. Let $N_\sigma$ be the null set of
Lemma~\ref{lem:maximal} for $g=w(\sigma)$; by \eqref{eq:ws-push} with $h=\mathbf 1_{N_\sigma}$, the sets $\{X_1(\sigma,\cdot)\in N_\sigma\}$ and
$\{X_2^\varepsilon(\sigma,\cdot)\in N_\sigma\}$ are $f_*$-null; here \eqref{eq:ws-push} applies to $h=\mathbf 1_{N_\sigma}$ because $N_\sigma$, the complement of
the Lebesgue points of the periodic function $w(\sigma)$, is periodic. By \eqref{eq:maxpw}, \eqref{eq:ws-push} and \eqref{eq:maxL2}, the third term is at most
$2C_dC_M\varrho_\sharp^{1/2}(\Rv_\sharp+\|u_1\|_{L^\infty})D_\varepsilon\|\nabla w\|_{L^2}$. By Young's inequality,
\begin{equation}\label{eq:ws-I}
|I_\varepsilon(\sigma)|\le\tfrac12\|\nabla w(\sigma)\|_{L^2}^2+C_2\,a(\sigma)\big(D_\varepsilon(\sigma)^2+\|w(\sigma)\|_{L^2}^2\big)\qquad\text{for a.e.\ }\sigma,
\end{equation}
with $C_2$ depending only on $d$, $\Rv_\sharp$ and $\varrho_\sharp$.

Finally, $|\psi_\sigma|\le(\Rv_\sharp+|u_1|)|w|$ for $|v|\le \Rv_\sharp$, which is integrable on $(t_0,T)\times\TT^d\times B_{\Rv_\sharp}$. By
\eqref{eq:ws-approx}, $\int_{t_0}^tI\dd\sigma=\lim_{k\to\infty}\int_{t_0}^tI_{\varepsilon_k}\dd\sigma$ for every $t$ and every sequence $\varepsilon_k\to0$.
Fix such a sequence. By \eqref{eq:ws-I}, Fatou's lemma applied to the nonnegative functions $a\,(\sup_{\varepsilon,t}D_\varepsilon^2-D_{\varepsilon_k}^2)$,
and \eqref{eq:ws-D} for the measurable function $\limsup_kD_{\varepsilon_k}(\sigma)^2$,
\begin{equation}\label{eq:ws-Ilim}
\int_{t_0}^tI\dd\sigma\le\tfrac12\int_{t_0}^t\|\nabla w\|_{L^2}^2\dd\sigma+C_2\int_{t_0}^ta(\sigma)\Big(C_1\int_{t_0}^\sigma\|w\|_{L^2}^2\dd\tau+\|w(\sigma)\|_{L^2}^2\Big)\dd\sigma .
\end{equation}

\proofstep{Step 5: Gronwall}
Let $y(t):=\sup_{t_0\le\sigma\le t}\|w(\sigma)\|_{L^2}^2$, which is finite and nondecreasing. Inserting \eqref{eq:ws-Ilim} into \eqref{eq:ws-rel}, and
using $T-t_0\le1$ and $\|\nabla u_1\|_{L^\infty}\le a$, we get $\frac12\|w(t)\|_{L^2}^2\le C_3\int_{t_0}^ta(\sigma)y(\sigma)\dd\sigma$ for every $t$, with
$C_3:=1+C_2(1+C_1)$. The right-hand side is nondecreasing in $t$, so $y(t)\le2C_3\int_{t_0}^ta\,y\dd\sigma$, and Gronwall's inequality gives
$y\equiv0$. Thus $u_1=u_2$, and then $f_1=f_2$ by Lemma~\ref{lem:DPL}.
\end{proof}

\subsection{Leray solutions with the strong energy inequality}

\begin{theo}[Existence]\label{thm:SEI}
Let $d\in\{2,3\}$, let $u_{\rm in}\in L^2(\TT^d)$ be divergence-free, and let $f_0\in L^\infty(\TT^d\times\RR^d)$ be nonnegative and vanish for $|v|\ge \Rv_0$. There is
an energy solution $(u,f)$ on $[0,1]$ with data $(u_{\rm in},f_0)$ which satisfies the strong energy inequality: for $s=0$ and for almost
every $s\in(0,1)$, and for all $t\in[s,1]$, \eqref{eq:SEI} holds. In particular, by Lemma~\ref{lem:L1Linf} if $d=3$ and by
Lemma~\ref{lem:L2Linf2D} if $d=2$, $(u,f)$ is a weak solution in the sense of Definition~\ref{def:weak}, and a Leray solution with the strong
energy inequality in the sense of Section~\ref{sec:structure}.
\end{theo}

The proof uses a Galerkin approximation of the fluid equation in which the kinetic equation is solved exactly, as in the rest of
the paper. The approximate solutions then satisfy the energy \emph{equality} on every time interval, and the strong energy inequality
follows from the strong convergence of the velocities in $L^2_{t,x}$ and the convergence of the kinetic energy at every time.

\begin{proof}
Let $\mathcal E(0)=\frac12\|u_{\rm in}\|_{L^2}^2+\frac12M_2(f_0)$ and $F_\infty:=e^d\|f_0\|_{L^\infty}$. Constants $C$ depend only on $\mathcal E(0)$,
$\|f_0\|_{L^1}$, $\|f_0\|_{L^\infty}$ and $\Rv_0$.

\proofstep{Step 1: Galerkin approximation}
Let $V_n$ and $\Pi_n$ be as in the proof of Lemma~\ref{lem:adjoint}. For $U\in C([0,1];V_n)$, the flow of $(v,U-v)$ is classical, and
$f[U]$ denotes the push-forward of $f_0$ along it, as in Section~\ref{sec:rough}. We look for $u_n\in C^1([0,1];V_n)$ such that
\begin{equation}\label{eq:galerkin}
u_n'=\Delta u_n-\Pi_n(u_n\cdot\nabla u_n)+\Pi_nF[f_n,u_n],\qquad u_n(0)=\Pi_nu_{\rm in},\qquad f_n:=f[u_n].
\end{equation}
Since $V_n$ is finite-dimensional, there is $C_n$ with $\|\phi\|_{W^{1,\infty}}+\|\Delta\phi\|_{L^2}\le C_n\|\phi\|_{L^2}$ for $\phi\in V_n$. Let
$A:=(2\mathcal E(0))^{1/2}+1$, and let $P_A(\phi):=\phi\min\{1,A/\|\phi\|_{L^2}\}$ be the retraction of $V_n$ onto its ball of radius $A$, which
is $1$-Lipschitz in $L^2$. For $U\in C([0,1];V_n)$ let $\tilde U(t):=P_A(U(t))$ and
\[
\mathcal G(U)(t):=\Big(\Delta\tilde U-\Pi_n(\tilde U\cdot\nabla\tilde U)+\Pi_nF\big[f[\tilde U],\tilde U\big]\Big)(t)\in V_n .
\]
Since $\|\tilde U\|_{L^\infty(0,1;W^{1,\infty})}\le C_nA$, Lemma~\ref{lem:kinstrong}(a) shows that $f[\tilde U]$ vanishes for $|v|\ge \Rv_0+C_nA$ and that
$\rho[\tilde U]$ is bounded by a constant depending on $n$. For $\phi\in V_n$, the function
$t\mapsto\int F[f[\tilde U],\tilde U](t)\cdot\phi\dd x=\int f_0(z)\big(V(t)-\tilde U(t,X(t))\big)\cdot\phi(X(t))\dd z$, where $(X,V)(t)=Z_{\tilde U}(t;0,z)$, is
continuous; hence $\mathcal G(U)\in C([0,1];V_n)$.

Let $U,U'\in C([0,1];V_n)$. For $\phi\in V_n$ with $\|\phi\|_{L^2}\le1$ and $\psi(x,v):=(v-\tilde U'(t,x))\cdot\phi(x)$,
\[
\int\Big(F\big[f[\tilde U],\tilde U\big]-F\big[f[\tilde U'],\tilde U'\big]\Big)(t)\cdot\phi\dd x=\iint\big(f[\tilde U]-f[\tilde U']\big)(t)\,\psi\dd z-\int\rho[\tilde U]\,(\tilde U-\tilde U')(t)\cdot\phi\dd x .
\]
By Lemma~\ref{lem:kinstrong}(b), with base field $\tilde U'$ and perturbation $\tilde U-\tilde U'$, the first term on the right is bounded by a
constant depending on $n$ times $\int_0^t\|(\tilde U-\tilde U')(s)\|_{L^\infty}\dd s$. Together with the elementary bounds for the linear and the
quadratic terms on $V_n$, this gives
\[
\|\mathcal G(U)(t)-\mathcal G(U')(t)\|_{L^2}\le L_n\sup_{0\le s\le t}\|U(s)-U'(s)\|_{L^2},
\]
with $L_n$ depending on $n$, $A$ and $f_0$. As in Step~2 of the proof of Lemma~\ref{lem:mean}, the map
$U\mapsto\Pi_nu_{\rm in}+\int_0^t\mathcal G(U)\dd s$ is a contraction of $C([0,1];V_n)$ for the norm $\sup_te^{-2L_nt}\|U(t)\|_{L^2}$. Its fixed
point $u_n$ belongs to $C^1([0,1];V_n)$ and solves $u_n'=\mathcal G(u_n)$, $u_n(0)=\Pi_nu_{\rm in}$.

It remains to remove the retraction. Let $T_n:=\sup\{t\in[0,1]:\ \|u_n\|_{L^2}<A\text{ on }[0,t]\}$, which is positive because
$\|\Pi_nu_{\rm in}\|_{L^2}\le\|u_{\rm in}\|_{L^2}<A$. On $[0,T_n)$ we have $\tilde u_n=u_n$, so \eqref{eq:galerkin} holds there. Testing with
$u_n$, and using $\int\Pi_n(u_n\cdot\nabla u_n)\cdot u_n=\int(u_n\cdot\nabla u_n)\cdot u_n=0$,
\[
\frac{\mathrm d}{\mathrm dt}\frac12\|u_n\|_{L^2}^2+\|\nabla u_n\|_{L^2}^2=\int F[f_n,u_n]\cdot u_n\dd x=\int j_n\cdot u_n\dd x-\int\rho_n|u_n|^2\dd x,
\]
where $\rho_n:=\rho[u_n]$ and $j_n:=j[u_n]$. For the particles, with $(X^{(n)},V^{(n)})(t,z):=Z_{u_n}(t;0,z)$,
\[
\frac{\mathrm d}{\mathrm dt}\frac12M_2(f_n)=\frac{\mathrm d}{\mathrm dt}\frac12\int f_0\,|V^{(n)}|^2\dd z=\int f_0\,V^{(n)}\cdot\big(u_n(t,X^{(n)})-V^{(n)}\big)\dd z=\int j_n\cdot u_n\dd x-M_2(f_n).
\]
Adding the two identities, with $\mathcal E_n(t):=\frac12\|u_n(t)\|_{L^2}^2+\frac12M_2(f_n(t))$,
\begin{equation}\label{eq:galen}
\mathcal E_n(t)+\int_s^t\Big(\|\nabla u_n\|_{L^2}^2+\iint|v-u_n|^2f_n\dd z\Big)\dd\sigma=\mathcal E_n(s)\qquad(0\le s\le t<T_n).
\end{equation}
In particular $\|u_n(t)\|_{L^2}^2\le2\mathcal E_n(0)\le2\mathcal E(0)=(A-1)^2$ on $[0,T_n)$. By continuity, $T_n=1$, and \eqref{eq:galerkin} and
\eqref{eq:galen} hold on $[0,1]$.

\proofstep{Step 2: uniform bounds}
By \eqref{eq:galen} and Lemma~\ref{lem:kinstrong}(a), for all $n$ and $t$,
\begin{equation}\label{eq:unif}
\|u_n(t)\|_{L^2}^2+M_2(f_n(t))\le2\mathcal E(0),\qquad\int_0^1\|\nabla u_n\|_{L^2}^2\dd t\le\mathcal E(0),\qquad 0\le f_n\le F_\infty,\qquad\|f_n(t)\|_{L^1}=\|f_0\|_{L^1}.
\end{equation}
In particular $\|u_n\|_{L^2(0,1;L^6)}\le C$. The function $M_3(f_n(t))=\int f_0|V^{(n)}(t,z)|^3\dd z$ is of class $C^1$, and by \eqref{eq:interp} with
$(k,m)=(2,3)$, together with $\frac{d+3}{d+2}\ge\frac65$ for $d\le3$,
\[
\frac{\mathrm d}{\mathrm dt}M_3(f_n)=3\int f_0\,|V^{(n)}|\,V^{(n)}\cdot\big(u_n(t,X^{(n)})-V^{(n)}\big)\dd z\le3\|u_n\|_{L^6}\|m_2f_n\|_{L^{\frac{d+3}{d+2}}}\le C\|u_n\|_{L^6}M_3(f_n)^{\frac{d+2}{d+3}}.
\]
Hence $M_3(f_n(t))^{\frac1{d+3}}\le(\Rv_0^3\|f_0\|_{L^1})^{\frac1{d+3}}+C\|u_n\|_{L^1L^6}\le C$, and by \eqref{eq:interp} with $m=3$, which gives
$\rho_n\in L^\infty L^{(d+3)/d}\subset L^\infty L^2$ and $j_n\in L^\infty L^{(d+3)/(d+1)}\subset L^\infty L^{3/2}$ for $d\le3$,
\begin{equation}\label{eq:unif3}
\sup_tM_3(f_n(t))+\|\rho_n\|_{L^\infty(0,1;L^2)}+\|j_n\|_{L^\infty(0,1;L^{3/2})}+\|F_n\|_{L^2(0,1;L^{3/2})}\le C,\qquad F_n:=F[f_n,u_n]=j_n-\rho_nu_n .
\end{equation}
For $\phi\in H^3(\TT^d;\RR^d)$, since $u_n'\in V_n$ and $\Pi_n$ is self-adjoint and commutes with derivatives,
\[
\int u_n'\cdot\phi\dd x=\int\Big(-\nabla u_n:\nabla\Pi_n\phi+u_n\otimes u_n:\nabla\Pi_n\phi+F_n\cdot\Pi_n\phi\Big)\dd x,
\]
and $\|\Pi_n\phi\|_{H^3}\le\|\phi\|_{H^3}$, $H^3\hookrightarrow W^{1,\infty}$ for $d\le3$. Hence
$\|u_n'\|_{H^{-3}}\le C\big(\|u_n\|_{H^1}+\|u_n\|_{L^2}^2+\|F_n\|_{L^{3/2}}\big)$, which is bounded in $L^2(0,1)$.

\proofstep{Step 3: compactness}
By \eqref{eq:unif}, the bound on $u_n'$, the Aubin--Lions lemma and the Arzel\`a--Ascoli theorem, there are a subsequence, not relabeled,
and a weakly divergence-free $u\in C_w([0,1];L^2)\cap L^2(0,1;H^1)$ such that
\[
u_n\to u\ \text{ in }L^2((0,1)\times\TT^d),\qquad u_n(t)\rightharpoonup u(t)\ \text{ in }L^2\text{ for every }t,\qquad u_n(s)\to u(s)\ \text{ in }L^2\text{ for a.e.\ }s,
\]
and $\nabla u_n\rightharpoonup\nabla u$ in $L^2_{t,x}$. Interpolating with the bound in $L^2(0,1;L^6)$, $u_n\to u$ in $L^2(0,1;L^p)$ for every $p<6$.

For $\psi\in C^1_c(\TT^d\times\RR^d)$, the function $t\mapsto\iint f_n(t)\psi\dd z=\int f_0\,\psi(Z_{u_n}(t;0,z))\dd z$ has derivative
$\iint f_n\big(v\cdot\nabla_x\psi+(u_n-v)\cdot\nabla_v\psi\big)\dd z$, which is bounded by $C_\psi(1+\|u_n(t)\|_{L^2})\le C_\psi$, uniformly in $n$ and $t$.
By the Arzel\`a--Ascoli theorem, a diagonal argument over a countable subset of $C^1_c$ that is dense in $L^1$, and the bound
$0\le f_n\le F_\infty$, a further subsequence satisfies
\[
f_n(t)\rightharpoonup^*f(t)\quad\text{in }L^\infty(\TT^d\times\RR^d)\quad\text{for every }t\in[0,1],
\]
where $0\le f(t)\le F_\infty$ and $t\mapsto\iint f(t)\psi\dd z$ is continuous. As in the proof of Lemma~\ref{lem:convergence}(b),
$f_n\rightharpoonup^*f$ in $L^\infty((0,1)\times\TT^d\times\RR^d)$. By \eqref{eq:unif3}, $\iint_{\{|v|\ge L\}}f_n(t)(1+|v|^2)\dd z\le CL^{-1}$ for $L\ge1$, uniformly
in $n$ and $t$, and the same holds for $f(t)$ by lower semicontinuity. Hence
\begin{equation}\label{eq:tight}
\iint f_n(t)\,\Psi\dd z\to\iint f(t)\,\Psi\dd z\qquad\text{for every }t\text{ and every }\Psi\in C(\TT^d\times\RR^d)\text{ with }|\Psi|\le C(1+|v|^2).
\end{equation}
In particular $M_2(f_n(t))\to M_2(f(t))$ and $\iint f(t)\dd z=\|f_0\|_{L^1}$ for every $t$. By \eqref{eq:tight}, dominated convergence in time
and \eqref{eq:unif3},
\begin{equation}\label{eq:rhoj}
\rho_n\rightharpoonup^*\rho_f\ \text{ in }L^\infty(0,1;L^2),\qquad j_n\rightharpoonup^*j_f\ \text{ in }L^\infty(0,1;L^{3/2}).
\end{equation}

\proofstep{Step 4: the limit equations}
Each $f_n$ satisfies Definition~\ref{def:weak}(iii) with field $u_n$, by Lemma~\ref{lem:kinstrong}(a). For a test function $\varphi$,
$\big|\iiint f_n(u_n-u)\cdot\nabla_v\varphi\big|\le F_\infty C_\varphi\|u_n-u\|_{L^1_{t,x}}\to0$, and the remaining terms pass to the limit because their
multipliers are integrable. Hence $f$ satisfies Definition~\ref{def:weak}(iii) with field $u$.

Let $\varphi$ be a smooth divergence-free field vanishing near $t=1$. Testing \eqref{eq:galerkin} with $\varphi$ and integrating by parts,
\[
\int\Pi_nu_{\rm in}\cdot\varphi(0)\dd x+\int_0^1\!\!\int\Big(u_n\cdot(\partial_t\varphi+\Delta\varphi)+u_n\otimes u_n:\nabla\Pi_n\varphi+F_n\cdot\Pi_n\varphi\Big)\dd x\dd t=0 .
\]
Since $\varphi$ is smooth, $\Pi_n\varphi\to\varphi$ in $C([0,1];W^{1,\infty})$. The first three terms converge because $u_n\to u$ in $L^2_{t,x}$. Moreover
$\iint j_n\cdot\Pi_n\varphi\to\iint j_f\cdot\varphi$ by \eqref{eq:rhoj}, and $\iint\rho_nu_n\cdot\Pi_n\varphi\to\iint\rho_fu\cdot\varphi$ by \eqref{eq:rhoj} and the strong
convergence of $u_n\cdot\Pi_n\varphi$ in $L^1(0,1;L^2)$. Hence Definition~\ref{def:weak}(iv) holds for $(u,f)$.

\proofstep{Step 5: the energy inequalities}
Fix $0\le s\le t\le1$ and pass to the limit in \eqref{eq:galen}.
\begin{itemize}
\item \emph{At time $t$.} $\|u(t)\|_{L^2}\le\liminf\|u_n(t)\|_{L^2}$ and $M_2(f_n(t))\to M_2(f(t))$, so $\mathcal E(t)\le\liminf\mathcal E_n(t)$.
\item \emph{Fluid dissipation.} $\int_s^t\|\nabla u\|_{L^2}^2\le\liminf\int_s^t\|\nabla u_n\|_{L^2}^2$.
\item \emph{Drag dissipation.} We have $\iint|v-u_n|^2f_n\dd z=M_2(f_n)-2\int j_n\cdot u_n\dd x+\int\rho_n|u_n|^2\dd x$. The first term is
bounded and converges at every time. Since $u_n\to u$ in $L^1(0,1;L^3)$ and $|u_n|^2\to|u|^2$ in $L^1(0,1;L^2)$, the other two converge
after integration in time, by \eqref{eq:rhoj}. Hence $\int_s^t\iint|v-u_n|^2f_n\to\int_s^t\iint|v-u|^2f$.
\item \emph{At time $s$.} $\mathcal E_n(s)\to\mathcal E(s)$ whenever $u_n(s)\to u(s)$ in $L^2$. This holds for a.e.\ $s$, and for $s=0$
because $u_n(0)=\Pi_nu_{\rm in}\to u_{\rm in}$.
\end{itemize}
Therefore \eqref{eq:SEI} holds for $s=0$ and a.e.\ $s\in(0,1)$, and all $t\in[s,1]$.

\proofstep{Step 6: conclusion}
By Steps~3--5, $(u,f)$ satisfies (E1)--(E4) of Definition~\ref{def:energy}: (E2) follows from $0\le f\le F_\infty$, the conservation of
mass and $M_2(f(t))\le2\mathcal E(0)$, and (E4) is \eqref{eq:SEI} with $s=0$. The last statement follows from Lemma~\ref{lem:L1Linf} or
Lemma~\ref{lem:L2Linf2D}.

Finally, \eqref{eq:SEI} is an every-time statement in $t$, so we fix the representative. In Steps~3--5, $f(t)$ is the weak-$*$ limit of
$f_n(t)$ for every $t$, and $t\mapsto\iint f(t)\psi\dd z$ is continuous for $\psi\in C^1_c$. The function $f$ has a representative in
$C([0,1];L^1)$, by Lemma~\ref{lem:L1Linf}(c) or Lemma~\ref{lem:L2Linf2D}(c). The two agree for a.e.\ $t$; both are weakly continuous in time against $C^1_c$ test functions, so they agree for
every $t$. In particular, $M_2(f(t))$ in \eqref{eq:SEI} is the same for both, and the statement holds for the representative in
$C([0,1];L^1)$ used in Sections~\ref{sec:local} and~\ref{sec:structure}.
\end{proof}

\subsection*{The two-dimensional torus}

In this paragraph $d=2$. The proof of the $L^2_tL^\infty$ bound of \cite[Lemma~3]{HMMM20} is written on $\RR^2$, for solutions in the
sense of \cite[Definition~1]{HMMM20}, and it uses the energy identity of \cite[Lemma~2]{HMMM20}. Its whole-space ingredients are maximal
regularity for the Stokes problem, a Gagliardo--Nirenberg inequality, and the estimate of the free heat flow through the embedding
$L^2(\RR^2)\hookrightarrow B^{-1}_{\infty,2}$. We prove the bound on $\TT^2$ for the energy solutions of Definition~\ref{def:energy}, with periodic
versions of these three ingredients. The energy inequality (E4) replaces \cite[Lemmas~1 and~2]{HMMM20}, so that no result of
\cite{HMMM20} is used.

\begin{lem}[Velocity bound in two dimensions; after {\cite[Lemma~3]{HMMM20}}]\label{lem:L2Linf2D}
Let $d=2$, and let $(u,f)$ be an energy solution on $[0,T]$, $T\le1$, with data $(u_{\rm in},f_0)$, in the sense of
Definition~\ref{def:energy}. Then:
\begin{enumerate}[label=(\alph*)]
\item $u\in L^2(0,T;L^\infty)$, with a bound depending only on $T$, $\mathcal E(0)$, $\|f_0\|_{L^1}$, $\|f_0\|_{L^\infty}$ and $\Rv_0$;
\item $0\le f\le e^2\|f_0\|_{L^\infty}$, and $f(t)$ vanishes for $|v|\ge \Rv_0+\|u\|_{L^1(0,T;L^\infty)}$;
\item $(u,f)$ is a weak solution on $[0,T]$ in the sense of Definition~\ref{def:weak}; $f$ is the unique bounded (hence
renormalized) solution of the kinetic equation with field $u$, it belongs to $C([0,T];L^1)$, and $\iint f(t)\dd z=\|f_0\|_{L^1}$ for
every $t$.
\end{enumerate}
\end{lem}

\begin{proof}
Let $F_\infty:=\|f\|_{L^\infty}$. Until the last step, constants $C$ depend on the listed quantities and on $F_\infty$; in the last step we
show that $F_\infty\le e^2\|f_0\|_{L^\infty}$; that $F_\infty<\infty$ is part of (E2). By (E1) and (E4),
$\|u\|_{L^\infty L^2}^2+\|u\|_{L^2H^1}^2\le C$ and $M_2(f(t))\le2\mathcal E(0)$ for a.e.\ $t$.

\proofstep{The forcing}
By \eqref{eq:interp} with $d=2$, $m=2$ and $k=0,1$, we have
$\rho_f\in L^\infty(0,T;L^2)$ and $j_f\in L^\infty(0,T;L^{4/3})$. By Ladyzhenskaya's inequality $\|u\|_{L^4}^2\le C\|u\|_{L^2}\|u\|_{H^1}$, so
$u\in L^4(0,T;L^4)$. Hence
\[
G:=j_f-\rho_fu-u\cdot\nabla u\in L^{4/3}(0,T;L^{4/3}),
\]
since $\|\rho_fu\|_{L^{4/3}}\le\|\rho_f\|_{L^2}\|u\|_{L^4}$ and $\|u\cdot\nabla u\|_{L^{4/3}}\le\|u\|_{L^4}\|\nabla u\|_{L^2}$.

\proofstep{Decomposition}
The mean $P_0u$ is bounded by $\|u\|_{L^\infty L^2}$. Let $W$ be the mean-free solution of $\partial_tW-\Delta W=(1-P_0)\QQ G$ with $W(0)=0$ given
by maximal $L^{4/3}$-regularity of the heat equation on $\TT^2$ \cite[Chapter~III]{Ama95}, so that
$\|\nabla^2W\|_{L^{4/3}(0,T;L^{4/3})}\le C\|G\|_{L^{4/3}L^{4/3}}$. We claim that
\[
(1-P_0)u=e^{t\Delta}(1-P_0)u_{\rm in}+W .
\]
We argue on Fourier coefficients, which avoids having to identify a function space for $\partial_tu$; we write $\hat g_n:=\int_{\TT^2}g\,\mathrm e^{-2\pi\mathrm in\cdot x}\dd x$. Fix
$n\in\ZZ^2\setminus\{0\}$ and a unit vector $e\perp n$. Testing Definition~\ref{def:weak}(iv) with the real and imaginary parts of
$\varphi=\eta(t)\,e\,\mathrm e^{-2\pi\mathrm in\cdot x}$, $\eta\in C^\infty_c([0,T))$, which are smooth and divergence-free since $e\cdot n=0$, shows that
$y(t):=e\cdot\hat u_n(t)$ satisfies $y'+4\pi^2|n|^2y=e\cdot\hat G_n$ in $\mathcal D'(0,T)$, with $y\in L^\infty(0,T)$ and $y(0)=e\cdot(\hat u_{\rm in})_n$ in the sense
of the weak formulation; here we used $\QQ\D(u\otimes u)=\QQ(u\cdot\nabla u)$ and $e\cdot\widehat{(\QQ g)}_n=e\cdot\hat g_n$. The right-hand side of the
claim has coefficients solving the same problem, and this scalar problem has a unique solution in $L^1(0,T)$: the difference of two
solutions is $c\,\mathrm e^{-4\pi^2|n|^2t}$, and testing with $\eta(0)=1$ gives $c=0$. The $n$-th coefficients of both sides of the claim
are orthogonal to $n$: for the left side because $u$ is divergence-free; for the right side because $u_{\rm in}$ is divergence-free and
$n\cdot\hat W_n$ solves the same scalar problem with zero datum and zero forcing, since the forcing $(1-P_0)\QQ G$ has coefficients
orthogonal to $n$. Hence the two sides agree for a.e.\ $t$. The left side is weakly continuous with values in $L^2$ and the right side
is continuous with values in $L^{4/3}$, so they agree for every $t$, and the claim follows. In particular
$W=(1-P_0)u-e^{t\Delta}(1-P_0)u_{\rm in}\in L^\infty(0,T;L^2)$. The Gagliardo--Nirenberg inequality for mean-free functions on $\TT^2$,
$\|\phi\|_{L^\infty}\le C\|\phi\|_{L^2}^{1/3}\|\nabla^2\phi\|_{L^{4/3}}^{2/3}$, gives
$\|W\|_{L^\infty}^2\le C\|W\|_{L^2}^{2/3}\|\nabla^2W\|_{L^{4/3}}^{4/3}\in L^1(0,T)$.

\proofstep{The free heat flow}
Let $g=(1-P_0)u_{\rm in}$ and $a_j=\|\Delta_jg\|_{L^2}$, $j\ge0$. By Bernstein's inequality and the decay of the heat semigroup on dyadic
annuli, $\|e^{t\Delta}\Delta_jg\|_{L^\infty}\le C2^je^{-ct4^j}a_j$. Hence
\[
\int_0^T\|e^{t\Delta}g\|_{L^\infty}^2\dd t\le C\sum_{j,k\ge0}2^{j+k}a_ja_k\int_0^\infty e^{-ct(4^j+4^k)}\dd t
\le C\sum_{j,k\ge0}2^{-|j-k|}a_ja_k\le C\sum_ja_j^2\le C\|u_{\rm in}\|_{L^2}^2,
\]
by Schur's test. Together with the previous step, $u\in L^2(0,T;L^\infty)$.

\proofstep{Velocity support and the bound on $f$}
We now know that $u\in L^1(0,T;W^{1,1})\cap L^1(0,T;L^\infty)$, and $f$ is a bounded weak solution of the kinetic equation with field $u$.
Step~8 of the proof of Lemma~\ref{lem:L1Linf}, which uses only these two facts, Lemma~\ref{lem:kinstrong}(a) and Lemma~\ref{lem:DPL}, applies in
dimension two with $e^2$ in place of $e^3$: $f$ is the weak-$*$ limit of the push-forwards $f[u_\varepsilon]$ of $f_0$. This proves (b) and
shows that $(u,f)$ satisfies Definition~\ref{def:weak}(i)--(ii). The remaining statements of (c) follow from Lemma~\ref{lem:DPL} and from the
conservation of mass, which is obtained as in Step~1 of the proof of Lemma~\ref{lem:L1Linf}. In fact the removal of the cutoff is now
trivial: by (b), for a.e.\ $t$ the function $f(t)$ vanishes for $|v|\ge \Rv_0+\|u\|_{L^1L^\infty}$, so the term containing $\nabla_v\chi_L$ is zero
for $L$ large.
\end{proof}

\section*{Acknowledgements}
The author thanks Xiaoyutao Luo for carefully reading an earlier version of the manuscript and for pointing out several errors
in Section~\ref{sec:regint}. The author is supported by the PRIN 2022 project ``Turbulent effects vs Stability in Equations from
Oceanography'' (TESEO), project number 2022HSSYPN.

\subsection*{Declaration on the use of AI tools}
The author used an AI language model (Claude, Anthropic) to check parts of the proofs and of the text,
and to harmonize language and notation. The author verified all mathematical content and takes full responsibility for it.

\appendix

\section{The convex integration step}\label{sec:CI}

In this appendix and the next we prove Proposition~\ref{prop:NSstep}, following \cite[Sections~4--5]{CL22};
Section~\ref{sec:NSdiff} summarizes the construction and the differences from \cite{CL22}.
In this section we construct $w$, $P$ and $R^{\rm CL}$ and verify the algebraic properties (i) and (ii) of
Proposition~\ref{prop:NSstep} (Proposition~\ref{prop:CIidentity}). The parameters are fixed and the quantitative
bounds (iii)--(v) are proved in Appendix~\ref{sec:estimates}.

\paragraph{Standing assumptions.} Throughout this appendix, $d\ge2$, $1\le p<2$, $1\le q<\infty$ and $\delta>0$ are
fixed, and $(u,R)$, $0<s<t_*\le1$ are as in Proposition~\ref{prop:NSstep}:
\[
u\in C^\infty([0,1]\times\TT^d;\RR^d),\quad \D u=0,\qquad R\in C^\infty([0,1]\times\TT^d;\Sz),\qquad R=0\ \text{on }[0,t_*]\times\TT^d .
\]
If $R\equiv0$, Proposition~\ref{prop:NSstep} holds with $w=0$, $P=0$, $R^{\rm CL}=0$. We therefore assume
\[
\|R\|_1:=\|R\|_{L^1([0,1]\times\TT^d)}>0 .
\]
The perturbation depends on four parameters,
\[
\begin{aligned}
&\nu\in\NN\ \text{(temporal oscillation)},&&\kappa\ge1\ \text{(temporal concentration)},\\
&\sigma\in\NN\ \text{(spatial oscillation)},&&\mu\ge1\ \text{(spatial concentration)},
\end{aligned}
\]
which are taken large and are fixed in Appendix~\ref{sec:estimates} as powers of a single frequency $\lambda$. All objects in this section
are well defined for every admissible choice of these parameters with $\mu\ge\mu_0$, where $\mu_0$ is the geometric
constant of Section~\ref{sec:mikado}.

\paragraph{Conventions.} For matrices we use the Frobenius norm $|A|=(A:A)^{1/2}$. For a vector $v$ and a matrix
$A$ we write $(v\cdot A)_i=v_jA_{ij}$, so that $\D(fA)=\nabla f\cdot A+f\D A$ with $(\D A)_i=\partial_jA_{ij}$. For
$f$ on $\TT^d$ and $\sigma\in\NN$, $f(\sigma\cdot)$ denotes the $\sigma^{-1}$-periodic function $x\mapsto f(\sigma x)$;
recall
\begin{equation}\label{eq:rescale}
\|f(\sigma\cdot)\|_{L^p(\TT^d)}=\|f\|_{L^p(\TT^d)},\qquad \nabla^m\big(f(\sigma\cdot)\big)=\sigma^m(\nabla^mf)(\sigma\cdot),\qquad
\fint f(\sigma\cdot)=\fint f .
\end{equation}

\subsection{Antidivergence operators}\label{sec:antidiv}

Besides the operator $\Rs$ of Section~\ref{sec:notation}, we use the bilinear antidivergence of
\cite[Appendix~B.3]{CL22},
\[
\Bs(v,A):=v\cdot\Rs A-\Rs\big(\nabla v\cdot\Rs A\big),
\]
defined for $v\in C^\infty(\TT^d;\RR^d)$ and $A\in C^\infty(\TT^d;\RR^{d\times d})$ with $\fint A=0$. Here $\Rs$ acts on
$A$ row by row, and the contractions are those of \cite[Appendix~B.3]{CL22}. We only use the following facts.

\begin{lem}[{(a), (c): \cite[Theorems~B.3, B.4]{CL22}}]\label{lem:antidiv}
\begin{enumerate}[label=(\alph*)]
\item $\Rs$ maps smooth vector fields to smooth $\Sz$-valued fields, $\D\Rs g=g-\fint g$, and
$\|\Rs g\|_{L^p}\lesssim\|g\|_{L^p}$ for $1\le p\le\infty$. If $\fint g=0$ and $\sigma\in\NN$, then
$\|\Rs(g(\sigma\cdot))\|_{L^p}\lesssim\sigma^{-1}\|g\|_{L^p}$.
\item For $1<p<\infty$, the operator $\Rs\D$ is bounded on $L^p(\TT^d;\RR^{d\times d})$.
\item Let $A$ be smooth, symmetric and of zero mean. Then $\Bs(v,A)$ is smooth and $\Sz$-valued,
\[
\D\Bs(v,A)=v\cdot A-\fint v\cdot A,\qquad \|\Bs(v,A)\|_{L^p}\lesssim\|v\|_{C^1}\|\Rs A\|_{L^p}\quad(1\le p\le\infty).
\]
\end{enumerate}
\end{lem}

Part (b) holds because $\Rs\D$ is a matrix of zeroth-order Fourier multipliers of Calder\'on--Zygmund type (compositions of
Riesz transforms), smooth away from the origin.

\subsection{Stationary Mikado flows}\label{sec:mikado}

The building blocks are the stationary Mikado flows of Daneri and Sz\'ekelyhidi \cite{DS17}, with spatial
concentration as in \cite{MSz18,CL22}. The following geometric lemma goes back to Nash; see \cite[Lemma~3.2]{DLS13}.

\begin{lem}[Geometric lemma]\label{lem:geometric}
Let $\mathcal N:=\{A\in\mathcal S^{d\times d}:\ |A-\mathrm{Id}|\le\frac12\}$, a compact subset of the positive definite
symmetric matrices. There are a finite set $\Xi\subset\ZZ^d\setminus\{0\}$ and smooth functions
$\Gamma_k\in C^\infty(\mathcal N;\RR)$, $k\in\Xi$, such that
\[
A=\sum_{k\in\Xi}\Gamma_k^2(A)\,e_k\otimes e_k\qquad\text{for all }A\in\mathcal N,\qquad e_k:=\frac k{|k|}.
\]
\end{lem}

Since only the directions $e_k$ enter, we may replace each $k\in\Xi$ by $k/\gcd(k_1,\dots,k_d)$; we therefore assume that every
$k\in\Xi$ is \emph{primitive}. After this replacement two elements of $\Xi$ may coincide (for instance $k$ and $2k$). We do not merge
them, since $(\Gamma_k^2+\Gamma_{k'}^2)^{1/2}$ need not be smooth where both functions vanish. Instead we regard $\Xi$ as a finite
\emph{indexed family}: in what follows, $k\ne k'$ refers to distinct indices, which may carry the same vector, and each index has its own
function $\Gamma_k$, its own base point $p_k$ and its own Mikado flow. The family $\Xi$ and the functions $\Gamma_k$ are fixed once and for
all, and depend only on $d$.

\paragraph{Periodic lines.} For $k\in\Xi$ and $p_k\in\TT^d$ let $l_k=\{p_k+tk:\ t\in\RR\}\subset\TT^d$. Since $k$ is primitive, it is a
closed geodesic of length $|k|$, traversed once as $t$ runs over $[0,1]$. Two periodic lines with parallel directions either coincide or are disjoint, and two lines with
non-parallel directions meet in finitely many points. We choose the points $p_k$, $k\in\Xi$, one after the other, so that
\begin{equation}\label{eq:lines}
l_k\cap l_{k'}=\emptyset\qquad\text{whenever }k\ne k'\text{ and }k\parallel k' .
\end{equation}
This is possible because, for each $k$, the forbidden set of $p_k$ is a finite union of periodic lines. We make this choice for
definiteness; it also covers parallel directions such as $k'=-k$, and two indices carrying the same vector. By compactness and transversality, there are $\mu_0\ge1$ and $M_\Xi\ge1$,
depending only on $\Xi$ and the points $p_k$, such that for $\mu\ge\mu_0$ the function $x\mapsto\dist(x,l_k)$ is
smooth on $\{0<\dist(x,l_k)\le\mu^{-1}\}$, and for $k\ne k'$
\begin{equation}\label{eq:intersect}
\{\dist(x,l_k)\le\mu^{-1}\}\cap\{\dist(x,l_{k'})\le\mu^{-1}\}\subset\bigcup_{y\in l_k\cap l_{k'}}B\big(y,M_\Xi\mu^{-1}\big),
\end{equation}
where the right-hand side is empty if $k\parallel k'$.

\paragraph{Profiles.} Fix a radial function $\phi\in C^\infty_c(\RR^{d-1})$, not identically zero, supported in the
annulus $\{\frac12\le|y|\le1\}$, and let $\psi=\Delta_{\RR^{d-1}}\phi$, which is radial, supported in the same annulus
and not identically zero. Write $\phi(|y|)$, $\psi(|y|)$ for their radial profiles. For $\mu\ge\mu_0$ and $k\in\Xi$ set
\begin{equation}\label{eq:psiphi}
\psi_k(x):=c_k\mu^{\frac{d-1}2}\psi\big(\mu\dist(x,l_k)\big),\qquad
\phi_k(x):=c_k\mu^{\frac{d-1}2-2}\phi\big(\mu\dist(x,l_k)\big),
\end{equation}
where $c_k=\big(|k|\int_{\RR^{d-1}}\psi^2\big)^{-1/2}$. In coordinates adapted to $l_k$ (arclength along $l_k$ and
Euclidean coordinates $y$ in the orthogonal $(d-1)$-plane), $\psi_k$ and $\phi_k$ depend only on $y$, so
\begin{equation}\label{eq:lap}
\Delta\phi_k=\psi_k,\qquad e_k\cdot\nabla\psi_k=e_k\cdot\nabla\phi_k=0,\qquad \int_{\TT^d}\psi_k^2\dd x=1 .
\end{equation}
The Mikado flows and their potentials are
\begin{equation}\label{eq:W}
W_k:=\psi_k\,e_k,\qquad \Omega_k:=e_k\otimes\nabla\phi_k-\nabla\phi_k\otimes e_k .
\end{equation}

\begin{lem}[Properties of Mikado flows]\label{lem:mikado}
For $\mu\ge\mu_0$ and $k,k'\in\Xi$:
\begin{enumerate}[label=(\alph*)]
\item $W_k,\Omega_k\in C^\infty(\TT^d)$, $\Omega_k$ is skew-symmetric, and
\[
\D\Omega_k=W_k,\qquad \D W_k=0,\qquad \D(W_k\otimes W_k)=0 .
\]
\item For $m\ge0$ and $1\le p\le\infty$, uniformly in $\mu$,
\[
\|\nabla^mW_k\|_{L^p}\lesssim_m\mu^{m+\frac{d-1}2-\frac{d-1}p},\qquad
\|\nabla^m\Omega_k\|_{L^p}\lesssim_m\mu^{m-1+\frac{d-1}2-\frac{d-1}p}.
\]
\item $\fint W_k\otimes W_k=e_k\otimes e_k$, $\fint W_k=0$, and, for $k\ne k'$ and $1\le p\le\infty$,
$\|W_k\otimes W_{k'}\|_{L^p}\lesssim\mu^{d-1-\frac dp}$.
\end{enumerate}
\end{lem}

\begin{proof}
(a) Smoothness follows from the smoothness of $\dist(\cdot,l_k)$ on the support of the profiles. $\Omega_k$ is skew by
definition, and by \eqref{eq:lap}
\[
(\D\Omega_k)_i=\partial_j\big(e_{k,i}\partial_j\phi_k-\partial_i\phi_k\,e_{k,j}\big)=e_{k,i}\Delta\phi_k-\partial_i(e_k\cdot\nabla\phi_k)=\psi_ke_{k,i}.
\]
Hence $\D W_k=\D\D\Omega_k=0$, since $\Omega_k$ is skew. Finally,
$\big(\D(W_k\otimes W_k)\big)_i=e_{k,i}\,e_k\cdot\nabla(\psi_k^2)=0$.

(b) In adapted coordinates, $\nabla^m$ falls on $\psi(\mu|y|)$ and produces $\mu^m$, while the $L^p$ norm of
$\mu^{\frac{d-1}2}\psi(\mu|y|)$ over a tube of length $|k|$ is $\simeq\mu^{\frac{d-1}2-\frac{d-1}p}$. The same applies to
$\nabla\phi_k$, which carries the factor $\mu^{\frac{d-1}2-2}\cdot\mu$.

(c) By \eqref{eq:lap} and $|\TT^d|=1$, $\fint W_k\otimes W_k=e_k\otimes e_k\int\psi_k^2=e_k\otimes e_k$. Also
$\fint W_k=\fint\D\Omega_k=0$. For $k\ne k'$, $W_k\otimes W_{k'}$ is supported in the set \eqref{eq:intersect}, which
has measure $\lesssim\mu^{-d}$ and is empty if $k\parallel k'$. Since $\|W_k\|_{L^\infty}\lesssim\mu^{\frac{d-1}2}$,
$\|W_k\otimes W_{k'}\|_{L^p}\lesssim\mu^{d-1}\mu^{-\frac dp}$.
\end{proof}

\subsection{Temporal intermittency}\label{sec:temporal}

Fix a nonnegative $g\in C^\infty_c((0,1))$ with $\int_0^1g^2=1$. For $\kappa\ge1$ let $g_\kappa$ be the $1$-periodic extension
of $t\mapsto\kappa^{1/2}g(\kappa t)$ from $[0,1)$ to $\RR$, and set
\[
h_\kappa(t):=\int_0^t\big(g_\kappa^2(\tau)-1\big)\dd\tau .
\]

\begin{lem}\label{lem:temporal}
For $\kappa\ge1$ and $\nu\in\NN$:
\begin{enumerate}[label=(\alph*)]
\item $g_\kappa$ is smooth and $1$-periodic, and for $n\ge0$, $1\le p\le\infty$,
\[
\big\|\partial_t^n\big(g_\kappa(\nu\cdot)\big)\big\|_{L^p(0,1)}\lesssim_n(\nu\kappa)^n\kappa^{\frac12-\frac1p}.
\]
\item $h_\kappa$ is smooth, $1$-periodic, $\|h_\kappa\|_{L^\infty}\le1$, and
\begin{equation}\label{eq:hprime}
\partial_t\big(\nu^{-1}h_\kappa(\nu t)\big)=g_\kappa^2(\nu t)-1 .
\end{equation}
\item (Temporal decorrelation.) For every $F\in C^1([0,1])$,
\[
\Big|\int_0^1g_\kappa^2(\nu t)F(t)\dd t-\int_0^1F(t)\dd t\Big|\le2\nu^{-1}\|F'\|_{L^\infty(0,1)} .
\]
\end{enumerate}
\end{lem}

\begin{proof}
(a) $\kappa^{1/2}g(\kappa\cdot)$ is supported in $(0,\kappa^{-1})\subset(0,1)$, so its periodic extension is smooth. On
$[0,1]$ the function $g_\kappa(\nu\cdot)$ consists of $\nu$ translated copies of $\kappa^{1/2}g(\nu\kappa\,\cdot)$, which gives
$\|\partial_t^n(g_\kappa(\nu\cdot))\|_{L^p(0,1)}=(\nu\kappa)^n\kappa^{\frac12-\frac1p}\|g^{(n)}\|_{L^p}$.

(b) Since $\int_0^1(g_\kappa^2-1)=0$, $h_\kappa$ is $1$-periodic. For $t\in[0,1]$, $h_\kappa(t)=\int_0^tg_\kappa^2-t$ with both
terms in $[0,1]$, so $|h_\kappa(t)|\le1$. The identity \eqref{eq:hprime} is the chain rule.

(c) Let $J_j=[j/\nu,(j+1)/\nu]$, $0\le j<\nu$. By periodicity, $\int_{J_j}(g_\kappa^2(\nu t)-1)\dd t=\nu^{-1}\int_0^1(g_\kappa^2-1)=0$, hence
\[
\begin{aligned}
\Big|\int_{J_j}\big(g_\kappa^2(\nu t)-1\big)F(t)\dd t\Big|&=\Big|\int_{J_j}\big(g_\kappa^2(\nu t)-1\big)\big(F(t)-F(j/\nu)\big)\dd t\Big|\\
&\le\nu^{-1}\|F'\|_{L^\infty}\int_{J_j}\big(g_\kappa^2(\nu t)+1\big)\dd t=2\nu^{-2}\|F'\|_{L^\infty}.
\end{aligned}
\]
Summing over $j$ gives the claim.
\end{proof}

Part (c) replaces the use of \cite[Lemma~B.1]{CL22} in time. That lemma is stated for periodic functions. In \cite{CL22} the stress
vanishes near $t=0$ and $t=1$, so it applies to $t\mapsto\int\rho\dd x$; here $R$ need not vanish near $t=1$, and the corresponding
function of time is not periodic.

\subsection{Cutoffs}\label{sec:cutoffs}

\paragraph{Stress cutoff.} Let $\tilde\chi\in C^\infty([0,\infty))$ be nondecreasing with
\[
\tilde\chi(r)=4\|R\|_1\ \ (0\le r\le\|R\|_1),\qquad \tilde\chi(r)=4r\ \ (r\ge2\|R\|_1).
\]
Such a function exists since $4\|R\|_1<8\|R\|_1=4\cdot2\|R\|_1$. Set $\chi(A):=\tilde\chi(|A|)$ for symmetric $A$, which is
smooth because $\tilde\chi$ is constant near $0$, and
\begin{equation}\label{eq:rho}
\rho:=\chi(R)\in C^\infty([0,1]\times\TT^d).
\end{equation}

\begin{lem}\label{lem:rho}
On $[0,1]\times\TT^d$,
\[
4\|R\|_1\le\rho\le8\|R\|_1+4|R|,\qquad \mathrm{Id}-\frac R\rho\in\mathcal N .
\]
In particular $\rho^{1/2}$ and $\Gamma_k(\mathrm{Id}-R/\rho)$ are smooth, and $\int_0^1\!\int_{\TT^d}\rho\dd x\dd t\le12\|R\|_1$.
\end{lem}

\begin{proof}
The bounds on $\rho$ follow from the monotonicity of $\tilde\chi$. If $|R|\le2\|R\|_1$, then
$|R|/\rho\le2\|R\|_1/(4\|R\|_1)=\frac12$; if $|R|\ge2\|R\|_1$, then $|R|/\rho=\frac14$. Hence
$|(\mathrm{Id}-R/\rho)-\mathrm{Id}|\le\frac12$, and $\mathrm{Id}-R/\rho$ is symmetric.
\end{proof}

\paragraph{Temporal cutoff.} Fix $\theta\in C^\infty(\RR;[0,1])$ with
\begin{equation}\label{eq:theta}
\theta=0\ \text{on }(-\infty,s],\qquad \theta=1\ \text{on }[t_*,\infty).
\end{equation}
Since $R=0$ on $[0,t_*]$,
\begin{equation}\label{eq:thetaR}
\theta R=\theta^2R=R\qquad\text{on }[0,1]\times\TT^d .
\end{equation}
This is the only property of $\theta$ used in the algebra below. In \cite{CL22} the cutoff (4.12) is two-sided and
adapted to the well-prepared set $I$; we do not use the concentration step, and a one-sided cutoff suffices.

\subsection{The velocity perturbation}\label{sec:perturbation}

\paragraph{Amplitudes.} For $k\in\Xi$ let
\begin{equation}\label{eq:ak}
a_k(t,x):=\theta(t)\,g_\kappa(\nu t)\,\rho^{1/2}(t,x)\,\Gamma_k\Big(\mathrm{Id}-\frac{R(t,x)}{\rho(t,x)}\Big).
\end{equation}
By Lemma~\ref{lem:rho}, $a_k\in C^\infty([0,1]\times\TT^d)$.

\paragraph{Principal part and divergence corrector.}
\begin{equation}\label{eq:wpc}
w^{(p)}:=\sum_{k\in\Xi}a_k\,W_k(\sigma x),\qquad
w^{(c)}:=\sigma^{-1}\sum_{k\in\Xi}\nabla a_k\cdot\Omega_k(\sigma x),
\end{equation}
where $(\nabla a\cdot\Omega)_i=\partial_ja\,\Omega_{ij}$. Since $\D(\Omega_k(\sigma x))=\sigma W_k(\sigma x)$,
\begin{equation}\label{eq:divform}
w^{(p)}+w^{(c)}=\sigma^{-1}\D\Big(\sum_{k\in\Xi}a_k\,\Omega_k(\sigma x)\Big).
\end{equation}

\paragraph{Temporal corrector.}
\begin{equation}\label{eq:wt}
w^{(t)}:=\nu^{-1}h_\kappa(\nu t)\big(\D R-\nabla\Delta^{-1}\D\D R\big),
\end{equation}
which is well defined since $\D\D R=\partial_i\partial_jR_{ij}$ has zero mean. Finally,
\begin{equation}\label{eq:w}
w:=w^{(p)}+w^{(c)}+w^{(t)} .
\end{equation}

\begin{lem}\label{lem:w}
\begin{enumerate}[label=(\alph*)]
\item $w\in C^\infty([0,1]\times\TT^d;\RR^d)$, $\D w=0$, and $\fint w(t)=0$ for all $t$.
\item $a_k=0$ on $[0,s]\times\TT^d$ and $w^{(t)}=0$ on $[0,t_*]\times\TT^d$. In particular $w=0$ on $[0,s]\times\TT^d$.
\item On $[0,1]\times\TT^d$,
\begin{equation}\label{eq:algebra}
\sum_{k\in\Xi}a_k^2\fint_{\TT^d}W_k\otimes W_k\dd x=\theta^2g_\kappa^2(\nu t)\,\rho\,\mathrm{Id}-g_\kappa^2(\nu t)R .
\end{equation}
\end{enumerate}
\end{lem}

\begin{proof}
(a) Smoothness is clear. By \eqref{eq:divform}, $w^{(p)}+w^{(c)}$ is the divergence of a skew-symmetric matrix field,
hence divergence-free and of zero mean. For \eqref{eq:wt}, $\D w^{(t)}=\nu^{-1}h_\kappa(\nu t)(\D\D R-\Delta\Delta^{-1}\D\D R)=0$,
since $\D\D R$ has zero mean, and $w^{(t)}$ has zero mean as the sum of a divergence and a gradient.

(b) $\theta=0$ on $[0,s]$ gives $a_k=0$, hence $w^{(p)}=w^{(c)}=0$ there. $R=0$ on $[0,t_*]$ gives $w^{(t)}=0$ there.

(c) By Lemma~\ref{lem:mikado}(c), Lemma~\ref{lem:geometric} and Lemma~\ref{lem:rho},
\[
\sum_ka_k^2e_k\otimes e_k=\theta^2g_\kappa^2(\nu t)\rho\sum_k\Gamma_k^2\Big(\mathrm{Id}-\frac R\rho\Big)e_k\otimes e_k
=\theta^2g_\kappa^2(\nu t)\rho\,\mathrm{Id}-\theta^2g_\kappa^2(\nu t)R,
\]
and $\theta^2R=R$ by \eqref{eq:thetaR}.
\end{proof}

\subsection{The new stress}\label{sec:newstress}

\paragraph{Oscillation.} Let
\[
T_k:=W_k\otimes W_k-\fint W_k\otimes W_k,
\]
which is smooth, symmetric, of zero mean, and satisfies $\D T_k=0$ by Lemma~\ref{lem:mikado}(a). Expanding the square
and separating the diagonal terms,
\begin{equation}\label{eq:square}
w^{(p)}\otimes w^{(p)}=\sum_ka_k^2\,T_k(\sigma x)+\sum_ka_k^2\fint W_k\otimes W_k+R_{\rm far}^{\rm full},\qquad
R_{\rm far}^{\rm full}:=\sum_{k\ne k'}a_ka_{k'}\,W_k(\sigma x)\otimes W_{k'}(\sigma x).
\end{equation}
The sum in $R_{\rm far}^{\rm full}$ runs over ordered pairs, so $R_{\rm far}^{\rm full}$ is symmetric; it is not traceless. Let
\[
R_{\rm far}:=R_{\rm far}^{\rm full}-\tfrac1d\big(\mathrm{tr}\,R_{\rm far}^{\rm full}\big)\mathrm{Id}\in\Sz .
\]
Since $\D(T_k(\sigma\cdot))=\sigma(\D T_k)(\sigma\cdot)=0$, we have
\[
\D\Big(\sum_ka_k^2T_k(\sigma x)\Big)=\sum_k\nabla(a_k^2)\cdot T_k(\sigma x),
\]
which has zero mean, being a divergence. Define
\begin{equation}\label{eq:Roscx}
R_{{\rm osc},x}:=\sum_{k\in\Xi}\Bs\big(\nabla(a_k^2),T_k(\sigma\cdot)\big),\qquad
R_{{\rm osc},t}:=\nu^{-1}h_\kappa(\nu t)\,\partial_tR .
\end{equation}
By Lemma~\ref{lem:antidiv}(c), $R_{{\rm osc},x}$ is $\Sz$-valued and $\D R_{{\rm osc},x}=\sum_k\nabla(a_k^2)\cdot T_k(\sigma x)$.

\begin{lem}[Oscillation identity]\label{lem:osc}
Let $R_{\rm osc}:=R_{{\rm osc},x}+R_{{\rm osc},t}+R_{\rm far}$ and
\begin{equation}\label{eq:P}
P:=-\theta^2g_\kappa^2(\nu t)\rho-\tfrac1d\,\mathrm{tr}\,R_{\rm far}^{\rm full}
+\big(g_\kappa^2(\nu t)-1\big)\Delta^{-1}\D\D R+\nu^{-1}h_\kappa(\nu t)\,\Delta^{-1}\D\D\,\partial_tR .
\end{equation}
Then
\begin{equation}\label{eq:oscid}
\partial_tw^{(t)}+\D\big(w^{(p)}\otimes w^{(p)}+R\big)+\nabla P=\D R_{\rm osc}.
\end{equation}
\end{lem}

\begin{proof}
By \eqref{eq:square} and \eqref{eq:algebra},
\[
\D\big(w^{(p)}\otimes w^{(p)}+R\big)=\D R_{{\rm osc},x}+\D R_{\rm far}+\tfrac1d\nabla\,\mathrm{tr}\,R_{\rm far}^{\rm full}
+\nabla\big(\theta^2g_\kappa^2(\nu t)\rho\big)+\big(1-g_\kappa^2(\nu t)\big)\D R ,
\]
where we used that $g_\kappa(\nu t)$ does not depend on $x$. By \eqref{eq:hprime},
\[
\partial_tw^{(t)}=\big(g_\kappa^2(\nu t)-1\big)\big(\D R-\nabla\Delta^{-1}\D\D R\big)
+\nu^{-1}h_\kappa(\nu t)\big(\D\partial_tR-\nabla\Delta^{-1}\D\D\partial_tR\big).
\]
The terms $\pm(g_\kappa^2(\nu t)-1)\D R$ cancel, and $\nu^{-1}h_\kappa(\nu t)\D\partial_tR=\D R_{{\rm osc},t}$. All remaining terms
are gradients, and they are cancelled by $\nabla P$.
\end{proof}

\begin{remark}
By \eqref{eq:hprime}, $\nu^{-1}\partial_t\big(R\,h_\kappa(\nu t)\big)=(g_\kappa^2(\nu t)-1)R+\nu^{-1}h_\kappa(\nu t)\partial_tR$, so
\[
P=-\theta^2g_\kappa^2(\nu t)\rho-\tfrac1d\,\mathrm{tr}\,R_{\rm far}^{\rm full}+\nu^{-1}\Delta^{-1}\D\D\,\partial_t\big(R\,h_\kappa(\nu t)\big).
\]
Thus \eqref{eq:P} differs from \cite[(4.25)]{CL22} only by the trace part of the far-field interaction and by the sign with which the
last term appears. These terms do not affect any estimate of the stress. Including the trace in the pressure makes
the new stress traceless, as required by the definition of the Reynolds system.
\end{remark}

\paragraph{Linear and correction errors.} As in \cite[(4.26)--(4.27)]{CL22}, define
\begin{align}
R_{\rm lin}&:=\Rs\Big(\partial_t\big(w^{(p)}+w^{(c)}\big)-\Delta w+\D(u\otimes w+w\otimes u)\Big),\label{eq:Rlin}\\
R_{\rm cor}&:=\Rs\,\D\Big(\big(w^{(c)}+w^{(t)}\big)\otimes w+w^{(p)}\otimes\big(w^{(c)}+w^{(t)}\big)\Big).\label{eq:Rcor}
\end{align}
The arguments of $\Rs$ have zero mean: $\partial_t(w^{(p)}+w^{(c)})$ by \eqref{eq:divform}, $\Delta w$ as a Laplacian, and the
other terms as divergences. Hence
\begin{equation}\label{eq:divlincor}
\begin{aligned}
\D R_{\rm lin}&=\partial_t\big(w^{(p)}+w^{(c)}\big)-\Delta w+\D(u\otimes w+w\otimes u),\\
\D R_{\rm cor}&=\D\Big(\big(w^{(c)}+w^{(t)}\big)\otimes w+w^{(p)}\otimes\big(w^{(c)}+w^{(t)}\big)\Big).
\end{aligned}
\end{equation}

\begin{prop}[The algebraic part of Proposition~\ref{prop:NSstep}]\label{prop:CIidentity}
Let $w$ be given by \eqref{eq:w}, $P$ by \eqref{eq:P}, and
\begin{equation}\label{eq:RCL}
R^{\rm CL}:=R_{\rm lin}+R_{\rm cor}+R_{\rm osc}=R_{\rm lin}+R_{\rm cor}+R_{{\rm osc},x}+R_{{\rm osc},t}+R_{\rm far}.
\end{equation}
Then $w$, $P$ and $R^{\rm CL}$ are smooth on $[0,1]\times\TT^d$, $\D w=0$, $R^{\rm CL}$ takes values in $\Sz$, and properties
(i) and (ii) of Proposition~\ref{prop:NSstep} hold:
\begin{enumerate}[label=(\roman*)]
\item $\partial_tw-\Delta w+\D(u\otimes w+w\otimes u+w\otimes w)+\D R+\nabla P=\D R^{\rm CL}$;
\item $\fint w(t)=0$ for all $t$, and $w=0$, $R^{\rm CL}=0$ on $[0,s]\times\TT^d$.
\end{enumerate}
This holds for every $\nu,\sigma\in\NN$, $\kappa\ge1$ and $\mu\ge\mu_0$.
\end{prop}

\begin{proof}
All ingredients are smooth in $(t,x)$. $R_{\rm lin}$, $R_{\rm cor}$ are $\Sz$-valued by Lemma~\ref{lem:antidiv}(a),
$R_{{\rm osc},x}$ by Lemma~\ref{lem:antidiv}(c), $R_{{\rm osc},t}$ because $R$ is, and $R_{\rm far}$ by construction.

For the identity (i), since $w=w^{(p)}+(w^{(c)}+w^{(t)})$,
\[
w\otimes w=w^{(p)}\otimes w^{(p)}+\big(w^{(c)}+w^{(t)}\big)\otimes w+w^{(p)}\otimes\big(w^{(c)}+w^{(t)}\big).
\]
Hence, using $\partial_tw=\partial_t(w^{(p)}+w^{(c)})+\partial_tw^{(t)}$,
\[
\begin{aligned}
&\partial_tw-\Delta w+\D(u\otimes w+w\otimes u+w\otimes w)+\D R+\nabla P\\
&\quad=\Big[\partial_t\big(w^{(p)}+w^{(c)}\big)-\Delta w+\D(u\otimes w+w\otimes u)\Big]
+\D\Big[\big(w^{(c)}+w^{(t)}\big)\otimes w+w^{(p)}\otimes\big(w^{(c)}+w^{(t)}\big)\Big]\\
&\qquad+\Big[\partial_tw^{(t)}+\D\big(w^{(p)}\otimes w^{(p)}+R\big)+\nabla P\Big]
=\D R_{\rm lin}+\D R_{\rm cor}+\D R_{\rm osc},
\end{aligned}
\]
by \eqref{eq:divlincor} and Lemma~\ref{lem:osc}.

For (ii), the mean and the vanishing of $w$ on $[0,s]$ are Lemma~\ref{lem:w}(a),(b). Since $w$ and the $a_k$ vanish identically on
$[0,s]\times\TT^d$ and are smooth, all their derivatives vanish there as well (one-sided at $t=s$). Hence the arguments of
$\Rs$ in \eqref{eq:Rlin}--\eqref{eq:Rcor} vanish on $[0,s]$, and so do $R_{{\rm osc},x}$ (which is bilinear in
$\nabla(a_k^2)$ and acts in space only) and $R_{\rm far}$ (which contains $a_ka_{k'}$). Finally, $R=0$ on $[0,t_*]$ gives
$\partial_tR=0$ there, so $R_{{\rm osc},t}=0$ on $[0,s]$.
\end{proof}

\begin{remark}\label{rem:remains}
Proposition~\ref{prop:CIidentity} contains no smallness. In Appendix~\ref{sec:estimates} the parameters are fixed as powers
of a single frequency $\lambda$, and the bounds (iii)--(v) of Proposition~\ref{prop:NSstep} are proved for $\lambda$ large.
\end{remark}

\section{Estimates and proof of Proposition~\ref{prop:NSstep}}\label{sec:estimates}

In this appendix we fix the parameters of Appendix~\ref{sec:CI} and prove the bounds (iii)--(v) of
Proposition~\ref{prop:NSstep} for the perturbation and stress of Proposition~\ref{prop:CIidentity}. We follow
\cite[Section~5]{CL22}. The standing assumptions of Appendix~\ref{sec:CI} are in force; in particular $\|R\|_1>0$.

\paragraph{Constants.} We write $A\lesssim B$ if $A\le CB$ with $C$ depending only on $d,p,q$ and the fixed choices
$\Xi,\Gamma_k,p_k,\phi,g$, which depend only on $d$. We write $C_u$ for constants that may in addition depend on $u$,
$R$, $s$, $t_*$ and $\theta$. Neither kind of constant depends on $\lambda$, and both may change from line to line. On
$[0,1]$ and on $\TT^d$, which have unit measure, $L^p$ norms are nondecreasing in $p$; we use this without comment.

\paragraph{Reduction.} Since $\|\cdot\|_{W^{1,q}}\le\|\cdot\|_{W^{1,q'}}$ for $q\le q'$, it suffices to prove
Proposition~\ref{prop:NSstep} with $q$ replaced by $\max\{q,2\}$. We therefore assume $2\le q<\infty$.

\subsection{Choice of parameters}\label{sec:parameters}

Fix $m\in\NN$ so large that $\gamma:=1/m$ satisfies
\begin{equation}\label{eq:gamma}
10d\gamma\le\min\Big\{\frac1p-\frac12,\ \frac1q\Big\},
\end{equation}
and set
\begin{equation}\label{eq:r}
r:=\frac d{d-\gamma},\qquad\text{so that}\qquad d-\frac dr=\gamma,\qquad 1<r<\frac d{d-1}\le2\le q .
\end{equation}
Both $\gamma$ and $r$ depend only on $d,p,q$. Since $q\ge2$, \eqref{eq:gamma} gives $\gamma\le\frac1{20d}$. For $N\in\NN$
let $\lambda:=N^m$ and
\begin{equation}\label{eq:params}
\nu:=\lambda^\gamma=N,\qquad \sigma:=\lambda^{1/\gamma}=N^{m^2},\qquad \mu:=\lambda,\qquad
\kappa:=\lambda^{K},\quad K:=\frac2\gamma+d+1-6\gamma .
\end{equation}
Then $\nu,\sigma,\mu\in\NN$, $\kappa\ge1$, and, since $K/2-1/\gamma=\frac{d+1}2-3\gamma>0$,
\[
\nu\ll\mu\ll\sigma\ll\kappa^{1/2}.
\]
In the end $N$, and hence $\lambda$, is chosen large depending on $u,R,s,t_*,\delta$. In particular we always assume
$\mu=\lambda\ge\mu_0$.

\begin{lem}[Scaling inequalities]\label{lem:scaling}
For every $\lambda\ge1$:
\begin{align}
&\nu\kappa^{\frac12}\sigma^{-1}\mu^{-1+\frac{d-1}2-\frac{d-1}r}\le\lambda^{-\gamma},\label{eq:S1}\\
&\kappa^{\frac12-\frac1p}\mu^{\frac{d-1}2}\le\lambda^{-\gamma},\label{eq:S2}\\
&\kappa^{-\frac12}\sigma\mu^{1+\frac{d-1}2-\frac{d-1}q}\le\lambda^{-\gamma},\label{eq:S3}\\
&\sigma^{-1}\mu^{-1+\frac{d-1}2}\le\lambda^{-\gamma},\qquad
\sigma^{-1}\mu^{(d-1)(1-\frac1r)}\le\lambda^{-\gamma},\qquad
\mu^{d-1-\frac dr}\le\lambda^{-\gamma},\label{eq:S4}\\
&\nu^{-1}=\lambda^{-\gamma},\qquad \sigma^{-1/2}\le\lambda^{-\gamma},\qquad \nu^{-1/2}=\lambda^{-\gamma/2}.\label{eq:S5}
\end{align}
\end{lem}

\begin{proof}
We compare exponents of $\lambda$, using \eqref{eq:gamma}, \eqref{eq:r} and $\gamma\le\frac1{20d}$.

\eqref{eq:S1}: the exponent is
$\gamma+\big(\frac1\gamma+\frac{d+1}2-3\gamma\big)-\frac1\gamma-1+\frac{d-1}2-\frac{d-1}r=-2\gamma+(d-1)\big(1-\frac1r\big)$, and
$(d-1)(1-\frac1r)\le d-\frac dr=\gamma$.

\eqref{eq:S2}: since $\frac12-\frac1p\le-10d\gamma<0$ and $K\ge\frac2\gamma$, the exponent is at most
$-10d\gamma\cdot\frac2\gamma+\frac{d-1}2=-20d+\frac{d-1}2\le-\gamma$.

\eqref{eq:S3}: the exponent is
$-\frac1\gamma-\frac{d+1}2+3\gamma+\frac1\gamma+1+\frac{d-1}2-\frac{d-1}q=3\gamma-\frac{d-1}q$, and
$\frac{d-1}q\ge10d(d-1)\gamma\ge20\gamma$.

\eqref{eq:S4}: the exponents are $-\frac1\gamma-1+\frac{d-1}2\le-20d+\frac{d-1}2\le-\gamma$; then
$-\frac1\gamma+(d-1)(1-\frac1r)\le-\frac1\gamma+\gamma\le-\gamma$; and $d-1-\frac dr=\gamma-1\le-\gamma$.

\eqref{eq:S5} is immediate from $\frac1{2\gamma}\ge\gamma$.
\end{proof}

\subsection{Amplitudes}

Write $a_k=\theta(t)g_\kappa(\nu t)b_k(t,x)$ with $b_k:=\rho^{1/2}\Gamma_k(\mathrm{Id}-R/\rho)$. By Lemma~\ref{lem:rho},
$b_k$ is smooth on $[0,1]\times\TT^d$, and it does not depend on $\lambda$. Hence
$\|\partial_t^n\nabla^m(\theta b_k)\|_{L^\infty_{t,x}}\le C_{u,m,n}$.

\begin{lem}[Amplitudes]\label{lem:ak}
For $k\in\Xi$, $m,n\ge0$ and $1\le p'\le\infty$,
\begin{equation}\label{eq:ak-est}
\big\|\partial_t^n\nabla^ma_k\big\|_{L^{p'}(0,1;L^\infty(\TT^d))}\le C_{u,m,n}(\nu\kappa)^n\kappa^{\frac12-\frac1{p'}},\qquad
\big\|\nabla^m(a_k^2)\big\|_{L^1(0,1;L^\infty(\TT^d))}\le C_{u,m}.
\end{equation}
Moreover, for every $t\in[0,1]$,
\begin{equation}\label{eq:akL2}
\|a_k(t)\|_{L^2(\TT^d)}^2\le C_\Gamma^2\,\theta^2(t)g_\kappa^2(\nu t)\int_{\TT^d}\rho(t,x)\dd x,\qquad
C_\Gamma:=\max_{k\in\Xi}\sup_{\mathcal N}|\Gamma_k| .
\end{equation}
\end{lem}

\begin{proof}
By the Leibniz rule, $\partial_t^n\nabla^ma_k=\sum_{i\le n}\binom ni\partial_t^i\big(g_\kappa(\nu t)\big)\,\partial_t^{n-i}\nabla^m(\theta b_k)$.
By Lemma~\ref{lem:temporal}(a) and $\nu\kappa\ge1$, $\|\partial_t^i(g_\kappa(\nu\cdot))\|_{L^{p'}(0,1)}\lesssim(\nu\kappa)^n\kappa^{\frac12-\frac1{p'}}$ for
$i\le n$, which gives the first bound. For the second, $a_k^2=g_\kappa^2(\nu t)\,\theta^2b_k^2$ and $\|g_\kappa^2(\nu\cdot)\|_{L^1(0,1)}=1$.
The last bound follows from $a_k^2=\theta^2g_\kappa^2(\nu t)\rho\,\Gamma_k^2(\mathrm{Id}-R/\rho)$ and Lemma~\ref{lem:rho}.
\end{proof}

\subsection{The velocity perturbation}

We use the improved H\"older inequality of Modena and Sz\'ekelyhidi \cite[Lemma~2.1]{MSz18}, see also
\cite[Lemma~B.1]{CL22}: for $1\le p'\le\infty$, smooth $a,f$ on $\TT^d$ and $\sigma\in\NN$,
\begin{equation}\label{eq:holder}
\big|\,\|af(\sigma\cdot)\|_{L^{p'}}-\|a\|_{L^{p'}}\|f\|_{L^{p'}}\big|\lesssim\sigma^{-1/p'}\|a\|_{C^1}\|f\|_{L^{p'}} .
\end{equation}

\begin{prop}[Principal part]\label{prop:wp}
There is $C_0\ge1$ depending only on $d$ such that
\begin{align}
\|w^{(p)}\|_{L^2([0,1]\times\TT^d)}&\le C_0\|R\|_1^{1/2}+C_u\lambda^{-\gamma/2},\label{eq:wpL2}\\
\|w^{(p)}\|_{L^p(0,1;L^\infty)}+\|w^{(p)}\|_{L^1(0,1;W^{1,q})}&\le C_u\lambda^{-\gamma}.\label{eq:wpLp}
\end{align}
\end{prop}

\begin{proof}
\proofstep{$L^2$ bound}
By \eqref{eq:holder} with $p'=2$ and $\|W_k\|_{L^2}=1$ (Lemma~\ref{lem:mikado}(c)), for every $t$,
\[
\|w^{(p)}(t)\|_{L^2}\le\sum_{k\in\Xi}\Big(\|a_k(t)\|_{L^2}+C\sigma^{-1/2}\|a_k(t)\|_{C^1}\Big).
\]
Taking the $L^2$ norm in time, by \eqref{eq:akL2} and \eqref{eq:ak-est} (with $p'=2$),
\[
\|w^{(p)}\|_{L^2_{t,x}}\le|\Xi|C_\Gamma\Big(\int_0^1g_\kappa^2(\nu t)F(t)\dd t\Big)^{1/2}+C_u\sigma^{-1/2},\qquad
F(t):=\theta^2(t)\int_{\TT^d}\rho(t,x)\dd x .
\]
$F$ is smooth and independent of $\lambda$, so $\|F'\|_{L^\infty}\le C_u$. By Lemma~\ref{lem:temporal}(c) and
Lemma~\ref{lem:rho},
\[
\int_0^1g_\kappa^2(\nu t)F(t)\dd t\le\int_0^1F(t)\dd t+2\nu^{-1}C_u\le12\|R\|_1+C_u\nu^{-1}.
\]
Using $\sqrt{A+B}\le\sqrt A+\sqrt B$ and \eqref{eq:S5}, we obtain \eqref{eq:wpL2} with $C_0:=\sqrt{12}\,|\Xi|C_\Gamma$.

\proofstep{$L^p_tL^\infty$ bound}
By Lemma~\ref{lem:mikado}(b), \eqref{eq:ak-est} with $p'=p$, and \eqref{eq:S2},
\[
\|w^{(p)}\|_{L^p_tL^\infty}\le\sum_k\|a_k\|_{L^p_tL^\infty}\|W_k\|_{L^\infty}\le C_u\kappa^{\frac12-\frac1p}\mu^{\frac{d-1}2}\le C_u\lambda^{-\gamma}.
\]

\proofstep{$L^1_tW^{1,q}$ bound}
By \eqref{eq:rescale} and Lemma~\ref{lem:mikado}(b), using $\sigma\mu\ge1$,
\[
\|W_k(\sigma\cdot)\|_{W^{1,q}}=\|W_k\|_{L^q}+\sigma\|\nabla W_k\|_{L^q}\lesssim\sigma\mu\,\mu^{\frac{d-1}2-\frac{d-1}q}.
\]
Since $\|aW_k(\sigma\cdot)\|_{W^{1,q}}\lesssim\|a\|_{C^1}\|W_k(\sigma\cdot)\|_{W^{1,q}}$, \eqref{eq:ak-est} with $p'=1$ and \eqref{eq:S3} give
\[
\|w^{(p)}\|_{L^1_tW^{1,q}}\le\sum_k\|a_k\|_{L^1_tC^1}\|W_k(\sigma\cdot)\|_{W^{1,q}}\le C_u\kappa^{-\frac12}\sigma\mu^{1+\frac{d-1}2-\frac{d-1}q}\le C_u\lambda^{-\gamma}.
\qedhere
\]
\end{proof}

\begin{prop}[Correctors]\label{prop:wct}
\begin{align}
\|w^{(c)}\|_{L^2(0,1;L^\infty)}+\|w^{(c)}\|_{L^1(0,1;W^{1,q})}&\le C_u\lambda^{-\gamma},\label{eq:wc}\\
\|w^{(t)}\|_{L^\infty(0,1;W^{1,\infty})}&\le C_u\lambda^{-\gamma}.\label{eq:wt-est}
\end{align}
\end{prop}

\begin{proof}
For \eqref{eq:wc}: by Lemma~\ref{lem:mikado}(b), \eqref{eq:ak-est} with $p'=2$, and \eqref{eq:S4},
\[
\|w^{(c)}\|_{L^2_tL^\infty}\le\sigma^{-1}\sum_k\|\nabla a_k\|_{L^2_tL^\infty}\|\Omega_k\|_{L^\infty}\le C_u\sigma^{-1}\mu^{-1+\frac{d-1}2}\le C_u\lambda^{-\gamma}.
\]
Similarly, $\|\Omega_k(\sigma\cdot)\|_{W^{1,q}}\lesssim\sigma\mu\,\mu^{-1+\frac{d-1}2-\frac{d-1}q}$, so by \eqref{eq:ak-est} with $p'=1$,
\[
\|w^{(c)}\|_{L^1_tW^{1,q}}\lesssim\sigma^{-1}\sum_k\|a_k\|_{L^1_tC^2}\|\Omega_k(\sigma\cdot)\|_{W^{1,q}}
\le C_u\kappa^{-\frac12}\mu^{\frac{d-1}2-\frac{d-1}q}=C_u(\sigma\mu)^{-1}\kappa^{-\frac12}\sigma\mu^{1+\frac{d-1}2-\frac{d-1}q},
\]
which is at most $C_u\lambda^{-\gamma}$ by \eqref{eq:S3}.

For \eqref{eq:wt-est}: since $|h_\kappa|\le1$ (Lemma~\ref{lem:temporal}(b)), and since $\D R-\nabla\Delta^{-1}\D\D R$ is smooth and
independent of $\lambda$, we get $\|w^{(t)}\|_{L^\infty_tW^{1,\infty}}\le C_u\nu^{-1}=C_u\lambda^{-\gamma}$.
\end{proof}

\begin{cor}[The perturbation]\label{cor:w}
Let $M:=\max\{1,2C_0\}$, which depends only on $d$. There is $\lambda_1$, depending on $u,R,s,t_*,\delta$,
such that, for $\lambda\ge\lambda_1$,
\[
\|w\|_{L^2_{t,x}}\le M\|R\|_1^{1/2},\qquad \|w\|_{L^p_tL^\infty}+\|w\|_{L^1_tW^{1,q}}\le C_u\lambda^{-\gamma}\le\delta .
\]
Moreover $\|w^{(c)}\|_{L^2_tL^\infty}+\|w^{(t)}\|_{L^2_tL^\infty}\le C_u\lambda^{-\gamma}$ and $\|w\|_{L^2_{t,x}}+\|w^{(p)}\|_{L^2_{t,x}}\le C_u$.
\end{cor}

\begin{proof}
Add Propositions~\ref{prop:wp} and~\ref{prop:wct}, using $\|\cdot\|_{L^2_{t,x}}\le\|\cdot\|_{L^2_tL^\infty}$,
$\|\cdot\|_{L^p_tL^\infty}\le\|\cdot\|_{L^2_tL^\infty}$ (since $p<2$) and $\|\cdot\|_{L^1_tW^{1,q}}\le\|\cdot\|_{L^\infty_tW^{1,\infty}}$. This gives
$\|w\|_{L^2_{t,x}}\le C_0\|R\|_1^{1/2}+C_u\lambda^{-\gamma/2}$. Since $\|R\|_1>0$, the last term is at most $C_0\|R\|_1^{1/2}$ for
$\lambda$ large.
\end{proof}

\subsection{The new stress}

Recall from Lemma~\ref{lem:antidiv} that $\Rs$ is bounded on $L^r$, that $\Rs\D$ is bounded on $L^r$ since $1<r<\infty$, and
that $\|\Rs(g(\sigma\cdot))\|_{L^r}\lesssim\sigma^{-1}\|g\|_{L^r}$ for $g$ of zero mean.

\begin{lem}[Linear error]\label{lem:lin}
$\|R_{\rm lin}\|_{L^1(0,1;L^r)}\le C_u\lambda^{-\gamma}$.
\end{lem}

\begin{proof}
We treat the three terms of \eqref{eq:Rlin} separately.

\textit{Dissipative term.} $\Rs\Delta w=\Rs\D(\nabla w)$, so $\|\Rs\Delta w\|_{L^1_tL^r}\lesssim\|\nabla w\|_{L^1_tL^r}\le\|w\|_{L^1_tW^{1,q}}\le C_u\lambda^{-\gamma}$,
since $r<q$.

\textit{Time derivative.} By \eqref{eq:divform},
$\partial_t(w^{(p)}+w^{(c)})=\sigma^{-1}\D\big(\sum_k\partial_ta_k\,\Omega_k(\sigma x)\big)$. Hence, by Lemma~\ref{lem:mikado}(b),
\eqref{eq:ak-est} with $n=1$, $p'=1$, and \eqref{eq:S1},
\[
\big\|\Rs\partial_t(w^{(p)}+w^{(c)})\big\|_{L^1_tL^r}\lesssim\sigma^{-1}\sum_k\|\partial_ta_k\|_{L^1_tL^\infty}\|\Omega_k\|_{L^r}
\le C_u\sigma^{-1}\nu\kappa^{\frac12}\mu^{-1+\frac{d-1}2-\frac{d-1}r}\le C_u\lambda^{-\gamma}.
\]

\textit{Transport term.} $\|\Rs\D(u\otimes w+w\otimes u)\|_{L^1_tL^r}\lesssim\|u\|_{L^\infty_{t,x}}\|w\|_{L^1_tL^r}\le C_u\|w\|_{L^p_tL^\infty}\le C_u\lambda^{-\gamma}$.
\end{proof}

\begin{lem}[Correction error]\label{lem:cor}
$\|R_{\rm cor}\|_{L^1(0,1;L^r)}\le C_u\lambda^{-\gamma}$.
\end{lem}

\begin{proof}
By the boundedness of $\Rs\D$ on $L^r$ and H\"older's inequality in $t$ and $x$, using $r<2$,
\[
\|R_{\rm cor}\|_{L^1_tL^r}\lesssim\big(\|w^{(c)}\|_{L^2_tL^\infty}+\|w^{(t)}\|_{L^2_tL^\infty}\big)\big(\|w\|_{L^2_{t,x}}+\|w^{(p)}\|_{L^2_{t,x}}\big)\le C_u\lambda^{-\gamma},
\]
by Corollary~\ref{cor:w}.
\end{proof}

\begin{lem}[Oscillation error]\label{lem:osc-est}
$\|R_{\rm osc}\|_{L^1(0,1;L^r)}\le C_u\lambda^{-\gamma}$.
\end{lem}

\begin{proof}
\textit{Spatial oscillation.} $T_k$ has zero mean, so by Lemma~\ref{lem:antidiv}(a),(c),
\[
\big\|\Bs\big(\nabla(a_k^2),T_k(\sigma\cdot)\big)\big\|_{L^r}\lesssim\|\nabla(a_k^2)\|_{C^1}\|\Rs(T_k(\sigma\cdot))\|_{L^r}\lesssim\sigma^{-1}\|a_k^2\|_{C^2}\|T_k\|_{L^r}.
\]
By Lemma~\ref{lem:mikado}(b),(c), $\|T_k\|_{L^r}\le\|W_k\|_{L^{2r}}^2+1\lesssim\mu^{(d-1)(1-\frac1r)}$. With \eqref{eq:ak-est} and \eqref{eq:S4},
$\|R_{{\rm osc},x}\|_{L^1_tL^r}\le C_u\sigma^{-1}\mu^{(d-1)(1-\frac1r)}\le C_u\lambda^{-\gamma}$.

\textit{Temporal oscillation.} Since $|h_\kappa|\le1$,
$\|R_{{\rm osc},t}\|_{L^1_tL^r}\le\nu^{-1}\|\partial_tR\|_{L^1_tL^r}\le C_u\nu^{-1}=C_u\lambda^{-\gamma}$.

\textit{Far-field interaction.} The traceless part is an orthogonal projection, so $|R_{\rm far}|\le|R_{\rm far}^{\rm full}|$. By
\eqref{eq:rescale}, $W_k(\sigma\cdot)\otimes W_{k'}(\sigma\cdot)=(W_k\otimes W_{k'})(\sigma\cdot)$ has the same $L^r$ norm as
$W_k\otimes W_{k'}$. Hence, by Lemma~\ref{lem:mikado}(c), \eqref{eq:ak-est} with $p'=2$, and \eqref{eq:S4},
\[
\|R_{\rm far}\|_{L^1_tL^r}\le\sum_{k\ne k'}\|a_k\|_{L^2_tL^\infty}\|a_{k'}\|_{L^2_tL^\infty}\|W_k\otimes W_{k'}\|_{L^r}
\le C_u\mu^{d-1-\frac dr}\le C_u\lambda^{-\gamma}.
\qedhere
\]
\end{proof}

\subsection{Proof of Proposition~\ref{prop:NSstep}}\label{sec:proofNS}

Let $d,p,q$ be given, and assume $q\ge2$ as explained at the beginning of this section. Let $r=r(d,p,q)$ be defined by
\eqref{eq:r}; then $1<r<\frac d{d-1}$. Let $M:=\max\{1,2C_0\}$, which depends only on $d$. Let $(u,R)$, $0<s<t_*\le1$ and $\delta>0$
be given. If $R\equiv0$, take $w=0$, $P=0$, $R^{\rm CL}=0$. Otherwise, let $w,P,R^{\rm CL}$ be given by
Proposition~\ref{prop:CIidentity} with the parameters \eqref{eq:params}. Then:
\begin{itemize}
\item smoothness, $\D w=0$, the $\Sz$-valuedness of $R^{\rm CL}$, and properties (i) and (ii) hold for every $\lambda$, by
Proposition~\ref{prop:CIidentity};
\item properties (iii) and (iv) hold for $\lambda\ge\lambda_1$, by Corollary~\ref{cor:w};
\item by \eqref{eq:RCL} and Lemmas~\ref{lem:lin}--\ref{lem:osc-est}, $\|R^{\rm CL}\|_{L^1_tL^r}\le C_u\lambda^{-\gamma}$, which is at most $\delta$
for $\lambda$ large; this is (v).
\end{itemize}
Choosing $N$, and hence $\lambda=N^m$, large enough (also so that $\lambda\ge\mu_0$) completes the proof.
\qed

\begin{remark}\label{rem:final}
\begin{enumerate}
\item In particular $\|w\|_{L^1_tL^\infty}\le\|w\|_{L^p_tL^\infty}\le C_u\lambda^{-\gamma}$. This is the only norm of $w$ that enters the
kinetic estimates of Section~\ref{sec:kinetic}.
\item The Laplacian plays no role in the construction. Proposition~\ref{prop:NSstep} therefore holds without the term
$-\Delta w$, the dissipative term in Lemma~\ref{lem:lin} being simply absent. This is the starting point for the inviscid
Vlasov--Euler system.
\end{enumerate}
\end{remark}

\end{document}